\documentclass{amsart}

\usepackage{amssymb}
\usepackage{eepic}
\usepackage{color}
\usepackage[linktocpage=true]{hyperref}
\usepackage{xypic}
\usepackage{enumitem}
\usepackage{mdwlist}
\usepackage{tabularx}

\allowdisplaybreaks

\numberwithin{equation}{section}

\newtheorem{theorem}{Theorem}[section]
\newtheorem{lemma}[theorem]{Lemma}
\newtheorem{proposition}[theorem]{Proposition}
\newtheorem{corollary}[theorem]{Corollary}
\newtheorem{remark}[theorem]{Remark}

\newtheorem{Atheorem}{Theorem}

\newtheorem{Aremark}[Atheorem]{Remark}

\newtheorem{Btheorem}{Theorem}
\newtheorem{Bproposition}[Btheorem]{Proposition}
\newtheorem{Bcorollary}[Btheorem]{Corollary}
\newtheorem{Blemma}[Btheorem]{Lemma}
\newtheorem{Bremark}[Btheorem]{Remark}

\newtheorem{TheoA}{Theorem A}

\newtheorem{TheoB}{Theorem B}

\newcommand{\N}{\mathbb{N}}
\newcommand{\Z}{\mathbb{Z}}
\newcommand{\R}{\mathbb{R}}
\newcommand{\C}{\mathbb{C}}
\newcommand{\F}{\mathbb{F}}

\newcommand{\summ}{\sum\nolimits}

\def\G{\mathrm{G}}
\def\1{\mathbf{1}}
\def\H{\mathcal{H}}
\def\Q{\mathcal{Q}}
\def\M{\mathcal{M}}
\def\wast{\mathrm{w}^\ast}
\def\A{\mathcal{A}}
\def\B{\mathcal{B}}
\def\k{\mathrm{k}}
\def\RR{\mathcal{R}}

\def\S{\mathcal{S}}

\def\op{\mathrm{op}}

\def\dTh{\mathrm{d}_\Theta}
\def\BMO{\mathrm{BMO}}

\def\CB{\mathcal{CB}}
\def\Hardy{\mathrm{H}}

\newcommand{\dem}{\noindent {\bf Proof. }}

\newcommand{\demA}{\noindent {\bf Proof of Theorem \ref{ThmBTorus}} i) {\bf and} ii){\bf.} }
\newcommand{\demAA}{\noindent {\bf Proof of Theorem \ref{ThmBTorus}} iii){\bf.} }

\newcommand{\fin}{\hspace*{\fill} $\square$ \vskip0.2cm}

\def\mean{- \hskip-11.2pt \int}

\def\meannn{- \hskip-9.5pt \int}

\addtocounter{tocdepth}{0}

\begin{document}

\null

\vskip-45pt

\null

\title[Singular integrals in quantum Euclidean spaces]{Singular integrals \\ in \hskip1pt quantum Euclidean spaces}

\author[Gonz\'alez-P\'erez, Junge, Parcet]
{Adri\'an M. Gonz\'alez-P\'erez, \\ Marius Junge and Javier Parcet}

\maketitle

\null

\vskip-45pt

\null

\begin{abstract}
In this paper, we establish the core of singular integral theory and pseudodifferential calculus over the archetypal algebras of noncommutative geometry: quantum forms of Euclidean spaces and tori. Our results go beyond Connes' pseudodifferential calculus for rotation algebras, thanks to a new form of Calder\'on-Zygmund theory over these spaces which crucially incorporates nonconvolution kernels. We deduce $L_p$-boundedness and Sobolev $p$-estimates for regular, exotic and forbidden symbols in the expected ranks. In the $L_2$ level both Calder\'on-Vaillancourt and Bourdaud theorems for exotic and forbidden symbols are also generalized to the quantum setting. As a basic application of our methods, we prove $L_p$-regularity of solutions for elliptic PDEs.  
\end{abstract}

\tableofcontents

\addtolength{\parskip}{+1ex}

\null

\vskip-60pt

\null

\section*{\bf Introduction}

Harmonic analysis and PDEs over Riemannian manifolds are paramount for the solution of many important problems in differential geometry, fluid mechanics or theoretical physics. In this paper, we establish the core of singular integral theory and pseudodifferential calculus over the archetypal algebras of noncommutative geometry. This includes the Heisenberg-Weyl algebra, quantum tori and other noncommutative deformations of Euclidean spaces of great interest in quantum field theory and quantum probability. Our approach crucially relies on a quantum form of the fruitful interplay \[ {\small \begin{array}{ccc} 
&\textsc{Kernels}& \\ [-3pt]
& \swarrow \hskip40pt \searrow & \\
& \hskip7pt \textsc{Symbols} \ \longrightarrow \ \textsc{Operators} &
\end{array}} \] at the interface of analysis and geometry. Strong reasons to develop such a program over matrix algebras and other noncommutative manifolds are also in connection to string theory, where several PDEs arise naturally over quantum spaces. We obtain optimal smoothness conditions for $L_p$-boundedness of singular integrals and corresponding Sobolev $p$-estimates for pseudodifferential operators. This is crucial for applications to PDEs, which we shall briefly discuss. In the line of the harmonic analysis school, a key point has been a profound analysis of the associated kernels which is specially challenging for noncommutative algebras. 

Let $\Theta$ be an anti-symmetric real $n \times n$ matrix. Roughly speaking, the quantum Euclidean space $\mathcal{R}_\Theta$ is the von Neumann algebra generated by certain family of unitaries $\{u_j(s) : 1 \le j \le n, \, s \in \R\}$ satisfying $$u_j(s)u_j(t) = u_j(s+t),$$ $$u_j(s) u_k(t) = e^{2\pi i \Theta_{jk} st} u_k(t) u_j(s).$$ Set $\lambda_\Theta(\xi) = u_1(\xi_1) u_2(\xi_2) \cdots u_n(\xi_n)$ for $\xi \in \R^n$ and $$\lambda_\Theta: \mathcal{C}_c(\R^n) \ni f \longmapsto \int_{\R^n} f(\xi) \lambda_\Theta(\xi) \, d\xi \in \mathcal{R}_\Theta.$$ Consider the trace determined by $\tau_\Theta(\lambda_\Theta(f)) = f(0)$ and the corresponding $L_p$ spaces $L_p(\RR_\Theta,\tau_\Theta)$ \cite{PX2}. Of course, $\Theta=0$ yields the Euclidean $L_p$-space in $\R^n$ with the Lebesgue measure and $(\RR_\Theta,\tau_\Theta)$ should be understood as a noncommutative deformation of it. Section \ref{Prelim} includes a careful presentation of $(\RR_\Theta,\tau_\Theta)$ for those potential readers not familiar with them. Our approach also contains a key Poincar\'e type inequality and a few more crucial results, maybe some known to experts. The lack of appropriate literature justifies a self-contained presentation.

The algebraic structure of these operator (type I) algebras is quite simple, but the connection to Euclidean spaces make them indispensable in a great variety of scenarios. If $\hbar$ stands for Planck's constant, the choice $$\Theta =  2 \pi \hbar \ id_{M_m} \otimes \Big( \hskip-5pt \begin{array}{rl} 0 & 1 \\ -1 & 0 \end{array} \Big)$$ yields the Heisenberg-Weyl algebra in quantum mechanics. Another description arises from the unbounded generators $x_{\Theta,j}$ of $u_j$ |a quantum analogue of the Euclidean variables| which satisfy $2 \pi i [x_{\Theta,j}, x_{\Theta,k}] = \Theta_{jk}$ and provide additional insight in our kernel manipulations below. Considering a Fock space representation $\mathcal{R}_\Theta$ becomes the CCR algebra associated to the symplectic form $\Theta$, thoroughly studied in quantum probability and quantum field theory. In this setting, it is simple to find a (nontracial) gaussian state with respect to which these $x_{\Theta,j}$ admit a gaussian distribution. In the physics literature, higher dimensional deformations are usually referred to as Moyal deformations of $\R^n$. An important instance in string theory is given by the noncommutative deformation of $\R^4$ associated to an invertible symbol $\Theta$, which leads to instantons on a noncommutative space in the influential papers \cite{CL,NS,SW}. In view of so many names for the same object, we have decided to rebaptize these algebras as \emph{quantum Euclidean spaces}, in consonance with quantum tori $\A_\Theta$ |also known in the literature as noncommutative tori or rotation algebras| which appear in turn as the subalgebra generated by $\lambda_\Theta(\xi)$ with $\xi$ running along $\Z^n$ or any other lattice of $\R^n$. Our main results in this paper about pseudodifferential operators hold for $\A_\Theta$ and $\RR_\Theta$.  

\begin{center} 
\textit{Calder\'on-Zygmund extrapolation}
\end{center}

In harmonic analysis, integral kernel representations play a central role to study the most relevant operators. In this particular form, pseudodifferential operators become well-behaved singular integrals, which admit a fruitful $L_p$-theory \cite{St}. A singular integral operator in a Riemannian manifold $(\mathrm{X},\mathrm{d}, \mu)$ admits the kernel representation $$T_kf(x) = \int_\mathrm{X} k(x,y) f(y) \, d\mu(y) \quad \mbox{for} \quad x \notin \mbox{supp} f.$$ Namely, $T_k$ is only assumed a priori to send test functions into distributions, so that it admits a distributional kernel in $\mathrm{X} \times \mathrm{X}$ which coincides in turn with a locally integrable function $k$ away from the diagonal $x=y$, where the kernel presents certain singularity. This already justifies the assumption $x \notin \mbox{supp} f$ in the kernel representation. The paradigm of singular integral theory is the Hilbert transform in $\R$, paramount to study the convergence of Fourier series and integrals. The challenge in higher dimensions required new real variable methods which culminated in the celebrated theorem of Calder\'on and Zygmund \cite{CZ}, who established sufficient conditions on a singular integral operator in $\R^n$ for its $L_p$-boundedness:
\begin{itemize}
\item[i)] Cancellation $$\big\| T_k: L_2(\R^n) \to L_2(\R^n) \big\| \le \mathrm{A}_1.$$

\item[ii)] Kernel smoothness $$\big| \nabla_x \, k (x,y) \big| + \big| \nabla_y \, k (x,y) \big| \, \le \, \frac{\mathrm{A}_2}{|x-y|^{n+1}}.$$
\end{itemize}
The same holds in Riemannian manifolds with nonnegative Ricci curvature \cite{B}.  

Noncommutative $L_p$ methods in harmonic analysis have gained a considerable momentum in recent years. The fast development of Fourier $L_p$ multiplier theory on group von Neumann algebras \cite{CPPR, GJP, JMP, JMP3, dLdlS, LdlS, PR} has been possible in part due to a deeper comprehension of the involved kernels. In spite of this, the validity of Calder\'on-Zygmund extrapolation principle over noncommutative manifolds is still widely open. Noncommutative martingale methods were used in \cite{Pa1} to establish endpoint estimates for singular  integrals over tensor product von Neumann algebras with an Euclidean factor, which have been the key for the recent solution in \cite{CPSZ} of the Nazarov-Peller conjecture. Other results in this direction include a CZ theory for group algebras over orthogonal crossed products $\R^n \rtimes \G$,  operator-valued kernels acting by left/right or Schur multiplication, other BMO spaces in a new approach towards nondoubling CZ theory, Littlewood-Paley estimates, H\"ormander-Mihlin multipliers or directional Hilbert transforms \cite{CMP,CP,HLMP,JMP,MP}. 

The Calder\'on-Zygmund theory presented below is the first form over a \lq fully noncommutative\rq${}$ von Neumann algebra. In other words, the singular integral acts on the whole algebra $\M$, not just over copies of $\R^n$ as tensor or crossed product factors in $\M$. A major challenge for such a von Neumann algebra $(\M,\tau)$ is to understand what it means to be a singular kernel. One has to identify the diagonal where the kernel singularity should be located, the quantum metric which measures the distance to it and its relation to the trace. A crucial point, undistinguishable in abelian algebras or the work cited so far, is to define kernels over $\M \bar\otimes \M_{\mathrm{op}}$ with the op-structure (reversed product law) in the second copy, see also \cite{JMP2}. In the case of $\RR_\Theta$, this is justified from the important map $$\pi_\Theta: L_\infty(\R^n) \to \mathcal{R}_\Theta \bar\otimes \mathcal{R}_\Theta^{\mathrm{op}},$$ $$\exp(2\pi i \langle \xi, \cdot \rangle) \longmapsto \lambda_\Theta(\xi) \otimes \lambda_\Theta(\xi)^*,$$ which extends to a normal $*$-homomorphism, for which the op-structure is strictly necessary. Note that $\pi_\Theta(f)(x,y) = f(x-y)$ for $\Theta=0$. In particular, if $| \cdot |$ stands for the Euclidean distance to $0$, the operator $$\mathrm{d}_\Theta = \pi_\Theta(| \cdot |)$$ is affiliated to $\mathcal{R}_\Theta \bar\otimes \mathcal{R}_\Theta^{\mathrm{op}}$ and implements the \emph{distance to the diagonal} as an unbounded operator. Similarly, the diagonal bands $\mathrm{b}_\Theta(\mathrm{R}) = \pi_\Theta(\chi_{|\cdot| \le \mathrm{R}})$ or smoothings of them will be indispensable to produce kernel truncations. An integral representation in $\RR_\Theta$ is formally given by $$T_k(\lambda_\Theta(f)) = (id \otimes \tau_\Theta) \big( k (\mathbf{1} \otimes \lambda_\Theta(f)) \big)$$ for some kernel $k$ affiliated to $\RR_\Theta \bar\otimes \RR_\Theta^{\mathrm{op}}$. We shall work with more general singular kernels which lead to $T_k \in \mathcal{L}(\S_\Theta, \S_\Theta')$, a map which sends the quantum Schwartz class $\S_\Theta = \lambda_\Theta(\S(\R^n))$ in $\RR_\Theta$ into its tempered distribution class $\S_\Theta'$. We shall also use the \lq free gradient\rq $$\nabla_{\hskip-1pt \Theta} = \sum_{j=1}^n s_j \otimes \partial_{\Theta}^j$$ associated to the partial derivatives $\partial_{\Theta}^j(\lambda_\Theta(\xi)) = 2 \pi i \xi_j \lambda_\Theta(\xi)$ and a free family $s_1, s_2, \ldots, s_n$ of semicircular random variables living in the free group algebra $\mathcal{L}(\mathbb{F}_n)$. 

\begin{TheoA} 
Let $T_k \in \mathcal{L}(\S_\Theta, \S_\Theta')$ and assume$\hskip1pt :$  
\begin{itemize}
\item[\emph{i)}] \emph{Cancellation} $$ \big\| T_k: L_2(\RR_\Theta) \to L_2(\RR_\Theta) \big\| \le \mathrm{A}_1.$$

\item[\emph{ii)}] \emph{Kernel smoothness} $$\Big| \mathrm{d}_\Theta^\alpha \bullet (\nabla_\Theta \otimes id) (k) \bullet \mathrm{d}_\Theta^\beta \Big| + \Big| \mathrm{d}_\Theta^{\alpha} \bullet (id \otimes \nabla_\Theta) (k) \bullet \mathrm{d}_\Theta^{\beta} \Big| \le \mathrm{A}_2,$$ for $(\alpha,\beta) = (n+1,0)$, $(\alpha,\beta) = (0,n+1)$ and $(\alpha,\beta) = (\frac{n+1}{2}, \frac{n+1}{2})$.
\end{itemize}
Then, $T_k: L_p(\RR_\Theta) \to L_p(\RR_\Theta)$ is completely bounded for every $1 < p < \infty$. 
\end{TheoA}

A more general statement is proved in Theorem \ref{CZExt4}. Our argument establishes $L_\infty \to \mathrm{BMO}$ endpoint estimates for a suitable noncommutative BMO. Interpolation with $L_p$ spaces is deduced in \cite{JM} from the theory of noncommutative martingales with continuous index set and a theory of Markov dilations. The convolution kernel case |in other words, quantum Fourier multipliers| is much easier to prove by transference methods stablished in \cite{CXY,Ri}. In the terminology of pseudodifferential operators, Fourier multipliers correspond to differential operators with constant coefficients. Of course, we aim to include nonconstant coefficients which leads to the analysis of the harder nonconvolution quantum kernels. Our statement above is very satisfactory and crucial for applications to pseudodifferential operator theory below. We shall also use other methods to justify that every CZ operator differs from its principal value by a left/right pointwise multiplier. This is fundamental in classical CZ theory and therefore of independent interest. 

\begin{center} 
\textit{The $L_p$ pseudodifferential calculus}
\end{center}

\renewcommand{\theequation}{$S_{\rho,\delta}^m$}
\addtocounter{equation}{-1}

The theory of pseudodifferential operators goes back to the mid 1960s with the work of Kohn, Nirenberg and H\"ormander. The basic idea is to exploit properties of the Fourier transform to produce a suitable representation $\Psi_L$ of partial differential operators $L = \sum_{|\alpha| \le m} a_\alpha(x) \partial_x^\alpha$ which can be inverted up to a controllable error term. This representation looks like $$\Psi_af(x) = \int_{\R^n} a(x,\xi) \widehat{f}(\xi) e^{2\pi i \langle x,\xi \rangle} \, d\xi$$ for a smooth symbol $a: \R^n \times \R^n \to \C$ satisfying 
\begin{equation}
\big| \partial_x^\beta \partial_\xi^\alpha a(x,\xi) \big| \le C_{\alpha\beta} \big( 1 + |\xi| \big)^{m - \rho |\alpha| + \delta |\beta|} \quad \mbox{for all} \quad \alpha, \beta \in \Z_+^n, 
\end{equation}
some $m \in \R$ and some $0 \le \delta \le \rho \le 1$. The realization of $\Psi_a$ as singular integral is given by partial Fourier inversion $k(x,y) = (id \otimes \mathcal{F}^{-1})(a)(x,x-y)$, which opens a door to CZ theory for Sobolev $p$-estimates of parametrices and error terms.  

\numberwithin{equation}{section}

In the noncommutative setting, this line took off in 1980 with Connes' work on pseudodifferential calculus for $\mathrm{C}^*$-dynamical systems \cite{Connes2}, originally conceived to extend the Atiyah-Singer index theorem for Lie group actions on $\mathrm{C}^*$-algebras, see also \cite{Ba1,Ba2,LM} for related results. Other applications in the context of quantum tori include a well-established elliptic operator theory  \cite{Connes3}, the Gauss-Bonnet theorem for 2D quantum tori \cite{CT,FK} and recent results on the local differential geometry of non-flat noncommutative tori \cite{BM, CM}. Unfortunately, the work of Connes and his collaborators does not include $L_p$ estimates for parametrices and error terms, which are paramount in harmonic analysis and partial differential equations. On the other hand, the only approach \cite{CXY,Ri,XXX} to harmonic $L_p$-analysis in quantum tori does not include pseudodifferential calculus, which requires Calder\'on-Zygmund  estimates in $\RR_\Theta$. In comparison with Connes' work |which focuses on the smallest H\"ormander class $S_{\rho,\delta}^m$ with $(\rho,\delta) = (1,0)$| our main contributions in this direction include all classes of symbols and $L_p$-estimates: 
\begin{itemize}
\item[\textbf{i)}] An $L_2$-theory for exotic and forbidden symbols $0 < \delta = \rho \le 1$.

\item[\textbf{ii)}] An $L_p$-theory for arbitrary H\"ormander classes and $1 < p < \infty$.
\end{itemize}
We refer to \cite{St,T1,T2} for the applications of these results in the Euclidean context. 

Pseudodifferential operators over quantum Euclidean spaces are easier to define than Calder\'on-Zygmund operators. The symbol $a(x,\xi)$ is now understood as an smooth function $a: \R^n \to \RR_\Theta$ since $\xi$ is still (dual) Euclidean, while $x$ becomes its $\Theta$-deformed analog $x_\Theta = (x_{\Theta,1}, x_{\Theta,2}, \ldots, x_{\Theta,n})$ as introduced above. We shall deal in this paper with two quantum forms of the H\"ormander classes:
\begin{itemize}
\item We say that $a \in S_{\rho,\delta}^m(\RR_\Theta)$ when $$\big| \partial_\Theta^\beta \partial_\xi^\alpha a(\xi) \big| \le C_{\alpha\beta} \big( 1 + |\xi| \big)^{m - \rho |\alpha| + \delta |\beta|}.$$ This is probably the most natural definition that comes to mind. 

\vskip5pt

\item We say that $a \in \Sigma_{\rho,\delta}^m(\RR_\Theta)$ when $$\hskip10pt \big| \partial_\Theta^\beta \partial_{\Theta,\xi}^{\alpha_1} \partial_\xi^{\alpha_2} a(\xi) \big| \le C_{\alpha_1\alpha_2\beta} \big( 1 + |\xi| \big)^{m - \rho |\alpha_1 + \alpha_2| + \delta |\beta|}.$$ Here $\partial_{\Theta,\xi}$ is a $\Theta$-deformation of $\partial_\xi$ by $\partial_{\Theta}$'s. More precisely, we have 
\begin{eqnarray*}
\partial_{\Theta, \xi}^j a (\xi) & = & \partial_\xi^j a(\xi) + 2 \pi i \big[ x_{\Theta,j}, a(\xi) \big] \\
& = & \partial_\xi^j a(\xi) + \frac{1}{2\pi i} \sum_{k=1}^n \Theta_{jk} \hskip1pt \partial_\Theta^k a(\xi) \\
& = & \lambda_\Theta(\xi)^\ast \partial_\xi^j \big\{ \lambda_\Theta(\xi) a(\xi) \lambda_\Theta(\xi)^\ast \big\} \lambda_\Theta(\xi).
\end{eqnarray*}
\end{itemize}
We clearly have $\Sigma_{\rho,\delta}^m(\RR_\Theta) \subset S_{\rho,\delta}^m(\RR_\Theta)$. It is very important to recall that both classes collapse into H\"ormander classical set of symbols $S_{\rho,\delta}^m$ when $\Theta=0$, so that both definitions above are a priori valid to generalize the Euclidean theory. It turns out that the $L_2$-theory holds for $S_{\rho,\delta}^m(\RR_\Theta)$, while the more involved class $\Sigma_{\rho,\delta}^m(\RR_\Theta)$ makes the $L_p$-theory valid. The reason has to do with the link to CZ theory and the two-sided nature of our Calder\'on-Zygmund conditions. Indeed, in all our past experiences with noncommutative Calder\'on-Zygmund theory certain amount of modularity is required. In this case, the bilateral form of our kernel conditions in Theorem A ultimately imposes the mixed quantum-classical derivatives $\partial_{\Theta,\xi}$. 

\noindent The pseudodifferential operator associated to $a: \R^n \to \RR_\Theta$ has the form
\begin{eqnarray*}
\Psi_a(\lambda_\Theta(f)) & = & \int_{\R^n} a(\xi) f(\xi) \lambda_\Theta(\xi) \, d\xi \\ & = & \big(id \otimes \tau_\Theta \big) \Big[ \Big( \underbrace{\int_{\R^n} (a(\xi) \otimes \1) (\lambda_\Theta(\xi) \otimes \lambda_\Theta(\xi)^*) d\xi}_{\mathrm{The \ kernel \ } k} \Big) \big( \1 \otimes \lambda_\Theta(f) \big) \Big].
\end{eqnarray*}
The algebra of pseudodifferential operators is formally generated by the derivatives $\partial_{\Theta}^j$ and the left multiplication maps $\lambda_\Theta(f) \mapsto x_{\Theta,j} \lambda_\Theta(f)$. The kernels affiliated to $\pi_\Theta(L_\infty(\R^n)) \subset \RR_\Theta \bar\otimes \RR_\Theta^{\mathrm{op}}$ implement Fourier multipliers $\lambda_\Theta(\xi) \mapsto m(\xi) \lambda_\Theta(\xi)$ in this setting, which correspond to the closure of pseudodifferential operators $\sum_\alpha a_\alpha \partial^\alpha$ with constant coefficients $a_\alpha$. A very subtle transference method |which avoids properly supported symbols| is required to obtain adjoint and product formulae in Section \ref{CompAdj}. Our main $L_p$ results are collected in the following statement. 


\begin{TheoB}
Let $a: \R^n \to \RR_\Theta$ and $1 < p < \infty$$\hskip1pt :$
\begin{itemize}
\item[\emph{i)}] If $a \in S_{\rho,\rho}^0(\RR_\Theta)$ with $0 \le \rho < 1$, $\Psi_a: L_2(\RR_\Theta) \to L_2(\RR_\Theta)$.

\vskip2pt

\item[\emph{ii)}] If $a \in S_{1,1}^0 \hskip1pt (\RR_\Theta) \cap S_{1,1}^0 \hskip1pt (\RR_\Theta)^*$, then $\Psi_a: L_2(\RR_\Theta) \to L_2(\RR_\Theta)$. 

\vskip2pt

\item[\emph{iii)}] If $a \in \Sigma_{1,1}^0(\RR_\Theta) \cap \Sigma_{1,1}^0(\RR_\Theta)^*$, then $\Psi_a: L_p(\RR_\Theta) \to L_p(\RR_\Theta)$. 
\end{itemize}
\end{TheoB}

\vskip3pt

Using $S_{\rho,\delta}^m(\RR_\Theta) \subset S_{\delta,\delta}^m(\RR_\Theta) \cap S_{\rho,\rho}^m(\RR_\Theta)$ for $0 \le \delta \le \rho \le 1$ |same inclusions for $\Sigma$-classes| we get $L_p$-estimates for regular, exotic and forbidden symbols in the expected ranks and Theorem B opens the core of the pseudodifferential $L_p$-calculus \cite{St,T1} to the context of quantum Euclidean spaces: 
\begin{itemize}
\item[$\bullet$] Theorem B i). Calder\'on-Vaillancourt theorem \cite{CV} on $L_2$-boundedness for exotic symbols quickly obtained a spectacular application \cite{BF} for $\rho=1/2$. Our proof of its quantum form for $\rho=0$ requires a careful approach due to the presence of a $\Theta$-phase. The case $\rho > 0$ also imposes an unexpected dilation argument among different deformed algebras $\RR_\Theta$.  

\vskip3pt

\item[$\bullet$] Theorem B ii). Bourdaud's theorem \cite{Bo} yields a form of the $T(1)$-theorem for pseudodifferential operators when $\rho = \delta = 1$: $\Psi_a$ is $L_2$-bounded iff the symbol $a_\dag^*$ of $\Psi_a^*$ remains in the same H\"ormander class. Our proof follows the classical one by showing that $\Psi_a$ is bounded in the Sobolev space $\mathrm{W}_{2,s}(\RR_\Theta)$ under a minimal amount of regularity $s>0$.

\vskip3pt

\item[$\bullet$] Theorem B iii). Our $L_p$-results follow by showing that any such symbol is a Calder\'on-Zygmund operator which fulfills all the hypotheses of Theorem A, the $L_2$-boundedness being assured by Theorem B ii). It is our CZ kernel condition what imposes the mixed quantum-classical derivatives $\partial_{\Theta,\xi}$ and the corresponding \lq\lq forbidden\rq\rq${}$ H\"ormander symbol classes $\Sigma_{\rho,\delta}^m(\RR_\Theta)$.

\vskip3pt

\item[$\bullet$] Related estimates. Our $L_p$-inequalities give rise to Sobolev $p$-estimates for symbols of arbitrary order $m$, we shall recollect these estimates in the body of the paper. On the other hand, the $L_p$-theory for symbols with $\rho < 1$ requires a negative degree to compensate lack of regularity. Fefferman proved in \cite{Fe} the $L_p$-bounds for the critical index $m = - (1-\rho) \frac{n}{2}$. We shall obtain nonoptimal $L_p$-estimates of this kind in $\RR_\Theta$. Interpolation yields even finer results for intermediate values of $m$.   
\end{itemize}


The analogue of Theorem B for quantum tori $\A_\Theta$ is proved by transference in Appendix A. The H\"ormander classes $S_{\rho,\delta}^m(\A_\Theta)$ and $\Sigma_{\rho,\delta}^m(\A_\Theta)$ involve discrete derivations over $\Z^n$ in the dual variable. In the line of Connes definition, we could also proceed by restriction to $\Z^n$ of symbols $\R^n \to \A_\Theta \subset \RR_\Theta$ in the corresponding H\"ormander classes. As in $\mathbb{T}^n$ both definitions turn out to be equivalent and this will be the source of our transference approach. The discrete form of difference operators has the advantage of being easier to be calculated with computers. 

\vskip5pt

\begin{center} 
\textit{An illustration for elliptic PDEs}
\end{center}

Pseudodifferential operators are a very powerful tool for linear and nonlinear partial differential equations \cite{T1,T2}. The existence, uniqueness and qualitative behavior of solutions for many PDEs are frequently understood by application of these methods. After the announced results so far, the potential applications for PDEs over quantum Euclidean spaces and tori are vast and beyond the scope of this paper. As a small but basic illustration, we prove in Theorem \ref{Lp-Regularity} the $L_p$ regularity for solutions of elliptic PDEs over quantum Euclidean spaces. We do not include this statement in the Introduction to avoid more terminology at this point. A profound analysis of partial differential equations over quantum spaces |$\A_\Theta$, $\RR_\Theta$ or even more general noncommutative manifolds| constitutes a long term program with conceivable implications for the geometry of such objects. 

\section{{\bf Quantum Euclidean spaces}} \label{Prelim}

Given an integer $n \ge 1$, fix an anti-symmetric $\R$-valued $n \times n$ matrix $\Theta$. Let us write $A_n(\R)$ to denote this class of matrices. We define $\mathrm{A}_\Theta$ as the universal $\mathrm{C}^*$-algebra generated by a family $u_1(s), u_2(s), \ldots, u_n(s)$ of one-parameter unitary groups in $s \in \R^n$ which are strongly continuous and satisfy the $\Theta$-commutation relations below $$u_j(s) u_k(t) = e^{2\pi i \Theta_{jk} st} u_k(t) u_j(s).$$  
If $\Theta=0$ and by Stone's theorem we may take $u_j(s) = \exp(2\pi i s \langle e_j, \cdot \rangle)$ and $\mathrm{A}_\Theta$ is the space of bounded continuous functions $\R^n \to \C$. Moreover, since $\R^n$ is amenable $\mathrm{A}_\Theta$ may be described as an spatial crossed product $\C \rtimes_\Theta \R^n \subset \mathcal{B}(L_2(\R^n))$ twisted by the $2$-cocycle determined by $\Theta$, as introduced by Zeller-Meier \cite{Z}. In general given $\xi \in \R^n$ we shall extensively use the unitaries $\lambda_\Theta(\xi) = u_1(\xi_1) u_2(\xi_2) \cdots u_n(\xi_n)$ and we define $\mathrm{E}_\Theta$ as the closure in $\mathrm{A}_\Theta$ of $\lambda_\Theta(L_1(\R^n))$ with $$\lambda_\Theta(f) = \int_{\R^n} f(\xi) \lambda_\Theta(\xi) \, d\xi.$$ If $\Theta=0$, we find $\mathrm{E}_\Theta = \mathcal{C}_0(\R^n)$. Given any $\Theta$, we easily see that
\begin{itemize}
\item $\lambda_\Theta(\xi)^* = e^{2\pi i \sum_{j > k} \Theta_{jk} \xi_j \xi_k} \lambda_\Theta(-\xi)$,

\vskip2pt

\item $\lambda_\Theta(\xi) \lambda_\Theta(\eta) = e^{2\pi i \langle \xi, \Theta \eta \rangle} \lambda_\Theta(\eta) \lambda_\Theta(\xi)$,

\vskip2pt

\item $\lambda_\Theta(\xi) \lambda_\Theta(\eta) = e^{2\pi i \sum_{j > k} \Theta_{jk} \xi_j \eta_k} \lambda_\Theta(\xi+\eta)$,

\vskip3pt

\item $\lambda_\Theta(f_1) \lambda_\Theta(f_2) = \lambda_\Theta(f_1 *_\Theta f_2)$ with $\Theta$-convolution given by $$\hskip20pt f_1 *_\Theta f_2(\xi) = \int_{\R^n} f_1(\xi - \eta) f_2(\eta) e^{2\pi i \sum_{j > k} \Theta_{jk} (\xi_j - \eta_j) \eta_k} \, d\eta.$$
\end{itemize}
Note that $\sum_{j > k} \Theta_{jk} \xi_j \eta_k = \langle \xi, \Theta_{\downarrow} \eta \rangle$ for the lower triangular truncation $\Theta_\downarrow$ of $\Theta$.

\subsection{Crossed product form} 

Let $$\tau_\Theta(\lambda_\Theta(f)) = \tau_\Theta \Big( \int_{\R^n} f(\xi) \lambda_\Theta(\xi)\, d\xi\Big) = f(0)$$ for $f: \R^n \to \C$ smooth and integrable. As we shall see, $\tau_\Theta$ extends to a normal faithful semifinite trace on $\mathrm{E}_\Theta$. Let $\mathcal{R}_\Theta = \mathrm{A}_\Theta'' = \mathrm{E}_\Theta''$ be the von Neumann algebra generated by $\mathrm{E}_\Theta$ in the GNS representation of $\tau_\Theta$. We obtain $\RR_\Theta = L_\infty(\R^n)$ for $\Theta = 0$. In general, we call the \emph{$\Theta$-deformation} $\mathcal{R}_\Theta$ a quantum Euclidean space. 

\subsubsection{Crossed products and trace} A $\mathrm{C}^*$-dynamical system is a triple formed by a $\mathrm{C}^*$-algebra $\mathrm{A}$, a locally compact group $\G$ and a continuous action $\beta: \G \to \mathrm{Aut}(\mathrm{A})$ by $*$-automorphisms. The reduced crossed product $\mathrm{A} \rtimes_{\beta,\mathrm{red}} \G$ is the norm closure in $\mathrm{A} \bar\otimes \mathcal{B}(L_2(\G))$ of the $*$-algebra generated by the representations $\rho: \mathrm{A} \to L_\infty(\G; \mathrm{A})$ and $\lambda: \G \to \mathcal{U}(L_2(\G))$, given by $$\rho(a)(g) = \beta_{g^{-1}}(a) \quad \mbox{and} \quad [\lambda(g)f] (h) = f(g^{-1}h).$$ The full crossed product $\mathrm{A} \rtimes_{\beta,\mathrm{full}} \G$ is the $\mathrm{C}^*$-algebra generated by all covariant representations $\gamma: \mathrm{A} \to \mathcal{B}(\mathcal{H})$ and $u: \G \to \mathcal{U}(\mathcal{H})$ over some Hilbert space $\mathcal{H}$: $u(g) \gamma(a) u(g)^* = \gamma (\beta_g(a))$. Given $f: \G \to \mathrm{A}$ continuous and integrable $$\Big\| \int_\G f_g \rtimes g \, d\mu(g) \Big\|_{\mathrm{A} \rtimes_{\beta, \mathrm{full}} \G} \, = \, \sup_{\begin{subarray}{c} \gamma, \, u \\ \mathrm{covariant} \end{subarray}} \Big\| \int_\G \gamma(f_g) u(g) \, d\mu(g) \Big\|_{\mathcal{B}(\mathcal{H})}.$$ It is a very well-known result \cite{BO} that $\mathrm{A} \rtimes_{\beta,\mathrm{full}} \G = \mathrm{A} \rtimes_{\beta,\mathrm{red}} \G$ when $\G$ is amenable. 

Given a pair $(\mathcal{M},\tau)$ formed by a von Neumann algebra $\M$ equipped with a normal faithful semifinite trace $\tau$ ---noncommutative measure space--- and a locally compact unimodular group $\G$ acting on $(\mathcal{M},\tau)$ by trace preserving automorphisms $\beta: \G \to \mathrm{Aut}(\M,\tau)$, the crossed product von Neumann algebra $\M \rtimes_\beta \G$ is the von Neumann subalgebra of $\mathcal{M} \bar\otimes \mathcal{B}(L_2(\G))$ generated by $\rho(\mathcal{M})$ and $\lambda(\G)$, defined as above. In other words, $\mathcal{M} \rtimes_\beta \G$ is the weak-$*$ closure of $\M \rtimes_{\beta, \mathrm{red}} \G$. 

Given $f: \G \to \M$ continuous and integrable, set $$\tau_\rtimes \Big( \int_\G f_g \rtimes \lambda(g) \, d\mu(g) \Big) = \tau(f_e)$$ where $\mu$ and $e$ stand for the Haar measure and the identity in the unimodular group $\G$. This determines a normal faithful semifinite trace which extends to the crossed product von Neumann algebra $\M \rtimes_\beta \G$, see Takesaki \cite{Ta}. In the following result, we provide an iterated crossed product characterization of quantum Euclidean spaces and construct a normal faithful semifinite trace on them.  

\begin{proposition} \label{Trace}
The following results hold$\hskip1pt :$
\begin{itemize}
\item[\emph{i)}] If $n=2$ and $\Theta \neq 0$, we have $$\hskip20pt \mathrm{E}_\Theta \hskip0.5pt \simeq \mathcal{C}_0(\R) \rtimes \R.$$ In this case, the crossed product action is given by $\R$-translations.

\vskip3pt 

\item[\emph{ii)}] If $n \ge 2$,  let us define $$\hskip20pt \tau_\Theta(\lambda_\Theta(f)) = \tau_\Theta \Big( \int_{\R^n} f(\xi) \lambda_\Theta(\xi)\, d\xi\Big) = f(0)$$ for $f: \R^n \to \C$ smooth and integrable. Then, $\tau_\Theta$ extends to a normal faithful semifinite trace on $\mathrm{E}_\Theta$. Moreover, let $\Xi$ denote the $(n-1) \times (n-1)$ upper left corner of $\Theta \in A_n(\R)$. Then there exists a continuous group action $\beta_{n-1}: \R \to \mathrm{Aut}(\mathrm{E}_\Xi)$ satisfying $$\hskip20pt \mathrm{E}_\Theta \simeq \mathrm{E}_\Xi \rtimes_{\beta_{n-1}} \R.$$

\vskip3pt 

\item[\emph{iii)}] Let $\mathcal{R}_\Theta = \mathrm{E}_\Theta''$ be the von Neumann algebra generated by $\mathrm{E}_\Theta$ in the \emph{GNS} representation determined by $\tau_\Theta$. We have $\RR_\Theta \simeq L_\infty(\R) \rtimes \R \simeq \mathcal{B}(L_2(\R))$ when $n=2$ and $\Theta \neq 0$, with $\rtimes$-action given by $\R$-translations. Moreover $\tau_\Theta$ extends to a n.f.s. trace on $\RR_\Theta$, and the action $\beta_{n-1}$ is trace preserving on $(\RR_\Xi,\tau_\Xi)$. Induction on $n$ and iteration give $$\hskip20pt \begin{array}{c} \RR_\Theta \simeq \RR_\Xi \rtimes_{\beta_{n-1}} \R, \\ \mathcal{R}_\Theta \simeq \Big( \big( L_\infty(\R) \rtimes_{\beta_1} \R \big) \cdots \rtimes_{\beta_{n-1}} \R \Big). \end{array}$$ 
\end{itemize}
\end{proposition}

\dem Given $\Theta \in A_n(\R)$ and a Hilbert space $\H_\pi$, every set of one-parameter unitary groups $\{u_{\pi j}(s) : 1 \le j \le n, s \in \R\}$ in $\mathcal{B}(\H_\pi)$ satisfying the $\Theta$-relations yields a $*$-representation $\pi: \mathrm{E}_\Theta \to \mathcal{B}(\H_\pi)$. In order to work with a concrete set of unitary groups $u_j(s)$, we consider the universal representation in the direct sum $\bigoplus_\pi \H_\pi$ over the set of all cyclic representations $\pi$ as above. We shall use in what follows |with no further reference| that $\mathcal{C}_0(\R)$ is the closure of $\mathcal{F}(L_1(\R))$, which can also be understood replacing the characters $\exp_s = \exp(2 \pi \langle s, \cdot \rangle)$ in the Fourier transform $\mathcal{F}$ by $u_j(s)$ for any fixed $1 \le j \le n$, since $\{ u_j(s): s \in \R \}$ forms a non-trivial one-parameter group of unitaries.

i) If $n=2$ and $\Theta \neq 0$, there must exist $\delta \neq 0$ with $\Theta = \delta (e_{12} - e_{21})$. We may rescale $u_1(s), u_2(t)$ and assume without loss of generality that $\delta=1$. Now, consider the map $$\mathrm{E}_\Theta \ni z = \int_{\R^2} z(s,t) u_1(s)u_2(t) \, dsdt \mapsto \int_\R f_t \rtimes t \, dt = f \in \mathcal{C}_0(\R) \rtimes \R$$ with $f_t \in \mathcal{C}_0(\R)$ given by $$f_t = \int_\R z(s,t) e^{2\pi i s \cdot} \, ds \simeq \int_\R z(s,t) u_1(s) \, ds.$$ If $\H = \oplus_\pi \H_\pi$, define $u: \R \to \mathcal{U}(\mathcal{H})$ and $\gamma: \mathcal{C}_0(\R) \to \mathcal{B}(\mathcal{H})$ by $$u(t) = u_2(t) \quad \mbox{and} \quad \gamma \Big( \int_\R \widehat{a}(s) e^{2\pi i s \cdot} \, ds \Big) = \int_\R \widehat{a}(s) u_1(s) \, ds.$$ The pair $(\gamma,u)$ forms a covariant representation since we have 
\begin{eqnarray*}
u(t) \gamma(a) u(t)^* & = & \int_\R \widehat{a}(s) u_2(t) u_1(s) u_2(-t) \, ds \\ & = & \int_\R \widehat{a}(s) e^{-2\pi i s t} u_1(s) \, ds \ = \ \gamma \Big( \int_\R \widehat{a}(s) e^{2\pi i s (\cdot - t)} \, ds \Big) \ = \ \gamma(\beta_t(a))
\end{eqnarray*}
where $\beta_t(a) = \lambda(t)[a]$ is the left regular representation at $t$ acting on $a$. This gives $$\|f\|_{\mathcal{C}_0(\R) \rtimes \R} \ge \Big\| \int_\R \gamma(f_t) u(t) \, dt \Big\|_{\mathcal{B}(\H)} = \Big\| \int_{\R^2} z(s,t) u_1(s) u_2(t) \, dsdt \Big\|_{\mathcal{B}(\H)} = \|z\|_{\mathrm{E}_\Theta}.$$ The reverse inequality is proved similarly. Indeed, let us consider the following map $$\mathcal{C}_0(\R) \rtimes \R \ni f = \int_\R f_t \rtimes t \, dt \mapsto \int_{\R^2} z(s,t) u_1(s) u_2(t) \, dsdt = z \in \mathrm{E}_\Theta$$ with $z(s,t) = \widehat{f}_t(s)$. Fix a Hilbert space $\mathcal{K}_{\gamma u}$ and a covariant representation $(\gamma,u)$ of the pair $(\mathcal{C}_0(\R),\R)$ in $\mathcal{B}(\mathcal{K}_{\gamma u})$. Define $w_1(s) = \gamma(e^{2\pi i s \cdot})$ and $w_2(t) = u(t)$. This shows that covariant representations of $(\mathcal{C}_0(\R),\R)$ with action given by translations are in one-to-one correspondence with $*$-representations of $\mathrm{E}_\Theta$ for the deformation $\Theta = e_{12} -e_{21}$. Indeed, $w_1(s)$ and $w_2(t)$ are one-parameter groups of unitaries since $\gamma$ is a $*$-representation and $u$ a unitary representation. Moreover, the commutation relations hold as a consequence of the covariant property $$w_1(s)w_2(t) = \gamma(\exp_s) u(t) = u(t) \gamma(\beta_{-t}(\exp_s)) = e^{2\pi i st} w_2(t) w_1(s)$$ for $\exp(s) = \exp(2 \pi i s \hskip1pt \cdot)$. In particular $$\|z\|_{\mathrm{E}_\Theta} \ge \Big\| \int_{\R^2} z(s,t) w_1(s) w_2(t) \, ds dt \Big\|_{\mathcal{B}(\mathcal{K}_{\gamma u})} = \Big\| \int_\R \gamma(f_t) u(t) \, dt \Big\|_{\mathcal{B}(\mathcal{K}_{\gamma u})}.$$ Taking the supremum over $(\gamma,u)$ covariant, we see that $\|z\|_{\mathrm{E}_\Theta} \ge \|f\|_{\mathcal{C}_0(\R) \rtimes \R}$. 

ii) When $n \ge 2$ we proceed by induction. To prove ii) for $n=2$, it suffices from i) to justify that $\tau_\Theta$ extends to a n.s.f. trace on $\mathrm{E}_\Theta$. Note that $\mathcal{C}_0(\R) \rtimes \R$ is generated by $\exp_\eta \rtimes \lambda(\zeta)$ for $(\eta,\zeta) \in \R \times \R$ where $\exp_\eta(x) = \exp (2 \pi i x \eta)$ and $\lambda(\zeta)f(x) = f(x-\zeta)$. According to i), this gives $$\lambda_\Theta(f) = \int_{\R^2} f(\xi) \lambda_\Theta(\xi) \, d\xi = \int_\R \Big[ \underbrace{\int_\R f(\eta,\zeta) \exp_\eta \, d\eta}_{\varphi_\zeta} \Big] \rtimes \lambda(\zeta) \, d\zeta.$$ This means that the crossed product trace $$\tau_\rtimes \Big( \int_\R \varphi_\zeta \rtimes \lambda(\zeta) \, d\zeta \Big) = \int_\R \varphi_0(x) \, dx = \int_\R \Big[ \int_\R f(\eta,0) \exp_\eta(x) \, d\eta \Big] dx = f(0)$$ coincides with $\tau_\Theta$ in $\mathrm{E}_\Theta$. Since $\mathcal{C}_0(\R) \rtimes \R$ embeds faithfully in $L_\infty(\R) \rtimes \R$, it turns out that $\tau_\Theta$ is a n.f.s. trace and completes the argument for the case $n=2$. Once this is settled, consider $\Theta \in A_n(\R)$ whose upper left $(n-1) \times (n-1)$ corner is denoted by $\Xi$. Assume ii) holds for any dimension smaller than $n$, and set $$\beta_{n-1}(s) \Big( \int_{\R^{n-1}} \varphi(z) \lambda_\Xi(z) \, dz \Big) = \int_{\R^{n-1}} \varphi(z) e^{-2\pi i \sum_{j < n} \Theta_{jn} z_js} \lambda_\Xi(z) \, dz.$$ Then, $\beta_{n-1}$ trivially yields a $\tau_\Xi$-preserving action on $(\mathrm{E}_\Xi, \tau_\Xi)$. Moreover, the map $\lambda_\Theta(\xi) \mapsto \lambda_\Xi(\xi_1, \xi_2, \ldots, \xi_{n-1}) \rtimes \lambda(\xi_n)$ also gives rise to $\mathrm{E}_\Theta \simeq \mathrm{E}_\Xi \rtimes_{\beta_{n-1}} \R$ and $\tau_\Theta = {\tau_{\rtimes}}{\mid_{\mathrm{E}_\Theta}}$ by arguing as above for $n=2$, details are left to the reader. 

iii) Now, for $n=2$ and $\Theta \neq 0$ we get $$\RR_\Theta = \mathrm{E}_\Theta'' = \big( \mathcal{C}_0(\R) \rtimes \R \big)'' = \mathcal{C}_0(\R)'' \rtimes \R = L_\infty(\R) \rtimes \R,$$ $\tau_\Theta = \tau_\rtimes$ on $\RR_\Theta$ and $\mathrm{E}_\Theta$ sits faithfully in $\RR_\Theta$. Moreover, $L_\infty(\R) \rtimes \R \subset \mathcal{B}(L_2(\R))$ acts on $L_2(\R)$ by modulation and translation, which implies $\RR_\Theta \simeq \mathcal{B}(L_2(\R))$ since only constant multiples of the identity map commute with all modulations and translations. When $n > 2$ we proceed by induction one more time to conclude that $\beta_{n-1}$ is $\tau_\Xi$-preserving, $\RR_\Theta \simeq \RR_\Xi \rtimes_{\beta_{n-1}} \R$, $\tau_\Theta = \tau_\rtimes$ and $\mathrm{E}_\Theta \subset \RR_\Theta$ faithfully. The last assertion follows trivially by iteration. This completes the proof. \fin 

\begin{remark} \label{Plancherel}
\emph{The map $$\lambda_\Theta: L_2(\R^n) \to L_2(\RR_\Theta,\tau_\Theta)$$ is an isometric isomorphism, extending Plancherel theorem for $\Theta=0$. Indeed, once we know $\tau_\Theta$ is a trace, it follows from the density of the quantum Schwartz class $\S_\Theta = \lambda_\Theta(\S(\R^n))$ in $L_2(\RR_\Theta)$ and the identity $\lambda_\Theta(f_1) \lambda_\Theta(f_2) = \lambda_\Theta(f_1 *_\Theta f_2)$.}
\end{remark}

\begin{remark}
\emph{When $n=1$, $\RR_\Theta = L_\infty(\R)$ generated by $u(s) = \exp(2 \pi i s \, \cdot)$. In the $2$D case, we find one more time $\RR_\Theta = L_\infty(\R^2)$ for $\Theta=0$. Otherwise, there exists $\delta \neq 0$ such that $\Theta = \delta (e_{12} - e_{21})$. Rescaling $\delta = 1$ and arguing as in Proposition \ref{Trace} iii) gives 
$$\Theta = \Big( \begin{array}{rl} 0 & 1 \\ - 1 & 0 \end{array} \Big) \ \Rightarrow \ \mathcal{R}_\Theta \simeq \mathcal{B}(L_2(\R)) \simeq L_\infty(\R) \rtimes \R$$ generated by modulations $\exp_\eta \rtimes \1$ and translations $\1 \rtimes \lambda(\zeta)$ for $\eta, \zeta \in \R$. These are the standard time/frequency unitaries in Fourier analysis. If we set $\Xi$ to be the $n \times n$ matrix with all its entries equal to 1, then the analogous space in dimension $2n$ is given by $\Theta = \Xi \otimes (e_{12} - e_{21})$ with $\RR_\Theta \simeq \mathcal{B}(L_2(\R^n)) \simeq L_\infty(\R^n) \rtimes \R^n$. The $3$D case admits other models. By Proposition \ref{Trace} iii)}
\begin{itemize}
\item[$\bullet$] If $\Theta=0$, then $\RR_\Theta = L_\infty(\R^3) \simeq L_\infty(\R) \bar\otimes L_\infty(\R) \bar\otimes L_\infty(\R)$,

\vskip12pt

\item[$\bullet$] If $\displaystyle \Theta = \left( \hskip4.3pt \begin{array}{rrl} 0 & 0 & 0 \\ 0 & 0 & \alpha \\ 0 & \hskip-7.6pt -\alpha & 0 \hskip1pt \end{array} \right) \ \Rightarrow \ \RR_\Theta \simeq L_\infty(\R) \bar\otimes \big( L_\infty(\R) \rtimes \R \big)$,

\vskip2pt

\item[$\bullet$] If $\displaystyle \Theta = \left( \hskip-4.5pt \begin{array}{rrl} 0 & 0 & \beta \\ 0 & 0 & \alpha \\ -\beta & \hskip-8pt -\alpha & 0 \hskip1pt \end{array} \right) \ \Rightarrow \ \RR_\Theta \simeq \big( L_\infty(\R) \bar\otimes L_\infty(\R) \big) \rtimes \R$,

\vskip2pt

\item[$\bullet$] If $\displaystyle \Theta = \left( \hskip-5pt  \begin{array}{rrl} 0 & \gamma & \beta \\ -\gamma & 0 & \alpha \\ \hskip-2pt -\beta & \hskip-8pt - \alpha & 0 \hskip2.5pt \end{array} \right) \ \Rightarrow \ \RR_\Theta \simeq \big( L_\infty(\R) \rtimes \R \big) \rtimes \R$,
\end{itemize}
\emph{for $\alpha, \beta, \gamma \neq 0$. Higher dimensions are treated similarly. When $\alpha \neq 0 = \beta = \gamma$, the $\rtimes$-action is $t \cdot f(s) = f(s - \alpha t)$. In the second case $\alpha, \beta \neq 0 = \gamma$, the $\rtimes$-action in $L_\infty(\R^2)$ is $t \cdot f(x,y) = f(x - \beta t, y - \alpha t)$. In the third case $\alpha, \beta, \gamma \neq 0$, both actions yield $t \cdot_\mathrm{ext} \big( (s \cdot_\mathrm{int} f) (r) \big) = \big( (s - \alpha t) \cdot_\mathrm{int} f \big) (r - \beta t) = f \big( r - \gamma s + (\alpha \gamma - \beta) t \big).$ In the particular case $\alpha = \beta = \gamma = 1$, we have full symmetry under the action of the permutation group. In this case, the time/frequency dichotomy described above for $n=2$ is replaced by three indistinguishable sets of unitaries.}
\end{remark} 

\subsubsection{The corepresentation map} \label{CorepSect}

Let us now recall a very useful consequence of the crossed product characterization of $\RR_\Theta$ given above, the normality (weak-$*$ continuity) of the corepresentation map $\sigma_\Theta: \lambda_\Theta(\xi) \mapsto \exp_\xi \otimes \lambda_\Theta(\xi)$, where $\exp_\xi$ stands for the character $x \mapsto \exp(2 \pi i \langle x, \xi \rangle)$ in $L_\infty(\R^n)$. This will be the source of several metric and differentiability considerations over quantum Euclidean spaces.

\begin{corollary} \label{Planc+Corep} 
$\sigma_\Theta: \RR_\Theta \to L_\infty(\R^n) \bar\otimes \RR_\Theta$ is a normal injective $*$-homomorphism.
\end{corollary}

\dem The assertion is a simple exercise in the $\mathrm{C}^*$-algebra level, so that we shall only justify normality. We proceed by induction on $n$, the case $n=1$ is nothing but comultiplication in $L_\infty(\R^n)$. In higher dimensions, $\sigma_\Theta$ factorizes as follows
\begin{equation*}
  \xymatrix{
    \RR_\Theta \ar@{-}[r]^{\simeq} \ar[d]^{\sigma_\Theta}      & \RR_\Xi \rtimes_{\beta_{n-1}} \R \ar[d]^{\widetilde{\sigma}_{\Xi} = \sigma_\Xi \rtimes id_\R}\\
    L_\infty(\R^n) \bar\otimes \RR_\Theta \ar@{-}[dr]^{\simeq} & L_\infty(\R^{n - 1}) \bar\otimes \RR_\Xi \rtimes_{\widehat{\beta}_{n - 1}} \R \ar[d]^{\Omega}\\
                                                               & L_\infty(\R) \bar\otimes \big( (L_\infty(\R^{n-1}) \bar\otimes \RR_\Xi) \rtimes_{\widehat{\beta}_{n-1}} \R \big)
  }
\end{equation*}

\noindent where $\widehat{\beta}_{n-1} = id_{\R^{n-1}} \otimes \beta_{n-1}$ and the map $\Omega$ is given by
\[
  \M \rtimes_\beta \R \ni \int_\R f_s \rtimes \lambda(s) \, ds
  \stackrel{\Omega}{\longrightarrow}
  \int_\R \exp_s \otimes \big( f_s \rtimes \lambda(s) \big) \, ds \in L_\infty(\R) \bar\otimes (\M \rtimes_\beta \R).
\]
By such factorization, it suffices to justify the normality of $\widetilde{\sigma}_\Xi$ and $\Omega$:
\begin{itemize}
\item The map $\sigma_\Xi$ is equivariant $$\hskip30pt \sigma_\Xi(\beta_{n-1}(s)(a)) = \widehat{\beta}_{n-1}(s) (\sigma_\Xi(a)).$$ Let $j = \rho \rtimes \lambda: \RR_\Xi \rtimes \R \to \RR_\Xi \bar\otimes \mathcal{B}(L_2(\R))$ be the natural injection. By the above equivariance, $j$ intertwines $\sigma_\Xi \rtimes id$ and $\sigma_\Xi \otimes id$ 
\[
  \hskip25pt \widetilde{\sigma}_\Xi = \sigma_\Xi \rtimes id = (id_{L_\infty(\R^{n-1})} \otimes j)^{-1} \circ (\sigma_\Xi \otimes id_{\mathcal{B}(L_2(\R))}) \circ j.
\]
Since $\sigma_\Xi \otimes id$ is normal by induction hypothesis, the same holds for $\widetilde{\sigma}_\Xi$.

\vskip5pt

\item The fundamental unitary on $\R^2$ $$\hskip30pt \mathrm{W}f(x,y) = f(x+y,y)$$ satisfies $\mathrm{W}^* (\1 \otimes \lambda(s)) \mathrm{W} = \lambda(s) \otimes \lambda(s)$. Using the isometric isomorphism $\Lambda: \mathcal{L}(\R) \ni \lambda(s) \mapsto \exp_s \in L_\infty(\R)$, we get 
\begin{eqnarray*}
\hskip20pt \Omega(f) & = & \int_\R \exp_s \otimes \big( f_s \rtimes \lambda(s) \big) \, ds \\ 
& = & (\Lambda \otimes id_{\M \rtimes \R}) (\1 \rtimes \mathrm{W}^*) \Big( \int_\R \1 \otimes \big( f_s \rtimes \lambda(s) \big) \, ds \Big) (\1 \rtimes \mathrm{W}).
\end{eqnarray*}
Thus $\Omega(f) = (\Lambda \otimes id_{\M \rtimes \R}) (\1 \rtimes \mathrm{W}^*) (\1 \otimes f) (\1 \rtimes \mathrm{W})$ and $\Omega$ is normal. \hfill $\square$
\end{itemize}
 
\subsection{Metrics and derivations}

In this paragraph, we exploit the corepresentation $\sigma_\Theta$ to introduce some other auxiliary operators which will help us to equip $\RR_\Theta$ with an induced metric, a natural $\mathrm{BMO}$ space and a differential structure. 

\subsubsection{A metric in $\RR_\Theta$ and $\mathrm{BMO}$} \label{MetricandBMO} 

Given a von Neumann algebra $\M$, its opposite algebra $\M_\mathrm{op}$ is obtained by preserving linear and adjoint structures but reversing the product $a_1 \cdot a_2 = a_2 a_1$. Several reasons justify why noncommutative singular integral operators require to understand the singular kernels as operators affiliated to $\M \bar\otimes \M_\mathrm{op}$, see \cite{JMP2} and Section \ref{SectKernels} below. We shall use from now on $\cdot$ for the $\M_\mathrm{op}$-product, as well as $\bullet$ for the product in $\M \bar\otimes \M_\mathrm{op}$, so that $$(a_1 \otimes a_2) \bullet (a_1' \otimes a_2') = (a_1 a_1') \otimes (a_2 \cdot a_2') = (a_1 a_1') \otimes (a_2' a_2).$$ Let us consider the linear map $\pi_\Theta$, determined by $$\exp_\xi \stackrel{\pi_\Theta}{\longmapsto} \lambda_\Theta(\xi) \otimes \lambda_\Theta(\xi)^*$$ where, as usual, we write $\exp_\xi$ for the Fourier characters $\exp(2\pi i \langle \xi, \cdot \rangle)$ in $\R^n$. As an illustration, recall that for $\Theta=0$ we may expect to get the following identity for any (say) Schwartz function $f: \R^n \to \C$ $$\pi_0(f)(x,y) = \pi_0 \Big( \int_{\R^n} \widehat{f}(\xi) \exp_\xi d\xi \Big) (x,y) = \int_{\R^n} \widehat{f}(\xi) \exp_\xi(x-y) \, d\xi = f(x-y).$$ Of course, this requires to justify the continuity properties of the map $\pi_\Theta$ which we shall do in Lemma \ref{QDistance} below. The quantum analogue of this map is particularly useful to identify the \emph{diagonal} in $\RR_\Theta \bar\otimes \RR_\Theta^{\mathrm{op}}$, where the kernel singularities of our operators are expected to live. Of particular relevance is the induced metric which we define by $$\mathrm{d}_\Theta = \pi_\Theta (| \cdot |)$$ for the Euclidean norm $|\cdot|$ or the bands around the diagonal $\mathrm{b}_\Theta(\mathrm{R}) = \pi_\Theta(\chi_{|\cdot| \le \mathrm{R}})$.

It is worth recalling that both $\sigma_\Theta$ and $\pi_\Theta$ take $L_2(\R^n)$ into $L_\infty(\R^n; L_2(\R^n))$ when $\Theta=0$. The quantum analogue for $\Theta \neq 0$ requires noncommutative forms of mixed-norm $L_\infty(L_2)$-spaces, whose construction we briefly recall. Given a Hilbert space $\H$ and $x = \sum_{j} m_j \otimes h_j \in \M \otimes_{\mathrm{alg}} \H$, we define 
\[
  \begin{array}{rcll>{\displaystyle}l}
    \| x \|_{\M \bar\otimes \H^r} & = & \big\| \langle x, x \rangle_r \big\|_{\M}^\frac12 & = & \Big\| \sum_{j,k} m_j m_k^\ast \langle h_j, h_k \rangle_\H \Big\|_{\M}^\frac12 \\
    \| x \|_{\M \bar\otimes \H^c} & = & \big\| \langle x, x \rangle_c \big\|_{\M}^\frac12 & = & \Big\| \sum_{j,k} m_j^\ast m_k \langle h_j, h_k \rangle_\H \Big\|_{\M}^\frac12.
  \end{array}
\]
Given $\dag \in \{r,c\}$, the space $\M \bar\otimes \H^\dag$ |also denoted by $\H^\dag \bar\otimes \M$ or $L_\infty(\M;\H^\dagger)$| is defined as the closure of $\M \otimes_{\mathrm{alg}} \H$ with respect to the weak topology generated by the functionals $$p_\omega(x) = \omega \big( \langle x, x \rangle_\dag^{\frac12} \big) \quad \mbox{for every $\omega \in \M_*$}.$$ Alternatively, $\M \bar\otimes \H^\dag$ is the weak-$\ast$ closed tensor product of the dual operator spaces $\M$ and $\H^{\dagger}$, the latter space representing the row or column operator space structure on $\H$. Indeed, if $\mathrm{X}$ and $\mathrm{Y}$ are dual operator spaces, there are completely isometric and weak-$\ast$ continuous injections $\pi_\mathrm{X}: \mathrm{X} \to \B(\H_\mathrm{X})$, $\pi_\mathrm{Y}: \mathrm{Y} \to \B(\H_\mathrm{Y})$ and we define $\mathrm{X} \bar\otimes \mathrm{Y}$ as $$\overline{\pi_\mathrm{X}[\mathrm{X}] \otimes_{\mathrm{alg}} \pi_\mathrm{Y}[\mathrm{Y}]^{\wast}} \subset \B(\H_\mathrm{X} \otimes_2 \H_\mathrm{X}).$$ It is well-known that such construction is representation-independent and when one of the tensor components is a von Neumann algebra, the predual is given by the projective tensor product $\mathrm{X}_* \widehat{\otimes} \mathrm{Y}_*$, see \cite{ERHopf,P3} for further details. Noncommutative mixed-norm spaces have also been studied in \cite{JP,P2}.

\begin{lemma} \label{QDistance}
$\pi_\Theta$ extends to a normal $*$-homomorphism $$\pi_\Theta: L_\infty(\R^n) \to \RR_\Theta \bar\otimes \RR_\Theta^{\mathrm{op}} \quad \mbox{satisfying} \quad (\sigma_\Theta \otimes id_{\RR_\Theta^{\mathrm{op}}}) \circ \pi_\Theta = (id_{\R^n} \otimes \pi_\Theta) \circ \Delta_{\R^n}$$ where $\Delta_{\R^n}(\exp_\xi) = \exp_\xi \otimes \exp_\xi$ is the comultiplication map in $\R^n$. This shows in particular that $\mathrm{d}_\Theta = \pi_\Theta(|\cdot|)$ is a well-defined operator affiliated to $\RR_\Theta \bar\otimes \RR_\Theta^\mathrm{op}$ as an increasing limit of the bounded operators $\mathrm{d}_\Theta(\mathrm{R}) = \pi_\Theta (\chi_{|\cdot| \le \mathrm{R}} |\cdot|)$. Moreover, the map $\pi_\Theta$ also extends to a complete isometry $\pi_\Theta: L_2^c(\R^n) \to L_2^c(\RR_\Theta) \bar\otimes \RR_\Theta^{\mathrm{op}}$.
\end{lemma}

\dem That $\pi_\Theta: \exp_\xi \mapsto \lambda_\Theta(\xi) \otimes \lambda_\Theta(\xi)^*$ extends to a $*$-homomorphism is a simple consequence of the product in $\RR_\Theta \bar\otimes \RR_\Theta^{\mathrm{op}}$, details are left to the reader. Let us then prove that $\pi_\Theta$ is weak-$*$ continuous. It is tedious but straightforward to check that $$\Lambda_\Theta: h \mapsto \int_{\R^n \times \R^n} h(\xi_1, \xi_2) \lambda_\Theta(\xi_1) \lambda_\Theta(\xi_2) \otimes \lambda_\Theta(\xi_1)^* \, d\xi_1 d\xi_2$$ yields an isometry $L_2(\R^n) \otimes_2 L_2(\R^n) \to L_2(\RR_\Theta) \otimes_2 L_2(\RR_\Theta^\mathrm{op})$. Indeed, by density it suffices to expand $\tau_\Theta \otimes \tau_\Theta ( \Lambda_\Theta(h)^* \Lambda_\Theta(h))$ for $h$ smooth, then calculate the trace applying twice the simple identity $\tau_\Theta (\lambda_\Theta(f) \lambda_\Theta(\xi)^*) = f(\xi)$ for a smooth integrable function $f$ in $\R^n$. Moreover, given any $z = \Lambda_\Theta(h) \in L_2(\RR_\Theta) \otimes_2 L_2(\RR_\Theta^\mathrm{op})$ it turns out that 
\begin{eqnarray*}
\pi_\Theta(\exp_\xi) (z) & = & \int_{\R^n \times \R^n} h(\xi_1,\xi_2) \, \lambda_\Theta(\xi) \lambda_\Theta(\xi_1) \lambda_\Theta(\xi_2) \otimes \lambda_\Theta(\xi)^* \hskip-2pt \cdot \hskip-2pt  \lambda_\Theta(\xi_1)^* \, d\xi_1 d\xi_2 \\ & = & \int_{\R^n \times \R^n} h(\xi_1,\xi_2) \, \lambda_\Theta(\xi_1+\xi) \lambda_\Theta(\xi_2) \otimes \lambda_\Theta(\xi_1+\xi)^* \, d\xi_1 d\xi_2 \\ & = & \int_{\R^n \times \R^n} h(\xi_1-\xi,\xi_2) \, \lambda_\Theta(\xi_1) \lambda_\Theta(\xi_2) \otimes \lambda_\Theta(\xi_1)^* \, d\xi_1 d\xi_2 \\ [5pt] & = & \big(\Lambda_\Theta \circ (\lambda_{\R^n}(\xi) \otimes id_{\R^n}) \circ \Lambda_\Theta^{-1}\big)(z),
\end{eqnarray*}
where $\lambda_{\R^n}$ denotes the left regular representation on $\R^n$. This shows that $\pi_\Theta$ is weak-$*$ continuous and satisfies the identity $\pi_\Theta(f) = \Lambda_\Theta \circ (f \otimes id_{\R^n}) \circ \Lambda_\Theta^{-1}$ for all $f \in L_\infty(\R^n)$, after identifying $\exp_\xi$ with $\lambda_{\R^n}(\xi)$. Once we have justified the weak-$*$ continuity, the relation $(\sigma_\Theta \otimes id_{\RR_\Theta^{\mathrm{op}}}) \circ \pi_\Theta = (id_{\R^n} \otimes \pi_\Theta) \circ \Delta_{\R^n}$ follows since it trivially holds when acting on $\exp_\xi$ for any $\xi \in \R^n$. Also, it implies that $\mathrm{d}_\Theta$ is affiliated to $\RR_\Theta \bar\otimes \RR_\Theta^{\mathrm{op}}$ and arises as an increasing limit of bounded operators $\mathrm{d}_\Theta(\mathrm{R})$ for $\mathrm{R} > 0$. It remains to show that $\pi_\Theta: L_2^c(\R^n) \to L_2^c(\RR_\Theta) \bar\otimes \RR_\Theta^\mathrm{op}$. Recall that the norm in $L_2^c(\M) \bar\otimes \M_\mathrm{op}$ is given by $$a \mapsto \Big\| \big( \tau \otimes id_{\M_{\mathrm{op}}} \big) (a^*a) \Big\|^{\frac12}_{\M_\mathrm{op}}.$$ When $\M = \RR_\Theta$ and $f \in L_2(\R^n)$ is smooth we find $$\big( \tau_\Theta \otimes id_{\RR_\Theta^{\mathrm{op}}} \big) (\pi_\Theta(f)) = \big( \tau_\Theta \otimes id_{\RR_\Theta^{\mathrm{op}}} \big) \Big( \int_{\R^n} \widehat{f} (\xi) \hskip2pt \lambda_\Theta(\xi) \otimes \lambda_\Theta(\xi)^* \, d\xi \Big) = \widehat{f}(0) \1_{\RR_\Theta^\mathrm{op}}.$$ Therefore, taking $f = \sum_{jk} f_{jk} \otimes e_{jk} \in M_n(L_2^c(\R^n))$ smooth, we obtain 
\begin{eqnarray*}
\big\| \big( \pi_\Theta(f_{jk}) \big)_{jk} \big\|_{M_n(L_2^c(\RR_\Theta) \bar\otimes \RR_\Theta)} & = & \Big\| \Big( \big( \tau_\Theta \otimes id_{\RR_\Theta^\mathrm{op}} \big) \big( \pi_\Theta (f^*f) \big) \Big)_{jk} \Big\|_{M_n \otimes_\mathrm{min} \RR_\Theta^\mathrm{op}}^\frac12 \\ & = & \Big\| \Big( \int_{\R^n} (f^*f)_{jk}(\xi) \, d\xi \Big) \Big\|_{M_n}^\frac12 \ = \ \|f\|_{M_n(L_2^c(\R^n))}.
\end{eqnarray*}
By density, we see that $\pi_\Theta \hskip-1pt : \hskip-1pt L_2^c(\R^n) \to L_2^c(\RR_\Theta) \bar\otimes \RR_\Theta^{\mathrm{op}}$ is a complete isometry. \fin

Let $$S_tf(x) = \int_{\R^n} \widehat{f}(\xi) e^{-t|\xi|^2} e^{2\pi i \langle x,\xi \rangle} \, d\xi$$ denote the heat semigroup acting on $f: \R^n \to \C$. Consider the induced semigroup $S_\Theta = (S_{\Theta,t})_{t \ge 0}$ on $\RR_\Theta$ determined by $\sigma_\Theta \circ S_{\Theta,t} = (S_t \otimes id_{\RR_\Theta}) \circ \sigma_\Theta$. This yields a Markov semigroup which formally acts as $$S_{\Theta,t}(\lambda_\Theta(f)) = \int_{\R^n} f(\xi) e^{-t|\xi|^2} \lambda_\Theta(\xi) \, d\xi.$$ Consider the corresponding column BMO norm 
\begin{eqnarray*}
\|a\|_{\mathrm{BMO}_c(\RR_\Theta)} & = & \hskip2pt \sup_{t > 0} \hskip1pt \Big\| \Big( S_{\Theta,t} (a^*a) - S_{\Theta,t}(a)^*S_{\Theta,t}(a) \Big)^\frac12 \Big\|_{\RR_\Theta} \\ & \sim & \sup_{Q \in \mathcal{Q}} \Big\| \Big( \mean_Q \big| \sigma_\Theta(a) - \sigma_\Theta(a)_Q \big|^2 \, d\mu \Big)^\frac12 \Big\|_{\RR_\Theta} \ = \ \|\sigma_\Theta(a)\|_{\mathrm{BMO}_c(\mathcal{Q}_\Theta)}
\end{eqnarray*}
where $\mathcal{Q}$ denotes the set of all Euclidean balls in $\R^n$, $\mu$ stands for the Lebesgue measure and $\sigma_\Theta(a)_Q$ is the average of $\sigma_\Theta(a)$ over the ball $Q$. The norm equivalence above ---which holds up to constants depending on the dimension $n$--- is a simple consequence of the intertwining identity $\sigma_\Theta \circ S_{\Theta,t} = (S_t \otimes id_{\RR_\Theta}) \circ \sigma_\Theta$ and the equivalence between the BMO norms respectively associated to the heat semigroup and the Euclidean metric in $\R^n$, see \cite[Section 1.2]{JMP} for further details. The space $\mathrm{BMO}_c(\mathcal{Q}_\Theta)$ is an illustration of the operator-valued spaces $\mathrm{BMO}_c(\R^n; \mathcal{B}(\H))$ which were extensively studied by Tao Mei in his PhD Thesis \cite{Mei}. 

We may use these latter spaces to properly define the column space $\mathrm{BMO}_c(\RR_\Theta)$. Indeed, we know from Corollary \ref{Planc+Corep} that $\sigma_\Theta(\RR_\Theta)$ is a subalgebra of $L_\infty(\R^n) \bar\otimes \RR_\Theta$, which in turn is included in $\mathrm{BMO}_c(\mathcal{Q}_\Theta)$. Since we know from \cite{Mei} that $\mathrm{BMO}_c(\mathcal{Q}_\Theta)$ admits a predual $\mathrm{H}^c_1(\mathcal{Q}_\Theta)$, we may define $$\mathrm{BMO}_c(\RR_\Theta) = \overline{\sigma_\Theta(\RR_\Theta)^{\wast}}$$ where the weak-$*$ closure is taken with respect to the pair $(\mathrm{H}^c_1(\mathcal{Q}_\Theta), \mathrm{BMO}_c(\mathcal{Q}_\Theta))$. This kind of BMO spaces over Markov semigroups have been deeply investigated in \cite{JM} for finite von Neumann algebras. The semifinite case is more subtle and we shall give in Appendix B a self-contained argument for $\RR_\Theta$. 

\subsubsection{A Poincar\'e type inequality}

Let $$\S_\Theta = \Big\{ \lambda_\Theta(f) \, : \, f \in \S(\R^n)=\R^n\mbox{-Schwartz class} \Big\}.$$ Define $\partial_{\Theta}^j$ as the linear extension of the map $$\partial_{\Theta}^j(\lambda_\Theta(\xi)) = 2 \pi i \xi_j \lambda_\Theta(\xi)$$ over the quantum Schwartz class $\S_\Theta$ for $1 \le j \le n$. Recall that $\S_\Theta$ is an $*$-algebra since $\lambda_\Theta(f_1) \lambda_\Theta(f_2) = \lambda_\Theta (f_1 *_\Theta f_2) \quad \mbox{and} \quad \lambda_\Theta(f)^* = \lambda_\Theta(f^*_\Theta)$ are stable in $\S_\Theta$. In what follows, we shall be working with this and other natural differential operators in $\RR_\Theta$. The following one is a free analogue of the gradient operator associated to the partial derivatives considered above. Let $\mathcal{L}(\F_n)$ denote the group von Neumann algebra associated to the free group over $n$ generators $\F_n$. It is well-known from (say) \cite{VDN} that $\mathcal{L}(\F_n)$ is generated by $n$ semicircular random variables $s_1, s_2, \ldots, s_n$. Let us consider the map $s: \R^n \to \mathcal{L}(\F_n)$ given by $$s(\xi) = \sum_{j=1}^n \langle \xi, e_j \rangle s(e_j) = \sum_{j=1}^n \xi_j s_j.$$ Then we introduce the $\Theta$-deformed free gradient $$\nabla_\Theta = \sum_{k=1}^n s_k \otimes \partial_{\Theta}^k: \S_\Theta \to \mathcal{L}(\F_n) \bar\otimes \RR_\Theta.$$ If $\nabla$ denotes the free gradient for $\Theta=0$, is easily checked that $$(id_{\mathcal{L}(\F_n)} \otimes \sigma_\Theta) \hskip-1pt \circ \hskip-1pt \nabla_\Theta = \sum_{k=1}^n s_k \otimes (\sigma_\Theta \hskip-1pt \circ \hskip-1pt \partial_{\Theta}^k) = \sum_{k=1}^n s_k \otimes (\partial_k \hskip-1pt \circ \hskip-1pt \sigma_\Theta) = (\nabla \otimes id_{\RR_\Theta}) \hskip-1pt \circ \hskip-1pt \sigma_\Theta.$$ Moreover, let us recall that $\nabla_\Theta(\lambda_\Theta(\xi)) = \displaystyle \summ_k s_k \otimes 2\pi i \xi_k \lambda_\Theta(\xi) = 2\pi i s(\xi) \otimes \lambda_\Theta(\xi)$.

\begin{proposition} \label{Prop-Poincare}
Let $\mathrm{B_R}$ and $q_\mathrm{R}$ stand for any ball of radius $\mathrm{R}$ in $\R^n$ and the characteristic function of it. Given a noncommutative measure space $(\M,\tau)$ and $\varphi: \mathrm{B_R} \to \M$ smooth with $\mathrm{B_R}$-average denoted by $\varphi_{\mathrm{B_R}}$, the following inequality holds for the free gradient $\nabla$ in $\R^n$ $$\Big\| \mean_{\mathrm{B}_\mathrm{R}} \big| \varphi - \varphi_{\mathrm{B_R}} \big|^2 d\mu \Big\|_{\M}^\frac12 \le 2 \sqrt{2} \mathrm{R} \Big\| (\1 \otimes q_\mathrm{R} \otimes \1) (\nabla \otimes id_{\M})(\varphi) \Big\|_{\mathcal{L}(\F_n) \bar\otimes L_\infty(\R^n) \bar\otimes \M}.$$
\end{proposition}

\dem Consider the derivation map $\delta(f) = f \otimes \1 - \1 \otimes f$. We shall use the following straightforward algebraic identity, which is valid for any normal state $\phi$ on any von Neumann algebra $$\phi \big( (f - \phi(f))^*(f - \phi(f)) \big) = \frac12 \phi \otimes \phi \big( \delta(f)^* \delta(f) \big).$$ Applying it for $\varphi = \sum_j f_j \otimes y_j$ we obtain $$(\phi \otimes id_\M) \Big( \big| \varphi - (\phi \otimes id_\M)(\varphi) \big|^2 \Big) = \frac12 (\phi \otimes \phi \otimes id_\M) \Big( \big| (\delta \otimes id_\M)(\varphi) \big|^2 \Big).$$ If $c_\mathrm{R}$ denotes the center of $\mathrm{B_R}$, we observe that $$(\delta \otimes id_\M)(\varphi) = (\varphi - \varphi(c_\mathrm{R})) \otimes \1_{\R^n} - \1_{\R^n} \otimes (\varphi - \varphi(c_\mathrm{R})).$$ Then, letting $\phi$ be the average over $\mathrm{B_R}$ we deduce the following inequality
\begin{eqnarray*}
\Big\| \mean_{\mathrm{B}_\mathrm{R}} \big| \varphi - \varphi_{\mathrm{B_R}} \big|^2 d\mu \Big\|_{\M}^\frac12 & = & \frac{1}{\sqrt{2}} \Big\| (\phi \otimes \phi \otimes id_\M) \Big( \big| (\delta \otimes id_\M)(\varphi) \big|^2 \Big) \Big\|_{\M}^\frac12 \\ & \le & \frac{1}{\sqrt{2}} \Big\| (\phi \otimes id_{\R^n} \otimes id_\M) \Big( \big| (\varphi - \varphi(c_\mathrm{R})) \otimes \1_{\R^n} \big|^2 \Big) \Big\|_{\M}^\frac12 \\ & + & \frac{1}{\sqrt{2}} \Big\| (id_{\R^n} \otimes \phi \otimes id_\M) \Big( \big| \1_{\R^n} \otimes (\varphi - \varphi(c_\mathrm{R})) \big|^2 \Big) \Big\|_{\M}^\frac12 \\ & = & \sqrt{2} \, \Big\| (\phi \otimes id_\M) \Big( \big| \varphi - \varphi(c_\mathrm{R}) \big|^2 \Big) \Big\|_{\M}^\frac12. 
\end{eqnarray*}
In order to estimate the latter term, we use integration by parts to obtain 
\begin{eqnarray*}
\varphi(x) - \varphi(c_\mathrm{R}) \!\!\! & = & \!\!\! \int_0^1 \sum_{k=1}^n \partial_k \varphi \big( t(x-c_\mathrm{R}) + c_{\mathrm{R}} \big) \langle x-c_\mathrm{R}, e_k \rangle \hskip1pt dt \\ \!\!\! & = & \!\!\! \int_0^1 (\underbrace{\tau_{\mathcal{L}(\F_n)} \otimes id_{\R^n}}_{\mathsf{E}_{\R^n}} ) \big( \underbrace{q_\mathrm{R}(x)  \nabla \varphi(t(x-c_\mathrm{R}) + c_{\mathrm{R}})}_{\mathrm{A}(t)} \underbrace{q_\mathrm{R}(x)  s(x-c_\mathrm{R})}_{\mathrm{B}} \big) \hskip1pt dt
\end{eqnarray*}
for $x \in \mathrm{B_R}$. By the operator-convexity of $| \cdot |^2$, we find the inequality below 
\begin{eqnarray*}
\Big\| \mean_{\mathrm{B}_\mathrm{R}} \big| \varphi - \varphi_{\mathrm{B_R}} \big|^2 d\mu \Big\|_{\M}^\frac12 \!\!\!\! & \le & \!\!\!\! \sqrt{2} \Big( \int_0^1 \big\| (\phi \otimes id_{\M}) \big( |\mathsf{E}_{\R^n}(\mathrm{A}(t) \mathrm{B}) |^2 \big) \big\|_\M \hskip1pt dt \Big)^\frac12 \\ \!\!\!\! & \le & \!\!\!\! \sqrt{2} \Big( \int_0^1 \| \mathrm{A}(t) \|^2_{\mathcal{L}(\F_n) \bar\otimes L_\infty(\R^n) \bar\otimes \M} \hskip1pt dt \Big)^\frac12 \|\mathrm{B}\|_{\mathcal{L}(\F_n) \bar\otimes L_\infty(\R^n)}.
\end{eqnarray*}
Now we observe that $\|\mathrm{A}(t)\| \le \|\mathrm{A}(1)\|$ for all $0 \le t \le 1$, so we conclude that
$$\Big\| \mean_{\mathrm{B}_\mathrm{R}} \big| \varphi - \varphi_{\mathrm{B_R}} \big|^2 d\mu \Big\|_{\M}^\frac12 \le \sqrt{2} \|\mathrm{B}\|_{\mathcal{L}(\F_n)} \Big\| (\1 \otimes q_\mathrm{R} \otimes \1) (\nabla \otimes id_{\M})(\varphi) \Big\|_{\mathcal{L}(\F_n) \bar\otimes L_\infty(\R^n) \bar\otimes \M}.$$ Finally, Voiculescu's inequality \cite{VDN} claims $$\|s(h)\|_{\mathcal{L}(\F_n)} = 2 \|h\|_{\R^n},$$ so that $\displaystyle \|\mathrm{B}\|_{\mathcal{L}(\F_n)} \le 2 \sup_{x \in \mathrm{B_R}} \|x-c_\mathrm{R}\| = 2 \mathrm{R}$ and the proof is complete. \fin

\begin{remark} \label{Poincare-BMO}
\emph{Recall that }
\begin{eqnarray*}
\|a\|_{\mathrm{BMO}_c(\RR_\Theta)} & \sim & \|\sigma_\Theta(a)\|_{\mathrm{BMO}_c(\mathcal{Q}_\Theta)} \\ & = & \sup_{\mathrm{R}>0} \Big\| \Big( \mean_{\mathrm{B_R}} \big| \sigma_\Theta(a) - \sigma_\Theta(a)_{\mathrm{B_R}} \big|^2 \, d\mu \Big)^\frac12 \Big\|_{\RR_\Theta}. 
\end{eqnarray*}
\emph{According to Proposition \ref{Prop-Poincare} for $\M = \RR_\Theta$, we deduce 
\begin{eqnarray*}
\|a\|_{\mathrm{BMO}_c(\RR_\Theta)} & \lesssim & \sup_{\mathrm{R}>0} \mathrm{R} \Big\| (\1 \otimes q_\mathrm{R} \otimes \1) (\nabla \otimes id_{\RR_\Theta}) \circ \sigma_\Theta(a) \Big\|_{\mathcal{L}(\F_n) \bar\otimes L_\infty(\R^n) \bar\otimes \RR_\Theta} \\ & = & \sup_{\mathrm{R}>0} \mathrm{R} \Big\| (\1 \otimes q_\mathrm{R} \otimes \1) (id_{\mathcal{L}(\F_n)} \otimes \sigma_\Theta) \circ \nabla_\Theta (a) \Big\|_{\mathcal{L}(\F_n) \bar\otimes L_\infty(\R^n) \bar\otimes \RR_\Theta}.
\end{eqnarray*}}
\end{remark}

\begin{remark}
\emph{Given $2 \le p \le \infty$, we have} 
\begin{eqnarray*}
\big\| \nabla_\Theta (a) \big\|_{L_p(\mathcal{L}(\F_n) \bar\otimes \RR_\Theta)} & \sim & \Big\| \Big( \sum_{k=1}^n (\partial_{\Theta}^ka)(\partial_{\Theta}^k a)^* \Big)^\frac12 \Big\|_{L_p(\RR_\Theta)} \\ & + & \Big\| \Big( \sum_{k=1}^n (\partial_{\Theta}^k a)^*(\partial_{\Theta}^k a) \Big)^\frac12 \Big\|_{L_p(\RR_\Theta)}
\end{eqnarray*}
\emph{from the operator-valued form of Voiculescu's inequality \cite{JPX,V}. Let us recall in passing that this norm equivalence holds in the category of operator spaces and moreover, the constants do not depend on the dimension $n$. This justifies our choice of free generators in the definition of $\nabla_\Theta$. An alternative choice would have been to work with Rademacher variables or matrix units, but the former does not lead to the same norm equivalences for $p=\infty$. If $\Theta=0$ we get} $$\big\| \nabla (f) \big\|_{L_p(\mathcal{L}(\F_n) \bar\otimes L_\infty(\R^n))} \sim \Big\| \Big( \sum_{k=1}^n |\partial_x^k f|^2 \Big)^\frac12 \Big\|_{L_p(\R^n)}.$$
\end{remark}

\subsection{Quantum Euclidean variables}

We consider other characterizations of $\RR_\Theta$ in terms of the infinitesimal generators of $u_j(s)$. These unbounded operators play the role in $\RR_\Theta$ of the Euclidean variables $x_j$ in $\R^n$. We will use them to study the quantum analogue of the Schwartz class, to give an intrinsic characterization of the quantum distance $\mathrm{d}_\Theta$ and to deduce the algebraic structure of $\RR_\Theta$.  

\subsubsection{Another approach towards $\RR_\Theta$}

Define $$x_{\Theta,j} = \frac{1}{2\pi i} \, \frac{d}{ds} \Big|_{s=0} (u_j(s)) \quad \mbox{for} \quad 1 \le j \le n,$$ with $u_j(s)$ the generating unitaries of the quantum Euclidean space $\RR_\Theta$. These are the (self-adjoint) infinitesimal generators of the one-parameter groups of unitaries $(u_j(s))_{s \in \R}$ given by Stone's theorem and may be regarded as quantum forms of the Euclidean variables. Namely, when $\Theta=0$ the one-parameter unitary group $u_j(s)$ is composed of multiplication operators by the Fourier characters $x \mapsto \exp(2 \pi i s x_j)$ and $$2\pi i x_j = \partial_s (e^{2\pi i s x_j})_{\mid_{s=0}}.$$ The operators $x_{\Theta,j}$ enjoy some fundamental properties of the Euclidean variables. 

\begin{proposition} \label{Quantumx}
The following results hold$\hskip1pt :$
\begin{itemize}
\item[\emph{i)}] The generators $x_{\Theta,j}$ satisfy for $1 \le j,k \le n$ $$\hskip25pt [x_{\Theta,j}, x_{\Theta,k}] = \frac{1}{2\pi i} \, \Theta_{jk} \ \Leftrightarrow \ u_j(s) u_k(t) = e^{2\pi i \Theta_{jk} st} u_k(t) u_j(s).$$

\vskip3pt

\item[\emph{ii)}] Recall the definition of the quantum Schwartz class $$\hskip25pt \S_\Theta = \Big\{ \lambda_\Theta(f) : f \in \S(\R^n)\Big\}.$$ The infinitesimal generators $x_{\Theta,j}$ are densely defined unbounded operators affiliated to $\RR_\Theta$. Moreover, in the \emph{GNS} representation on $L_2(\RR_\Theta)$ we find $\S_\Theta \subset \mathrm{dom}(x_{\Theta,j})$ and $x_{\Theta,j} \S_\Theta, \S_\Theta x_{\Theta,j} \subset \S_\Theta$. More precisely 
\begin{eqnarray*}
\hskip30pt x_{\Theta,j} \lambda_\Theta(f) = \lambda_\Theta(D_{\Theta,j}^\ell f) & \mbox{where} & D_{\Theta,j}^\ell = \hskip4pt \sum_{k=1}^{j-1} \Theta_{jk} M_{\xi_k} - \frac{1}{2\pi i} \partial_\xi^j, \\ \hskip30pt \lambda_\Theta(f) x_{\Theta,j} = \lambda_\Theta(D_{\Theta,j}^r f) & \mbox{where} & D_{\Theta,j}^r = \sum_{i=j+1}^{n} \Theta_{ij} M_{\xi_i} - \frac{1}{2\pi i} \partial_\xi^j,
\end{eqnarray*} 
for $f \in \S(\R^n)$ and $M_{\xi_k}f(\xi) = \xi_k f(\xi)$. In addition we get $[D_{\Theta,j}^\ell, D_{\Theta,k}^r]=0$.

\vskip6pt

\item[\emph{iii)}] Let $(x_{\Theta,j})_j$ and $(y_{\Theta,j})_j$ be the infinitesimal generators associated to $\RR_\Theta \otimes \1$ and $\1 \otimes \RR_\Theta^\mathrm{op}$ respectively. Then, we may relate the quantum distance $\mathrm{d}_\Theta$ with these quantum variables as follows $$\hskip25pt \mathrm{d}_\Theta = \Big( \sum_{j=1}^n \big( x_{\Theta,j} - y_{\Theta,j} \big)^2 \Big)^{\frac{1}{2}}.$$
\end{itemize}
\end{proposition}

\dem All the assertions are quite standard. Assume $[x_{\Theta,j}, x_{\Theta,k}] = \frac{\Theta_{jk}}{2\pi i}$ and define the maps $\varphi_j(s)[\mathrm{X}] = \exp(2\pi i s x_{\Theta,j}) \mathrm{X} \exp(-2\pi i s x_{\Theta,j}) = u_j(s) \mathrm{X} u_j(-s).$ Recalling that $\varphi_j(0)$ is the identity map and noticing that $$\frac{d}{ds} \big( \varphi_j(s) \big) = 2\pi i \Big( x_{\Theta,j} \varphi_j(s) - \varphi_j(s) x_{\Theta,j} \Big),$$ we deduce $\varphi_j(s)[\mathrm{X}] = \exp \big( 2\pi i s \delta_j(\mathrm{X}) \big)$ for the derivation $\delta_j(\mathrm{X}) = [x_{\Theta,j}, \mathrm{X}]$, so that $$u_j(s) u_k(t) u_j(-s) = \sum_{n \ge 0} \frac{(2 \pi i s)^n}{n!} \delta_j^n(u_k(t)).$$ Consider the map given by $\Theta_{jk} z_{jk} = 2\pi i x_{\Theta,k}$, so that $\delta_j(z_{jk}) = [x_{\Theta,j}, z_{jk}] = 1$ and $\delta_j(h(z_{jk})) = h'(z_{jk})$. Since $u_k(t) = \exp (2\pi i t x_{\Theta,k}) = \exp(t \Theta_{jk} z_{jk}) = h(z_{jk})$ we deduce the identity $$\delta_j^n(u_k(t)) = h^{(n)}(z_{jk}) = (t \Theta_{jk})^n u_k(t)$$ which yields $u_j(s)u_k(t) = e^{2\pi i \Theta_{jk} st} u_k(t) u_j(s)$ as desired. Reciprocally, assuming this commutation relation we express the commutator of $x_{\Theta,j}$ and $x_{\Theta,k}$ as follows 
\begin{eqnarray*}
[x_{\Theta,j}, x_{\Theta,k}] & = & \frac{1}{- 4 \pi^2} \frac{d^2}{dsdt} {\Big|_{s=t=0}} \big( u_j(s)u_k(t) - u_k(t)u_j(s) \big) \\ & = & \frac{1}{- 4 \pi^2} \frac{d^2}{dsdt} {\Big|_{s=t=0}} \big( e^{2\pi i \Theta_{jk}st} - 1 \big) u_k(t) u_j(s) \ = \ \frac{2\pi i}{- 4 \pi^2} \Theta_{jk} \ = \ \frac{\Theta_{jk}}{2 \pi i}.
\end{eqnarray*}
Regarding the second assertion ii), the first identity can be justified as follows
\begin{eqnarray*}
x_{\Theta,j} \lambda_\Theta(f) & = & \frac{1}{2\pi i} \frac{d}{ds} {\Big|_{s=0}} (u_j(s)) \lambda_\Theta(f) \\ & = & \frac{1}{2\pi i} \frac{d}{ds} {\Big|_{s=0}} \int_{\R^n} f(\xi) u_j(s) \lambda_ \Theta(\xi) d\xi \\ & = & \frac{1}{2\pi i} \frac{d}{ds} {\Big|_{s=0}} \int_{\R^n} f(\xi) e^{2\pi i \sum_{k < j} \Theta_{jk} s \xi_k} \lambda_ \Theta(\xi + se_j) d\xi \\ & = & \frac{1}{2\pi i} \frac{d}{ds} {\Big|_{s=0}} \int_{\R^n} f(\xi - s e_j) e^{2\pi i \sum_{k < j} \Theta_{jk} s \xi_k} \lambda_ \Theta(\xi) d\xi \\ & = & \frac{1}{2\pi i} \Big[ \int_{\R^n} \Big( 2\pi i \sum_{k < j} \Theta_{jk} \xi_k \Big) f(\xi) \lambda_ \Theta(\xi) d\xi - \int_{\R^n} \partial_\xi^j f(\xi) \lambda_ \Theta(\xi) d\xi \Big].
\end{eqnarray*}
The second identity is proved similarly. This shows that $\mathcal{S}_\Theta$ is a common core of the $x_{\Theta,j}$ for $1 \le j \le n$. Thus, it just remains to show that $[D_{\Theta,j}^\ell, D_{\Theta,k}^r]=0$ to complete the proof of ii). This is clear for $j \le k$, as for $j>k$   
\begin{eqnarray*}
[D_{\Theta,j}^\ell, D_{\Theta,k}^r] & = & -\frac{\Theta_{jk}}{2\pi i} \Big(  [\partial_\xi^j, M_{\xi_j}] + [M_{\xi_k}, \partial_\xi^k] \Big) \ = \ 0.
\end{eqnarray*}
Finally, since $\mathrm{d}_\Theta = \pi_\Theta(| \cdot |)$ and $\pi_\Theta$ is a $*$-homomorphism, assertion iii) reduces to show that $\pi_\Theta ( x_j ) = x_{\Theta,j} - y_{\Theta,j}$ for $1 \le j \le n$. This can be proved again with a differentiation argument as follows
\begin{eqnarray*}
\pi_\Theta (x_j) & = & \frac{1}{2\pi i} \frac{d}{ds} {\Big|_{s=0}} \pi_\Theta \big( e^{2\pi i s \langle \cdot, e_j \rangle} \big) \\ & = & \frac{1}{2\pi i} \frac{d}{ds} {\Big|_{s=0}} \big( \lambda_\Theta(se_j) \otimes \lambda_\Theta(se_j)^* \big) \\ & = & \Big( \frac{1}{2\pi i} \frac{d}{ds} {\Big|_{s=0}} u_j(s) \Big) \otimes \1 + \1  \otimes \Big( \frac{1}{2\pi i} \frac{d}{ds} {\Big|_{s=0}} u_j(-s) \Big) = x_{\Theta,j} - y_{\Theta,j},
\end{eqnarray*}
according to our definition of $x_{\Theta,j}$ and $y_{\Theta,j}$ in $\RR_\Theta \otimes \1$ and $\1 \otimes \RR_\Theta^{\mathrm{op}}$ respectively. \fin

\begin{remark} \label{RemQuantumx}
\emph{A few comments are in order:}
\begin{itemize}
\item[$\bullet$] \emph{$\RR_\Theta$ is generated by the spectral projections of the quantum variables $x_{\Theta,j}$.}

\vskip5pt

\item[$\bullet$] \emph{The Euclidean Schwartz $\S(\R^n)$ class is the space of infinitely differentiable functions $f: \R^n \to \C$ which satisfy that $f$ and its derivatives decay at infinity faster than polynomials. In $\RR_\Theta$, we find 
\begin{eqnarray*}
\lefteqn{\hskip-30pt \Big( \prod_{1 \le r \le m_\ell}^{\rightarrow} x_{\Theta,j_r} \Big) \partial_\Theta^\beta \big( \lambda_\Theta(f) \big) \Big( \prod_{1 \le s \le m_r}^{\rightarrow} x_{\Theta,k_s} \Big)} \\ \hskip60pt & = & \lambda_\Theta \Big[ \Big( \prod_{1 \le r \le m_\ell}^{\rightarrow} D_{\Theta,j_r}^\ell \Big) \Big( \prod_{1 \le s \le m_r}^{\rightarrow} D_{\Theta, k_s}^r \Big) M_{(2\pi i \xi)^\beta} f \Big] \ \in \ \S_\Theta,
\end{eqnarray*}
which admits other representations since $D_{\Theta,j_r}^\ell$ and $D_{\Theta, k_s}^r$ commute. It shows that the quantum Schwartz class is also closed under differentiation and left/right multiplication by quantum polynomials.} 

\vskip5pt

\item[$\bullet$] \emph{Proposition \ref{Quantumx} iii) establishes a canonical Pythagorean formula for the quantum Euclidean distance $\mathrm{d}_\Theta$ in terms of  quantum variables. This shows that the metric $\mathrm{d}_\Theta$ that we shall be using along the rest of the paper is not \emph{induced} but somehow \emph{intrinsic} to $\RR_\Theta$. This gives some evidence that our main results in this paper are formulated in their most natural way.} 

\vskip5pt

\item[$\bullet$] \emph{We also note in passing that the quantum variables $z_{\Theta,j} = x_{\Theta,j} - y_{\Theta,j}$ from Proposition \ref{Quantumx} iii) are pairwise commuting for different values of $1 \le j \le n$ since}
\begin{eqnarray*}
\hskip20pt [z_{\Theta,j}, z_{\Theta,k}] & = & z_{\Theta,j} \bullet z_{\Theta,k} - z_{\Theta,k} \bullet z_{\Theta,j} \\ & = & [x_{\Theta,j}, x_{\Theta,k}]_{\RR_\Theta} \otimes \1 + \1 \otimes [y_{\Theta,j}, y_{\Theta,k}]_{\RR_\Theta^{\mathrm{op}}} \ = \ 0.
\end{eqnarray*}  
\end{itemize}
\end{remark}

\subsubsection{On the quantum Schwartz class} \label{QSch}

Using quantum variables, we are ready to prove some fundamental properties of the quantum Schwartz class. The analogues in the commutative case $\Theta=0$ are rather easy to prove. Let us consider the map $j_\Theta: \S(\R^n) \to \S_\Theta$ given by $$j_\Theta \Big( \int_{\R^n} f(\xi) e^{2 \pi i \langle \cdot, \xi \rangle} \, d\xi \Big) = \int_{\R^n} f(\xi) \lambda_\Theta(\xi) \, d\xi,$$ so that $j_\Theta(f) = \lambda_\Theta(\widehat{f} \hskip2pt)$. By Remark \ref{Plancherel} and Plancherel theorem, $j_\Theta$ extends to an isometric isomorphism $L_2(\R^n) \to L_2(\RR_\Theta)$. We shall also need the space $\S_\Theta'$ of continuous linear functionals on $\S_\Theta$, tempered quantum distributions. Finally recalling that the quantum variables $x_{\Theta,j}$ are affiliated to $\RR_\Theta$, we set for $1 \le j \le n$ $$\RR_{\Theta,j} = \Big\langle \mathrm{spectral \ projections \ of \ } x_{\Theta,j} \Big\rangle'' \subset \RR_\Theta.$$ We write $\RR_j$ for $\RR_{\Theta,j}$ with $\Theta = 0$. We begin with an elementary auxiliary result. 

\begin{lemma} \label{j-spectral}
We have$\hskip1pt:$
\begin{itemize}
\item[\emph{i)}] $j_\Theta(x_j^k) = x_{\Theta,j}^k$ in the sense of distributions.

\item[\emph{ii)}] $j_\Theta: \RR_j \to \RR_{\Theta,j}$ is a normal $*$-homomorphism.
\end{itemize}
\end{lemma}

\dem Every element in $\S_\Theta$ may be represented in the form $j_\Theta(f)$ for some $f$ in the Schwartz class of $\R^n$. On the other hand, since $j_\Theta: L_2(\R^n) \to L_2(\RR_\Theta)$ is an isometric isomorphism, we define $j_\Theta(x_j^k) \in \S_\Theta'$ by $$\big\langle j_\Theta(x_j^k), j_\Theta(f) \big\rangle = \int_{\R^n} x_j^k f(x) \, dx.$$ Thus, it suffices to see that this quantity coincides with $$\tau_\Theta \big(x_{\Theta,j}^k j_\Theta(f) \big) = \tau_\Theta \big( x_{\Theta,j}^k \lambda_\Theta(\widehat{f} \hskip1pt ) \big) = \tau_\Theta \Big( \lambda_\Theta\big( (D_{\Theta,j}^\ell)^k \widehat{f} \hskip2pt \big) \Big) = (D_{\Theta,j}^\ell)^k \widehat{f} \hskip1pt (0)$$ for $D_{\Theta,j}^\ell = \sum_{s<j} \Theta_{js} M_{\xi_s} - \frac{1}{2\pi i} \partial_{j}^\xi$, by Proposition \ref{Quantumx}. A simple computation shows that this is indeed the case. Next, assertion ii) follows from the fact that $\RR_{\Theta,j} \simeq L_\infty(\R)$ for $1 \le j \le n$ no matter which is the deformation $\Theta$. Indeed, in order to use the same terminology as in the proof of Proposition \ref{Trace}, we shall assume for convenience that $j=n$. Then, we may identify $\RR_{\Theta,n}$ with the subalgebra $\1 \rtimes_{\beta_{n-1}} \R$ of the von Neumann algebra $\RR_{\Xi} \rtimes_{\beta_{n-1}} \R$, which in turn is isomorphic to $\RR_\Theta$. Then, it is a well-known fact that we have $\1 \rtimes_{\beta_{n-1}} \R \simeq L_\infty(\R)$. \fin

\begin{proposition} \label{quantumSchwarzt1}
If $\Theta \in A_n(\R)$ and $\gamma > \frac12$, we find $$\prod_{1 \le j \le n}^\rightarrow \big( \1 + |x_{\Theta,j}|^\gamma \big)^{-1} \in L_2(\RR_\Theta).$$ In particular, the quantum Schwartz class $\S_\Theta \subset L_p(\RR_\Theta)$ for all $p > 0$.
\end{proposition}

\dem According to Lemma \ref{j-spectral} $$\big( \1 + |x_{\Theta,j}|^\gamma \big)^{-1} = j_\Theta \Big( \frac{1}{1 + |x_j|^\gamma} \Big).$$ Let us proceed by induction on $n$, the case $n=1$ being trivial. According to Proposition \ref{Trace}, $\tau_\Theta$ coincides with the crossed product trace $\tau_\rtimes$ in $\RR_\Xi \rtimes \R$ which in turn factorizes for operators with separated variables. This means that 
\begin{eqnarray*}
\lefteqn{\hskip-30pt \tau_\Theta \Big( \Big| \prod_{1 \le j \le n}^\rightarrow \big( \1 + |x_{\Theta,j}|^\gamma \big)^{-1} \Big|^2 \Big)} \\ & = & \tau_\Xi \Big( \Big| \prod_{1 \le j \le n-1}^\rightarrow \big( \1 + |x_{\Theta,j}|^\gamma \big)^{-1} \Big|^2 \Big) \int_\R \frac{dx}{(1 + |x|^\gamma)^2}
\end{eqnarray*}
and we conclude by induction. To prove the last assertion, since $\S_\Theta \subset \RR_\Theta$ it clearly suffices to show that $\S_\Theta \subset L_p(\RR_\Theta)$ for $p$ small. Assume $p = 1/m$ for $m \in \Z_+$ and let $$\mathrm{Q} = \Big| \prod_{1 \le j \le n}^\rightarrow \big( \1 + |x_{\Theta,j}|^\gamma \big) \Big|^{2m}.$$ According to H\"older's inequality, we find for $f \in \mathcal{S}(\R^n)$ 
\begin{eqnarray*}
\big\| \lambda_\Theta(f) \big\|_p & \le & \big\| Q^{-1} \big\|_p \big\| Q \lambda_\Theta(f) \big\|_\infty \\ & = & \tau_\Theta \Big[ \Big| \prod_{1 \le j \le n}^\rightarrow \big( \1 + |x_{\Theta,j}|^\gamma \big)^{-1} \Big|^2 \Big]^{\frac{1}{p}} \big\| \lambda_\Theta( Q[D_{\Theta,j}^\ell] f) \big\|_\infty
\end{eqnarray*}
where $Q[D_{\Theta,j}^\ell]$ is the differential operator associated to $Q$ according to the second point of Remark \ref{RemQuantumx}. Since $Q[D_{\Theta,j}^\ell] f \in \mathcal{S}(\R^n)$, the finiteness of the quantity in the right hand side is guaranteed by the first assertion in the statement. \fin

\begin{proposition} \label{quantumSchwarzt2}
We have$\hskip1pt:$
\begin{itemize}
\item[\emph{i)}] $\S_\Theta$ is weak-$*$ dense in $\RR_\Theta$.

\item[\emph{ii)}] $\S_\Theta$ is dense in $L_p(\RR_\Theta)$ for all $p > 0$.
\end{itemize}
In particular, the same density results hold for $\lambda_\Theta(L_1(\R^n)) \subset \RR_\Theta$.
\end{proposition}

\dem Since finite sums of the elementary frequencies $\lambda_\Theta(\xi)$ are weak-$*$ dense in $\RR_\Theta$ by construction, it suffices to approximate $\lambda_\Theta(\xi)$ by elements of $\S_\Theta$ in the weak-$*$ topology. In other words, we need to find a family of functions $\phi_{\xi,\varepsilon} \in \mathcal{S}(\R^n)$ so that
$$\lim_{\varepsilon \to 0} \tau_\Theta \big( (\lambda_\Theta(\phi_{\xi,\varepsilon}) - \lambda_\Theta(\xi)) a \big) = 0 \quad \mbox{for all} \quad a \in L_1(\RR_\Theta).$$
If $\mathrm{B}_\varepsilon(\xi)$ denotes the Euclidean ball around $\xi$ of radius $\varepsilon$, let $\phi_{\xi,\varepsilon}$ be a smoothing of the function $|\mathrm{B}_\varepsilon(\xi)|^{-1} \chi_{\mathrm{B}_\varepsilon(\xi)}$, so that Lebesgue differentiation theorem holds for the family $\{ \phi_{\xi,\varepsilon}: \varepsilon > 0 \}$. Now, since $\S_\Theta$ is dense in $L_2(\RR_\Theta)$, the same holds for $\S_\Theta \S_\Theta \subset L_2(\RR_\Theta) L_2(\RR_\Theta) = L_1(\RR_\Theta)$ and we may approximate $a$ by a sequence $\lambda_\Theta(f_j) \in \S_\Theta$. Recall that $$\big\| \lambda_\Theta(\phi_{\xi,\varepsilon}) - \lambda_\Theta(\xi) \big\|_{\RR_\Theta} \le 1 + \int_{\R^n} \phi_{\xi,\varepsilon}(\zeta) d\zeta = 2.$$ Given $\delta > 0$, there exists $j_\delta \ge 1$ so that $\| a - \lambda_\Theta(f_{j_\delta})\|_1 < \delta/2$. Thus
\[
  \Big| \lim_{\varepsilon \to 0} \tau_\Theta \big( (\lambda_\Theta(\phi_{\xi,\varepsilon}) - \lambda_\Theta(\xi)) a \big) \Big|
  \le \delta + \Big| \lim_{\varepsilon \to 0} \tau_\Theta \big( (\lambda_\Theta(\phi_{\xi,\varepsilon}) - \lambda_\Theta(\xi)) \lambda_\Theta(f_{j_\delta}) \big) \Big|.
\]
On the other hand, since $\lambda_\Theta(\xi)^* = e^{2\pi i \langle \xi, \Theta_{\downarrow} \xi \rangle} \lambda_\Theta(-\xi)$ we find 
\begin{eqnarray*}
\lefteqn{\tau_\Theta \big( (\lambda_\Theta(\phi_{\xi,\varepsilon}) - \lambda_\Theta(\xi)) \lambda_\Theta(f_{j_\delta}) \big)} \\ [3pt]
& = & \phi_{\xi,\varepsilon} *_\Theta f_{j_\delta}(0) - e^{-2 \pi i \langle \xi, \Theta_{\downarrow} \xi \rangle} f_{j_\delta}(-\xi) \\
& = & \int_{\R^n} \phi_{\xi,\varepsilon}(\zeta) e^{-2 \pi i \langle \zeta, \Theta_{\downarrow} \zeta \rangle} f_{j_\delta}(- \zeta) d\zeta - e^{-2 \pi i \langle \xi, \Theta_{\downarrow} \xi \rangle} f_{j_\delta}(-\xi), 
\end{eqnarray*}
where $\Theta_\downarrow$ is the lower triangular part of $\Theta$. The expression above converges to $0$ as $\varepsilon \to 0$. Letting $\delta \to 0$ we conclude that $\S_\Theta$ is weak-$*$ dense in $\RR_\Theta$. Let us now prove that $\S_\Theta$ is norm dense in $L_p(\RR_\Theta)$ for all $p > 0$. Since $\S_\Theta \S_\Theta \subset \S_\Theta$, it suffices from H\"older inequality to prove norm density in the case $p > 2$. Given $a \in L_p(\RR_\Theta)$ for some $p > 2$, we may approximate it in the $L_p$-norm by another element in $\RR_\Theta$ which is left/right supported by a finite projection. In other words, we may assume that $a$ itself belongs to $\RR_\Theta$ and $a = qaq$ for some projection $q$ satisfying $\tau_\Theta (q) < \infty$. Pick two sequences $f_j, g_k \in \mathcal{S}(\R^n)$ satisfying that $$\mathop{\mathrm{w}^*\mbox{-\hskip1pt lim}}_{j \to \infty} \lambda_\Theta(f_j) = a \quad \mbox{and} \quad \mathop{\mathrm{w}^*\mbox{-\hskip1pt lim}}_{k \to \infty} \lambda_\Theta(g_k) = q.$$ By Kaplanski density theorem, we may also assume that $$\sup_{j,k \ge 1} \Big( \|\lambda_\Theta(f_j)\|_{\RR_\Theta} + \|\lambda_\Theta(g_k)\|_{\RR_\Theta} \Big) \le 1 + \|a\|_{\RR_\Theta} < \infty$$ and both convergences hold strongly. Therefore, since $a \in L_2(\RR_\Theta)$, given $\delta > 0$ there must exists $k_\delta$ satisfying $\|a (q - \lambda_\Theta(g_{k_\delta}))\|_2 < \frac{\delta}{2}$. Moreover, once the index $k_\delta$ is fixed and since $\lambda_\Theta (g_{k_\delta}) \in L_2(\RR_\Theta)$ there must exists an index $j_\delta$ satisfying the inequality $\|(a - \lambda_\Theta(f_{j_\delta})) \lambda_\Theta(g_{k_\delta})\|_2 < \frac{\delta}{2}$. Combining these estimates $$\big\| a - \lambda_\Theta(f_{j_\delta}) \lambda_\Theta(g_{k_\delta}) \big\|_2 \le \big\|a (q - \lambda_\Theta(g_{k_\delta}))\big\|_2 + \big\|(a - \lambda_\Theta(f_{j_\delta})) \lambda_\Theta(g_{k_\delta})\big\|_2 < \delta.$$ On the other hand, by the three lines lemma 
\begin{eqnarray*}
\big\| a - \lambda_\Theta(f_{j_\delta}) \lambda_\Theta(g_{k_\delta}) \big\|_p & \le & \big\| a - \lambda_\Theta(f_{j_\delta}) \lambda_\Theta(g_{k_\delta}) \big\|_\infty^{1-\frac{2}{p}} \big\| a - \lambda_\Theta(f_{j_\delta}) \lambda_\Theta(g_{k_\delta}) \big\|_2^{\frac{2}{p}} \\ & \le & \delta^{\frac{2}{p}} \big\| a - \lambda_\Theta(f_{j_\delta}) \lambda_\Theta(g_{k_\delta}) \big\|_\infty^{1-\frac{2}{p}} \ \le \ \big( 2 \|a\|_{\RR_\Theta} \big)^{1-\frac{2}{p}} \delta^{\frac{2}{p}}.
\end{eqnarray*}  
Taking $\delta \to 0$ we see that $\S_\Theta$ is norm dense in $L_p(\RR_\Theta)$ for $p > 2$. \fin

\subsubsection{Structure of $\RR_\Theta$}

We start by showing the very simple algebraic structure of quantum Euclidean spaces $\RR_\Theta$. Indeed, given $\Theta \in A_n(\R)$ and according to the spectral theorem, there exist $d_1, d_2 \in \Z_+$ with $d_1 + 2d_2 = n$ and $\kappa_1, \kappa_2, \ldots \kappa_{d_2} \in \R \setminus \{0\}$ satisfying the following relation for some orthogonal matrix $\mathrm{B} \in \mathrm{SO}(n)$ and for $\Phi$ the $d_1 \times d_1$ 0-matrix$$\Theta = \mathrm{B} \hskip1pt \Big[ \underbrace{\Phi \oplus \bigoplus_{j=1}^{d_2} \kappa_j \Big( \hskip-5pt \begin{array}{rl} 0 & 1 \\ - 1 & 0 \end{array} \Big)}_{\Delta} \Big] \hskip1pt \mathrm{B}^* = \mathrm{B} \hskip2pt \Delta \hskip2pt \mathrm{B}^*.$$ Arguing as in Proposition \ref{Trace} iii) for $n=2$, we find that 
\begin{eqnarray*}
\RR_\Delta & \simeq & L_\infty(\R^{d_1}) \hskip1pt \bar\otimes \hskip1pt \Big( \bigotimes_{j=1}^{d_2} \mathcal{B}(L_2(\R)) \Big) \ \simeq \ L_\infty \big( \R^{d_1}; \mathcal{B}(L_2(\R^{d_2})) \big) \\ & \simeq & L_\infty(\R^{d_1}) \hskip1pt \bar\otimes \hskip1pt \Big( \bigotimes_{j=1}^{d_2} L_\infty(\R) \rtimes \R \Big) \ \simeq \ L_\infty \big( \R^{d_1}; L_\infty(\R^{d_2}) \rtimes \R^{d_2} \big)
\end{eqnarray*}
is a type I von Neumann algebra. Since the commutation relations are determined by $\lambda_\Theta(\xi) \lambda_\Theta(\eta) = \exp (2\pi i \langle \xi, \Theta \eta \rangle) \lambda_\Theta(\eta) \lambda_\Theta(\xi)$ it is tempting to set $\lambda_\Delta(\xi) = \lambda_\Theta(\mathrm{B}\xi)$ to conclude that $\RR_\Theta \simeq \RR_\Delta$ is also a type I von Neumann  algebra. This choice of unitaries do not arise however from a family of one-parameter groups of unitaries as expected. The right change of variables is $x_{\Theta} \mapsto \mathrm{B} x_\Theta$, where $x_\Theta$ stands for $(x_{\Theta,1}, x_{\Theta,2}, ..., x_{\Theta,n})$, at the level of infinitesimal generators. If we want to take exponentials to generate one-parameter groups of unitaries $s \mapsto \exp(2 \pi i s (\mathrm{B} x_\Theta)_j)$ new extra terms appear due to nonvanishing commutators. 

\begin{proposition}
The unitaries $$\lambda_\Delta(\xi) = \exp \Big( \pi i \sum_{j < k} \big( \xi_j \xi_k \Delta_{jk} - (\mathrm{B}\xi)_j (\mathrm{B}\xi)_k \Theta_{jk} \big) \Big) \lambda_\Theta(\mathrm{B}\xi)$$ generate $\RR_\Delta$. In particular, $\RR_\Theta \simeq \RR_\Delta$ so that quantum Euclidean spaces $\RR_\Theta$ are always type \emph{I} von Neumann algebras which are invariant under conjugation by $\mathrm{SO}(n)$. Moreover, the traces coincide $\tau_\Theta = \tau_\Delta$ and the one-parameter unitary groups $w_j(s) = \exp(2 \pi i s x_{\Delta,j}) = \lambda_\Delta(s e_j)$ have the form $$w_j(s) = \exp \Big( - \pi i s^2 \sum_{\alpha < \beta} \mathrm{B}_{j \alpha}^* \Theta_{\alpha \beta} \mathrm{B}_{\beta j} \Big) \lambda_\Theta(s \mathrm{B}e_j).$$
\end{proposition}

\dem Consider the self-adjoint operators $$x_{\Delta,j} = \sum_{k=1}^n \mathrm{B}_{jk}^* \, x_{\Theta,k} = \sum_{k=1}^n \mathrm{B}_{kj} x_{\Theta,k}.$$ It follows from Proposition \ref{Quantumx} that the quantum Schwartz class $\S_\Theta$ is a common core for the family $x_{\Delta,j}$ with $1 \le j \le n$. In particular, these operators are densely defined in the Hilbert space $L_2(\RR_\Theta) \simeq L_2(\R^n) \simeq L_2(\RR_\Delta)$. On the other hand, the commutators are
$$[x_{\Delta,j}, x_{\Delta,k}] = \sum_{1 \le \alpha, \beta \le n} \mathrm{B}_{j\alpha}^* [x_{\Theta,\alpha}, x_{\Theta,\beta}] \mathrm{B}_{\beta k} = \frac{1}{2 \pi i} \sum_{1 \le \alpha, \beta \le n} \mathrm{B}_{j\alpha}^* \Theta_{\alpha \beta} \mathrm{B}_{\beta k} = \frac{1}{2\pi i} \Delta_{jk}.$$
Therefore, Proposition \ref{Quantumx} implies that $\RR_\Delta$ is the weak-$*$ closure of the $\mathrm{C}^*$-algebra generated by the one-parameter unitary groups $w_j(s) = \exp(2 \pi i s x_{\Delta,j})$ for $j \le n$ or equivalently by the products $$w_1(\xi_1) w_2(\xi_2) \cdots w_n(\xi_n) = \prod_{1 \le j \le n}^{\rightarrow} \exp(2 \pi i \xi_j x_{\Delta,j}).$$ Consequently, if we can justify the equality $$\prod_{1 \le j \le n}^{\rightarrow} \exp(2 \pi i \xi_j x_{\Delta,j}) = \exp \Big( \pi i \sum_{j < k} \big( \xi_j \xi_k \Delta_{jk} - (\mathrm{B}\xi)_j (\mathrm{B}\xi)_k \Theta_{jk} \big) \Big) \lambda_\Theta(\mathrm{B}\xi)$$ it will follow automatically that $\RR_\Theta \simeq \RR_\Delta$ as expected. The identity $\tau_\Theta = \tau_\Delta$ and the expression given for $w_j(s)$ also follow easily from the above equality. This is proved from the Baker-Campbell-Hausdorff formula. Namely, since we know that $[x_{\Delta,j}, x_{\Delta,k}] = \frac{1}{2\pi i} \Delta_{jk}$ we may use the simple identity below for operators $\mathrm{X}, \mathrm{Y}$ with vanishing iterated brackets $$\log \big( \exp \mathrm{X} \exp \mathrm{Y} \big) = \mathrm{X} + \mathrm{Y} + \frac12 [\mathrm{X}, \mathrm{Y}].$$ Taking $\mathrm{X}_j = 2 \pi i \xi_j x_{\Delta,j}$ we have $[\mathrm{X}_j, \mathrm{X}_k] = 2 \pi i \xi_j \xi_k \Delta_{jk}$, so that 
\begin{eqnarray*}
\prod_{1 \le j \le n}^{\rightarrow} \exp(2 \pi i \xi_j x_{\Delta,j}) & = & \prod_{1 \le j \le n}^{\rightarrow} \exp \mathrm{X}_j \\ & = & \exp \Big( \frac12 \sum_{j < k} [\mathrm{X}_j, \mathrm{X}_k] \Big) \exp \Big( \sum_{j=1}^n \mathrm{X}_j \Big) \\ & = & \exp \Big( \pi i \sum_{j < k} \xi_j \xi_k \Delta_{jk} \Big) \exp \Big( 2 \pi i \sum_{j=1}^n \xi_j x_{\Delta,j} \Big) \\ & = & \exp \Big( \pi i \sum_{j < k} \xi_j \xi_k \Delta_{jk} \Big) \exp \Big( 2 \pi i \sum_{k=1}^n (\mathrm{B} \xi)_k x_{\Theta,k} \Big).
\end{eqnarray*} 
Using the same formula for the family $\mathrm{Z}_j = 2 \pi i (\mathrm{B} \xi)_j x_{\Theta,j}$ we may conclude. \fin

\subsubsection{A $\Theta$-deformation of $\partial_\xi$}

We finish this section with another local operator acting on a given symbol $a: \R^n \to \RR_\Theta$. It plays a crucial role in the H\"ormander classes $\Sigma_{\rho,\delta}^m(\RR_\Theta)$ from the Introduction. The \emph{mixed classical-quantized derivative} is given by
$$\partial_{\Theta, \xi}^j a (\xi) = \lambda_\Theta(\xi)^\ast \partial_\xi^j \big\{ \lambda_\Theta(\xi) a(\xi) \lambda_\Theta(\xi)^\ast \big\} \lambda_\Theta(\xi).$$

\begin{lemma}
We have $$\partial_{\Theta, \xi}^j a (\xi) = \partial_\xi^j a(\xi) + 2 \pi i \big[ x_{\Theta,j}, a(\xi) \big].$$
\end{lemma}

\dem Note that 
\begin{eqnarray*}
\frac{d}{ds}{\Big|}_{s = 0} \lambda_\Theta(\xi + s e_j) & = & \hskip2pt \lambda_\Theta(\xi) \Big( \frac{d}{ds} {\Big|}_{s = 0} e^{-2 \pi i s \langle \xi, \Theta_\downarrow e_j \rangle} \lambda_\Theta(s e_j) \Big) \hskip3pt, \\
\frac{d}{ds}{\Big|}_{s = 0} \lambda_\Theta(\xi + s e_j)^* \hskip-4.8pt & = & \Big( \frac{d}{ds} {\Big|}_{s = 0} e^{2 \pi i s \langle \xi, \Theta_\downarrow e_j \rangle} \lambda_\Theta(-s e_j) \Big) \lambda_\Theta(\xi)^*.
\end{eqnarray*}
A simple calculation then yields that
\begin{eqnarray*}
\partial_{\Theta, \xi}^j a (\xi)
\!\!\! & = & \!\!\! \lambda_\Theta(\xi)^\ast \bigg( \frac{d}{ds} {\Big|}_{s = 0} \lambda_\Theta(\xi + s e_j) a(\xi + s e_j) \lambda_\Theta(\xi + s e_j)^\ast \bigg) \lambda_\Theta(\xi) \\
\!\!\! & = & \!\!\! \big( 2 \pi i x_{\Theta,j} - 2 \pi i \langle \xi, \Theta_\downarrow e_j \rangle \big) a(\xi) + \partial_\xi^j a(\xi)
+ a(\xi) \big( 2 \pi i \langle \xi, \Theta_\downarrow e_j \rangle - 2 \pi i x_{\Theta,j} \big).
\end{eqnarray*}
Eliminating vanishing terms and rearranging gives the desired identity. \fin

The commutator vanishes in the Euclidean setting $\Theta=0$. Therefore, we should understand $\partial_{\Theta,\xi}$ as a $\Theta$-deformation of the classical derivative $\partial_\xi$ which |as it is indicated by the result below| is also very much related to the quantum derivatives $\partial_\Theta$. Thus, we get a $\Theta$-deformation of $\partial_\xi$ by $\partial_\Theta$. 

\begin{lemma} \label{Mx.comIden}
Given $\varphi \in \S_\Theta \subset \RR_\Theta$ we have 
$$\big[x_{\Theta,j}, \varphi \big] = \frac{1}{2\pi i} \sum_{k=1}^n \Theta_{jk} \hskip1pt \partial_\Theta^k \varphi.$$
In particular, we obtain the following estimate $$\big\| \big[ x_{\Theta,j}, \varphi \big] \big\|_{\RR_\Theta} \leq \frac{1}{2\pi} \Big( \sum_{k=1}^n |\Theta_{jk}|^2 \Big)^\frac12 \Big\| \Big( \sum_{k=1}^n |\partial_\Theta^k \varphi|^2 \Big)^\frac12 \Big\|_{\RR_\Theta}.$$
\end{lemma}

\dem Observe that
$$\big[ x_{\Theta,j},\varphi \big] = \int_{\R^n} \widehat{\varphi}(\xi) \big[ x_{\Theta,j}, \lambda_\Theta(\xi) \big] \, d\xi = \int_{\R^n} \widehat{\varphi}(\xi) \frac{1}{2 \pi i} \frac{d}{ds} {\bigg|}_{s = 0} \big[ \lambda_\Theta(s e_j), \lambda_\Theta(\xi) \big] \, d\xi.$$
Using $[\lambda_\Theta(s e_j), \lambda_\Theta(\xi)] = \lambda_\Theta(s e_j) \lambda_\Theta(\xi) (1 - e^{2 \pi i s \langle \xi, \Theta e_j \rangle})$ and applying the Leibniz rule, we easily deduce the first assertion. The second one is straightforward. \fin

The last estimate above provides a uniform and linear bound that explicitly gives the convergence $\partial_{\Theta, \xi}^j \to \partial_{\xi}^j$ in the point-operator norm when $\Theta \to 0^+$. We also recall that the assumption $\varphi \in \S_\Theta$ is just needed a priori and can be extended to any $\varphi$ in the
weak-$\ast$ closed domain of $\nabla_\Theta$.

\section{{\bf Calder\'on-Zygmund $L_p$ theory}}
\label{Sect2}

We now introduce kernel representations of linear operators acting on a von Neumann algebra and develop a very satisfactory Calder\'on-Zygmund theory over quantum Euclidean spaces. We shall prove $L_p$-boundedness from $L_2$-boundedness and Calder\'on-Zygmund conditions for the kernel. $L_2$-boundedness will be analyzed below in this paper for pseudodifferential operators. Our kernel conditions are given in terms of the intrinsic metric and gradient introduced above and resemble very neatly the classical conditions. This is the first form of Calder\'on-Zygmund theory over a fully noncommutative von Neumann algebra. We refer to \cite{JMP,Pa1} for related results over tensor product and crossed product algebras containing an abelian factor. An algebraic/probabilistic approach |lacking the geometric aspects of the present one| will be presented in \cite{JMP2} for more general von Neumann algebras.

\subsection{Kernels and symbols}
\label{SectKernels}

Given a measure space $(\Omega,\mu)$ and a linear map $T$ acting on certain function space X over $\Omega$, a kernel representation has the following form for functions $f$ living in some dense domain in X $$T_kf(x) = \int_\Omega k(x,y) f(y) \, d\mu(y),$$ where the kernel $k: \Omega \times \Omega \to \C$ is only assumed a priori to be defined almost everywhere and measurable. Now, given a noncommutative measure space $(\M,\tau)$ composed by a semifinite von Neumann algebra $\M$ and normal faithful semifinite trace $\tau$, the kernel representation takes the analogous form $$T_k \varphi = \big( id \otimes \tau \big) \big( k ( \1 \otimes \varphi ) \big) = \big( id \otimes \tau \big) \big( ( \1 \otimes \varphi ) k \big)$$ with the only difference that $k$ is now an operator affiliated to $\M \bar\otimes \M_{\mathrm{op}}$, instead of $\M \bar\otimes \M$ as one could have expected. This novelty |undistinguishable in the abelian case, where $\M = \M_\mathrm{op}$| is crucial to develop a consistent theory. Let us begin by showing the fundamental properties of these kernel representations. This will simplify the task of justifying our choice of $\M \bar\otimes \M_{\mathrm{op}}$. Recall the products $\cdot$ and $\bullet$ in $\M_\mathrm{op}$ and $\M \bar\otimes \M_\mathrm{op}$ respectively from Section \ref{MetricandBMO} above. 

\begin{remark} \label{MvsMop}
\emph{Rigorously speaking, the map $T_k$ so defined should send operators $\varphi$ affiliated to $\M_\mathrm{op}$ to another operator $T_k\varphi$ affiliated to $\M$. Of course, this is not an obstruction since the set of affiliated operators coincides for $\M$ and $\M_\mathrm{op}$ and $\tau_\M = \tau_{\M_\mathrm{op}}$. We will regard $\varphi$ as affiliated to $\M$, so that $T_k$ becomes a linear map over the noncommutative measure space $(\M,\tau)$.}
\end{remark} 

\begin{lemma} \label{KernelL}
The following properties hold$\hskip1pt :$
\begin{itemize}
\item[\emph{i)}] Adjoints and composition $$\hskip15pt T_k^* = T_{\mathrm{flip}(k)^*} \quad \mbox{with} \quad \mathrm{flip}(a \otimes b) = b \otimes a,$$ $$\hskip20pt T_{k_1} \circ T_{k_2} = T_k \quad \mbox{with} \quad k = \big( id \otimes \tau \otimes id \big) \Big( \big( k_1 \otimes \1 \big) \big( \1 \otimes k_2 \big) \Big).$$

\item[\emph{ii)}] Schur lemma and factorization 
\begin{eqnarray*}
T_{k_\mathrm{A} \bullet k_\mathrm{B}}\varphi \!\!\! & = & \!\!\! \big( id \otimes \tau \big) \big( k_\mathrm{A} \big( \1 \otimes \varphi \big) k_\mathrm{B} \big), \\ \hskip40pt \big\| T_{k_\mathrm{A} \bullet k_\mathrm{B}}: L_2(\M) \to L_2(\M) \big\| \!\!\! & \le & \!\!\! \big\| (id \otimes \tau) (k_\mathrm{A}^{\null} k_\mathrm{A}^*) \big\|_\M^\frac12 \big\| (\tau \otimes id) (k_\mathrm{B}^{\null} k_\mathrm{B}^*) \big\|_\M^\frac12, \\ \big\| T_{k_\mathrm{A} \bullet k_\mathrm{B}}: L_\infty(\M) \to L_\infty(\M) \big\| \!\!\! & \le & \!\!\! \big\| (id \otimes \tau) (k_\mathrm{A}^{\null} k_\mathrm{A}^*) \big\|_\M^\frac12 \big\| (id \otimes \tau) (k_\mathrm{B}^* k_\mathrm{B}^{\null} ) \big\|_\M^\frac12.
\end{eqnarray*}

\vskip3pt

\item[\emph{iii)}] $T_k$ is completely positive if and only if $k$ is positive as affiliated to $\M \bar\otimes \M_\mathrm{op}$.
\end{itemize}
\end{lemma}

\dem We have 
\begin{eqnarray*}
\big\langle T_k\varphi_1, \varphi_2 \big\rangle & = & \tau \Big( (id \otimes \tau) \big( k (\1 \otimes \varphi_1) \big)^* \varphi_2 \Big) \\ & = & \tau \Big( (\tau \otimes id) \big( (\1 \otimes \varphi_1^*) k^* (\varphi_2 \otimes \1) \big) \Big) \\ & = & \tau \Big( \varphi_1^* (id \otimes \tau) \big( \mathrm{flip}(k)^* (\1 \otimes \varphi_2) \big) \Big) \ = \ \big\langle \varphi_1, T_k^*\varphi_2 \big\rangle,
\end{eqnarray*}
which proves the kernel formula for the adjoint. Regarding the composition   
\begin{eqnarray*}
T_{k_1} (T_{k_2} \varphi) & = & (id \otimes \tau) \Big( k_1 \big( \1 \otimes (id \otimes \tau)(k_2(\1 \otimes \varphi)) \Big) \\
& = & (id \otimes \tau \otimes \tau) \Big( (k_1 \otimes \1) (\1 \otimes k_2) (\1 \otimes \1 \otimes \varphi) \Big) \\
& = & (id \otimes \tau) \Big( \big[(id \otimes \tau \otimes id) (k_1 \otimes \1) (\1 \otimes k_2) \big] (\1 \otimes \varphi) \Big).
\end{eqnarray*}
In both cases |adjoints and compositions| we have regarded one more time the involved operators as affiliated to $\M$ or $\M_\mathrm{op}$ according to the context, as we explain in Remark \ref{MvsMop}. Next, the factorization identity in ii) uses in a fundamental way the product $\bullet$ in $\M \bar\otimes \M_\mathrm{op}$ $$T_{k_\mathrm{A} \bullet k_\mathrm{B}}\varphi = \big( id \otimes \tau \big) \big( k_\mathrm{A} \bullet k_\mathrm{B} \big( \1 \otimes \varphi \big) \big) = \big( id \otimes \tau \big) \big( k_\mathrm{A} \big( \1 \otimes \varphi \big) k_\mathrm{B} \big).$$ Namely, in the last identity above the first coordinate remains unchanged since $\1 \otimes \varphi$ does not affect the product in $\M \otimes \1$, whereas the second coordinate in $\1 \otimes \M_\mathrm{op}$ is explained using its product $\cdot$ as follows $$\tau(\alpha \cdot \beta \cdot \varphi) = \tau (\alpha \cdot \varphi \beta) = \tau (\alpha \varphi \beta).$$ Let us now prove the announced inequalities. By the Cauchy-Schwarz inequality for the operator-valued inner product $(x,y) = (id \otimes \tau)(x^*y)$ over the von Neumann algebra $\M \bar\otimes \M_\mathrm{op}$ we note that $$\big| (id \otimes \tau)(x^*y) \big|^2 \le \big\| (id \otimes \tau) (x^*x) \big\|_\M (id \otimes \tau)(y^*y),$$ see for instance \cite[Proposition 1.1]{La}. In particular 
\begin{eqnarray*}
\big| T_{k_\mathrm{A} \bullet k_\mathrm{B}}\varphi \big|^2 & = & \Big| \big( id \otimes \tau \big) \big( k_\mathrm{A} (\1 \otimes \varphi) k_\mathrm{B} \big)\Big|^2 \\ & \le & \big\| (id \otimes \tau) (k_\mathrm{A} k_\mathrm{A}^*) \big\|_\M \big( id \otimes \tau \big) \Big( k_\mathrm{B}^* (\1 \otimes \varphi^*\varphi) k_\mathrm{B} \Big).
\end{eqnarray*}
The $L_\infty$-estimate announced for $T_{k_\mathrm{A} \bullet k_\mathrm{B}} \varphi$ follows immediately from the inequality above. In order to prove the $L_2$-estimate, pick a unit vector $\varphi \in L_2(\M)$. We just need to take the trace and apply Fubini
\begin{eqnarray*}
\tau \big( \big| T_{k_\mathrm{A} \bullet k_\mathrm{B}} \varphi \big|^2 \big) & \le & \big\| (id \otimes \tau) (k_\mathrm{A} k_\mathrm{A}^*) \big\|_\M \big( \tau \otimes \tau \big) \Big( k_\mathrm{B}^* \bullet k_\mathrm{B} (\1 \otimes \varphi^*\varphi) \Big) \\ & = & \big\| (id \otimes \tau) (k_\mathrm{A} k_\mathrm{A}^*) \big\|_\M \hskip1pt \tau \Big( (\tau \otimes id) (k_\mathrm{B}^* \bullet k_\mathrm{B}) \varphi^*\varphi \Big) \\ & = & \big\| (id \otimes \tau) (k_\mathrm{A} k_\mathrm{A}^*) \big\|_\M \hskip1pt \tau \Big( (\tau \otimes id) (k_\mathrm{B} k_\mathrm{B}^*) \varphi^*\varphi \Big) \\ [3pt] & \le & \big\| (id \otimes \tau) (k_\mathrm{A} k_\mathrm{A}^*) \big\|_\M \big\| (\tau \otimes id) (k_\mathrm{B} k_\mathrm{B}^*) \big\|_\M. 
\end{eqnarray*}
It remains to prove the last assertion iii). As an operator affiliated to $\M \bar\otimes \M_\mathrm{op}$, the kernel $k$ is positive iff there exists $\kappa$ also affiliated to $\M \bar\otimes \M_\mathrm{op}$ so that $k = \kappa^* \bullet \kappa$ and the factorization identity above gives in that case $$T_kf = \big( id \otimes \tau \big) \big( \kappa^* (\1 \otimes \varphi) \kappa \big)$$ which is clearly a completely positive map. Reciprocally, let $T_k$ be completely positive. Assume for simplicity that $T_k$ is well-defined over projections in $\M$. Then, given any pair of projections $p,q \in \M$ we know from positivity of $T_k$ that $\tau(T_k(q)p) \ge 0$. However
\begin{eqnarray*}
\tau(T_k(q)p) & = & \tau \Big( (id \otimes \tau) \big( k (\1 \otimes q) \big) p \Big) \\ & = & \tau \Big( (id \otimes \tau) \big( k \bullet (\1 \otimes q) \big) p \Big) \\ & = & \tau \Big( (id \otimes \tau) \big( k \bullet (p \otimes q) \big) \Big) \ = \ \tau \otimes \tau \big( k \bullet (p \otimes q) \big).
\end{eqnarray*}
The positivity of the last term for arbitrary projections implies the assertion. \fin

\begin{remark}
\emph{Lemma \ref{KernelL} i) also holds for kernels affiliated to $\M \bar\otimes \M$, contrary to points ii) and iii). On the other hand, recall that the norm of a completely positive map is determined by its value at $\1$. The $L_\infty$-estimate in Lemma \ref{KernelL} ii) rephrases it in terms of kernels when $k \ge 0$ and $k_\mathrm{A} = k_\mathrm{B} = \sqrt{k}$, so this estimate provides a generalization for nonpositive maps. Also, the $L_2$-estimate generalizes a classical result for integral kernels known as Schur lemma \cite[Lemma in page 284]{St}. Finally, the use of kernels $k$ affiliated to $\M \bar\otimes \M_\mathrm{op}$ |essential for the properties in Lemma \ref{KernelL} ii) and iii) above| is consistent with the duality} $$L_1(\M)^* = \M_\mathrm{op}$$ \emph{via the pairing $\langle x, y \rangle = \tau(xy)$, we refer to Pisier's book \cite{P3} for further details.}
\end{remark}

\begin{remark}
\emph{Ignoring for the moment more general kernels which will arise as tempered distributions, let us assume that $k$ is affiliated to $\RR_\Theta \bar\otimes \RR_\Theta^{\mathrm{op}}$ and admits an expression $$k = \int_{\R^n \times \R^n} \widehat{k}(\xi,\eta) \lambda_\Theta(\xi) \otimes \lambda_\Theta(\eta)^* \, d\mu(\xi,\eta)$$ for some measure $\mu$ on $\R^n \times \R^n$. Noticing that $\tau_\Theta(\lambda_\Theta(f) \lambda_\Theta(\xi)^*) = f(\xi)$ for $f$ smooth we may interpret the kernel $k$ as a bilinear form where |regardless $\lambda_\Theta(\xi), \lambda_\Theta(\eta)$ are not in $L_2(\RR_\Theta)$| we put $\lambda_\Theta(\xi) \otimes \lambda_\Theta(\eta)^* \sim | \lambda_\Theta(\xi) \rangle \langle \lambda_\Theta(\eta)|$ following the bra-ket notation. This is easily checked for Dirac measures $\mu = \delta_{\xi_0,\eta_0}$ $$\big\langle T_k(\lambda_\Theta(f_1)), \lambda_\Theta(f_2) \big\rangle = \tau_\Theta \big( T_k(\lambda_\Theta(f_1))^* \lambda_\Theta(f_2) \big) = \widehat{k}(\xi_0,\eta_0) \overline{f_1(\eta_0)} f_2(\xi_0).$$ We will approximate general measures  as limits of finite sums of Dirac measures.}
\end{remark} 


This paper is devoted to investigate singular integral operators over quantum Euclidean spaces, both in terms of Calder\'on-Zygmund conditions for the kernel and H\"ormander smoothness for the associated symbol. Let us therefore briefly describe the kernels and symbols we will be working with. There exists a very well-known relation between kernels and symbols of classical pseudodifferential operators, the reader can look for instance in \cite{St,T1} or almost any textbook on pseudodifferential operators. Indeed, given 
$$\Psi_af(x) = \int_{\R^n} a(x,\xi) \widehat{f}(\xi) e^{2\pi i \langle x,\xi \rangle} \, d\xi,$$
it turns out that $\Psi_af = T_kf$ for $$k = (id \otimes \mathcal{F}^{-1})(a)(x,x-y) = \int_{\R^n} a(x,\xi) e^{2\pi i \langle x-y,\xi \rangle} \, d\xi.$$ Given $n \ge 2$, let us know consider a $n \times n$ deformation $\Theta$. As explained in the Introduction, noncommutative symbols over quantum Euclidean spaces are smooth functions $a: \R^n \to \S_\Theta$ and pseudodifferential operators look like  
\begin{eqnarray*}
\Psi_a(\lambda_\Theta(f)) & = & \int_{\R^n} a(\xi) f(\xi) \lambda_\Theta(\xi) \, d\xi \\ & = & \int_{\R^n} a(\xi) \tau_\Theta \big( \lambda_\Theta(f) \lambda_\Theta(\xi)^* \big) \lambda_\Theta(\xi) \, d\xi \\ & = & \big(id \otimes \tau_\Theta \big) \Big[ \Big( \underbrace{\int_{\R^n} (a(\xi) \otimes \1) (\lambda_\Theta(\xi) \otimes \lambda_\Theta(\xi)^*) d\xi}_{\mathrm{The \ kernel \ } k} \Big) \big( \1 \otimes \lambda_\Theta(f) \big) \Big].
\end{eqnarray*}
Thus, we find formally that $\Psi_a(\lambda_\Theta(f)) = T_k(\lambda_\Theta(f))$ for 
$$k = \int_{\R^n \times \R^n} \widehat{k}(\xi,\eta) \lambda_\Theta(\xi) \otimes \lambda_\Theta(\eta)^* \, d\mu(\xi,\eta) = \int_{\R^n} a(\xi) \lambda_\Theta(\xi) \otimes \lambda_\Theta(\xi)^* d\xi.$$
Reciprocally, we also have the following expression for $a$ $$a(\xi) = \int_{\R^n} \widehat{k}(\xi,\eta) \big( \1 \otimes \lambda_\Theta(\eta)^* \lambda_\Theta(\xi) \big) d\eta$$ when $\mu$ is the Lebesgue measure in $\R^n \times \R^n$. The algebra of pseudodifferential operators is formally generated by the derivatives $\partial_{\Theta}^j$ and the left multiplication maps $\lambda_\Theta(f) \mapsto x_{\Theta,j} \lambda_\Theta(f)$. 
\begin{remark}
\emph{Given $a:\mathbb{R}^n \to \RR_\Theta$ define $$\Psi^\mathrm{right}_a(\lambda_\Theta(f)) = \int_{\R^n} f(\xi) \lambda_\Theta(\xi) a(\xi) \, d\xi.$$ Both theories for right and left pseudodifferential operators are analogous. In fact, if we denote by $S$ the bounded operator in $L_p(\RR_\Theta)$ given by extension of $S(x) = x^\ast$ we have that $$S \Psi_a^{\mathrm{left}} = \Psi^\mathrm{right}_{b} S \quad \mbox{where} \quad b(\xi) = a(-\xi)^\ast.$$ The proof of such intertwining identity amounts to a straightforward calculation.}
\end{remark}

\subsection{CZ extrapolation: Model case}

We are ready to prove $L_p$-boundedness of operators associated to elementary kernels satisfying cancellation and smoothness conditions of Calder\'on-Zygmund type. Our kernels will belong along this paragraph to $\S_\Theta \otimes_\mathrm{alg} \S_\Theta^{\mathrm{op}}$, so that $$k = \summ_j \int_{\R^n} \int_{\R^n} \kappa_{1j}(\xi) \kappa_{2j}(\eta) \lambda_\Theta(\xi) \otimes \lambda_\Theta(\eta)^* d\xi d\eta,$$ where the sum above is finite and $\kappa_{ij} \in \mathcal{S}(\R^n)$. We will temporarily refer to these kernels as \emph{algebraic kernels}. Of course, in this case $T_k$ is $L_p$-bounded for $1 \le p \le \infty$ with constants a priori depending on the family $\kappa_{ij}$. Our goal is to provide $L_\infty \to \mathrm{BMO}$ estimates with constants which only depend on structural properties of the whose $k$ since this will allow us to include general singular kernels below. The following result is the basic core of this paper. We shall use the quantum metric $\mathrm{d}_\Theta$ defined in Section \ref{MetricandBMO}, the notation $\nabla_\Theta k$ to denote the operator $$\big( \nabla_\Theta \otimes id_{\RR_\Theta^\mathrm{op}} \big)(k) \in \mathcal{L}(\F_n) \bar\otimes \RR_\Theta \bar\otimes \RR_\Theta^\mathrm{op}$$ for $k \in \S_\Theta \otimes_\mathrm{alg} \S_\Theta^\mathrm{op}$ and the dimensional constant $\mathrm{K}_n = \frac12 (n+1)$.

\begin{theorem} \label{CZExt1}
Let $k \in \S_\Theta \otimes_\mathrm{alg} \S_\Theta^{\mathrm{op}}$ and assume$\hskip1pt :$
\begin{itemize}
\item[\emph{i)}] \emph{Cancellation} $$\hskip10pt \big\| T_k: L_2(\RR_\Theta) \to L_2(\RR_\Theta) \big\| \le \mathrm{A}_1.$$

\item[\emph{ii)}] \emph{Kernel smoothness.} There exists $$\hskip10pt \alpha < \mathrm{K}_n - \frac12 < \beta < \mathrm{K}_n + \frac12$$ satisfying the gradient conditions below for $\rho = \alpha, \beta$ $$\hskip10pt \Big| \mathrm{d}_\Theta^{\rho} \bullet (\nabla_\Theta \otimes id_{\RR_\Theta^\mathrm{op}}) (k) \bullet \mathrm{d}_\Theta^{n+1-\rho} \Big| \le \mathrm{A}_2.$$
\end{itemize}
Then, we find the following $L_\infty \to \mathrm{BMO}_c$ estimate $$\big\| T_k: L_\infty(\RR_\Theta) \to \mathrm{BMO}_c(\RR_\Theta) \big\|_{\mathrm{cb}} \le C_n(\alpha, \beta) \big( \mathrm{A}_1 + \mathrm{A}_2 \big).$$
\end{theorem}

\dem Since $k$ is an algebraic kernel, $T_k$ is bounded on $\RR_\Theta$. Moreover, $T_k^*$ is also bounded on $L_1(\RR_\Theta^\mathrm{op})$ and $T_k$ is a normal operator. On the other hand, the weak-$*$ topology in $\RR_\Theta$ is stronger than the inherited one from the weak-$*$ topology in $\mathrm{BMO}_c(\RR_\Theta)$. Therefore, by Kaplansky density theorem, it suffices to estimate the norm of $T_k$ on the weak-$*$ dense subalgebra $\S_\Theta$. Given $\varphi = \lambda_\Theta(f) \in \S_\Theta$ $$\|T_k \varphi\|_{\mathrm{BMO}_c(\RR_\Theta)} \sim \sup_{\mathrm{R}>0} \sup_{\mathrm{B_R} \in \Q_\mathrm{R}} \Big\| \Big( \mean_{\mathrm{B_R}} \big| \sigma_\Theta(T_k\varphi) - \sigma_\Theta(T_k\varphi)_{\mathrm{B_R}} \big|^2 \, d\mu \Big)^\frac12 \Big\|_{\RR_\Theta}$$ where the second supremum runs over the set $\Q_\mathrm{R}$ of Euclidean balls $\mathrm{B_R}$ of radius $\mathrm{R}$ and center $c_{\mathrm{B_R}}$. Recall that $\sigma_\Theta(T_k\varphi) = T_{k_\sigma}\varphi$ for $k_\sigma = (\sigma_\Theta \otimes id)(k)$. Now let $\psi: \R^n \to [0,1]$ be a Schwartz function which is identically 1 over the unit ball $\mathrm{B}_1(0)$ and identically 0 outside its concentric 2-dilation $\mathrm{B}_2(0)$. Define $$\psi_\mathrm{R}(x) = \psi \Big( \frac{x + c_{\mathrm{B_R}}}{2\mathrm{R}} \Big)$$ to decompose the kernel $k_\sigma = (\sigma_\Theta \otimes id)(k)$ as follows
\[
  k_\sigma = \underbrace{k_\sigma \bullet \pi_\Theta(\psi_\mathrm{R})}_{k_{\sigma  1}(\mathrm{R})} +
  \underbrace{k_\sigma \bullet \big( \1 - \pi_\Theta(\psi_\mathrm{R}) \big)}_{k_{\sigma  2}(\mathrm{R})}.
\]
Note here that $\pi_\Theta(\psi_\mathrm{R}) \simeq \1 \otimes \pi_\Theta(\psi_\mathrm{R})$ is an element of $\RR_\Theta \bar\otimes \RR_\Theta^\mathrm{op}$ and only $k_\sigma$ has a component in $L_\infty(\R^n)$. We claim that the following inequality for $k_{\sigma 1}(\mathrm{R})$ holds up to constants independent of the ball $\mathrm{B_R} \in \Q_\mathrm{R}$ and the radius $\mathrm{R}$
\begin{equation} \label{Est1}
\Big\| \Big( \mean_{\mathrm{B_R}} \big| T_{k_{\sigma  1}(\mathrm{R})}\varphi - (T_{k_{\sigma  1}(\mathrm{R})}\varphi)_{\mathrm{B_R}} \big|^2 \, d\mu \Big)^\frac12 \Big\|_{\RR_\Theta}
\le C_n \mathrm{A}_1 \|\varphi\|_{\RR_\Theta}.
\end{equation}
Before proving this first claim, let us continue with the argument. Of course, it would suffice to give a similar estimate for $k_{\sigma  2}(\mathrm{R})$. To do so and according to the Poincar\'e type inequality in Proposition \ref{Prop-Poincare} and its relation to $\mathrm{BMO}_c(\RR_\Theta)$ outlined in Remark \ref{Poincare-BMO}, it suffices to estimate
\[
    \mathrm{R} \Big\| \underbrace{(\1 \otimes q_\mathrm{R} \otimes \1) (\nabla \otimes id_{\RR_\Theta}) (T_{k_{\sigma  2}(\mathrm{R})}\varphi)}_{T_\mathbf{K}(\varphi)} \Big\|_{\mathcal{L}(\F_n) \bar\otimes L_\infty(\R^n) \bar\otimes \RR_\Theta}
\]
for $q_\mathrm{R} = 1_{\mathrm{B_R}}$ with constants independent of R. Since $$(\nabla \otimes id_{\RR_\Theta}) \circ \sigma_\Theta = (id_{\mathcal{L}(\mathbb{F}_n)} \otimes \sigma_\Theta) \circ \nabla_\Theta$$ we may rewrite $T_\mathbf{K}(\varphi)$ as follows
\begin{eqnarray*}
T_\mathbf{K}(\varphi) \!\! & = & \!\! (\1 \otimes q_\mathrm{R} \otimes \1) (\nabla \otimes id) \Big( (id^{\otimes 2} \otimes \tau_\Theta) \big( k_{\sigma  2}(\mathrm{R}) (\1^{\otimes 2} \otimes \varphi) \big) \Big) \\ \!\! & = & \!\! (\1 \otimes q_\mathrm{R} \otimes \1) \big( id^{\otimes 3} \otimes \tau_\Theta \big) \Big( (\nabla \otimes id^{\otimes 2}) (\sigma_\Theta \otimes id) [k] \bullet (\1 - \Psi_\mathrm{R}) (\1^{\otimes 3} \otimes \varphi) \Big) \\ \!\! & = & \!\! \big( id^{\otimes 3} \otimes \tau_\Theta \big) \Big( \underbrace{(\1 \otimes q_\mathrm{R} \otimes \1^{\otimes 2}) (id \otimes \sigma_\Theta \otimes id) [\nabla_\Theta k] \bullet (\1 - \Psi_\mathrm{R})}_{\mathbf{K}} (\1^{\otimes 3} \otimes \varphi) \Big)
\end{eqnarray*}
with $\mathbf{K} \in \mathcal{L}(\F_n) \bar\otimes L_\infty(\R^n) \bar\otimes \RR_\Theta \bar\otimes \RR_\Theta^{\mathrm{op}}$ and $\Psi_\mathrm{R} = \pi_\Theta(\psi_\mathrm{R})$. For simplicity, we shall use a more compact terminology and write $\mathbf{K} = q_\mathrm{R} \nabla_\Theta^\sigma k \bullet (\1 - \Psi_\mathrm{R})$. We may now decompose $\mathbf{K}$ as follows $$\mathbf{K} = \underbrace{q_\mathrm{R} \Psi_\mathrm{R} \bullet \nabla_\Theta^\sigma k \bullet (\1 - \Psi_\mathrm{R})}_{\mathbf{K}_1} + \underbrace{q_\mathrm{R} (\1 - \Psi_\mathrm{R}) \bullet \nabla_\Theta^\sigma k \bullet (\1 - \Psi_\mathrm{R})}_{\mathbf{K}_2}.$$ We claim that the following inequality holds for $j=1,2$
\begin{equation} \label{Est2}
\big\| T_{\mathbf{K}_j} (\varphi) \big\|_{\mathcal{L}(\F_n) \bar\otimes L_\infty(\R^n) \bar\otimes \RR_\Theta} \le C_n(\alpha, \beta) \frac{\mathrm{A}_2}{\mathrm{R}} \|\varphi\|_{\RR_\Theta}.
\end{equation}
Our discussion so far has reduced the proof to justifying \eqref{Est1} and \eqref{Est2}:

\noindent \textbf{Proof of \eqref{Est1}.} It is clear that $$\Big\| \Big( \meannn_{\mathrm{B_R}} \big| T_{k_{\sigma  1}(\mathrm{R})}\varphi - (T_{k_{\sigma  1}(\mathrm{R})}\varphi)_{\mathrm{B_R}} \big|^2 \, d\mu \Big)^\frac12 \Big\|_{\RR_\Theta} \!\! \le C_n  \mathrm{R}^{-\frac{n}{2}} \Big\| \Big( \int_{\R^n} \big| T_{k_{\sigma  1}(\mathrm{R})}\varphi \big|^2 \, d\mu \Big)^\frac12 \Big\|_{\RR_\Theta} \hskip-1pt .$$ On the other hand, given $T: \S_\Theta \to L_2^c(\R^n) \bar\otimes \RR_\Theta$ we need to introduce its module extension $T': \S_\Theta \bar\otimes \RR_\Theta \ni \lambda_\Theta(f) \otimes \varphi \mapsto T(\lambda_\Theta(f)) (\1 \otimes \varphi) \in L_2^c(\R^n) \bar\otimes \RR_\Theta$. Recall the following (elementary) algebraic identity $$(id \otimes \tau_\Theta) \big( k_\sigma \bullet (a \otimes b) \bullet (\1 \otimes \varphi) \big) = T_{k_\sigma}(\varphi b) (\1 \otimes a) \ = \ T_{k_\sigma}' \big( (\varphi \otimes \1) \mathrm{flip}(a \otimes b) \big).$$ Then, noticing that $k$ is assumed to be an algebraic kernel in $\S_\Theta \otimes_{\mathrm{alg}} \S_\Theta^{\mathrm{op}}$, it is not difficult to check that the above formula extends from elementary tensors $a \otimes b$ to arbitrary elements in $\RR_\Theta \bar\otimes \RR_\Theta^\mathrm{op}$. This yields $$T_{k_{\sigma  1}(\mathrm{R})} \varphi = T_{k_\sigma \bullet \Psi_\mathrm{R}} \varphi = T_{k_\sigma}' \big( (\varphi \otimes \1) \mathrm{flip}(\pi_\Theta(\psi_\mathrm{R}))\big).$$ In particular, we easily obtain the following estimate 
\begin{eqnarray*}
\lefteqn{\Big\| \Big( \int_{\R^n} \big| T_{k_{\sigma  1}(\mathrm{R})}\varphi \big|^2 \, d\mu \Big)^\frac12 \Big\|_{\RR_\Theta}} \\ & \le & \Big\| T_{k_\sigma}': L_2^c(\RR_\Theta) \bar\otimes \RR_\Theta \to L_2^c(\R^n) \bar\otimes \RR_\Theta \Big\|_{\mathrm{cb}} \, \big\| \mathrm{flip}(\pi_\Theta(\psi_\mathrm{R})) \big\|_{L_2^c(\RR_\Theta) \bar\otimes \RR_\Theta} \, \|\varphi \|_{\RR_\Theta}.
\end{eqnarray*}
According to Lemma \ref{QDistance} $$\big\| \mathrm{flip}(\pi_\Theta(\psi_\mathrm{R})) \big\|_{L_2^c(\RR_\Theta) \bar\otimes \RR_\Theta} \le \big\| \psi_\mathrm{R} \big\|_{L_2(\R^n)} \le C_n \mathrm{R}^{\frac{n}{2}}$$ since the argument given there for $\pi_\Theta$ also holds for $\mathrm{flip} \circ \pi_\Theta$. Therefore, it remains to estimate the cb-norm of $T_{k_\sigma}'$. We claim that it is bounded by $\mathrm{A}_1$, the $L_2$-norm of $T_k$. To justify it, we introduce the map $$W: L_2^c(\R^n) \bar\otimes \RR_\Theta \ni \int_{\R^n} \exp_\xi \otimes \, a(\xi) \, d\xi \mapsto \int_{\R^n} \exp_\xi \otimes \, \lambda_\Theta(\xi) a(\xi) \, d\xi \in L_2^c(\R^n) \bar\otimes \RR_\Theta.$$ It is straightforward to show that $W$ extends to an isometry. On the other hand, letting $j_\Theta: \exp_\xi \mapsto \lambda_\Theta(\xi)$ be the $L_2$-isometry introduced in Section \ref{QSch}, we observe that 
\begin{eqnarray*}
W(j_\Theta^* \otimes id) (\lambda_\Theta(f) \otimes a) & = & W \Big( \int_{\R^n} f(\xi) \exp_\xi \otimes \, a \, d\xi \Big) \\ & = & \int_{\R^n} f(\xi) \exp_\xi \otimes \, \lambda_\Theta(\xi) a \, d\xi \ = \ \sigma_\Theta(\lambda_\Theta(f)) (\1 \otimes a).
\end{eqnarray*}
Identifying $T_k\varphi$ with $\lambda_\Theta(f)$ for some smooth $f: \R^n \to \C$ we obtain the identity $W(j_\Theta^*T_k \otimes id) (\varphi \otimes a) = T_{k_\sigma}'(\varphi \otimes a)$. This a fortiori implies that the cb-norm of the map $T_{k_\sigma}': L_2^c(\RR_\Theta) \bar\otimes \RR_\Theta \to L_2^c(\R^n) \bar\otimes \RR_\Theta$ is dominated by the $L_2(\RR_\Theta) \to L_2(\RR_\Theta)$ norm of $T_k$, as desired. This completes the proof of claim \eqref{Est1}.

\noindent \textbf{Proof of \eqref{Est2}.} Let $$\mathrm{D}_\Theta = (\sigma_\Theta \otimes id) (\mathrm{d}_\Theta)$$ and decompose the kernels $\mathbf{K}_j$ for $j=1,2$ as follows
\[
  \mathbf{K}_1 = \Big( q_\mathrm{R} \, \Psi_\mathrm{R} \bullet \mathrm{D}_\Theta^{-\alpha} \Big) \bullet
  \Big( \mathrm{D}_\Theta^{\alpha} \bullet \nabla_\Theta^\sigma k \bullet \mathrm{D}_\Theta^{\alpha'} \Big) \bullet
  \Big( \mathrm{D}_\Theta^{-\alpha'} \bullet (\1 - \Psi_\mathrm{R}) q_\mathrm{R} \Big),
\]
\[
  \mathbf{K}_2 = \Big( q_\mathrm{R} \, (\1 - \Psi_\mathrm{R}) \bullet \mathrm{D}_\Theta^{-\beta} \Big) \bullet
  \Big( \mathrm{D}_\Theta^{\beta} \bullet \nabla_\Theta^\sigma k \bullet \mathrm{D}_\Theta^{\beta'} \Big) \bullet
  \Big( \mathrm{D}_\Theta^{-\beta'} \bullet (\1 - \Psi_\mathrm{R}) q_\mathrm{R} \Big),
\]
with $\alpha' = n+1-\alpha$ and $\beta' = n+1-\beta$. Using the terminology $\mathbf{K}_j = \mathbf{A}_j \bullet \mathbf{B}_j \bullet \mathbf{C}_j$ for the brackets above and according to the operator-valued (trivial) extension of Lemma \ref{KernelL}, we find for $\B_\Theta = \mathcal{L}(\mathbb{F}_n) \bar\otimes L_\infty(\R^n) \bar\otimes \RR_\Theta$ that $$\big\| T_{\mathbf{K}_j}: \RR_\Theta \to \B_\Theta \big\| \le \big\| (id \otimes \tau_\Theta)(\mathbf{A}_j \mathbf{A}_j^*) \big\|_{\B_\Theta}^\frac12 \big\| \mathbf{B}_j \big\|_{\B_\Theta \bar\otimes \RR_\Theta^{\mathrm{op}}} \big\| (id \otimes \tau_\Theta)(\mathbf{C}_j^* \mathbf{C}_j) \big\|_{\B_\Theta}^\frac12.$$ Since $\mathbf{B}_j = (\sigma_\Theta \otimes id) (\mathrm{d}_\Theta^{\rho} \bullet \nabla_\Theta k \bullet \mathrm{d}_\Theta^{n+1-\rho})$ for $\rho = \alpha, \beta$ and $\sigma_\Theta$ is a $*$-homomorphism we deduce from the hypotheses that $\|\mathbf{B}_1\| + \|\mathbf{B}_2\| \le \mathrm{A}_2$. Therefore, recalling that we have $\alpha + \alpha' = \beta + \beta' = n+1$, it suffices to prove the following inequalities for the terms associated to $\mathbf{A}_j$ and $\mathbf{C}_j$ 
\begin{eqnarray*}
\big\| (id \otimes \tau_\Theta)(\mathbf{A}_1 \mathbf{A}_1^*) \big\|_{\B_\Theta}^\frac12 & \le & C_n(\hskip1.5pt \alpha \hskip1.5pt) \mathrm{R}^{\frac{n}{2} - \alpha}, \\ 
\big\| \hskip1pt (id \otimes \tau_\Theta)(\mathbf{C}_1^* \mathbf{C}_1) \big\|_{\B_\Theta}^\frac12 & \le & C_n(\alpha') \mathrm{R}^{\frac{n}{2} - \alpha'}, \\ 
\big\| (id \otimes \tau_\Theta)(\mathbf{A}_2 \mathbf{A}_2^*) \big\|_{\B_\Theta}^\frac12 & \le & C_n(\hskip1.5pt \beta \hskip1.5pt) \mathrm{R}^{\frac{n}{2} - \beta}, \\ 
\big\| \hskip1pt (id \otimes \tau_\Theta)(\mathbf{C}_2^* \mathbf{C}_2) \big\|_{\B_\Theta}^\frac12 & \le & C_n(\beta') \mathrm{R}^{\frac{n}{2} - \beta'}, \\ 
\end{eqnarray*}
for any $\alpha < \mathrm{K}_n - \frac12 < \beta < \mathrm{K}_n + \frac12$. We will justify the above estimates for $\mathbf{A}_1$ and $\mathbf{C}_1$, the proof of the others is very similar. Let $\Delta_{\R^n} f (x,y) = f(x+y)$ be the comultiplication map on $\R^n$. According to Lemma \ref{QDistance} $$(\sigma_\Theta \otimes id) \circ \pi_\Theta = (id \otimes \pi_\Theta) \circ \Delta_{\R^n}.$$ In particular, $\mathbf{A}_1^* = (id \otimes \pi_\Theta) \big( \Delta_{\R^n} ( | \cdot |^{-\alpha} ) q_\mathrm{R} (\1 \otimes \psi_\mathrm{R}) \big)$ and we find 
\begin{eqnarray*}
  \big\| (id \otimes \tau_\Theta)(\mathbf{A}_1 \mathbf{A}_1^*) \big\|_{\B_\Theta}^\frac12 & = & \big\| \mathbf{A}_1^* \big\|_{L_\infty(\R^n) \bar\otimes \RR_\Theta \bar\otimes L_2^c(\RR_\Theta^{\mathrm{op}})} \\
  & \le & \big\| id \otimes \pi_\Theta \big\| \, \big\| \Delta_{\R^n} ( | \cdot |^{-\alpha} ) q_\mathrm{R} (\1 \otimes \psi_\mathrm{R}) \big\|_{L_\infty(\R^n) \bar\otimes L_2^c(\R^n)},
\end{eqnarray*}
where $id \otimes \pi_\Theta: L_\infty(\R^n) \bar\otimes L_2^c(\R^n) \to L_\infty(\R^n) \bar\otimes \RR_\Theta \bar\otimes L_2^c(\RR_\Theta^{\mathrm{op}})$. Its norm is dominated by the cb-norm of $\pi_\Theta: L_2^c(\R^n) \to \RR_\Theta \bar\otimes L_2^c(\RR_\Theta^{\mathrm{op}})$. We already proved in Lemma \ref{QDistance} the cb-contractivity of $\pi_\Theta: L_2^c(\R^n) \to L_2^c(\RR_\Theta) \bar\otimes \RR_\Theta^{\mathrm{op}}$ and the exact same argument can be trivially adapted to show that $\|id \otimes \pi_\Theta\| \le 1$ in the right hand side of the above inequality. It then remains to estimate 
$$\sup_{x \in \R^n} q_\mathrm{R}(x) \Big( \int_{\R^n} |x+y|^{-2 \alpha} \psi_{\mathrm{R}}^2(y) \, dy \Big)^\frac12 \le \Big( \int_{\mathrm{B}_{5\mathrm{R}}(0)} |y|^{-2\alpha} dy \Big)^\frac12 \lesssim \mathrm{R}^{\frac{n}{2} - \alpha}$$ for $\alpha < \frac{n}{2}$. In the case of $\mathbf{C}_1 = (id \otimes \pi_\Theta) \big( \Delta_{\R^n} (| \cdot |^{-\alpha'}) q_\mathrm{R} \otimes (1 - \psi_\mathrm{R}) \big)$ we get 
$$\sup_{x \in \R^n} q_\mathrm{R}(x) \Big( \int_{\R^n} |x+y|^{-2 \alpha'} (1-\psi_{\mathrm{R}}(y))^2 \, dy \Big)^\frac12 \le \Big( \int_{\mathrm{B}_{\mathrm{R}}^\mathrm{c}(0)} |y|^{-2\alpha'} dy \Big)^\frac12 \lesssim \mathrm{R}^{\frac{n}{2} - \alpha'}$$
for $\alpha' > \frac{n}{2}$. The same argument applies for $\mathbf{A}_2$ and $\mathbf{C}_2$. This proves claim \eqref{Est2}. 

\noindent \textbf{Conclusion.} The argument above proves that $T_k: L_\infty(\RR_\Theta) \to \mathrm{BMO}_c(\RR_\Theta)$ defines a bounded operator with norm dominated by $C_n(\alpha, \beta) (\mathrm{A}_1 + \mathrm{A}_2)$. The exact same argument can be used after matrix amplification to prove that the cb-norm of $T_k$ satisfies the same upper bound. This completes the proof of the theorem. \fin

\begin{remark}
\emph{Our kernel conditions in Theorem \ref{CZExt1} are natural extensions of the classical ones \cite{G2,St}. The price for noncommutativity is a concrete and balanced left/right location of the exponents $\alpha, \beta$, which of course is meaningless in the commutative setting. These surprisingly transparent Calder\'on-Zygmund kernel conditions are possible due to the very precise geometric information on $\RR_\Theta$ that we collected in Section \ref{Prelim}. Our results below for pseudodifferential operators crucially rely on these conditions. For more general von Neumann algebras, the resulting conditions are necessarily less transparent \cite{JMP2}.}
\end{remark}

\subsection{CZ extrapolation: General case} \label{CZExtGral}

If $\S_\Theta'$ denotes the space of continuous linear functionals on $\S_\Theta$ |tempered distributions| the aim of this section is to generalize our $L_\infty \to \mathrm{BMO}_c$ estimate in Theorem \ref{CZExt1} for continuous linear maps $T \in \mathcal{L}(\S_\Theta, \S_\Theta')$ to incorporate actual Calder\'on-Zygmund kernels. This imposes a careful analysis of tempered $\Theta$-distributions and how this affects our former kernel manipulations. By symmetrization, interpolation and duality, we shall conclude by proving the $L_p$-boundedness of quantum Calder\'on-Zygmund operators. 

\subsubsection{Tempered distributions}

The Schwartz class $\S(\R^n)$ comes equipped with the locally convex topology determined by the seminorms $p_{\alpha,\beta}(f) = \sup_{x \in \R^n} | x^\alpha \partial_x^\beta f(x) |$ for $f: \R^n \to \C$ and all $\alpha, \beta \in \N^n$. The quantum analogue for $\varphi = \lambda_\Theta(f)$ in $\S_\Theta$ was described in Remark \ref{RemQuantumx}. In particular, a sequence $\varphi_j = \lambda_\Theta(f_j)$ converges to 0 in $\S_\Theta$ as $j \to \infty$ when $$\lim_{j \to \infty} \big\| P(x_\Theta) \partial_\Theta^\beta \big( \lambda_\Theta(f_j) \big) Q(x_\Theta) \big\|_{\RR_\Theta} = 0$$ for all $\beta \in \N^n$ and all quantum monomials $$P(x_\Theta) = \prod_{1 \le r \le m_\ell}^{\rightarrow} x_{\Theta,j_r} \quad \mbox{and} \quad Q(x_\Theta) = \prod_{1 \le s \le m_r}^{\rightarrow} x_{\Theta,k_s}.$$ By Remark \ref{RemQuantumx}, this holds iff $$\lim_{j \to \infty} \Big\| \lambda_\Theta \Big[ \Big( \prod_{1 \le r \le m_\ell}^{\rightarrow} D_{\Theta,j_r}^\ell \Big) \Big( \prod_{1 \le s \le m_r}^{\rightarrow} D_{\Theta, k_s}^r \Big) M_{(2\pi i \xi)^\beta} f_j \Big] \Big\|_{\RR_\Theta} = 0.$$ 

\begin{lemma} \label{ContinuitySchwartz}
If $\Theta_1, \Theta_2 \in A_n(\R)$ we find that $$\lim_{j \to \infty} \lambda_{\Theta_1}(f_j) = 0 \mbox{ in $\S_{\Theta_1}$} \ \Leftrightarrow \ \lim_{j \to \infty} \lambda_{\Theta_2}(f_j) = 0 \mbox{ in $\S_{\Theta_2}$}.$$
\end{lemma}

\dem Let us set $$p_{P,Q,\beta}^{\Theta} \big( \lambda_{\Theta}(f) \big) = \big\| P(x_\Theta) \partial_\Theta^\beta \big( \lambda_\Theta(f) \big) Q(x_\Theta) \big\|_{\RR_\Theta}$$ and assume that $\lim_{j \to \infty} p_{P_1,Q_1,\beta_1}^{\Theta_1} \big( \lambda_{\Theta_1}(f_j) \big) = 0$ for all $\beta_1 \in \N^n$ and all quantum $\Theta_1$-monomials $P_1, Q_1$. Given $\beta_2 \in \N^n$ and quantum $\Theta_2$-monomials $P_2, Q_2$ it then suffices to show that $$\lim_{j \to \infty} p_{P_2,Q_2,\beta_2}^{\Theta_2} \big( \lambda_{\Theta_2}(f_j) \big) = 0.$$ According to Remark \ref{RemQuantumx}, we may find two commuting operators $P_2[D_{\Theta_2,j}^\ell]$ and $Q_2[D_{\Theta_2,j}^r]$ |sums of modulations and derivations| satisfying the following identity $$p_{P_2, Q_2, \beta_2}^{\Theta_2}(\lambda_{\Theta_2}(f_j)) = \Big\| \lambda_{\Theta_2} \big( P_2[D_{\Theta_2,j}^\ell] Q_2[D_{\Theta_2,j}^r] M_{(2\pi i \xi)^{\beta_2}} f_j \big) \Big\|_{\RR_{\Theta_2}}.$$ Applying Riemann-Lebesgue and Cauchy-Schwartz with $B(\xi) = 1 + |\xi|^n$ and $B^{-1}(\xi)$
\begin{eqnarray*}
p_{P_2, Q_2, \beta_2}^{\Theta_2}(\lambda_{\Theta_2}(f_j)) \!\! & \le & \!\! \Big\| P_2[D_{\Theta_2,j}^\ell] Q_2[D_{\Theta_2,j}^r] M_{(2\pi i \xi)^{\beta_2}} f_j \Big\|_{L_1(\R^n)} \\ \!\! & \lesssim & \!\! \Big\| M_{B(\xi)} P_2[D_{\Theta_2,j}^\ell] Q_2[D_{\Theta_2,j}^r] M_{(2\pi i \xi)^{\beta_2}} f_j \Big\|_{L_2(\R^n)} \\ \!\! & = & \!\! \Big\| \lambda_{\Theta_1} \big( M_{B(\xi)} P_2[D_{\Theta_2,j}^\ell] Q_2[D_{\Theta_2,j}^r] M_{(2\pi i \xi)^{\beta_2}} f_j \big) \Big\|_{L_2(\RR_{\Theta_1})}
\end{eqnarray*}
By Proposition \ref{quantumSchwarzt1}, there exists a quantum $\Theta_1$-polynomial $R_1(x_{\Theta_1})$, whose inverse lives in $L_2(\RR_{\Theta_1})$. Therefore, multiplying and dividing by it on the left hand side yields   
\begin{eqnarray*}
p_{P_2, Q_2, \beta_2}^{\Theta_2}(\lambda_{\Theta_2}(f_j)) \!\! & \lesssim & \!\! \Big\| R_1(x_{\Theta_1}) \lambda_{\Theta_1} \big( M_{B(\xi)} \mathbf{D}_{\Theta_2}(P_2) \mathbf{D}_{\Theta_2}(Q_2) M_{(2\pi i \xi)^{\beta_2}} f_j \big) \Big\|_{\RR_{\Theta_1}} \\ \!\! & = & \!\! \Big\| \lambda_{\Theta_1} \big( \underbrace{R_1[D_{\Theta_1,j}^\ell] M_{B(\xi)} P_2[D_{\Theta_2,j}^\ell] Q_2[D_{\Theta_2,j}^r] M_{(2\pi i \xi)^{\beta_2}} f_j}_{\mathbf{A}(f_j)} \big) \Big\|_{\RR_{\Theta_1}}.
\end{eqnarray*}
Using standard commutation relations, $\mathbf{A}(f_j)$ may be written as $$\mathbf{A}(f_j) = \summ_k P_{1k}[D_{\Theta_1,j}^\ell] Q_{1k}[D_{\Theta_1,j}^r] M_{(2\pi i \xi)^{\beta_{1k}}} f_j$$ for finitely many $\beta_{1k} \in \N^n$ and quantum $\Theta_1$-monomials $P_{1k}$ and $Q_{1k}$. \fin 

A linear functional $L: \S_\Theta \to \C$ is in $\S_\Theta'$ when it satisfies that $\lim_{j \to \infty} \langle L, \varphi_j \rangle = 0$ for any sequence $\varphi_j \in \S_\Theta$ converging to $0$. Using the unitaries $j_\Theta: \exp_\xi \mapsto \lambda_\Theta(\xi)$ we construct $j_{\Theta_1 \Theta_2} = j_{\Theta_2} \circ j_{\Theta_1}^*: \lambda_{\Theta_1}(\xi) \mapsto \lambda_{\Theta_2}(\xi)$. According to Lemma \ref{ContinuitySchwartz} and given $L \in \S_{\Theta_1}'$, this means that $\langle j_{\Theta_1 \Theta_2}L, \lambda_{\Theta_2}(f) \rangle := \langle L, \lambda_{\Theta_1}(f) \rangle$ defines a tempered distribution in $\S_{\Theta_2}'$. Since this process is invertible, it turns out that the theory of tempered distributions in $\RR_\Theta$ is formally equivalent to the classical theory.

Let us now consider continuous linear operators $T \in \mathcal{L}(\S_\Theta, \S_\Theta')$. Of course, since the topology in $\S_\Theta'$ is that of pointwise convergence, a linear map $T: \S_\Theta \to \S_\Theta'$ is continuous whenever $\lim_j \langle T(\lambda_\Theta(f_j)), \lambda_\Theta(g) \rangle = 0$ for any family $\lambda_\Theta(f_j)$ which converges to $0$ in $\S_\Theta$ and any $\lambda_\Theta(g) \in \S_\Theta$. To identify the kernel of $T \in \mathcal{L}(\S_\Theta, \S_\Theta')$ consider $j_\Theta^* T j_\Theta \in \mathcal{L}(\S(\R^n), \S(\R^n)')$ where $j_\Theta: \S(\R^n) \to \S_\Theta$ and $j_\Theta^*: \S_\Theta' \to \S(\R^n)'$ by our discussion above. Then, by a well-known result of Schwartz, there exists a unique kernel $m \in \mathcal{S}'(\R^{2n}) = (\S(\R^n) \otimes_\pi \S(\R^n))'$ satisfying $$\big\langle j_\Theta^* T j_\Theta f, g \big\rangle = \big\langle m, g \otimes f \big\rangle \quad \mbox{for all} \quad f, g \in \mathcal{S}(\R^n).$$ Therefore, given $T \in \mathcal{L}(\S_\Theta, \S_\Theta')$ we find its associated kernel $$k = j_\Theta \otimes j_\Theta (m) \in \S_{\Theta \oplus \Theta}' \simeq (\S_\Theta \otimes_\pi \S_\Theta)'$$ such that 
\begin{eqnarray*}
\big\langle T(\lambda_\Theta(f)), \lambda_\Theta(g) \big\rangle & = & \big\langle j_\Theta^* T j_\Theta (\lambda_0(f)), \lambda_0(g) \big\rangle \\ & = & \big\langle m, \lambda_0(g) \otimes \lambda_0(f) \big\rangle \ = \ \big\langle k, \lambda_\Theta(g) \otimes \lambda_\Theta(f) \big\rangle.
\end{eqnarray*}
Now, according to the density of the quantum Schwartz class $\S_{\Theta \oplus \Theta}$ in $\S_{\Theta \oplus \Theta}'$ |since the same result also holds in the commutative case| we easily conclude the density of the algebraic tensor product $\mathcal{S}_\Theta \otimes_{\mathrm{alg}} \mathcal{S}_\Theta$ in $\S_{\Theta \oplus \Theta}'$. This proves that the family of algebraic kernels we considered for the proof of Theorem \ref{CZExt1} are dense in the space $\S_{\Theta \oplus \Theta}'$ of arbitrary kernels for maps in $\mathcal{L}(\S_\Theta, \S_\Theta')$. Moreover, by the weak-$*$ density of trigonometric polynomials in $\RR_{\Theta \oplus \Theta}$, we may also approximate $\S_{\Theta \oplus \Theta}'$ by finite sums of the form $$k = \summ_j w(\xi_j,\eta_j) \lambda_\Theta(\xi_j)^* \otimes \lambda_\Theta(\eta_j)^*.$$ According to our expression above |Section \ref{SectKernels}| for kernels affiliated to $\RR_\Theta \bar\otimes \RR_\Theta^{\mathrm{op}}$ $$k = \int_{\R^n \times \R^n} \widehat{k}(\xi,\eta) \lambda_\Theta(\xi) \otimes \lambda_\Theta(\eta)^* \, d\mu(\xi,\eta),$$ this result amounts to say (not surprising) that kernels associated to finite sums of Dirac deltas are dense. Note that identity $\big\langle T(\lambda_\Theta(f)), \lambda_\Theta(g) \big\rangle = \big\langle k, \lambda_\Theta(g) \otimes \lambda_\Theta(f) \big\rangle$ takes the following form for finite sums of Dirac deltas $$\big\langle T(\lambda_\Theta(f)), \lambda_\Theta(g) \big\rangle = \summ_j w(\xi_j,\eta_j) f(\eta_j) g(\xi_j) = \big\langle k, \lambda_\Theta(g) \otimes \lambda_\Theta(f) \big\rangle.$$

\subsubsection{Kernel manipulations and derivations}

In the model case treated previously we decomposed the kernel $k_\sigma = (\sigma_\Theta \otimes id)(k)$ as $k_\sigma \bullet \pi_\Theta(\psi_\mathrm{R}) + k_\sigma \bullet (\1 - \pi_\Theta(\psi_\mathrm{R}))$ and, after applying our Poincar\'e type inequality to the second term, we further decomposed the resulting kernel $\mathbf{K}$ as $\pi_\Theta(\psi_\mathrm{R}) \bullet \mathbf{K} + (\1 - \pi_\Theta(\psi_\mathrm{R})) \bullet \mathbf{K}$. This leads us to understand the same operation for general kernels in $\S_{\Theta \oplus \Theta}'$. To that end, we introduce the following operations for $L \in \S_\Theta'$ $$\big\langle \lambda_\Theta(\xi) L, \varphi \big\rangle = \big\langle L,  \varphi \lambda_\Theta(\xi) \big\rangle \quad \mbox{and} \quad \big\langle L \lambda_\Theta(\xi), \varphi \big\rangle = \big\langle L, \lambda_\Theta(\xi) \varphi \big\rangle.$$

\begin{lemma} \label{LRDLemma}
Given $\psi \in \S(\R^n)$ and $T \in \mathcal{L}(\S_\Theta, \S_\Theta')$, the maps 
\begin{eqnarray*}
M_\psi^\ell(T)(\lambda_\Theta(f)) & = & \int_{\R^n} \widehat{\psi}(\xi) \lambda_\Theta(\xi) T(\lambda_\Theta(\xi)^* \lambda_\Theta(f)) \, d\xi, \\
M_\psi^r(T)(\lambda_\Theta(f)) & = & \int_{\R^n} \widehat{\psi}(\xi) T(\lambda_\Theta(f) \lambda_\Theta(\xi)^*) \lambda_\Theta(\xi) \, d\xi,
\end{eqnarray*}
belong to $\mathcal{L}(\S_\Theta, \S_\Theta')$ and their kernels extend $\pi_\Theta(\psi) \bullet k$ and $k \bullet \pi_\Theta(\psi)$ respectively.
\end{lemma}

\dem We shall prove the assertion only for $M_\psi^\ell(T)$, since both operators can be handled similarly. In order to prove continuity, assume that $\lambda_\Theta(f_j) \to 0$ in $\S_\Theta$ as $j \to \infty$. Then we need to show that 
\begin{eqnarray*}
\lim_{j \to \infty} \big\langle M_\psi^\ell(T)(\lambda_\Theta(f_j)), \lambda_\Theta(g) \big\rangle \!\! & = & \!\! \lim_{j \to \infty} \int_{\R^n} \widehat{\psi}(\xi) \big\langle T(\lambda_\Theta(\xi)^* \lambda_\Theta(f_j)), \lambda_\Theta(g) \lambda_\Theta(\xi) \big\rangle d\xi \\ \!\! & = & \!\! \lim_{j \to \infty} \int_{\R^n} \widehat{\psi}(\xi) \big\langle k, \lambda_\Theta(g) \lambda_\Theta(\xi) \otimes \lambda_\Theta(\xi)^* \lambda_\Theta(f_j) \big\rangle d\xi
\end{eqnarray*}
vanishes for all $g \in \S(\R^n)$. We first note that $$\lim_{j \to \infty} \big\langle T(\lambda_\Theta(\xi)^* \lambda_\Theta(f_j)), \lambda_\Theta(g) \lambda_\Theta(\xi) \big\rangle = 0$$ for all $\xi \in \R^n$, since $\lambda_\Theta(\xi)^* \lambda_\Theta(f_j) \to 0$ in $\S_\Theta$ as $j \to \infty$. Indeed, we have $\lambda_\Theta(\xi)^* \lambda_\Theta(f_j) = \lambda_\Theta(f_{j\xi})$ for $f_{j\xi}(\eta) = f_j(\eta + \xi) e^{- 2\pi i \langle \xi, \Theta_{\downarrow} \eta \rangle}$ and we then use Lemma \ref{ContinuitySchwartz} with $(\Theta_1,\Theta_2) = (\Theta,0)$. Once this is known, we use the dominated convergence theorem, for which we need an integrable upper bound of $$\Phi(\xi) = \sup_{j \ge 1} \big| \widehat{\psi}(\xi) \big\langle k, \lambda_\Theta(g) \lambda_\Theta(\xi) \otimes \lambda_\Theta(\xi)^* \lambda_\Theta(f_j) \big\rangle \big|.$$ Since $\lambda_\Theta(g) \lambda_\Theta(\xi) = \lambda_\Theta(g_\xi)$ for $g_\xi(\eta) = g(\eta-\xi) e^{2\pi i \langle \eta - \xi, \Theta_{\downarrow} \xi \rangle}$, we have $$\big\langle k, \lambda_\Theta(g) \lambda_\Theta(\xi) \otimes \lambda_\Theta(\xi)^* \lambda_\Theta(f_j) \big\rangle = \big\langle j_\Theta^* \otimes j_\Theta^*(k), \widehat{g}_\xi \otimes \widehat{f}_{j\xi} \big\rangle$$ with $j_\Theta^* \otimes j_\Theta^*(k) \in \S(\R^{2n})'$. According to \cite[Proposition 2.3.4]{G1}, there exists an absolute constant $C$ and positive integers $k, m$ such that the following inequality holds 
\begin{eqnarray*}
\Phi(\xi) & \le & C |\widehat{\psi}(\xi)| \sup_{j \ge 1} \sum_{\begin{subarray}{c} |\alpha| \le k \\ |\beta| \le m \end{subarray}} p_{\alpha,\beta} \big( \widehat{g}_\xi \otimes \widehat{f}_{j\xi} \big) \\ & \le & C |\widehat{\psi}(\xi) P(\xi)| \sup_{j \ge 1} \sum_{\begin{subarray}{c} |\alpha| \le k \\ |\beta| \le m \end{subarray}} p_{\alpha,\beta} \big( \widehat{g} \otimes \widehat{f}_{j} \big) \ \lesssim \ |\widehat{\psi}(\xi) P(\xi)| \ \in \ \S(\R^n)
\end{eqnarray*}
for certain polynomial $P$. On the other hand, for integrable kernels 
\begin{eqnarray*}
M_\psi^\ell(T)(\lambda_\Theta(f)) & = & \int_{\R^n} \widehat{\psi}(\xi) \lambda_\Theta(\xi) (id \otimes \tau_\Theta) \big( (\1 \otimes \lambda_\Theta(\xi)^* \lambda_\Theta(f)) \bullet k \big) \, d\xi \\ & = & \int_{\R^n} \widehat{\psi}(\xi) (id \otimes \tau_\Theta) \big( (\1 \otimes \lambda_\Theta(f)) \bullet (\lambda_\Theta(\xi) \otimes \lambda_\Theta(\xi)^*) \bullet k \big) \, d\xi \\ [6pt] & = & (id \otimes \tau_\Theta) \big( (\1 \otimes \lambda_\Theta(f)) \bullet \pi_\Theta(\psi) \bullet k \big) \ = \ T_{\pi_\Theta(\psi) \bullet k}(\lambda_\Theta(f)). 
\end{eqnarray*}  
Interchanging trace and integral is justified for finite tensors by evaluation against a test function, and a fortiori by density of these kernels in $\mathcal{L}(\S_\Theta, \S_\Theta') \simeq \S_{\Theta \oplus \Theta}'$. \fin

\begin{remark} \label{L2BoundPsi}
\emph{Given $\dag \in \{\ell,r\}$, it is also clear that
  \[
    \big\| M_\psi^\dag(T): L_2(\RR_\Theta) \to L_2(\RR_\Theta) \big\| \le \Big( \int_{\R^n} |\widehat{\psi}(\xi)| \, d\xi \Big) \, \big\| T: L_2(\RR_\Theta) \to L_2(\RR_\Theta) \big\|.
  \]}
\end{remark}

\begin{remark} \label{Kernel}
\emph{It will also be relevant to observe that 
\begin{eqnarray*}
\Big\langle M_\psi^\ell(T)(\lambda_\Theta(f)), \lambda_\Theta(g) \Big\rangle & = & \Big\langle k, \big( \lambda_\Theta(g) \otimes \lambda_\Theta(f) \big) \bullet \pi_\Theta(\psi) \Big\rangle, \\ \Big\langle M_\psi^r(T)(\lambda_\Theta(f)), \lambda_\Theta(g) \Big\rangle & = & \Big\langle k, \pi_\Theta(\psi) \bullet \big( \lambda_\Theta(g) \otimes \lambda_\Theta(f) \big) \Big\rangle,
\end{eqnarray*}
for the kernel $k \in \S_{\Theta \oplus \Theta}'$ associated to $T$ and any $\psi \in \S(\R^n)$. Indeed, we have
\begin{eqnarray*}
\Big\langle M_\psi^\ell(T)(\lambda_\Theta(f)), \lambda_\Theta(g) \Big\rangle \!\!\! & = & \!\!\! \int_{\R^n} \widehat{\psi}(\xi) \Big\langle T(\lambda_\Theta(\xi)^* \lambda_\Theta(f)), \lambda_\Theta(g) \lambda_\Theta(\xi) \Big\rangle \, d\xi \\ \!\!\! & = & \!\!\! \int_{\R^n} \widehat{\psi}(\xi) \Big\langle k, \lambda_\Theta(g) \lambda_\Theta(\xi) \otimes \lambda_\Theta(\xi)^* \lambda_\Theta(f) \Big\rangle \, d\xi \\ \!\!\! & = & \!\!\! \int_{\R^n} \widehat{\psi}(\xi) \Big\langle k, \big( \lambda_\Theta(g) \otimes \lambda_\Theta(f) \big) \bullet \big( \lambda_\Theta(\xi) \otimes \lambda_\Theta(\xi)^* \big) \Big\rangle \, d\xi
\end{eqnarray*} which gives the desired identity. Moreover, these identities hold for any function $\psi$ for which both $( \lambda_\Theta(g) \otimes \lambda_\Theta(f) ) \bullet \pi_\Theta(\psi)$ and $\pi_\Theta(\psi) \bullet ( \lambda_\Theta(g) \otimes \lambda_\Theta(f) )$ stay in $\S_{\Theta \oplus \Theta}$.}  
\end{remark}

Again as we did in the model case above, we shall need to operate with module extensions. Given a linear map $T: \S_\Theta \to (\S(\R^n) \otimes_\pi \S_\Theta)'$ we will use  its module extension $T': \S_\Theta \otimes_\pi \S_\Theta \to (\S(\R^n) \otimes_\pi \S_\Theta)'$ given by $$T': \lambda_\Theta(f) \otimes \lambda_\Theta(g) \mapsto T(\lambda_\Theta(f)) (\1 \otimes \lambda_\Theta(g)),$$ where $\big\langle T(\lambda_\Theta(f)) (\1 \otimes \lambda_\Theta(g)), (a \otimes b) \big\rangle = \big\langle T(\lambda_\Theta(f)), (\1 \otimes \lambda_\Theta(g)) (a \otimes b) \big\rangle.$

\begin{lemma} \label{SigmaLD}
There exists a continuous map $\sigma_\Theta: \S_\Theta' \to (\S(\R^n) \otimes_\pi \S_\Theta)'$ which extends the corepresentation $\sigma_\Theta: \RR_\Theta \ni \lambda_\Theta(\xi) \mapsto \exp_\xi \otimes \lambda_\Theta(\xi) \in L_\infty(\R^n) \bar\otimes \RR_\Theta$. In particular, given $T_k \in \mathcal{L}(\S_\Theta, \S_\Theta')$ with kernel $k \in \S_{\Theta \oplus \Theta}'$ the composition $T_{k_\sigma} = \sigma_\Theta T_k$ belongs to $\mathcal{L}(\S_\Theta, (\S(\R^n) \otimes_\pi \S_\Theta)')$ with kernel $k_\sigma = (\sigma_\Theta \otimes id) k$. Moreover, $T_{k_\sigma}'$ extends to a continuous right modular map $\S_\Theta \otimes_\pi \S_\Theta \to (\S(\R^n) \otimes_\pi \S_\Theta)'$ with  $$T_{k_\sigma}' \big( (\lambda_\Theta(f) \otimes \1) \mathrm{flip} (\pi_\Theta(\psi)) \big) = M_\psi^r(T_{k_\sigma})(\lambda_\Theta(f)).$$
\end{lemma}

\dem The map $q(\widehat{f},\widehat{g}) = \widehat{fg}$ arises as the conjugation of the multiplication map $(f,g) \mapsto fg$ by the Fourier transform. It thus follows from the Leibniz rule that it defines a continuous map $\S(\R^n) \otimes_\pi \S(\R^n) \to \S(\R^n)$. Letting $q_\Theta = j_\Theta \circ q \circ (id \otimes j_\Theta^*)$ we define for $L \in \S_\Theta'$ $$\big\langle \sigma_\Theta L, \lambda_0(f) \otimes \lambda_\Theta(g) \big\rangle = \big\langle L, q_{\Theta} (\lambda_0(f) \otimes \lambda_{\Theta}(g)) \big\rangle.$$ It is clear that $\sigma_\Theta: \S_\Theta' \to (\S(\R^n) \otimes_\pi \S_\Theta)'$ is continuous and we find 
\begin{eqnarray*}
\big\langle \sigma_\Theta(\lambda_\Theta(\xi)^*), \lambda_0(f) \otimes \lambda_\Theta(g) \big\rangle & = & \tau_\Theta \big( \lambda_\Theta(\xi)^* \lambda_\Theta(fg) \big) \ = \ f(\xi) g(\xi) \\ & = & \big\langle \exp_{-\xi} \otimes \lambda_\Theta(\xi)^*, \lambda_0(f) \otimes \lambda_\Theta(g) \big\rangle
\end{eqnarray*}
and thus $\sigma_\Theta$ so defined extends the corepresentation $\sigma_\Theta$ introduced in Section \ref{CorepSect}. This immediately implies that $T_{k_\sigma} = \sigma_\Theta T_k$ belongs to $\mathcal{L}(\S_\Theta, (\S(\R^n) \otimes_\pi \S_\Theta)')$ and its kernel $k_\sigma = (\sigma_\Theta \otimes id) (k)$. Let us now justify the continuity of the module extension $T_{k_\sigma}'$. Indeed, the module extension of $\sigma_\Theta j_\Theta$, defines a continuous linear map $$W: \S(\R^n)' \otimes \S_\Theta \to (\S(\R^n) \otimes \S_\Theta)'$$ satisfying $W(\exp_\xi \otimes \lambda_\Theta(\eta)) = \exp_\xi \otimes \lambda_\Theta(\xi) \lambda_\Theta(\eta)$. Its continuity follows easily from the continuity of $\sigma_\Theta$. Next, observe that $T_{k_\sigma}' = W \circ (j_\Theta^* T_k \otimes id)$ since it trivially holds for the dense class of finite sums $k = \sum_j w(\xi_j, \eta_j) \lambda_\Theta(\xi_j)^* \otimes \lambda_\Theta(\eta_j)^*$. This automatically implies the continuity of the module extension $T_{k_\sigma}'$. It remains to justify the given identity for $T_{k_\sigma}'$ 
\begin{eqnarray*}
T_{k_\sigma}' \big( (\lambda_\Theta(f) \otimes \1) \mathrm{flip} (\pi_\Theta(\psi)) \big) & = & \int_{\R^n} \widehat{\psi}(\xi) T_{k_\sigma}' \big( \lambda_\Theta(f) \lambda_\Theta(\xi)^* \otimes \lambda_\Theta(\xi) \big) \, d\xi \\ & = & \int_{\R^n} \widehat{\psi}(\xi) T_{k_\sigma} \big( \lambda_\Theta(f) \lambda_\Theta(\xi)^* \big) (\1 \otimes \lambda_\Theta(\xi)) \, d\xi,
\end{eqnarray*}
which is the definition of $M_\psi^r(T_{k_\sigma})(\lambda_\Theta(f))$. This completes the proof. \fin

Our next goal is to generalize the Poincar\'e type inequality in Proposition \ref{Prop-Poincare} to the context of tempered distributions. Of course, the free $\Theta$-gradient can be understood as a map $\nabla_\Theta: \S_\Theta' \to \mathcal{L}(\mathbb{F}_n) \otimes \S_\Theta'$ in the canonical way $$\nabla_\Theta L = \sum_{k=1}^n s_k \otimes \partial_{\Theta}^k L \quad \mbox{for} \quad L \in \S_\Theta',$$ where $\langle \partial_{\Theta}^k L, \lambda_\Theta(f) \rangle = - \langle L,  \partial_{\Theta}^k \lambda_\Theta(f) \rangle = - 2 \pi i \langle L, \lambda_\Theta (f_{[k]}) \rangle$ and $f_{[k]}(\xi) = \xi_k f(\xi)$. Now, given a $\R^n$-ball $\mathrm{B_R}$ of radius $\mathrm{R}$ with characteristic function $q_\mathrm{R}$, Proposition \ref{Prop-Poincare} gives an upper bound for $$\hskip5pt \big\| \varphi - \varphi_{\mathrm{B_R}} \big\|_{\RR_\Theta \bar\otimes L_2^c(\phi)} = \Big\| \Big( \mean_{\mathrm{B}_\mathrm{R}} \big( \varphi - \varphi_{\mathrm{B_R}} \big)^* \big( \varphi - \varphi_{\mathrm{B_R}} \big) d\mu \Big)^\frac12 \Big\|_{\RR_\Theta}$$ in terms of the operator norm of the gradient of $\varphi$ localized at $\mathrm{B}_\mathrm{R}$. Let us recall that the predual of $\RR_\Theta \bar\otimes L_2^c(\phi)$ with respect to the linear bracket is given by the space $\A_{\Theta}(\mathrm{B_R}) = L_1(\RR_\Theta^\mathrm{op}) \widehat\otimes L_2^r(\phi)$, whose norm is
\[
  \big\| \psi - \psi_{\mathrm{B_R}} \big\|_{L_1(\RR_\Theta^\mathrm{op}) \widehat\otimes L_2^r(\phi)} = \Big\| \Big( \mean_{\mathrm{B}_\mathrm{R}} \big( \psi - \psi_{\mathrm{B_R}} \big) \big( \psi - \psi_{\mathrm{B_R}} \big)^* d\mu \Big)^\frac12 \Big\|_{L_1(\RR_\Theta)}.
\] It is a simple exercise to show that $\S(\R^n) \otimes_\pi \S_\Theta$ is norm dense in $L_1(\RR_\Theta^\mathrm{op}) \widehat\otimes L_2^r(\phi)$. In particular, the following result gives an extension of our Poincar\'e type inequality.
 
\begin{proposition} \label{PoincareDist}
Given $L \in (\S(\R^n) \otimes_\pi \S_\Theta)'$, assume $$\big( \partial_k \otimes id_{\RR_\Theta} \big)(L) \in L_\infty(\mathrm{B_R}) \bar\otimes \RR_\Theta$$ for $1 \le k \le n$. Then, the following Poincar\'e type inequality holds 
\begin{eqnarray*}
\lefteqn{\hskip-30pt \sup_{\begin{subarray}{c} \psi \in \S(\R^n) \otimes_\pi \S_\Theta \\ \|\psi - \psi_{\mathrm{B_R}}\|_{\A_{\Theta}(\mathrm{B_R})} \le 1 \end{subarray}} \big| \big\langle q_\mathrm{R} L, \psi - \psi_{\mathrm{B_R}} \big\rangle \big|} \\ & \le & 2 \sqrt{2} \mathrm{R} \Big\| (\1 \otimes q_\mathrm{R} \otimes \1) (\nabla \otimes id_{\RR_\Theta})(L) \Big\|_{\mathcal{L}(\F_n) \bar\otimes L_\infty(\R^n) \bar\otimes \RR_\Theta}.
\end{eqnarray*}
\end{proposition}

\dem Assume for clarity that $\mathrm{B_R}$ is centered at the origin, see Proposition \ref{Prop-Poincare} for the standard modifications in the general case. Since $\partial_k L \in L_\infty(\R^n) \bar\otimes \RR_\Theta$, we may define 
\begin{eqnarray*}
\widetilde{L}(x) & = & \int_0^1 \sum_{k=1}^n \partial_k L(tx) x_k \, dt \\ & = & \int_0^1 (\underbrace{\tau_{\mathcal{L}(\mathbb{F}_n)} \otimes id_{\R^n}}_{\mathsf{E}_{\R^n}}) \big( \underbrace{q_\mathrm{R}(x) \nabla L(tx)}_{\mathrm{A}(t)} \underbrace{q_\mathrm{R}(x) s(x)}_{\mathrm{B}} \big) \, dt
\end{eqnarray*}
for $x \in \mathrm{B_R}$. Now let $\varphi_j \in \S(\R^n) \otimes_\pi \S_\Theta$ be an approximating sequence for $L$ and define the functions $\widetilde{\varphi}_j(x) = \varphi_j(x) - \varphi_j(0)$ accordingly. In particular, the following identity holds for every test function $\psi \in \S(\R^n) \otimes_\pi \S_\Theta$ $$\int_{\mathrm{B_R}} \tau_\Theta \big( \widetilde{\varphi}_j(x) (\psi(x) - \psi_{\mathrm{B_R}}) \big) \, dx = \int_{\mathrm{B_R}} \tau_\Theta \big( \varphi_j(x) (\psi(x) - \psi_{\mathrm{B_R}}) \big) \, dx.$$
By approximation we get 
\begin{eqnarray*}
\big| \big\langle q_\mathrm{R} L, \psi - \psi_{\mathrm{B_R}} \big\rangle \big| & = & \big| \big\langle q_\mathrm{R} \widetilde{L}, \psi - \psi_{\mathrm{B_R}} \big\rangle \big| \\ [5pt] & = & \Big| \int_0^1 \big\langle \mathsf{E}_{\R^n} (\mathrm{A}(t) \mathrm{B}), \psi - \psi_{\mathrm{B_R}} \big\rangle \, dt \Big| \\ & \le & \Big( \int_0^1 \big\| \mathsf{E}_{\R^n} (\mathrm{A}(t) \mathrm{B}) \big\|_{\RR_\Theta \bar\otimes L_2^c(\phi)} \, dt  \Big) \big\| \psi - \psi_{\mathrm{B_R}} \big\|_{\A_{\Theta}(\mathrm{B_R})}. 
\end{eqnarray*}
Now we may complete the argument as we did in the proof of Proposition \ref{Prop-Poincare}. \fin

\subsubsection{A Calder\'on-Zygmund extrapolation theorem}

We are finally ready to prove an estimate for general Calder\'on-Zygmund operators $T_k$. According to the classical theory, we impose cancellation and smoothness conditions on the kernel. To be more precise, let $T_k \in \mathcal{L}(\S_\Theta, \S_\Theta')$ admit a kernel $k \in \S_{\Theta \oplus \Theta}'$ with gradient $(\nabla_\Theta \otimes id_{\RR_\Theta^{\mathrm{op}}})(k)$ affiliated to $\mathcal{L}(\mathbb{F}_n) \bar\otimes \RR_\Theta \bar\otimes \RR_\Theta^{\mathrm{op}}$. Then, we shall call $T_k$ a \emph{column Calder\'on-Zygmund operator} with parameters $(\mathrm{A}_j, \alpha_j, \beta_j)$ when: 
\begin{itemize}
\item[i)] \emph{Cancellation} $$\hskip10pt \big\| T_k: L_2(\RR_\Theta) \to L_2(\RR_\Theta) \big\| \le \mathrm{A}_1.$$

\item[ii)] \emph{Kernel smoothness.} If $\mathrm{K}_n = \frac12 (n+1)$, there exists $$\hskip10pt \alpha < \mathrm{K}_n - \frac12 < \beta < \mathrm{K}_n + \frac12$$ satisfying the gradient conditions below for $\rho = \alpha, \beta$ $$\hskip30pt \Big| \mathrm{d}_\Theta^{\rho} \bullet (\nabla_\Theta \otimes id_{\RR_\Theta^\mathrm{op}}) (k) \bullet \mathrm{d}_\Theta^{n+1-\rho} \Big| \le \mathrm{A}_2.$$
\end{itemize}

\begin{remark} \label{GradientAffiliated}
\emph{We implicitly use that $$\Big\langle \mathrm{d}_\Theta^{\gamma} \bullet \nabla_\Theta k \bullet \mathrm{d}_\Theta^{\eta}, z \Big\rangle = \Big\langle \nabla_\Theta k, \mathrm{d}_\Theta^\eta \bullet z \bullet \mathrm{d}_\Theta^\gamma \Big\rangle,$$ $$\big\| \mathrm{d}_\Theta^{\gamma} \bullet \nabla_\Theta k \bullet \mathrm{d}_\Theta^{\eta} \big\|_{\RR_\Theta \bar\otimes \RR_\Theta^\mathrm{op}}^{\null} = \sup_{\begin{subarray}{c} z \in \mathcal{L}(\mathbb{F}_n) \bar\otimes \S_{\Theta \oplus \Theta} \\ \|z\|_{L_1(\mathcal{L}(\mathbb{F}_n) \bar\otimes \RR_\Theta \bar\otimes \RR_\Theta^{\mathrm{op}})} \le 1 \end{subarray}} \Big\langle \nabla_\Theta k, \mathrm{d}_\Theta^\eta \bullet z \bullet \mathrm{d}_\Theta^\gamma \Big\rangle.$$ As explained in Remark \ref{Kernel}, this is justified for any tempered distribution when $\gamma, \eta \in 2\Z$, but not for $\gamma, \eta < 2$ since $\mathrm{d}_\Theta^\eta \bullet z \bullet \mathrm{d}_\Theta^\gamma$ does not stay in the test space $\S_{\Theta \oplus \Theta}$. The necessity of using these values |only for $n=2$ in the simpler statement of Theorem A| forces us to impose that $\nabla_\Theta k$ is, in addition, affiliated to the algebra. Although our assumption is admissible in view of the classical theory we could have alternatively used an approximation argument $\mathrm{d}_\Theta^\gamma = \lim_{\varepsilon} ( \mathrm{d}_\Theta^2 + \varepsilon \1 )^{\gamma/2}$ to avoid it. On the other hand, the kernel $k$ |not its gradient| should be treated as a distribution since this allows certain Dirac deltas which do appear in the classical theory, see Paragraph \ref{PV} below for further details.}
\end{remark}

\begin{proposition} \label{CZExt2}
If $T_k$ is a column \emph{CZO} and $\varphi \in \S_\Theta$ $$\big\| T_k (\varphi) \big\|_{\mathrm{BMO}_c(\RR_\Theta)} \le C_n(\alpha, \beta) \big( \mathrm{A}_1 + \mathrm{A}_2 \big) \|\varphi\|_{\RR_\Theta}.$$
\end{proposition}

\dem We shall adapt our argument in the model case of Theorem \ref{CZExt1}. Given $\varphi = \lambda_\Theta(f) \in \S_\Theta$, this means that we need to control the operator-valued BMO norm of $\sigma_\Theta(T_k\varphi)$. According to Lemma \ref{SigmaLD}, we have $\sigma_\Theta T_k = T_{k_\sigma}$ for $k_\sigma = (\sigma_\Theta \otimes id)(k)$ and we may decompose it as follows $$k_\sigma = \underbrace{k_\sigma \bullet \pi_\Theta(\psi_\mathrm{R})}_{k_{\sigma1}(\mathrm{R})} + \underbrace{k_\sigma \bullet \big( \1 - \pi_\Theta(\psi_\mathrm{R}) \big)}_{k_{\sigma2}(\mathrm{R})},$$ where the decomposition uses Lemma \ref{LRDLemma} and Remark \ref{Kernel}. Next, we need to show the validity of \eqref{Est1}. To that end we follow the argument in Theorem \ref{CZExt1} by recalling the crucial identity $$T_{k_{\sigma}}' \big( (\lambda_\Theta(f) \otimes \1) \mathrm{flip} (\pi_\Theta(\psi)) \big) = M_\psi^r(T_{k_{\sigma}})(\lambda_\Theta(f)),$$ which was justified in Lemma \ref{SigmaLD} for general kernels. This is the part of the proof which requires $L_2$-boundedness of $T_k$. Once we have completed our argument for $k_{\sigma_1}(\mathrm{R})$, we apply the Poincar\'e type inequality in Proposition \ref{PoincareDist} to the term associated to $k_{\sigma2}(\mathrm{R})$. This gives 
\begin{eqnarray*}
\lefteqn{\hskip-25pt \Big\| \Big( \mean_{\mathrm{B_R}} \big| T_{k_{\sigma 2} (\mathrm{R})}\varphi - (T_{k_{\sigma 2} (\mathrm{R})}\varphi)_{\mathrm{B_R}} \big|^2 d\mu \Big)^\frac12 \Big\|_{\RR_\Theta}} \\ \hskip20pt & \lesssim & \mathrm{R} \Big\| \underbrace{(\1 \otimes q_\mathrm{R} \otimes \1) (\nabla \otimes id_{\RR_\Theta}) (T_{k_{\sigma2}(\mathrm{R})}\varphi)}_{T_\mathbf{K}(\varphi)} \Big\|_{\mathcal{L}(\F_n) \bar\otimes L_\infty(\R^n) \bar\otimes \RR_\Theta}.
\end{eqnarray*}
As in Theorem \ref{CZExt1}, the goal is to show that $\mathbf{K}$ is the distribution $$\mathbf{K} = q_\mathrm{R} \nabla_\Theta^\sigma k \bullet (\1 - \Psi_\mathrm{R}) \in \mathcal{L}(\mathbb{F}_n) \bar\otimes \big( \S(\R^n) \otimes_\pi \S_\Theta \otimes_\pi \S_\Theta \big)'$$ with $q_\mathrm{R} \nabla_\Theta^\sigma k = (\1 \otimes q_\mathrm{R} \otimes \1^{\otimes 2}) (id \otimes \sigma_\Theta \otimes id) (\nabla_\Theta \otimes id) [k]$ and $\Psi_\mathrm{R} = \pi_\Theta(\psi_\mathrm{R})$. By density, it suffices to justify it for elementary kernels $k = \lambda_\Theta(\xi) \otimes \lambda_\Theta(\eta)$. This is possible via the identity $(\nabla \otimes id_{\RR_\Theta}) \circ \sigma_\Theta = (id_{\mathcal{L}(\mathbb{F}_n)} \otimes \sigma_\Theta) \circ \nabla_\Theta$ due to our extensions of the maps $\sigma_\Theta$ and $\nabla_\Theta$ for tempered distributions above. We may now decompose $\mathbf{K}$ by means of Lemma \ref{LRDLemma} as follows $$\mathbf{K} = \underbrace{q_\mathrm{R} \Psi_\mathrm{R} \bullet \nabla_\Theta^\sigma k \bullet (\1 - \Psi_\mathrm{R})}_{\mathbf{K}_1} + \underbrace{q_\mathrm{R} (\1 - \Psi_\mathrm{R}) \bullet \nabla_\Theta^\sigma k \bullet (\1 - \Psi_\mathrm{R})}_{\mathbf{K}_2}.$$ At this point, the argument follows verbatim the proof of Theorem \ref{CZExt1}. Indeed, we further decompose $\mathbf{K}_j = \mathbf{A}_j \bullet \mathbf{B}_j \bullet \mathbf{C}_j$ as we did there and apply Lemma \ref{KernelL} |valid for affiliated kernels, as we assume for $\mathbf{K}_j$| to obtain
\begin{eqnarray} \label{FactDist}
\lefteqn{\hskip-10pt \big\| T_{\mathbf{K}_j}: \RR_\Theta \to \B_\Theta \big\|} \\ \hskip20pt \nonumber \!\!\! & \le & \!\!\! \big\| (id \otimes \tau_\Theta)(\mathbf{A}_j \mathbf{A}_j^*) \big\|_{\B_\Theta}^\frac12 \big\| \mathbf{B}_j \big\|_{\B_\Theta \bar\otimes \RR_\Theta^{\mathrm{op}}} \big\| (id \otimes \tau_\Theta)(\mathbf{C}_j^* \mathbf{C}_j) \big\|_{\B_\Theta}^\frac12
\end{eqnarray}
with $\B_\Theta = \mathcal{L}(\mathbb{F}_n) \bar\otimes L_\infty(\R^n) \bar\otimes \RR_\Theta$. The estimates for $\mathbf{A}_j, \mathbf{B}_j, \mathbf{C}_j$ also apply here. \fin  

\begin{remark} \label{AlternativePf}
\emph{Alternatively, if we do not want to assume that $\nabla_\Theta k$ is affiliated to $\mathcal{L}(\mathbb{F}_n) \bar\otimes \RR_\Theta \bar\otimes \RR_\Theta^\mathrm{op}$ and use the approximation argument indicated in Remark \ref{GradientAffiliated}, we should be able to generalize the inequality \eqref{FactDist} for tempered distributions. Recall that the norm of $T_{\mathbf{K}_j}: \RR_\Theta \to \B_\Theta$ can be expressed as the supremum |over Schwartz elements $\varphi, \phi$ respectively in the unit ball of $\RR_\Theta^{\mathrm{op}}$ and $L_1(\B_\Theta)$| of the linear brackets 
\begin{eqnarray*}
\big| \big\langle T_{\mathbf{K}_j} \varphi, \phi \big\rangle \big| & = & \big| \big\langle \mathbf{A}_j \bullet \mathbf{B}_j \bullet \mathbf{C}_j, \phi \otimes \varphi \big\rangle \big| \\ 
[5pt] & \le & \big\| \mathbf{A}_j \bullet \mathbf{B}_j \bullet \mathbf{C}_j \big\|_{\B_\Theta \bar\otimes L_1(\RR_\Theta^{\mathrm{op}})} \|\phi\|_{L_1(\B_\Theta)} \|\varphi\|_{\RR_\Theta^{\mathrm{op}}}.
\end{eqnarray*}
Now we use the following characterization of the norm in $\M \bar\otimes \mathrm{X}$ $$\big\| \mathrm{A} \big\|_{\M \bar\otimes \mathrm{X}} = \sup_{\alpha, \beta \in \mathrm{B}_1(L_2(\M))} \big\| (\alpha \otimes \1) \mathrm{A} (\beta \otimes \1) \big\|_{L_1(\M; \mathrm{X})},$$ which is due to Pisier \cite{P2} when $\M$ is hyperfinite and $\mathrm{X}$ is any operator space. It is also well-known that Pisier's identity still holds for non-hyperfinite von Neumann algebras |as in our case with $\M = \B_\Theta$| as long as $\mathrm{X}$ is a noncommutative $L_p$ space. In fact, Pisier's identity generalizes to arbitrary mixed $L_p(L_q)$-norms. In our case
\begin{eqnarray*}
\lefteqn{ \hskip-20pt \big\| \mathbf{A}_j \bullet \mathbf{B}_j \bullet \mathbf{C}_j \big\|_{\B_\Theta \bar\otimes L_1(\RR_\Theta^{\mathrm{op}})}} \\ & = & \sup_{\mathrm{a}, \mathrm{c} \in \mathrm{B}_1(L_2(\B_\Theta))} \big\| (\mathrm{a} \otimes \1) \bullet \mathbf{A}_j \bullet \mathbf{B}_j \bullet \mathbf{C}_j \bullet (\mathrm{c} \otimes \1) \big\|_{L_1(\B_\Theta \bar\otimes \RR_\Theta^{\mathrm{op}})}. 
\end{eqnarray*}
In particular, we find $$\big\| (\mathrm{a} \otimes \1) \bullet \mathbf{A}_j \bullet \mathbf{B}_j \bullet \mathbf{C}_j \bullet (\mathrm{c} \otimes \1) \big\|_1 \le \big\| (\mathrm{a} \otimes \1) \bullet \mathbf{A}_j \big\|_2 \big\| \mathbf{B}_j \big\|_\infty \big\| \mathbf{C}_j \bullet (\mathrm{c} \otimes \1) \big\|_2$$ and the elementary inequalities below complete the proof of \eqref{FactDist}
\begin{eqnarray*}
\big\| (\mathrm{a} \otimes \1) \bullet \mathbf{A}_j \big\|_2 & \le & \|\mathrm{a}\|_2 \big\| (id \otimes \tau_\Theta) (\mathbf{A}_j \mathbf{A}_j^*) \big\|_{\B_\Theta}^\frac12, \\
\big\| \mathbf{C}_j \bullet (\mathrm{c} \otimes \1) \big\|_2 & \le & \|\mathrm{c}\|_2 \big\| (id \otimes \tau_\Theta) (\mathbf{C}_j^* \mathbf{C}_j) \hskip1pt \big\|_{\B_\Theta}^\frac12. 
\end{eqnarray*}}
\end{remark}

\begin{proposition} \label{CZExt3}
Every column \emph{CZO} is normal. In particular $$\big\| T_k: \RR_\Theta \to \mathrm{BMO}_c(\RR_\Theta) \big\|_{\mathrm{cb}} \le C_n(\alpha, \beta) \big( \mathrm{A}_1 + \mathrm{A}_2 \big).$$
\end{proposition}

\dem Let $T_k^*: L_2(\RR_\Theta) \to L_2(\RR_\Theta)$ denote the adjoint of $T_k$, so that 
\begin{eqnarray*}
\big| \tau_\Theta \big( T_k^*(\lambda_\Theta(f)) \lambda_\Theta(g)^* \big) \big| \!\!\! & = & \!\!\! \big| \tau_\Theta \big( \lambda_\Theta(f) T_k(\lambda_\Theta(g))^* \big) \big| \\ \!\!\! & \le & \!\!\! C_n(\alpha, \beta) \big( \mathrm{A}_1 + \mathrm{A}_2 \big) \big\| \lambda_\Theta(f) \big\|_{\mathrm{H}_c^1(\RR_\Theta)} \big\| \lambda_\Theta(g) \big\|_{\RR_\Theta}
\end{eqnarray*}
for all $\lambda_\Theta(f), \lambda_\Theta(g) \in \S_\Theta$. Indeed, here $\mathrm{H}_c^1(\RR_\Theta)$ denotes the predual of $\mathrm{BMO}_c(\RR_\Theta)$ with respect to the antilinear duality bracket above, as described in Appendix B below. In particular, this inequality directly follows from Proposition \ref{CZExt2}. Now we claim that this implies $$\big\| T_k^* (\varphi) \big\|_{L_1(\RR_\Theta)} \le C_n(\alpha, \beta) \big( \mathrm{A}_1 + \mathrm{A}_2 \big) \|\varphi\|_{\mathrm{H}_c^1(\RR_\Theta)}$$ for all $\varphi = \lambda_\Theta(f) \in \S_\Theta$. Indeed, let us prove that $$z = \frac{T_k^*(\varphi)}{C_n(\alpha, \beta) \big( \mathrm{A}_1 + \mathrm{A}_2 \big) \|\varphi\|_{\mathrm{H}_c^1(\RR_\Theta)}}$$ belongs to the unit ball of $L_1(\RR_\Theta)$. To that end, it clearly suffices to prove that $|\tau_\Theta(qzqa)| \le 1$ for every contraction $a$ in $\RR_\Theta$ and every $\tau_\Theta$-finite projection $q$. Since $z \in L_2(\RR_\Theta)$, we have $zq \in L_1(\RR_\Theta)$ and Kaplansky density theorem provides a sequence $u_j \in \S_\Theta$ in the unit ball of $\RR_\Theta$ so that $$|\tau_\Theta(qzqa)| = \lim_{j \to \infty} |\tau_\Theta (u_jzq)|.$$ Moreover, since $u_j z \in L_1(\RR_\Theta)$, we also find $v_k \in \S_\Theta$ in the unit ball of $\RR_\Theta$ with 
$$|\tau_\Theta(qzqa)| \, = \, \lim_{j \to \infty} |\tau_\Theta (u_jzq)| \, = \, \lim_{j \to \infty} \lim_{k \to \infty} |\tau_\Theta (u_jzv_k)|.$$
Finally, since $|\tau_\Theta (z w) | \le 1$ for every $w \in \S_\Theta$ in the unit ball of $\RR_\Theta$ |as we recalled at the beginning of the proof| and the Schwartz class $\S_\Theta$ is a $*$-algebra we obtain that $|\tau_\Theta(qzqa)| \le 1$ as expected. This proves our claim. Next, we use the norm density of $\S_\Theta$ in $\mathrm{H}_c^1(\RR_\Theta)$ from Corollary \ref{Bdensity} in Appendix B below to  conclude that $T_k^*: \mathrm{H}_c^1(\RR_\Theta) \to L_1(\RR_\Theta)$ is bounded. The operator $T_k^*$ is the antilinear adjoint corresponding to the duality $$\overline{L_1(\RR_\Theta)}^* = \RR_\Theta$$ with respect to the antilinear duality bracket. Thus $$T_k: \RR_\Theta \to \overline{\mathrm{H}_c^1(\RR_\Theta)}^* \simeq \mathrm{BMO}_c(\RR_\Theta)$$ with the same constants, see Appendix B for further details on the duality $\mathrm{H}_1 - \mathrm{BMO}$ in this setting. This proves the $L_\infty \to \mathrm{BMO}_c$ boundedness of $T_k$. As in the model case proved in Theorem \ref{CZExt1}, the cb-boundedness follows similarly and it just requires a more involved notation to incorporate matrix amplifications. \fin

Once we have proved the complete $L_\infty \to \mathrm{BMO}_c$ boundedness of column CZOs, the general extrapolation theorem follows from additional assumptions of the same kind on the kernel, which makes them more symmetric. More precisely, we know that $T_k: \RR_\Theta \to \mathrm{BMO}_r(\RR_\Theta)$ is cb-bounded iff the operator $$T_k^\dag(\lambda_\Theta(f)) = T_k(\lambda_\Theta(f)^*)^*$$ defines a completely bounded map from $\RR_\Theta \to \mathrm{BMO}_c(\RR_\Theta)$. When this is the case we get a cb-map $T_k: \RR_\Theta \to \mathrm{BMO}(\RR_\Theta)$. Of course, the same assumptions for the adjoint $T_k^*$ trivially imply that $T_k$ also defines a cb-map $T_k: \mathrm{H}_1(\RR_\Theta) \to L_1(\RR_\Theta)$ and interpolation |see Appendix B for details| yields complete $L_p$-boundedness for $1 < p < \infty$. This means that we should impose that the maps $T_k^\dag, T_k^*, T_k^{*\dag}$ are column Calder\'on-Zygmund operators. It is clear that $L_2$-boundedness follows automatically from $T_k$. Therefore, we just need to impose new kernel smoothness conditions. We have $$\mathrm{kernel}(T_k^\dag) = k^*, \quad \mathrm{kernel}(T_k^*) = \mathrm{flip}(k)^*, \quad \mathrm{kernel}(T_k^{*\dag}) = \mathrm{flip}(k).$$ Therefore, the results so far imply the following extrapolation theorem for CZOs.
\begin{theorem} \label{CZExt4}
Let $T_k \in \mathcal{L}(\S_\Theta, \S_\Theta')$ and assume$\hskip1pt :$  
\begin{itemize}
\item[\emph{i)}] \emph{Cancellation} $$ \big\| T_k: L_2(\RR_\Theta) \to L_2(\RR_\Theta) \big\| \le \mathrm{A}_1.$$

\item[\emph{ii)}] \emph{Kernel smoothness.} There exists $$\alpha < \mathrm{K}_n - \frac12 < \beta < \mathrm{K}_n + \frac12 < \gamma$$ satisfying the gradient conditions below for $\rho = \alpha, \beta, \gamma$ $$\hskip20pt \Big| \mathrm{d}_\Theta^{\rho} \bullet (\nabla_\Theta \otimes id) (k) \bullet \mathrm{d}_\Theta^{n+1-\rho} \Big| + \Big| \mathrm{d}_\Theta^{\rho} \bullet (id \otimes \nabla_\Theta) (k) \bullet \mathrm{d}_\Theta^{n+1-\rho} \Big| \le \mathrm{A}_2.$$
\end{itemize}
Then, we find the following endpoint estimates for $T_k$ 
\begin{eqnarray*}
\big\| T_k: \mathrm{H}_1(\RR_\Theta) \to L_1(\RR_\Theta) \big\|_{\mathrm{cb}} & \le & C_n(\alpha, \beta, \gamma) \big( \mathrm{A}_1 + \mathrm{A}_2 \big), \\ \big\| T_k: L_\infty(\RR_\Theta) \to \mathrm{BMO}(\RR_\Theta) \big\|_{\mathrm{cb}} & \le & C_n(\alpha, \beta, \gamma) \big( \mathrm{A}_1 + \mathrm{A}_2 \big).
\end{eqnarray*}
In particular, $T_k: L_p(\RR_\Theta) \to L_p(\RR_\Theta)$ is completely bounded for every $1 < p < \infty$. 
\end{theorem}

In what follows, a \emph{Calder\'on-Zygmund operator} over the quantum Euclidean space $\RR_\Theta$ associated to the parameters $(\mathrm{A}_j, \alpha_j, \beta_j)$ will be any linear map $T_k \in \mathcal{L}(\S_\Theta, \S_\Theta')$ satisfying the hypotheses in Theorem \ref{CZExt4}. The kernel considerations for the adjoint also appear in commutative Calder\'on-Zygmund theory, whereas the $\dag$-operation is standard and arises from noncommutativity.

\subsubsection{The principal value of kernel truncations}
\label{PV}

As in classical Calder\'on-Zygmund theory, we want to understand how far is an operator $T_k \in \B(L_2(\RR_\Theta))$ from the principal value singular integral determined by its kernel truncations. Our aim is to show that the difference is a left/right multiplier. Let us be more precise. Consider a smooth function $\psi \in \S(\R^n)$ which is identically 1 over $\mathrm{B}_1(0)$ and vanishes over $\R^n \setminus \mathrm{B}_2(0)$. Define $$\Psi_{\Delta,\delta} = \pi_\Theta(\psi_{\Delta,\delta}) \quad \mbox{with} \quad \psi_{\Delta,\delta} (x) = \psi \big( \frac{x}{\Delta} \big) - \psi \big( \frac{x}{\delta} \big) = \psi_\Delta(x) - \psi_\delta(x)$$ for $0 < \delta << \Delta < \infty$. We shall study the kernel truncations $\Psi_{\Delta,\delta} \bullet k$ and $k \bullet \Psi_{\Delta,\delta}$ and how their limits are related to $T_k$. To that end, we introduce the notion of \emph{admissible projection}. A projection $p \in \RR_\Theta$ will be called admissible when the function $\R^n \to \mathrm{Proj}(\RR_\Theta)$ defined as $$\delta \longmapsto \bigvee_{s \in \mathrm{B}_\delta(0)} \sigma_\Theta^s(p)$$ is weak-$*$ continuous around $\delta = 0$. Here $\sigma_\Theta^s(\lambda_\Theta(\xi)) = \exp( 2 \pi i \langle s, \xi \rangle) \, \lambda_\Theta(\xi)$. 

\begin{remark} \label{Rem1DClosedProj}
\emph{Even in the Euclidean setting with $\Theta = 0$, not all projections are admissible. In that case, the projection-valued function defined above associates a measurable set $A$ with $\mathrm{B}_\delta[A]$, the union of all the balls of radius $\delta$ with center in $A$. If we take, for instance, a dense open subset of $[0,1]^n$ with measure strictly less than $1$, we will have that $[0,1] \subset \mathrm{B}_\delta[A]$ for every $\delta > 0$, which poses an obstruction to admissibility. This can be easily fixed in the Euclidean setting by considering measurable sets which are closed up to a null set.}
\end{remark}

\begin{lemma} \label{LemmAllAlgebra}
The bicommutant of admissible projections is the whole algebra $\RR_\Theta$.
\end{lemma}

\dem It suffices to observe that one-dimensional spectral projections of the form $\chi_{[a,b]}(x_{\Theta,j})$ are admissible for each of the quantum variables $x_{\Theta,j}$ by Remark \ref{Rem1DClosedProj} and also that this family trivially generates $\RR_\Theta$. This completes the proof. \fin 

\begin{remark}
\emph{Let us define a \emph{closed projection} $p \in \RR_\Theta$ as those whose complement $\1 - p$ is the left support of certain element $\varphi \in \mathrm{E}_\Theta$ as defined at beginning of Section \ref{Prelim}. By the $\ast$-stability of $\S_\Theta$ we could have replaced the left support $\ell(\varphi)$ by the right one $r(\varphi)$ or even by the full support $s(\varphi)$ of self-adjoint elements. We conjecture that all closed projections so defined are indeed admissible. At the time of this writing we have not been able to confirm this conjecture, but this will have no consequence in Theorem \ref{PrincipalValue} below.}
\end{remark}

\begin{lemma} \label{LocalizationL}
Given $\varphi \in \RR_\Theta$, there exist projections $p_\delta, q_\delta$ such that
\begin{eqnarray*}
\pi_\Theta (\psi_{\delta}) \bullet (\1 \otimes \varphi) & = & \pi_\Theta (\psi_{\delta}) \bullet (p_{\delta} \otimes \varphi), \\
\pi_\Theta (\psi_{\delta}) \bullet (\varphi \otimes \1) & = & \pi_\Theta (\psi_{\delta}) \bullet (\varphi \otimes q_{\delta}).
\end{eqnarray*}
If $r(\varphi)$ is admissible $\displaystyle \mathrm{w}^*\mbox{-}\lim_{\delta \to 0} p_{\delta} = r(\varphi)$, if $\ell(\varphi)$ is admissible $\displaystyle \mathrm{w}^*\mbox{-}\lim_{\delta \to 0} q_{\delta} = \ell(\varphi).$
\end{lemma}

\dem The assertions concerning $q_\delta$ follow from those for $p_\delta$ after applying the map $\mathrm{flip}^*: a \otimes b \mapsto b^* \otimes a^*$, details are left to the reader. Now, let us recall that the map $T: \RR_\Theta \bar\otimes \RR_\Theta^\mathrm{op} \to \mathcal{CB}(L_1(\RR_\Theta), \RR_\Theta)$ sending a kernel $k$ to the corresponding map $T_k$ is a complete isometry. Moreover, observe that 
$$T_{k \bullet (\1 \otimes \varphi)}(\phi) = T_k(\phi \varphi) \quad \mbox{and} \quad T_{k \bullet (\varphi \otimes \1)}(\phi) = T_k(\phi) \varphi.$$ 
Since we clearly have $$\| \pi_\Theta(\psi_\delta) \|_{\RR_\Theta \bar\otimes \RR_\Theta^\mathrm{op}} \leq \| \widehat{\psi}_\delta \|_{L_1(\R^n)} = \| \widehat{\psi} \hskip1pt \|_{L_1(\R^n)} < \infty,$$ we know that $T_{\pi_\Theta(\psi_\delta)}$ is uniformly in $\mathcal{CB}(L_1(\RR_\Theta),\RR_\Theta)$. Let us define
\[
  \mathcal{N}_{\delta} =
  \overline{\mathrm{span}^{\mathrm{w}^\ast}} \big\{ \varpi \, T_{\pi_\Theta (\psi_{\delta})}(\phi \varphi): \phi \in L_1(\RR_\Theta), \, \varpi \in \RR_\Theta \big\}
  \subset \RR_\Theta.
\]
Clearly $\mathcal{N}_{\delta}$ is a weak-$\ast$ closed left module. In particular, there must exist certain projection $p_\delta \in \RR_\Theta$ satisfying $\mathcal{N}_\delta = \RR_\Theta \, p_\delta$ and the following identity holds for every element $\phi \in L_1(\RR_\Theta)$
\[
  T_{\pi_\Theta (\psi_{\delta}) \bullet (\1 \otimes \varphi)}(\phi) =
  T_{\pi_\Theta (\psi_{\delta})}(\phi \varphi) =
  T_{\pi_\Theta (\psi_{\delta})}(\phi \varphi) \, p_\delta =
  T_{\pi_\Theta (\psi_{\delta}) \bullet (p_\delta \otimes \varphi)}(\phi).
\]
Since $T$ is (completely) isometric, $\pi_\Theta (\psi_{\delta}) \bullet (\1 \otimes \varphi) = \pi_\Theta (\psi_{\delta}) \bullet (p_\delta \otimes \varphi)$. It remains to show that the projections $p_\delta$ so defined converges weakly to $r(\varphi)$ as $\delta \to 0^+$. Given any $\phi \in L_1(\RR_\Theta)$, notice that
$$
T_{\pi_\Theta(\psi_\delta)}(\phi \varphi) = \int_{\R^n} \widehat{\psi}_\delta(\xi) \, \tau_\Theta \big( \lambda_\Theta(\xi)^\ast \phi \varphi \big) \lambda_\Theta(\xi) \, d\xi 
= \int_{\mathrm{B}_{2\delta}(0)} \psi_\delta(s) \sigma_\Theta^s(\phi \varphi) \, ds.
$$
Therefore, its right support satisfies
\[
  r\big(T_{\pi_\Theta(\psi_\delta)}(\phi \varphi) \big)
  \le
  \bigvee_{s \in \mathrm{B}_{2\delta}(0)} \sigma_\Theta^s(r(\varphi)) \quad \Rightarrow \quad p_\delta \le \bigvee_{s \in \mathrm{B}_{2\delta}(0)} \sigma_\Theta^s(r(\varphi)).
\]
Hence, since $p_\delta \ge r(\varphi)$, we conclude by admissibility that $\displaystyle \mathrm{w}^*\mbox{-}\lim_{\delta \to 0} p_\delta = r(\varphi)$. \fin

\noindent Given $T_k \in \mathcal{B}(L_2(\RR_\Theta)) \subset \mathcal{L}(\S_\Theta, \S_\Theta')$, we truncate it as follows $$T_{\Delta,\delta}^\ell = M_{\Psi_{\Delta,\delta}}^\ell(T_k) \quad \mbox{and} \quad T_{\Delta,\delta}^r = M_{\Psi_{\Delta,\delta}}^r(T_k).$$ According to Remark \ref{L2BoundPsi} both truncations $T_{\Delta, \delta}^\dag$ satisfy the $L_2$-estimate $$\big\| T_{\Delta, \delta}^\dag: L_2(\RR_\Theta) \to L_2(\RR_\Theta) \big\| \le 2 \, \big\| \widehat{\psi} \big\|_1 \big\| T_k: L_2(\RR_\Theta) \to L_2(\RR_\Theta) \big\|.$$ In particular, the Banach-Alaoglu theorem confirms that certain subfamily of our truncations $T_{\Delta, \delta}^\dag$ converges to some $L_2$-bounded operator $S_k^\dag: L_2(\RR_\Theta) \to L_2(\RR_\Theta)$ for $\dag \in \{\ell,r\}$. We shall assume for simplicity of notation that the whole family of truncations converges to $S_k^\dag$ as $\Delta \to \infty$ and $\delta \to 0$. 

\begin{theorem} \label{PrincipalValue}
  There exist $z_\dag \in \RR_\Theta$ such that
  \[
    (T_k - S_k^\ell) (a) = a z_\ell \quad \mbox{and} \quad (T_k - S_k^r) (a) = z_r a.
  \]
\end{theorem}

\dem Given an admissible projection $p$ and by Remark \ref{Kernel}
\[
  \Big\langle T_{\Delta,\delta}^r (\lambda_\Theta(f)p), \lambda_\Theta(g) \Big\rangle
  = \Big\langle \hskip1pt k, \hskip1pt \pi_\Theta(\Psi_{\Delta,\delta}) \bullet \big( \lambda_\Theta(g) \otimes \lambda_\Theta(f)p \big) \Big\rangle.
\]
Since $\pi_\Theta(\chi_{\mathrm{B}_{\mathrm{R}}(0)})$ converges to $\1$ in the strong operator topology, we can safely
assume that $\lambda_\Theta(g) \otimes \lambda_\Theta(f) p = (\lambda_\Theta(g) \otimes p) \bullet (\1 \otimes \lambda_\Theta(f))$ is left supported by $\pi_\Theta(\chi_{\mathrm{B}_{\mathrm{R}}(0)})$ for $\mathrm{R}$ large enough. Then we have
\begin{eqnarray*}
  \lefteqn{\Big\langle (T_k - T_{\Delta,\delta}^r) (\lambda_\Theta(f)p), \lambda_\Theta(g) \Big\rangle} \\
  \!\!\! & = & \!\!\! \Big\langle \hskip1pt k, \hskip1pt \pi_\Theta(\1 - \psi_{\Delta,\delta}) \bullet \big( \lambda_\Theta(g) \otimes \lambda_\Theta(f)p \big) \Big\rangle \\
  \!\!\! & = & \!\!\! \Big\langle \hskip1pt k, \hskip1pt \pi_\Theta(\psi_{\delta}) \bullet \big( \lambda_\Theta(g) \otimes \lambda_\Theta(f)p \big) \Big\rangle + \Big\langle \hskip1pt k, \hskip1pt \pi_\Theta(\1 - \psi_{\Delta}) \bullet \big( \lambda_\Theta(g) \otimes \lambda_\Theta(f) \big) \Big\rangle.
\end{eqnarray*}
Since $\ell(\lambda_\Theta(g) \otimes \lambda_\Theta(f)p) \leq \pi_\Theta(\psi_\Delta)$ for large $\Delta$, the second term vanishes when $\Delta$ is large. The identity $\lambda_\Theta(g) \otimes \lambda_\Theta(f) p = (\1 \otimes p) \bullet (\lambda_\Theta(g) \otimes \lambda_\Theta(f))$ allows to apply Lemma \ref{LocalizationL} to get $$\pi_\Theta(\psi_{\delta}) \bullet \big( \lambda_\Theta(g) \otimes \lambda_\Theta(f) p \big) = \pi_\Theta(\psi_{\delta}) \bullet \big( p_\delta \lambda_\Theta(g) \otimes \lambda_\Theta(f) p \big)$$ for some projection $p_\delta$ converging to $p$ in the weak$-*$ topology. This gives
\begin{eqnarray*}
\lefteqn{\Big\langle (T_k - T_{\Delta,\delta}^r) (\lambda_\Theta(f) p), \lambda_\Theta(g) \Big\rangle} \\ \!\!\! & = & \!\!\! \Big\langle \hskip1pt k, \hskip1pt \pi_\Theta(\psi_{\delta}) \bullet \big( p_{\delta} \lambda_\Theta(g) \otimes \lambda_\Theta(f) p \big) \Big\rangle = \Big\langle (T_k - T_{\Delta,\delta}^r) (\lambda_\Theta(f) p) p_{\delta}, \lambda_\Theta(g) \Big\rangle.
\end{eqnarray*}
Taking limits in $\Delta \to \infty$ and $\delta \to 0$, we get $(T_k - S_k^r) (\lambda_\Theta(f) p) = (T_k - S_k^r) (\lambda_\Theta(f))p$ for any admissible projection $p \in \RR_\Theta$. This readily implies that $T_k - S_k^r$ commutes with the von Neumann algebra generated by right multiplication with admissible projections and, by Lemma \ref{LemmAllAlgebra}, we conclude that $T_k - S_k^r$ belongs to the commutant in $\B(L_2(\RR_\Theta))$ of $\RR_\Theta$ acting by right multiplication. Such algebra is given by $\RR_\Theta$ acting on the left and so, there is a unique $z_r \in \RR_\Theta$ such that $(T_k - S_k^r)(a) = z_r a$. A symmetric argument works for $S_k^\ell$, which also satisfies the assertion. \fin

\begin{remark}
\emph{We may also consider two-sided principal values $T_{\Delta,\delta}^\ell T_{\Delta',\delta'}^r$. Taking first a weak-$\ast$ accumulation point in $(\Delta,\delta)$ and then another in $(\Delta',\delta')$ gives an element $S_k$ such that $S_k(\lambda_\Theta(f)) = z_r \lambda_\Theta(f) + \lambda_\Theta(f) z_\ell$, for certain $z_\dag \in \RR_\Theta$. This is the quantum analogue of a basic result in Calder\'on-Zygmund theory, further details can be found in \cite[Proposition 8.1.11]{G2}.}
\end{remark}

\section{\bf Pseudodifferential $L_p$ calculus}
\label{Sect3}

The aim of this section is to establish sufficient smoothness conditions on a given symbol $a: \R^n \to \RR_\Theta$ for the $L_p$-boundedness of the pseudodifferential operator $\Psi_a$ associated to it. This is the content of Theorem B in the Introduction. Sobolev $p$-estimates naturally follow from this analysis. Before that, subtle transference methods will be needed to extend the classical composition and adjoint formulae to the context of quantum Euclidean spaces. The proof of Theorem B is divided into several blocks. We begin with an analysis of $L_2$-boundedness, which includes the quantum forms of Calder\'on-Vaillancourt theorem and Bourdaud's condition stated in Theorem B i) and ii) respectively. Theorem B iii) follows from it and Theorem A, once we prove that $\Psi_a$ is a Calder\'on-Zygmund operator.

\subsection{Adjoint and product formulae}
\label{CompAdj}

Recall that a symbol over $\RR_\Theta$ must be understood as a smooth function $a: \R^n \to \RR_\Theta$ whose  associated pseudodifferential operator takes the form $$\Psi_a(\lambda_\Theta(f)) = \int_{\R^n} a(\xi) f(\xi) \lambda_\Theta(\xi) \, d \xi.$$ Given $m \in \R$ and $0 \le \delta \le \rho \le 1$, the H\"ormander classes $S_{\rho,\delta}^m(\RR_\Theta)$ are $$S_{\rho,\delta}^m(\RR_\Theta) = \Big\{ a: \R^n \to \RR_\Theta \, : \, \big| \partial_\Theta^\beta \, \partial_\xi^{\alpha} a(\xi) \big|\le C_{\alpha,\beta} \langle \xi \rangle^{m - \rho |\alpha| + \delta |\beta|} \ \mbox{for all} \ \alpha, \beta \in \Z_+^n \Big\}.$$ Here we follow standard notation $\langle \xi \rangle = (1 + |\xi|^2)^{1/2}$. Pseudodifferential operators are formally generated by Fourier multipliers and left multiplication operators. It is easy to see that these families of operators generate in turn the whole $\B(L_2(\RR_\Theta))$ as a von Neumann algebra. It is therefore reasonable to think that adjoints and composition of pseudodifferential operators are pseudodifferential operators. Our first goal is to develop asymptotic formulae for adjoints and compositions to justify that the adjoint of a regular $(\delta < \rho)$ pseudodifferential operator of degree $m$ is again a pseudodifferential operator of degree $m$ and that the composition of operators of degrees $m_1$ and $m_2$ yields a pseudodifferential operator of degree $m_1 + m_2$. 

We start by defining $\Psi_a$ in the distributional sense. First, we have $\Psi_a: \S_\Theta \to \S_\Theta$ continuously whenever $a \in \S(\R^n;\S_\Theta)$ is a Schwartz function itself. Indeed, recall that $a \in \S(\R^n;\S_\Theta)$ means that 
\begin{equation} \label{SchwartzSymbol}
a(\xi) = \int_{\R^n} \widehat{a}(z,\xi) \lambda_\Theta(z) \, d z
\end{equation}
for some $\widehat{a} \in \S(\R^n \times \R^n)$. This immediately gives $$\Psi_a(\lambda_\Theta(f)) = \int_{\R^n} \underbrace{\int_{\R^n} \widehat{a}(z-\xi, \xi) f(\xi) e^{2\pi i \langle z-\xi, \Theta_{\downarrow} \xi \rangle} d\xi}_{F(z)} \lambda_\Theta(z) \, d z$$ with $F \in \S(\R^n)$, which implies the assertion. The following lemma refines it.

\begin{lemma} \label{PSD.Apriori}
Given any $a \in S^m_{\rho,\delta}(\RR_\Theta)$, we have that $\Psi_a:\S_\Theta \to \S_\Theta$ continuously. 
\end{lemma}

\dem Note that
\begin{eqnarray*}
\big\| \Psi_a( \lambda_\Theta(f)) \big\|_{\RR_\Theta} & = & \Big\| \int_{\R^n} a(\xi) \langle \xi \rangle^{n+1} f(\xi) \lambda_\Theta(\xi) \, \frac{d\xi}{\langle \xi \rangle^{n+1}} \Big\|_{\RR_\Theta} \\ & \lesssim & \sup_{\xi \in \R^n} \Big\{ \big\| \langle \xi \rangle^{-m} a(\xi) \big\|_{\RR_\Theta} \big| \, \langle \xi \rangle^{n + m + 1} f(\xi) \big| \Big\}.
\end{eqnarray*}
The first term is bounded by the H\"ormander condition with $\alpha = \beta = 0$ and the second one since $f \in \S(\R^n)$. According to Remark \ref{RemQuantumx} and Lemma \ref{ContinuitySchwartz}, it suffices to see that the operators $$P(x_{\Theta}) \partial_\Theta^\beta \Psi_a(\lambda_\Theta(f)) Q(x_\Theta)$$ satisfy similar inequalities for arbitrary monomials $P,Q$ and $\beta \in \Z_+^n$. Recall that $$\partial_{\Theta}^j \Psi_a(\lambda_\Theta(f)) = \Psi_{\partial_{\Theta}^j a} (\lambda_\Theta(f)) + \Psi_a(\partial_{\Theta}^j (\lambda_\Theta(f))),$$ but $\partial_{\Theta}^j a \in S_{\rho,\delta}^{m + \delta}(\RR_\Theta)$ and $\partial_{\Theta}^j (\lambda_\Theta(f)) = \lambda_\Theta(2 \pi i \xi_j f)$. In particular, $\partial_\Theta^\beta \Psi_a(\lambda_\Theta(f))$ behaves as $\Psi_a(\lambda_\Theta(f))$ and we may ignore $\beta$. Thus, it will be enough to illustrate the argument for $(P,Q,\beta) = (1,x_{\Theta,j},0)$ and $(P,Q,\beta) = (x_{\Theta,j},1,0)$. In the first case, since our pseudodifferential operators act by left multiplication of the symbol $a$, the exact same argument given in the proof of Proposition \ref{Quantumx} gives the identity below, even for $a$ taking values in $\RR_\Theta$ as it is the case 
\begin{eqnarray*}
\Psi_a(\lambda_\Theta(f)) x_{\Theta, j} \!\!\! & = & \!\!\! \int_{\R^n} D_{\Theta,j}^r (af)(\xi) \lambda_\Theta(\xi) \, d\xi \\ \!\!\! & = & \!\!\! \int_{\R^n} \Big( a(\xi) D_{\Theta,j}^r (f)(\xi) - \frac{1}{2\pi i } \partial_\xi^ja(\xi) f(\xi) \Big) \lambda_\Theta(\xi) \, d\xi.
\end{eqnarray*}
Clearly $D_{\Theta,j}^r f \in \S(\R^n)$ and $\partial_\xi^j a \in S_{\rho,\delta}^{m - \rho}(\RR_\Theta)$, so we may proceed as above. We need a similar expression when $x_{\Theta,j}$ acts by left multiplication. In this second case we need to be a bit more careful
\begin{eqnarray*}
\lefteqn{x_{\Theta, j} \Psi_a(\lambda_\Theta(f))} \\ 
& = & \!\!\!\! x_{\Theta,j} \int_{\R^n} \Big( \int_{\R^n} \widehat{a}(z,\xi) \lambda_\Theta(z) \, dz \Big) f(\xi) \lambda_\Theta(\xi) \, d\xi \\ 
& = & \!\!\!\! \frac{1}{2\pi i}\frac{d}{ds}\Big|_{s=0}\int_{\R^n}\int_{\R^n}\widehat{a}(z,\xi) f(\xi) \lambda_\Theta(se_j)  \lambda_\Theta(z) \lambda_\Theta(\xi) \, dz d\xi \\
& = & \!\!\!\! \frac{1}{2\pi i} \frac{d}{ds} \Big|_{s=0} \int_{\R^n} \underbrace{\Big[ \int_{\R^n} \widehat{a}(z,\xi) e^{2\pi i s \langle e_j, \Theta z \rangle} \lambda_\Theta(z) dz \Big]}_{a_{js}(\xi)} f(\xi) e^{2\pi i s \langle e_j, \Theta_\downarrow \xi \rangle} \lambda_\Theta(\xi + se_j) d\xi.
\end{eqnarray*}
Equivalently, we may write it as follows 
$$x_{\Theta, j} \Psi_a(\lambda_\Theta(f)) 
\, = \, \frac{1}{2\pi i} \frac{d}{ds} \Big|_{s=0} \int_{\R^n} a_{js}(\xi - se_j) f(\xi - se_j) e^{2\pi i s \langle e_j, \Theta_\downarrow \xi \rangle} \lambda_\Theta(\xi) d\xi.$$ 
In particular, Leibniz rule and the argument in Proposition \ref{Quantumx} give
\begin{eqnarray*}
x_{\Theta, j} \Psi_a(\lambda_\Theta(f)) 
\!\!\! & = & \!\!\! \frac{1}{2\pi i} \int_{\R^n} \frac{d}{ds} \Big|_{s=0} a_{js}(\xi) f(\xi) \lambda_\Theta(\xi) d\xi \\
\!\!\! & + & \!\!\! \frac{1}{2\pi i} \int_{\R^n} \frac{d}{ds} \Big|_{s=0} \Big( a(\xi - se_j) f(\xi - se_j) e^{2\pi i s \langle e_j, \Theta_\downarrow \xi \rangle} \Big) \lambda_\Theta(\xi) d\xi \\
\!\!\! & = & \!\!\! \frac{1}{2\pi i} \int_{\R^n} \Big\{ \Big( \sum_{k=1}^n \Theta_{jk} \partial_\Theta^k a(\xi) \Big) f(\xi) + 2 \pi i D_{\Theta,j}^\ell (af)(\xi) \Big\} \lambda_\Theta(\xi) d\xi.
\end{eqnarray*}
Since the new terms $\partial_\Theta^k a \in S_{\rho,\delta}^{m+\delta}(\RR_\Theta)$, we may proceed as above once more. \fin

Consider a pair of symbols $a_1, a_2: \R^n \to \RR_\Theta$. In order to properly identify $\Psi_{a_j}$ with $a_j$, we need to confirm that $\Psi_{a_1} = \Psi_{a_2}$ implies that $a_1 = a_2$. This is the case when the symbols $a_j$ are of polynomial growth |there exists $k \ge 0$ such that $|a_j(\xi)| \le C_j \langle \xi \rangle^k$| and $\Psi_{a_1} = \Psi_{a_2}$ holds as operators in $\mathcal{B}(\S_\Theta,\S_\Theta)$. This result will be enough for our purposes and it follows by an elementary application of Fourier inversion for distributions, which we omit. 

\begin{lemma} \label{PSD.direct}
Given $a, a_1, a_2 \in \S(\R^n; \RR_\Theta)$, we find$\hskip1pt :$
\begin{itemize}
\item[\emph{i)}] $\Psi_a^{\ast} = \Psi_{a_\dagger^\ast}$ where $$a_\dagger(\xi) = \int_{\R^n} \widehat{a}(z,\xi - z) \lambda_\Theta(z) \, d z.$$

\item[\emph{ii)}] $\Psi_{a_1} \circ \Psi_{a_2} = \Psi_{a_1 \diamond a_2}$ where $$(a_1 \diamond a_2)(\xi) = \int_{\R^n} a_1(z) \widehat{a}_2(z-\xi,\xi) \lambda_\Theta(z - \xi) \, d z.$$
\end{itemize}
\end{lemma}

\dem By Lemma \ref{KernelL} i) $\Psi_a^\ast = T_{k_a}^\ast = T_{\mathrm{flip}(k_a)^\ast}$. By \eqref{SchwartzSymbol} 
\begin{eqnarray*}
\mathrm{flip}(k_a)^* & = & \int_{\R^n}  \lambda_\Theta(\xi) \otimes \lambda_\Theta(\xi)^\ast a(\xi)^\ast \, d \xi\\
                     & = & \int_{\R^n} \int_{\R^n} \overline{\widehat{a}(z,\xi)} \, e^{- 2 \pi i \langle z, \Theta_{\downarrow} \xi \rangle} \lambda_\Theta(\xi) \otimes \lambda_\Theta(z + \xi)^\ast \, dz \, d\xi\\
                     & = & \int_{\R^n} \Big( \int_{\R^n} \widehat{a}(z,\xi - z) \lambda_\Theta(z) dz \Big)^\ast \lambda_\Theta(\xi) \otimes \lambda_\Theta(\xi)^\ast \, d\xi = k_{a_\dagger^\ast},
\end{eqnarray*}
which implies $\Psi_a^\ast = \Psi_{a_\dagger^\ast}$. The composition formula is obtained similarly. \fin

The formulas above are difficult to treat directly. Following the classical setting we introduce double pseudodifferential operators. Namely, if $A: \R^n \to \RR_\Theta \bar\otimes \RR_\Theta^\op$ is a double symbol, its associated operator is given by $$\mathcal{D}_A(\varphi) = (id \otimes \tau_\Theta) \Big\{ \Big( \int_{\R^n} A(\xi) \bullet \big( \lambda_\Theta(\xi) \otimes \lambda_\Theta(\xi)^* \big) \, d \xi \Big) \, (\1 \otimes \varphi) \Big\}.$$ Observe that if the double symbol is of the form $A(\xi) = a(\xi) \otimes \1$ then $\mathcal{D}_{A} = \Psi_{a}$. The advantage of the above class of operators is that they admit simpler expressions for adjoints $\mathcal{D}_{a \otimes \1}^{\ast} = \mathcal{D}_{\1 \otimes a^\ast}$ and products $\mathcal{D}_{a_1 \otimes \1} \circ \mathcal{D}_{\1 \otimes a_2} = \mathcal{D}_{a_1 \otimes a_2}$. We now introduce extended H\"ormander classes for double symbols. To that end, we recall the definition of the Haagerup tensor product. Given $z \in \RR_\Theta \otimes_{\mathrm{alg}} \RR_\Theta$, let us define
$$\|z\|_{\RR_\Theta \otimes_h \RR_\Theta} = \inf \Big\{
             \Big\| \summ_{j} x_j x_j^{\ast} \Big\|_{\RR_\Theta}^\frac12 \, \Big\| \summ_{j} y_j^{\ast} y_j \Big\|_{\RR_\Theta}^\frac12 
               : \ z = \summ_{j} x_j \otimes y_j
             \Big\}.$$
The Haagerup tensor product $\RR_\Theta \otimes_h \RR_\Theta$ is defined by completion and also admits a natural operator space structure \cite{P3}. We will say that $A: \R^n \to \RR_\Theta \bar\otimes \RR_\Theta^\mathrm{op}$ belongs to $S^m_{\rho, \delta_1, \delta_2}(\RR_\Theta)$ when 
$$\Big\|(\partial_{\Theta}^{\beta_1} \otimes \partial_{\Theta}^{\beta_2}) \, \partial_{\xi}^\alpha \! A(\xi) \Big\|_{\RR_\Theta \otimes_h \RR_\Theta} \le C_{\alpha, \beta_1, \beta_2} \, \langle \xi \rangle^{m - \rho |\alpha| + \delta_1 |\beta_1| + \delta_2 |\beta_2|},$$
for all multindices $\alpha, \beta_1, \beta_2$. Our next result provides a compression map $$\mathfrak{B}: \S(\R^n; \S_\Theta \otimes_\pi \S_\Theta) \to \S(\R^n; \S_\Theta),$$ $$\mathfrak{B} (A)(\xi) = (id \otimes \tau_\Theta) \int_{\R^n} A(\eta) \bullet (\lambda_\Theta(\eta) \otimes \lambda_\Theta(\eta)^\ast) \bullet (\lambda_\Theta(\xi)^\ast \otimes \lambda_\Theta(\xi)) \, d\eta,$$ which sends double symbols into symbols inducing the same operators. This map involves in turn the map $m: \S_\Theta \otimes_\pi \S_\Theta \to \S_\Theta$ defined by linear extension of $\varphi_1 \otimes \varphi_2 \mapsto \varphi_1 \varphi_2$. If $\Theta = 0$, $m$ is the restriction to the diagonal $\varphi(x,y) \mapsto \varphi(x,x)$ which extends to a positive preserving contraction with the $\mathrm{C}^\ast$-norm. This fails in general for nonabelian algebras. Instead, the Haagerup tensor product can be understood as the smallest (operator space) tensor product making the operation $m$ continuous. This justifies our use of the Haagerup tensor product in the above definition of double H\"ormander classes. Define $$L_{\Theta,\xi} = \exp \Big( \frac{1}{2\pi i} \sum_{j=1}^n \partial_\xi^j \otimes id_{\RR_\Theta} \otimes \partial_{\Theta}^j \Big) \in \mathcal{B} \big( \S(\R^n; \S_\Theta \otimes_\pi \S_\Theta) \big).$$

\begin{theorem} \label{PSD.compress}
The compression map $\mathfrak{B}: \S(\R^n; \S_\Theta \otimes_\pi \S_\Theta) \to \S(\R^n; \S_\Theta)$ above satisfies the identity $\mathcal{D}_A = \Psi_{\mathfrak{B}(A)}$ as operators in $\mathcal{B}(\S_\Theta, \S_\Theta')$. In addition, the following identities hold 
\begin{eqnarray*}
\mathfrak{B}(A)(\xi) & = & \mathcal{D}_A(\lambda_\Theta(\xi)) \lambda_\Theta(\xi)^* \\ & = & m \big( L_{\Theta,\xi} A(\xi) \big) \ \sim \ \sum_{\gamma \in \Z_+^n} \frac{m \big((\partial_\xi^\gamma \otimes id_{\RR_\Theta} \otimes \partial_\Theta^\gamma)A(\xi) \big)}{(2\pi i)^{|\gamma|} \gamma!}.
\end{eqnarray*}
Moreover, given $A \in S_{\rho,\delta_1,\delta_2}^m(\RR_\Theta)$ with $\delta_2 < \rho$ and $\mathrm{N} \in \Z_+$ large we find $$\Big\| \mathfrak{B}(A)(\xi) - \sum_{|\gamma| < \mathrm{N}} \frac{m \big( (\partial_\xi^\gamma \otimes id \otimes \partial_\Theta^\gamma) A(\xi) \big)}{(2 \pi i)^{|\gamma|} \gamma!} \Big\|_{\RR_\Theta} \, \le \, C_\mathrm{N} \, \langle \xi \rangle^{m+n-(\rho - \delta_2) \mathrm{N}}.$$
In particular, $\mathfrak{B}: S^m_{\rho, \delta_1, \delta_2}(\RR_\Theta) \to S^m_{\rho, \delta}(\RR_\Theta)$ for $\delta = \max\{\delta_1, \delta_2\}$ whenever $\delta_2 < \rho$.
\end{theorem}

\dem The proof is divided into three blocks:

\noindent \textbf{A. Expressions for $\mathfrak{B}(A)$.} Clearly $\mathfrak{B}(A)(\xi) = \mathcal{D}_A(\lambda_\Theta(\xi)) \lambda_\Theta(\xi)^*$, so $$\mathcal{D}_A(\varphi) = \mathcal{D}_A \Big( \int_{\R^n} \widehat{\varphi}(\xi) \lambda_\Theta(\xi) \, d\xi \Big) = \int_{\R^n} \widehat{\varphi}(\xi) \mathfrak{B}(A)(\xi) \lambda_\Theta(\xi) \, d\xi = \Psi_{\mathfrak{B}(A)}(\varphi)$$ for any $\varphi \in \S_\Theta$. To prove the identity $\mathfrak{B}(A)(\xi) = m ( L_{\Theta,\xi} A(\xi) )$, we write $$A(\xi) = \int_{\R^n} \widehat{A}(u) e^{2\pi i \langle u,\xi \rangle} du = \int_{\R^n} \Big( \int_{\R^n} \widetilde{A}(u,v) \otimes \lambda_\Theta(v) \, dv \Big) e^{2\pi i \langle u,\xi \rangle}  \, du$$ 
where $\widehat{A}$ is the Euclidean Fourier transform of $A: \R^n \to \RR_\Theta \otimes \RR_\Theta^\mathrm{op}$ and $\widetilde{A}$ is the quantum (partial) Fourier transform of it in the second tensor. In other words, we have 
\begin{eqnarray*}
\widetilde{A}(u,v) & = & (id \otimes \tau_\Theta) \Big( \widehat{A}(u) \bullet (\1 \otimes \lambda_\Theta(v)^*) \Big) \\ & = & \int_{\R^n} (id \otimes \tau_\Theta) \Big( A(s) \bullet (\1 \otimes \lambda_\Theta(v)^*) \Big) e^{-2\pi i \langle u,s \rangle} \, ds.
\end{eqnarray*}
Now, using the Taylor series expansion $$L_{\Theta,\xi} = \sum_{k=0}^\infty \frac{(2\pi i)^k}{k!} \sum_{j_1, j_2, \ldots, j_k=1}^n \Big( \prod_{s=1}^k \frac{\partial^{j_s}_\xi}{2\pi i} \Big) \otimes id \otimes \Big( \prod_{s=1}^k \frac{\partial_{\Theta}^{j_s}}{2\pi i} \Big)$$ we easily get the following identity for $L_{\Theta,\xi}A$ $$L_{\Theta,\xi}A = \int_{\R^n} \Big( \int_{\R^n} \widetilde{A}(u,v) \otimes \lambda_\Theta(v) \, dv \Big) e^{2\pi i \langle u,\xi + v \rangle} \, du.$$
Applying $m$ to this expression and writing $\widetilde{A}$ in terms of $A$, we get
\begin{eqnarray*}
m(L_{\Theta,\xi}A) & = & (id \otimes \tau_\Theta) \int_{\R^n \times \R^n \times \R^n} A(s) \bullet \pi_\Theta(\exp_v) \, e^{2 \pi i \langle u, \xi + v -s \rangle} \, ds du dv \\
& = & (id \otimes \tau_\Theta) \int_{\R^n \times \R^n} A(s) \bullet \pi_\Theta(\exp_v) \, \delta_\xi(s-v) \, ds dv \\
& = & (id \otimes \tau_\Theta) \int_{\R^n} A(s) \bullet \pi_\Theta(\exp_{s-\xi}) \, ds.
\end{eqnarray*}
This proves $\mathfrak{B}(A) = m(L_{\Theta,\xi}A)$. On the other hand $$L_{\Theta,\xi} = \sum_{\gamma \in \Z_+^n} \frac{\partial_\xi^\gamma \otimes id_{\RR_\Theta} \otimes \partial_{\Theta}^\gamma}{(2\pi i)^{|\gamma|} \gamma!}$$ by standard modification of the Taylor series. This gives the formal series expansion. 

\noindent \textbf{B. Estimate for the remainder.} Thus, our next goal is to justify the Taylor remainder estimate in the statement. This requires yet another expression for $\mathfrak{B}(A)$. We begin by noticing that
\begin{eqnarray*}
\mathfrak{B}(A)(\xi) & = & m \Big\{ (id \otimes \tau_\Theta) \Big( \int_{\R^n} A(\eta) \bullet \pi_\Theta(\exp_{\eta-\xi}) \, d\eta \Big) \otimes \1 \Big\} \\ & = & m \Big\{ \int_{\R^n} \int_{\R^n} (id \otimes \sigma_\Theta^z) \big( A(\eta) \bullet \pi_\Theta(\exp_{\eta-\xi}) \big) \, d\eta dz \Big\} \\ & = & m \Big\{ \int_{\R^n} \int_{\R^n} (id \otimes \sigma_\Theta^z) \big( A(\eta) \big) \bullet \pi_\Theta(\exp_{\eta-\xi}) e^{-2\pi i \langle z,\eta-\xi \rangle} \, d\eta dz \Big\}.
\end{eqnarray*}
The first identity follows from \textbf{A} above. The second identity reduces to $$\int_{\R^n} \sigma_\Theta^z(\varphi) \, dz = \tau_\Theta(\varphi) \, \1 \quad \mbox{with} \quad \sigma_\Theta^z(\lambda_\Theta(\zeta)) = e^{2\pi i \langle z,\zeta \rangle} \lambda_\Theta(\zeta)$$ and the last one since $\sigma_\Theta^z$ is a $*$-homomorphism. Using $m(\mathsf{A} \bullet \pi_\Theta(\exp_\zeta)) = m(\mathsf{A})$ we get
\begin{equation} \label{PSD.OscilatoryB}
\mathfrak{B}(A)(\xi) = \int_{\R^n} \underbrace{\Big( \int_{\R^n} m \big( (id \otimes \sigma^z_\Theta) A(\xi + \eta) \big) e^{-2 \pi i \langle z, \eta \rangle} \, dz \Big)}_{\Omega_\eta(\xi + \eta)} \, d\eta.
\end{equation}
On the other hand, we use $\partial_y^\gamma \sigma_\Theta^y = \sigma_\Theta^y \partial_\Theta^\gamma$ to deduce
\begin{eqnarray*}
\lefteqn{\hskip-20pt m \big( (\partial_\xi^\gamma \otimes id \otimes \partial_\Theta^\gamma) A(\xi) \big)} \\ [5pt]
& = & \partial_\xi^\gamma \partial_y^\gamma \Big|_{y=0} m \big( (id \otimes \sigma_\Theta^y) A(\xi) \big) \\ [5pt]
& = & \partial_\xi^\gamma \partial_y^\gamma \Big|_{y=0} \Big( \int_{\R^n} \int_{\R^n} m \big( (id \otimes \sigma_\Theta^z) A(\xi) \big) e^{2\pi i \langle y-z,\eta \rangle} \, dz d\eta \Big) \\ 
& = & \partial_\xi^\gamma \Big( \int_{\R^n} \int_{\R^n} m \big( (id \otimes \sigma_\Theta^z) A(\xi) \big) e^{-2\pi i \langle z,\eta \rangle} (2\pi i \eta)^\gamma \, dz d\eta \Big) \\ 
& = & \partial_\xi^\gamma \int_{\R^n} \Omega_\eta(\xi) (2\pi i \eta)^\gamma \, d\eta.  
\end{eqnarray*}
This implies that 
$$\mathfrak{B}(A)(\xi) - \sum_{|\gamma| < \mathrm{N}} \frac{m \big( (\partial_\xi^\gamma \otimes id \otimes \partial_\Theta^\gamma) A(\xi) \big)}{(2 \pi i)^{|\gamma|} \gamma!} = \int_{\R^n} \underbrace{\Omega_\eta(\xi + \eta) - \sum_{|\gamma| < \mathrm{N}} \frac{1}{\gamma!} \partial_\xi^\gamma \Omega_\eta(\xi) \eta^\gamma}_{\mathrm{R}_\xi(\eta)} \, d\eta.$$
By Taylor remainder formula 
$$\mathrm{R}_\xi(\eta) = \sum_{|\gamma| = \mathrm{N}} \frac{\mathrm{N}}{\gamma!} \Big( \int_0^1 (1-t)^{\mathrm{N}-1} \partial_s^\gamma \Big|_{s=\eta+t\xi} \Omega_\eta(s) \, dt \Big) \eta^\gamma.$$
In particular, we obtain the following estimate
$$\Big\| \int_{\R^n} \mathrm{R}_\xi(\eta) \, d\eta \Big\|_{\RR_\Theta} \le C_\mathrm{N} \sup_{\begin{subarray}{c} |\gamma| = \mathrm{N} \\ 0 \le t \le 1 \end{subarray}} \Big\| \int_{\R^n} \eta^\gamma \partial_s^\gamma \Omega_\eta(\eta+t\xi) \, d\eta \Big\|_{\RR_\Theta}.$$
Since $\partial_\xi^\gamma$ commutes with $m$, we get the identity
$$\partial_s^\gamma \Omega_\eta(s) = \int_{\R^n} m \big( (id \otimes \sigma^z_\Theta) \partial_s^\gamma A(s) \big) e^{-2 \pi i \langle z, \eta \rangle} \, dz.$$
Next, we use the standard oscillatory integral trick
\begin{eqnarray*}
e^{-2 \pi i \langle z, \eta \rangle} & = & \frac{(- \Delta_\eta)^{n}}{(4 \pi^2 |z|^2)^{n}} e^{-2 \pi i \langle z, \eta \rangle}, \\
e^{-2 \pi i \langle z, \eta \rangle} & = & \frac{(\1 - \Delta_z)^{\frac{\mathrm{N}}{2}}}{(1+ 4 \pi^2 |\eta|^2)^{\frac{\mathrm{N}}{2}}} e^{-2 \pi i \langle z, \eta \rangle}. 
\end{eqnarray*} 
Taking $M(\xi, \eta, z ,t) = m \big( (id \otimes \sigma^z_\Theta) \partial_s^\gamma A(\eta+t\xi) \big)$ and integrating by parts
\begin{eqnarray*}
\lefteqn{\int_{\R^n} \eta^\gamma \partial_s^\gamma \Omega_\eta(\eta+t\xi) \, d\eta} \\
\hskip-5pt & = & \hskip-5pt \int_{\R^n} \int_{\R^n} \frac{\eta^\gamma (\1 - \Delta_z)^{\frac{\mathrm{N}}{2}}}{(1+ 4 \pi^2 |\eta|^2)^{\frac{\mathrm{N}}{2}}} \Big( M(\xi, \eta, z ,t) \Big) e^{-2 \pi i \langle z, \eta \rangle} \, dz d\eta \\
\hskip-5pt & = & \hskip-5pt \int_{\mathrm{B}_1(0)} \Big( \int_{\R^n} \frac{\eta^\gamma (\1 - \Delta_z)^{\frac{\mathrm{N}}{2}}}{(1+ 4 \pi^2 |\eta|^2)^{\frac{\mathrm{N}}{2}}} \Big( M(\xi, \eta, z ,t) \Big) e^{-2 \pi i \langle z, \eta \rangle} \, d\eta \Big) \, dz \\
\hskip-5pt & + & \hskip-5pt \int_{\mathrm{B}_1^c(0)} \Big( \int_{\R^n} \frac{(- \Delta_\eta)^{n}}{(4 \pi^2 |z|^2)^{n}} \Big[ \frac{\eta^\gamma(\1 - \Delta_z)^{\frac{\mathrm{N}}{2}}}{(1+ 4 \pi^2 |\eta|^2)^{\frac{\mathrm{N}}{2}}}\Big( M(\xi, \eta, z ,t) \Big)  \Big] e^{-2 \pi i \langle z, \eta \rangle} \, d\eta \Big) \, dz.
\end{eqnarray*}
Let us write $\Pi_1$ and $\Pi_2$ for the two terms in the right hand side. Then, we use one more time the identity $\partial_z^\gamma \sigma_\Theta^z = \sigma_\Theta^z \partial_\Theta^\gamma$ together with the contractivity of the map $m: \RR_\Theta \otimes_h \RR_\Theta \to \RR_\Theta$ and the $S_{\rho,\delta_1,\delta_2}^m$-- condition. This yields the following inequality for any $|\gamma| = \mathrm{N}$ and $0 \le t \le 1$ 
\begin{eqnarray*}
\|\Pi_1\|_{\RR_\Theta} & \lesssim & \int_{\R^n} \Big\| \big( id \otimes (\1 -\Delta_\Theta)^{\frac{\mathrm{N}}{2}} \big) \partial_s^\gamma A(\eta + t \xi) \Big\|_{\RR_\Theta \otimes_h \RR_\Theta} d\eta \\ & \lesssim & \int_{\R^n} \max \big\{ \langle \xi \rangle, \langle \eta \rangle \big\}^{m - (\rho - \delta_2) \mathrm{N}} \, d\eta \ \lesssim \ \langle \xi \rangle^{m+n - (\rho - \delta_2) \mathrm{N}}.
\end{eqnarray*}
Similarly, $\Pi_2$ is dominated by 
$$\sum_{|\nu_1 + \nu_2| = 2n} \int_{\mathrm{B}_1^c(0) \times \R^n} \Big| \partial_\eta^{\nu_1} \Big( \frac{\eta^\gamma}{\langle \eta \rangle^{\mathrm{N}}} \Big) \Big| \, \Big\| \big( id \otimes (\1 -\Delta_\Theta)^{\frac{\mathrm{N}}{2}} \big) \partial_s^{\gamma + \nu_2} A(\eta + t \xi) \Big\|_{\RR_\Theta \otimes_h \RR_\Theta} \frac{dz d\eta}{|z|^{2n}}$$
which is bounded by $\langle \xi \rangle^{m+n - (\rho - \delta_2) \mathrm{N}}$. This completes the estimate of the remainder.

\noindent \textbf{C. $\mathfrak{B}$ respects the H\"ormander classes.} It remains to show that $\mathfrak{B}(A)$ belongs to the H\"ormander class $S_{\rho,\delta}^m(\RR_\Theta)$ for $\delta = \max \{\delta_1, \delta_2\}$ whenever $A \in S^m_{\rho, \delta_1, \delta_2}(\RR_\Theta)$ and $\delta_2 < \rho$. Since we have $$\partial_\Theta^\beta \circ m = \sum_{\beta_1 + \beta_2 = \beta} \frac{\beta!}{\beta_1!\beta_2!} \, m \circ \big( \partial_\Theta^{\beta_1} \otimes \partial_\Theta^{\beta_2} \big),$$ it turns out that the following inequality holds for any $\gamma \in \Z_+^n$
\begin{eqnarray*}
\lefteqn{\hskip-30pt \Big\| \partial_\Theta^\beta \partial_\xi^\alpha m \big( (\partial_\xi^{\gamma} \otimes id \otimes \partial_\Theta^{\gamma}) A(\xi) \big) \Big\|_{\RR_\Theta}} \\ & \le & \sum_{\beta_1 + \beta_2 = \beta} \frac{\beta!}{\beta_1!\beta_2!} \Big\| \big( \partial_\Theta^{\beta_1} \otimes \partial_\Theta^{\gamma + \beta_2} \big) \partial_\xi^{\gamma + \alpha} A (\xi) \Big\|_{\RR_\Theta \otimes_h \RR_\Theta}. 
\end{eqnarray*}
Since the H\"ormander classes are nested in the degree $m$, this implies that $$\sum_{|\gamma| < \mathrm{N}}^{\null} m \big( (\partial_\xi^{\gamma} \otimes id \otimes \partial_\Theta^{\gamma}) A(\xi) \big) \in \bigcup_{|\gamma| < \mathrm{N}} S_{\rho,\delta}^{m - (\rho-\delta_2)|\gamma|}(\RR_\Theta) = S_{\rho,\delta}^m(\RR_\Theta)$$ as a consequence of $A \in S_{\rho,\delta_1,\delta_2}^m(\RR_\Theta)$, $\delta = \max \{\delta_1, \delta_2\}$ and $\delta_2 < \rho$. Therefore, the inclusion $\mathfrak{B}(A) \in S_{\rho,\delta}^m(\RR_\Theta)$ will follow if there exists a large enough $\mathrm{N} \in \Z_+$ satisfying the inequality $$\Big\| \partial_\Theta^\beta \partial_\xi^\alpha \Big( \mathfrak{B}(A) (\xi) - \sum_{|\gamma| < \mathrm{N}}^{\null} m \big( (\partial_\xi^{\gamma} \otimes id \otimes \partial_\Theta^{\gamma}) A(\xi) \big) \Big) \Big\|_{\RR_\Theta} \le C_{\mathrm{N},\alpha,\beta} \langle \xi \rangle^{m - \rho |\alpha| + \delta |\beta|}.$$ Our estimate for the Taylor remainder above shows that this is indeed the case when $\alpha = \beta = 0$. Using $\partial_\xi^\alpha m = m \partial_\xi^\alpha$ and the commutation formula for $\partial_\Theta^\beta \circ m$ given above, the exact same argument applies for general $\alpha, \beta$. This gives that any $\mathrm{N} \ge n/(\rho - \delta_2)$ works, details are left to the reader. \fin

\begin{corollary} \label{PSD.stabilityA}
The following stability results hold$\hskip1pt :$ 
\begin{itemize}
\item[\emph{i)}] If $a \in S_{\rho, \delta}^m(\RR_\Theta)$ and $\rho < \delta$, then $\Psi_a^* = \Psi_{a_\dag^*}$ with $$a_\dag^* \sim \sum_{\gamma \in \Z_+^n} \frac{\partial_\Theta^\gamma \partial_\xi^\gamma a^*(\xi)}{(2\pi i)^{|\gamma|} \gamma!} \in S_{\rho,\delta}^m(\RR_\Theta).$$

\item[\emph{ii)}] If $a_j \in S_{\rho_j, \delta_j}^{m_j}$, then $\Psi_{a_1} \circ \Psi_{a_2} = \Psi_{a_1 \diamond a_2}$ with $$a_1 \diamond a_2 \sim \sum_{\gamma \in \Z_+^n} \frac{\partial_\xi^\gamma a_1(\xi) \partial_\Theta^\gamma a_2(\xi)}{(2\pi i)^{|\gamma|} \gamma!} \in S_{\rho,\delta}^m(\RR_\Theta)$$ for $m= m_1+m_2$, $\rho = \min \{\rho_1,\rho_2\}$ and $\delta = \max \{\delta_1,\delta_2\}$ when $\delta_2 < \rho$.
\end{itemize}
\end{corollary}

\dem Recall that $$\Psi_a^\ast = \mathcal{D}_{a \otimes \1}^\ast = \mathcal{D}_{\1 \otimes a^\ast} = \Psi_{\mathfrak{B}(\1 \otimes a^\ast)}.$$ If $a \in S^m_{\rho,\delta}(\RR_\Theta)$ then $a \otimes \1 \in S^m_{\rho,\delta,0}(\RR_\Theta)$ and $\1 \otimes a^\ast \in S_{\rho,0,\delta}^m(\RR_\Theta)$. By Theorem \ref{PSD.compress} we have that $a_\dag^\ast = \mathfrak{B}(\1 \otimes a^\ast) \in S^m_{\rho,\delta}(\RR_\Theta)$. The second assertion follows similarly by recalling that 
\begin{eqnarray*}
\Psi_{a_1 \diamond a_2} & = & \Psi_{a_1} \circ \Psi_{a_2}^{\ast \ast} \ = \ \Psi_{a_1} \circ \Psi_{a_{2\dag}^\ast}^\ast \\
& = & \mathcal{D}_{a_1 \otimes \1} \circ \mathcal{D}_{\1 \otimes a_{2\dag}} = \mathcal{D}_{a_1 \otimes a_{2\dag}} \ = \ \Psi_{\mathfrak{B}(a_1 \otimes a_{2\dag})}.
\end{eqnarray*}
Indeed, according to the first assertion, we know that $a_1 \otimes a_{2\dag} \in S_{\rho,\delta_1,\delta_2}^m(\RR_\Theta)$. The asymptotic expansions also follow easily from Theorem \ref{PSD.compress} using the identities $a_\dag^* = \mathfrak{B}(\1 \otimes a^\ast)$ and $a_1 \diamond a_2 = \mathfrak{B}(a_1 \otimes a_{2\dag})$, see e.g. \cite{T1} for a similar approach. \fin

\begin{remark} \label{PSD.remarkCompSigma}
\emph{A natural question is whether the classes $\Sigma_{\rho,\delta}^m(\RR_\Theta)$ are closed under products and adjoints for $\delta < \rho$. This question is still open. Indeed, proceeding as for $S^m_{\rho,\delta}(\RR_\Theta)$ we may define a new class $\Sigma_{\rho, \delta_1, \delta_2}^m(\RR_\Theta)$ of mixed double symbols $A: \R^n \to \RR_\Theta \otimes_{h} \RR_\Theta$ satisfying the condition 
$$\Big\|(\partial_{\Theta}^{\beta_1} \otimes \partial_{\Theta}^{\beta_2}) \, \partial_{\Theta, \xi}^{\alpha_1} \, \partial_{\xi}^{\alpha_2} \! A(\xi) \Big\|_{\RR_\Theta \otimes_h \RR_\Theta} \le C_{\alpha_1, \alpha_2, \beta_1, \beta_2} \, \langle \xi \rangle^{m - \rho |\alpha_1+\alpha_2| + \delta_1 |\beta_1| + \delta_2 |\beta_2|},$$
where, abusing of notation, $\partial_{\Theta, \xi}^j$ acts on $\S(\R^n; S(\RR_\Theta) \otimes_\pi \S(\RR_\Theta))$ as follows
\begin{eqnarray*}
(\partial_{\Theta, \xi}^j A)(\xi) & = & \partial_\xi^j A (\xi) + 2 \pi i \big[ A(\xi), \mathrm{d}_{\Theta, j} \big] \\
& = & \pi_\Theta(\exp_\xi)^\ast \bullet \partial_\xi^j \big\{ \pi_\Theta(\exp_\xi) \bullet A(\xi) \bullet \pi_\Theta(\exp_\xi)^\ast \big\} \bullet \pi_\Theta(\exp_\xi).
\end{eqnarray*}
The operator $\mathrm{d}_{\Theta, j}$ is just $x_{\Theta, j} \otimes \1 - \1 \otimes x_{\Theta, j}$. We shall identify the first term with $x_{\Theta, j}$ and the second with $y_{\Theta, j}$. Of course, we expect that our contraction map satisfies $\mathfrak{B}: \Sigma_{\rho, \delta_1, \delta_2}^m(\RR_\Theta) \to \Sigma_{\rho, {\delta_1 \vee \delta_2}}^m(\RR_\Theta)$ for $\delta_2 < \rho$. Unfortunately, our argument above does not admit a direct generalization. The problem arises since the automorphism $\sigma_\Theta$ in the oscillatory integral \eqref{PSD.OscilatoryB} for $\mathfrak{B}$ does not commute with
$\partial_{\Theta, \xi}^j$. We refer to Lemmas \ref{LemmaNegativeIndex} and \ref{LemmaParametrix} and Remark \ref{RemParametrix} for the calculus of parametrices in this setting. On the other hand, a minimum stability for products |necessary for our Sobolev $p$-estimates, see the proof of Corollary \ref{PSD.lpBound2}| does hold. Namely, if $a_1 \in \Sigma^{m_1}_{\rho_1,0}(\RR_\Theta)$ takes values in $\C \1$ or, more generally, in the center of $\RR_\Theta$, we have that
$$a_1 \diamond a_2 \in \Sigma_{{\rho_1 \wedge \rho_2},\delta}^{m_1 + m_2}(\RR_\Theta) \quad \mbox{ whenever } \quad a_2 \in \Sigma_{\rho_2,\delta}^{m_2}(\RR_\Theta).$$
In particular, composition with polynomials of $\partial_\Theta^j$'s transforms degrees as expected.}
\end{remark}

\subsection{$L_2$-boundedness: Sufficient conditions}
\label{PSD-L2}

We now explore $L_2$-boundedness of pseudodifferential operators in $S_{\rho,\delta}^0(\RR_\Theta)$. Since $S_{\rho,\delta}^0(\RR_\Theta) \subset S_{\delta,\delta}^0(\RR_\Theta) \cap S_{\rho,\rho}^0(\RR_\Theta)$ it suffices to study $L_2$-boundedness for exotic $0 \le \delta = \rho < 1$ and forbidden $\delta = \rho = 1$ symbols. The first case $\rho < 1$ requires a quantum analogue of the celebrated Calder\'on-Vaillancourt theorem \cite{CV}. The second one also requires an additional assumption extending Bourdaud's condition \cite{Bo}, which can be regarded as a form of the $T(1)$ theorem for pseudodifferential operators. 

\subsubsection{The Calder\'on-Vaillancourt theorem in $\RR_\Theta$}
\label{PSD-CV}

As in the Euclidean setting, the hardest part in proving a quantum form of Calder\'on-Vaillancourt theorem is still the case $\rho = 0$. Our argument follows from a combination of \cite{C,T1} adapted to $\RR_\Theta$ which demands a careful argument due to the presence of a $\Theta$-phase. Given $a \in S_{0,0}^0(\RR_\Theta)$, the first step consists in decomposing the symbol as follows. The double Fourier transform of $a$ in the quantum and classical variables $(x_\Theta,\xi)$ is given by 
\begin{eqnarray*}
\widehat{\widehat{a}}(z,\zeta) & = & \int_{\R^n} \tau_\Theta \big( a(\xi) \lambda_\Theta(z)^*\big) e^{- 2\pi i \langle \xi, \zeta \rangle} d\xi \\ & = & \underbrace{\big( 1 + 4 \pi^2 |z|^2 \big)^{\mathrm{N}} \big( 1 + 4 \pi^2 |\zeta|^2 \big)^{\mathrm{N}} \hskip1pt \widehat{\widehat{a}}(z,\zeta)}_{\widehat{\widehat{b}}(z,\zeta)} \underbrace{\big( 1 + 4 \pi^2 |z|^2 \big)^{-\mathrm{N}}  \big( 1 + 4 \pi^2 |\zeta|^2 \big)^{-\mathrm{N}}}_{\widehat{\widehat{g}}(z,\zeta)}.
\end{eqnarray*}
Here we fix $\mathrm{N}$ large enough. We shall also use the terminology $$\widehat{a}(z,\xi) = \tau_\Theta \big( a(\xi) \lambda_\Theta(z)^* \big) = \int_{\R^n} \widehat{\widehat{a}}(z,\zeta) e^{2\pi i \langle \xi, \zeta \rangle} d\zeta$$ for $a,b$ and $g$. In order to express $\Psi_a$ in terms of $b$ and $g$ we need to introduce two auxiliary maps. The first one is a left-module extension $\Pi_\Theta: \RR_\Theta \bar\otimes \RR_\Theta^\mathrm{op} \to \RR_\Theta \bar\otimes \RR_\Theta^\mathrm{op}$ of the $*$-homomorphism $\pi_\Theta$ defined as follows $$\Pi_\Theta \big( \lambda_\Theta(\xi) \otimes \lambda_\Theta(\eta) \big) = \lambda_\Theta(\xi) \otimes \lambda_\Theta(\xi)^* \lambda_\Theta(\eta) = \big( \1 \otimes \lambda_\Theta(\eta) \big) \bullet \pi_\Theta(\exp_\xi).$$ $\Pi_\Theta^* \big( \lambda_\Theta(\xi) \otimes \lambda_\Theta(\eta) \big) = \lambda_\Theta(\xi) \otimes \lambda_\Theta(\xi) \lambda_\Theta(\eta) = ( \1 \otimes \lambda_\Theta(\eta) ) \bullet ( \lambda_\Theta(\xi) \otimes \lambda_\Theta(\xi) )$ gives the adjoint with respect to the module bracket $\langle \langle \alpha, \beta \rangle \rangle = (\tau_\Theta \otimes id) (\alpha \bullet \beta^*)$. The second one is the left-modulation map $M_\eta(\varphi_1 \otimes \varphi_2) = \lambda_\Theta(\eta) \varphi_1 \otimes \varphi_2$ with adjoint $M_\eta^*(\varphi_1 \otimes \varphi_2) = \lambda_\Theta(\eta)^* \varphi_1 \otimes \varphi_2$ with respect to the same bracket above. In the next result we shall use the following symbol 
\begin{eqnarray*}
g_\eta(\xi) & = & \int_{\R^n} \widehat{g}_\eta(z,\xi) \lambda_\Theta(z) \, dz \\ 
& = & \int_{\R^n} \widehat{g}(z,\xi) e^{-2\pi i \langle \xi, \Theta_\downarrow z \rangle} e^{2\pi i \langle \Theta \eta, z \rangle} \lambda_\Theta(z) \, dz.
\end{eqnarray*} 

\begin{lemma} \label{LemmaConvolution}
If $\Phi_\eta = \Pi_\Theta^* \circ M_\eta^*$, the following identity holds $$\Psi_a (\varphi) = (id \otimes \tau_\Theta) \int_{\R^n} (\1 \otimes b(\eta)) \Big( \Phi_\eta^* \circ \big( \Psi_{g_\eta} \otimes id \big) \circ \Phi_\eta \Big) (\varphi \otimes \1) \, d\eta.$$
\end{lemma}

\dem We first claim that $$a(\xi) = (id \otimes \tau_\Theta) \int_{\R^n} (\1 \otimes b(\eta)) \Big( \underbrace{\int_{\R^n} \widehat{g}(z,\xi-\eta) \pi_\Theta(\exp_z) \, dz}_{\Gamma(\xi - \eta)} \Big) \, d\eta.$$ Indeed, writing the symbol $a$ in terms of $b$ and $g$ we obtain 
\begin{eqnarray*}
a(\xi) & = & \int_{\R^n} \int_{\R^n} \widehat{\widehat{b}}(z,\zeta) \widehat{\widehat{g}}(z,\zeta) e^{2\pi i \langle \xi, \zeta \rangle} \lambda_\Theta(z) \, dz d\zeta \\
& = & \int_{\R^n} \int_{\R^n} \widehat{b}(z,\eta) \widehat{g}(z,\xi - \eta) \lambda_\Theta(z) \, dz d\eta.
\end{eqnarray*}
Now the claim follows from the quantum form of convolution via the identity $$\int_{\R^n} f_1(z) f_2(z) \lambda_\Theta(z) \, dz = (id \otimes \tau_\Theta) \Big\{ \Big( \1 \otimes \lambda_\Theta(f_1) \Big) \Big( \int_{\R^n} f_2(z) \pi_\Theta(\exp_z) \, dz \Big) \Big\}.$$ Next we use the claim to produce an expression for $\Psi_a(\varphi)$. Namely, we have 
\begin{eqnarray*}
\lefteqn{\Psi_a(\varphi) = (id \otimes \tau_\Theta) \int_{\R^n} (a(\xi) \otimes \varphi) \pi_\Theta(\exp_\xi) \, d\xi} \\
& = & (id \otimes \tau_\Theta \otimes \tau_\Theta) \int_{\R^n} \int_{\R^n} \Big[ (\1 \otimes b(\eta)) \Gamma(\xi-\eta) \otimes \varphi \Big] \pi_\Theta(\exp_\xi)_{[13]} \, d\xi d\eta \\
& = & (id \otimes \tau_\Theta) \int_{\R^n} (\1 \otimes b(\eta)) \Big\{ (id \otimes id \otimes \tau_\Theta) \int_{\R^n} (\Gamma(\xi-\eta) \otimes \varphi) \pi_\Theta(\exp_\xi)_{[13]} \, d\xi \Big\} \, d\eta
\end{eqnarray*} 
with $(a \otimes b)_{[13]} = a \otimes \1 \otimes b$. The assertion reduces to prove the following identity $$\mathrm{A} := (id \otimes id \otimes \tau_\Theta) \int_{\R^n} (\Gamma(\xi-\eta) \otimes \varphi) \pi_\Theta(\exp_\xi)_{[13]} \, d\xi = \Big( \Phi_\eta^* \circ \big( \Psi_{g_\eta} \otimes id \big) \circ \Phi_\eta \Big) (\varphi \otimes \1) =: \mathrm{B}.$$
Expanding $\Gamma(\xi - \eta)$ it is clear that $$\mathrm{A} = \int_{\R^n} \int_{\R^n} \widehat{\varphi}(\xi) \widehat{g}(z,\xi-\eta) \lambda_\Theta(z) \lambda_\Theta(\xi) \otimes \lambda_\Theta(z)^* \, dz d\xi.$$ On the other hand, we have the identity $$\mathrm{B} = (\lambda_\Theta(\eta) \otimes \1) \Pi_\Theta \Big\{ \int_{\R^n} \hskip-5pt \Big( g_\eta(\xi) \otimes \underbrace{(\tau_\Theta \otimes id) \big( \Phi_\eta(\varphi \otimes \1) (\lambda_\Theta(\xi)^* \otimes \1) \big)}_{\beta_\eta(\xi)} \Big) (\lambda_\Theta(\xi) \otimes \1) \, d\xi \Big\}$$ where it is easily checked that 
\begin{eqnarray*}
\Phi_\eta(\varphi \otimes \1) & = & \Pi_\Theta^* \Big( \int_{\R^n} \widehat{\varphi}(s) \lambda_\Theta(\eta)^* \lambda_\Theta(s) \otimes \1 \, ds \Big) \\ & = & \int_{\R^n} \widehat{\varphi}(s) e^{2 \pi i \langle \eta, \Theta_\downarrow(\eta-s) \rangle} \lambda_\Theta (s-\eta) \otimes \lambda_\Theta(s - \eta) \, ds,
\end{eqnarray*}
so that $\beta_\eta(\xi) = \widehat{\varphi}(\xi + \eta) e^{-2\pi i \langle \eta, \Theta_\downarrow \xi \rangle} \lambda_\Theta(\xi)$. This yields  
\begin{eqnarray*}
\mathrm{B} \!\!\! & = & \!\!\! (\lambda_\Theta(\eta) \otimes \1) \Pi_\Theta \Big\{ \int_{\R^n} \hskip-5pt \Big( g_\eta(\xi) \otimes \beta_\eta(\xi) \Big) (\lambda_\Theta(\xi) \otimes \1) \, d\xi \Big\} \\
\!\!\! & = & \!\!\! (\lambda_\Theta(\eta) \otimes \1) \Pi_\Theta \Big\{ \int_{\R^n} \int_{\R^n} \hskip-5pt \Big( \widehat{g}_\eta(z,\xi) \lambda_\Theta(z) \otimes \beta_\eta(\xi) \Big) (\lambda_\Theta(\xi) \otimes \1) \, dz d\xi \Big\} \\
\!\!\! & = & \!\!\! (\lambda_\Theta(\eta) \otimes \1) \int_{\R^n} \int_{\R^n}  e^{2\pi i \langle z, \Theta_\downarrow \xi \rangle} \widehat{g}_\eta(z,\xi) \lambda_\Theta(z+\xi) \otimes \lambda_\Theta(z+\xi)^* \beta_\eta(\xi) \, dz d\xi \\
\!\!\! & = & \!\!\! \int_{\R^n} \int_{\R^n} \widehat{\varphi}(\xi  + \eta) e^{2\pi i \langle z - \eta, \Theta_\downarrow \xi \rangle} \widehat{g}_\eta(z,\xi) \lambda_\Theta(\eta) \lambda_\Theta(z+\xi) \otimes \lambda_\Theta(z+\xi)^* \lambda_\Theta(\xi) \, dz d\xi.
\end{eqnarray*}
Rearranging and using $\widehat{g}_\eta(z,\xi) = \widehat{g}(z,\xi) e^{-2\pi i \langle \xi, \Theta_\downarrow z \rangle} e^{2\pi i \langle \Theta \eta, z \rangle}$ yields $\mathrm{A} = \mathrm{B}$. \fin

\begin{remark}
\emph{The above lemma may be regarded as the quantum analogue of the identity in \cite[Lemma XIII.1.1]{T1}, whose Euclidean proof is trivial. Unfortunately the quantum analogue gives an extra $\Theta$-phase which vanishes for $\Theta=0$. It is this phase what forces us to be very careful in adapting Cordes argument \cite{C} below.} 
\end{remark}

\begin{lemma} \label{LemmaTraceClass}
$\Psi_{g_\eta}$ admits the factorization $$\Psi_{g_\eta} = A_\eta^* \circ B \circ A_\eta \quad \mbox{with} \quad \sup_{\eta \in \R^n} \big\| A_\eta: L_2(\RR_\Theta) \to L_2(\RR_\Theta) \big\| \, \big\| B \big\|_{S_1(L_2(\RR_\Theta))} < \infty.$$
\end{lemma}

\dem Let $w_\eta(z,\xi) = e^{-2\pi i \langle \xi, \Theta_\downarrow z \rangle} e^{2\pi i \langle \Theta \eta, z \rangle}$, so that 
\begin{eqnarray*}
\Psi_{g_\eta}(\varphi) & = & \int_{\R^n} g_\eta(\xi) \widehat{\varphi}(\xi) \lambda_\Theta(\xi) \, d\xi \\
& = & \int_{\R^n} \int_{\R^n} \widehat{g}(z,\xi) w_\eta(z,\xi) \widehat{\varphi}(\xi) \lambda_\Theta(z) \lambda_\Theta(\xi) \, dz d\xi \\
& = & \int_{\R^n} \int_{\R^n} \int_{\R^n} \widehat{\widehat{g}}(z,\zeta) w_\eta(z,\xi) e^{2\pi i \langle \xi, \zeta \rangle} \widehat{\varphi}(\xi) \lambda_\Theta(z) \lambda_\Theta(\xi) \, dz d \zeta d\xi.
\end{eqnarray*}
Let us define $j_{\xi\eta}: \R^n \to \C$ and $m_\mathrm{N}: \R^n \to \C$ as follows 
\begin{eqnarray*}
j_{\xi\eta}(z) & = & \frac{w_\eta(z,\xi)}{(1 + 4 \pi^2 |z|^2)^{\mathrm{N}}}, \\
m_\mathrm{N}(\xi) & = & \int_{\R^n} \frac{e^{2\pi i \langle \xi, \zeta \rangle}}{(1 + 4 \pi^2 |\zeta|^2 )^{\mathrm{N}}} \, d\zeta \ = \ \widehat{j}_{00}(\xi),
\end{eqnarray*}
where $\widehat{j}_{00}$ stands for the Euclidean Fourier transform of $j_{\xi\eta}$ when $(\xi,\eta) = (0,0)$. Inserting our definition of $\widehat{\widehat{g}}(z,\zeta)$, we finally end up with the following factorization  
\begin{eqnarray*}
\Psi_{g_\eta}(\varphi) & = & \int_{\R^n} \int_{\R^n} j_{\xi\eta}(z) m_\mathrm{N}(\xi) \widehat{\varphi}(\xi) \lambda_\Theta(z) \lambda_\Theta(\xi) \, dz d\xi \\
& = & \int_{\R^n} \Big( \int_{\R^n} j_{\xi\eta}(z) m_\mathrm{N}(\xi - z) \widehat{\varphi}(\xi - z) e^{2\pi i \langle z, \Theta_\downarrow (\xi - z)\rangle} \, dz \Big) \lambda_\Theta(\xi) \, d\xi \\
& = & \int_{\R^n} \Big( \int_{\R^n} \underbrace{\big( j_{\xi\eta}(\xi - z) e^{2\pi i \langle \xi-z, \Theta_\downarrow z \rangle} \big) m_\mathrm{N}(z)}_{k_\eta(\xi,z)} \, \widehat{\varphi}(z) \, dz \Big) \lambda_\Theta(\xi) \, d\xi.  
\end{eqnarray*}
This gives $\Psi_{g_\eta} = \lambda_\Theta \circ T_{k_\eta} \circ \lambda_\Theta^{-1}$, which reduces our goal to justify the assertion for $T_{k_\eta}$ instead of $\Psi_{g_\eta}$. Indeed, assume $T_{k_\eta} = \mathsf{A}_\eta^* \circ \mathsf{B} \circ \mathsf{A}_\eta$ with $\mathsf{A}_\eta$ uniformly bounded in $\mathcal{B}(L_2(\R^n))$ and $\mathsf{B}$ a trace class operator on the Hilbert space $L_2(\R^n)$. Then we consider the maps $$A_\eta = \lambda_\Theta \circ \mathsf{A}_\eta \circ \lambda_\Theta^{-1} \quad \mbox{and} \quad B = \lambda_\Theta \circ \mathsf{B} \circ \lambda_\Theta^{-1},$$ which factorize $\Psi_{g_\eta}$ and satisfy
\begin{eqnarray*}
\|A_\eta\|_{\mathcal{B}(L_2(\RR_\Theta))} & = & \|\mathsf{A}_\eta\|_{\mathcal{B}(L_2(\R^n))}, \\
\| \hskip0.5pt B \hskip0.5pt \|_{S_1(L_2(\RR_\Theta))} & = & \| \hskip0.5pt \mathsf{B} \hskip1pt \|_{S_1(L_2(\R^n))}.
\end{eqnarray*}
The kernel $k_\eta$ can be written as follows $$k_\eta(x,y) = e^{2\pi i \langle x - y, \Theta \eta + \Theta_\downarrow y + \Theta_\uparrow x \rangle} j_{00}(x-y) \, m_\mathrm{N}(y) = e^{2 \pi i \langle x-y, \Theta \eta \rangle} k(x,y).$$ If $\mathsf{A}_\eta f(x) = e^{-2\pi i \langle x, \Theta \eta \rangle} f(x)$, we see that $T_{k_\eta} = \mathsf{A}_\eta^* \circ \mathsf{B} \circ \mathsf{A}_\eta$ with $\mathsf{B} = T_k$ and $\mathsf{A}_\eta$ unitaries. Thus, it suffices to show that $\mathsf{B}$ is trace class on $L_2(\R^n)$. Composing it with the Euclidean Fourier transform $\mathcal{F} = \lambda_0^{-1}$ as in the proof of \cite[Lemma 1]{C} we end up with $L = \mathcal{F} \circ T_k$, whose kernel is given by $$\ell(x,y) = e^{-2\pi i \langle x, y \rangle} \widehat{\alpha}(x - \Theta y) m_\mathrm{N}(y),$$ where $\widehat{\alpha}$ is the Euclidean Fourier transform of $\alpha(z) = j_{00}(z) e^{-2\pi i \langle z, \Theta_\downarrow z \rangle} = j_{z0}(z)$. This is very similar to the kernel in \cite[Lemma 1 - (1.25)]{C}, in fact we recover the same kernel for $\Theta=0$. Unfortunately, due to the $\Theta$-phase we are carrying, we do not have separated variables as in \cite{C}. However, a detailed analysis of Cordes argument shows that what really matters is that the $x$-factor of the kernel |$\psi_\tau(x)$ in \cite{C}| yields a pointwise multiplier by $m_\mathrm{N}$. We only have that in the $y$-variable. Taking the adjoint $L^* = T_{k_\eta}^* \circ \mathcal{F}^{-1}$ we get the kernel $$\ell^*(x,y) = \overline{\ell(y,x)} = e^{2 \pi i \langle x, y \rangle} \overline{m_\mathrm{N}(x) \widehat{\alpha}(y - \Theta x)}.$$ Then, Cordes factorization $m_\mathrm{N}(x) = \zeta(x) \kappa(x)$ with $\zeta(x) = \exp (-\frac12 \langle x \rangle)$ implies in turn that $L^* = R \circ S$ where their respective kernels $r(x,z)$ and $s(z,y)$ are given by 
\begin{eqnarray*}
r(x,z) & = & \frac{\overline{\kappa(x)}}{(1 + 4 \pi^2 |z|^2)^{\mathrm{M}}} \int_{\R^n} \frac{e^{2\pi i \langle x-z, s \rangle}}{(1 + 4 \pi^2 |s|^2 )^{\mathrm{M}}} \, ds, \\
[5pt] s(z,y) & = & (1 + 4 \pi^2 |z|^2)^{\mathrm{M}} (\1 - \Delta_z)^{\mathrm{M}} \Big( \zeta(z) e^{2 \pi i \langle z, y \rangle} \overline{\widehat{\alpha}(y - \Theta z)} \Big).
\end{eqnarray*}         
By \cite{C}, $R$ is Hilbert-Schmidt for $\mathrm{M}$ large enough. Since $\widehat{\alpha}$ is as smooth as $\widehat{j}_{00}$, it is $\mathcal{C}^k(\R^n)$ for $\mathrm{N} > \frac12 (n+k)$ and exponentially decreasing at $\infty$. In particular, $S$ is also Hilbert-Schmidt for $\mathrm{N}$ large enough. Thus $\mathsf{B} = \mathcal{F}^{-1} S^* R^* \in S_1(L_2(\R^n))$. \fin

\begin{theorem} \label{PSD.CalVaiThm00}
If $a \in S_{0,0}^0(\RR_\Theta)$, then $\Psi_a \hskip-3pt: L_2(\RR_\Theta) \to L_2(\RR_\Theta)$ is bounded.
\end{theorem}

\dem According to Lemmas \ref{LemmaConvolution} and \ref{LemmaTraceClass} we find 
$$\Psi_a (\varphi) = (id \otimes \tau_\Theta) \int_{\R^n} (\1 \otimes b(\eta)) \Big( \Phi_\eta^* \big( \underbrace{A_\eta^* B A_\eta \otimes id}_{\mathbf{A}_\eta^* \mathbf{B} \mathbf{A}_\eta} \big) \Phi_\eta \Big) (\varphi \otimes \1) \, d\eta.$$
Given $\varphi_1, \varphi_2$ in the unit ball of $L_2(\RR_\Theta)$, it suffices to get a uniform bound for 
\begin{eqnarray*}
\big\langle \Psi_a(\varphi_1), \varphi_2 \big\rangle & = & \int_{\R^n} (\tau_\Theta \otimes \tau_\Theta) \Big\{ \Phi_\eta^* \mathbf{A}_\eta^* \mathbf{B} \mathbf{A}_\eta \Phi_\eta (\varphi_1 \otimes \1) (\varphi_2 \otimes b(\eta)^*)^* \Big\} \, d\eta \\
& = & \int_{\R^n} (\tau_\Theta \otimes \tau_\Theta) \Big\{ \mathbf{B}_1 \mathbf{A}_\eta \Phi_\eta (\varphi_1 \otimes \1) \mathbf{B}_2 \mathbf{A}_\eta \Phi_\eta(\varphi_2 \otimes b(\eta)^*)^* \Big\} \, d\eta,
\end{eqnarray*}
where $\mathbf{B} = (u |\mathbf{B}|^\frac12) |\mathbf{B}|^\frac12 = \mathbf{B}_2^* \mathbf{B}_1$ from polar decomposition. By Cauchy-Schwarz
\begin{eqnarray*}
\lefteqn{\hskip-6pt \big| \big\langle \Psi_a(\varphi_1), \varphi_2 \big\rangle \big| \le \Big( \int_{\R^n} (\tau_\Theta \otimes \tau_\Theta) \Big\{ |\mathbf{B}_1|^2 \mathbf{A}_\eta \Phi_\eta (\varphi_1 \otimes \1) \mathbf{A}_\eta \Phi_\eta(\varphi_1 \otimes \1 )^* \Big\} \, d\eta \Big)^\frac12} \\ 
& \times & \Big( \int_{\R^n} (\tau_\Theta \otimes \tau_\Theta) \Big\{ |\mathbf{B}_2|^2 \mathbf{A}_\eta \Phi_\eta (\varphi_2 \otimes b(\eta)^*) \mathbf{A}_\eta \Phi_\eta(\varphi_2 \otimes b(\eta)^*)^* \Big\} \, d\eta \Big)^\frac12 = \alpha \beta.
\end{eqnarray*}
Writing $\mathbf{B}_2 = B_2 \otimes id$, we claim that the second term above $\beta$ is dominated by $$\sup_{\eta \in \R^n} \big\| b(\eta): L_2(\RR_\Theta) \to L_2(\RR_\Theta) \big\| \, \big\| |B_2|^2 \big\|_{S_1(L_2(\RR_\Theta))}^\frac12.$$ Note that the same estimate applies to the first term with $b(\eta)=\1$ and $(\varphi_2, \mathbf{B}_2)$ replaced by $(\varphi_1, \mathbf{B}_1)$. Moreover, since $|B_j|^2 \le |B| + u|B|u^*$ and $B$ is trace class, it suffices to check that $b(\xi) = (\1 - \Delta_\Theta)^\mathrm{N} (\1 - \Delta_\xi)^\mathrm{N} a(\xi)$ is uniformly bounded in $\RR_\Theta$ which follows from the fact that $a \in S_{00}^0(\RR_\Theta)$. Therefore, it only remains to justify our claim above. Since $|B_2|^2$ is trace class, let $s_j$ denote its singular numbers and consider the corresponding set $u_j$ of unit eigenvectors. This gives $$|\mathbf{B}_2|^2(h) = \summ_j s_j (\tau_\Theta \otimes id) \big(h (u_j \otimes \1)^*\big) (u_j \otimes \1).$$ In particular, using the module bracket $\langle \langle h_1, h_2 \rangle \rangle = (\tau_\Theta \otimes id) (h_1 \bullet h_2^*)$, we get
\begin{eqnarray*}
\beta^2 & = & \summ_j s_j \int_{\R^n} \tau_\Theta \Big\{ \big| \big\langle \big\langle \mathbf{A}_\eta \Phi_\eta(\varphi_2 \otimes b(\eta)^*), u_j \otimes \1 \big\rangle \big\rangle \big|^2 \Big\} \, d\eta \\
& = & \summ_j s_j \int_{\R^n} \tau_\Theta \Big\{ \big| \big\langle \big\langle M_\eta^*(\varphi_2 \otimes b(\eta)^*), \Pi_\Theta \mathbf{A}_\eta ^*(u_j \otimes \1) \big\rangle \big\rangle \big|^2 \Big\} \, d\eta.
\end{eqnarray*}
Now, recalling that $\mathbf{A}_\eta^*(u_j \otimes \1) = A_\eta^*(u_j) \otimes \1 = \lambda_\Theta \circ \mathsf{A}_\eta^* \circ \lambda_\Theta^{-1}(u_j) \otimes \1$, we get 
\begin{eqnarray*}
\lefteqn{\hskip-15pt \Big\langle \Big\langle M_\eta^*(\varphi_2 \otimes b(\eta)^*), \Pi_\Theta \mathbf{A}_\eta ^*(u_j \otimes \1) \Big\rangle \Big\rangle} \\ [3pt] & = &  \Big\langle \Big\langle \lambda_\Theta(\eta)^* \varphi_2 \otimes b(\eta)^*, \int_{\R^n} e^{2\pi i \langle \xi, \Theta \eta \rangle} \widehat{u}_j(\xi) \pi_\Theta(\exp_\xi) \, d\xi \Big\rangle \Big\rangle \\
& = & \Big( \int_{\R^n} e^{-2\pi i \langle \xi, \Theta \eta \rangle} \overline{\widehat{u}_j(\xi)} (\tau_\Theta \otimes id) \Big( \big( \lambda_\Theta(\eta)^* \varphi_2 \otimes \1 \big) \pi_\Theta(\exp_\xi)^* \Big) \, d\xi \Big) b(\eta)^* \\
& = & \Big( \int_{\R^n} \underbrace{e^{-2\pi i \langle \xi, \Theta \eta \rangle} e^{-2\pi i \langle \eta, \Theta_\downarrow \xi \rangle} \overline{\widehat{u}_j(\xi)} \widehat{\varphi}_2(\xi+ \eta)}_{\widehat{\phi}_{\eta j}(\xi)} \lambda_\Theta(\xi) \, d\xi \Big) b(\eta)^*.
\end{eqnarray*}
This gives 
\begin{eqnarray*}
\beta^2 & \le & \sup_{\eta \in \R^n} \big\| b(\eta) \big\|_{\RR_\Theta}^2 \summ_j s_j \int_{\R^n} \big\| \lambda_\Theta(\phi_{\eta_j}) \big\|_{L_2(\RR_\Theta)}^2 \, d\eta \\
& \le & \sup_{\eta \in \R^n} \big\| b(\eta) \big\|_{\RR_\Theta}^2 \summ_j s_j \int_{\R^n} \int_{\R^n} \big| \widehat{u}_j(\xi) \widehat{\varphi}_2(\xi+\eta) \big|^2 \, d\eta d\xi,
\end{eqnarray*}
which is exactly the estimate we were looking for. This completes the proof. \fin

\begin{remark} \label{RemOptimalCV}
\emph{A careful analysis of the function $\widehat{\alpha}$ in the proof of Lemma \ref{LemmaTraceClass} could lead as in \cite{C} to the sharp condition $\mathrm{N} > n/4$. This would imply that Theorem \ref{PSD.CalVaiThm00} holds under the optimal assumption} $$\big| \partial_\Theta^\beta \partial_\xi^\alpha a(\xi) \big| \le C_{\alpha\beta} \quad \mbox{for} \quad |\alpha|, |\beta| \le \Big[\frac{n}{2}\Big] + 1.$$ 
\end{remark}

Now we are ready to study the $L_2$-boundedness for exotic symbols in $S_{\rho,\rho}^0(\RR_\Theta)$ with $0 < \rho < 1$. A weak form of Cotlar's almost orthogonality lemma naturally plays a crucial role. Namely, given a family of operators $(T_j)_{j \ge 0} \subset \mathcal{B}(\H)$ and a summable sequence $(c_j)_{j \ge 0} \subset \R_+$ we find $$\Big\| \sum_{j \ge 0} T_j \Big\|_{\mathcal{B}(\H)} \, \lesssim \, \sum_{j \ge 0} c_j$$ provided that the following conditions hold for $j \neq k$ $$\sup_{j \ge 0} \big\| T_j \big\|_{\mathcal{B}(\H)} < \infty, \qquad \big\| T_j T_k^\ast \big\|_{\mathcal{B}(\H)} = 0, \qquad \big\| T_j^\ast T_k \big\|_{\mathcal{B}(\H)} \le c_jc_k.$$ The other ingredient is a dilation argument among different deformations $\RR_\Theta$. 

\begin{lemma} \label{LemmaDilation}
Given $\mathrm{R}>0$, the map $$\mathrm{D}_\mathrm{R}: \RR_\Theta \ni \lambda_\Theta(\xi) \mapsto \lambda_{\mathrm{R}^2 \Theta}\Big(\frac{\xi}{\mathrm{R}}\Big) \in \RR_{\mathrm{R}^2 \Theta}$$ is a $*$-homomorphism. Moreover, $\Psi_a = \mathrm{D}_{\mathrm{R}}^{-1} \Psi_{\widetilde{a}_\mathrm{R}} \mathrm{D}_{\mathrm{R}}$ for $$a: \R^n \to \RR_\Theta \quad \mbox{and} \quad \widetilde{a}_\mathrm{R}(\xi) = \int_{\R^n} \widehat{a}(z,\mathrm{R} \xi) \lambda_{\mathrm{R}^2 \Theta}(z/\mathrm{R}) \, dz = \mathrm{D}_\mathrm{R}(a(\mathrm{R}\xi)) \in \RR_{\mathrm{R}^2 \Theta}.$$ 
\end{lemma}

\dem To prove that $\mathrm{D}_\mathrm{R}$ is a $*$-homomorphism is straightforward. Now
\begin{eqnarray*}
\lefteqn{\mathrm{D_R} \Psi_a \mathrm{D}_{\mathrm{R}}^{-1} (\varphi) 
= \int_{\R^n} (\widehat{\Psi_a \mathrm{D}_{\mathrm{R}}^{-1} \varphi})(\xi) \lambda_{\mathrm{R}^2\Theta} (\xi/\mathrm{R}) \, d\xi} \\
\!\!\! & = & \!\!\! \int_{\R^n} \tau_\Theta \Big\{ \Big( \int_{\R^n} a(\eta) \widehat{\mathrm{D}_{\mathrm{R}}^{-1} \varphi}(\eta) \lambda_\Theta(\eta) \, d\eta  \Big) \lambda_\Theta(\xi)^* \Big\} \lambda_{\mathrm{R}^2\Theta}(\xi/\mathrm{R}) \, d\xi \\
\!\!\! & = & \!\!\! \int_{\R^n} \tau_\Theta \Big\{ \Big( \int_{\R^n} \int_{\R^n} \widehat{a}(z,\eta) \lambda_\Theta(z) \mathrm{R}^{-n} \widehat{\varphi}(\eta/\mathrm{R}) \lambda_\Theta(\eta) \, dz d\eta  \Big) \lambda_\Theta(\xi)^* \Big\} \lambda_{\mathrm{R}^2\Theta}(\xi/\mathrm{R}) \, d\xi \\
\!\!\! & = & \!\!\! \int_{\R^n} \int_{\R^n} e^{2\pi i \langle \xi - \eta, \Theta_\downarrow \eta \rangle} \widehat{a}(\xi-\eta,\eta) \mathrm{R}^{-n} \widehat{\varphi}(\eta/\mathrm{R}) \lambda_{\mathrm{R}^2\Theta} (\xi/\mathrm{R}) \, d\eta d\xi \\
\!\!\! & = & \!\!\! \int_{\R^n} \int_{\R^n} \widehat{a}(\xi,\eta) \mathrm{R}^{-n} \widehat{\varphi}(\eta/\mathrm{R}) \lambda_{\mathrm{R}^2\Theta} (\xi/\mathrm{R}) \lambda_{\mathrm{R}^2\Theta}(\eta/\mathrm{R}) \, d\eta d\xi \\
\!\!\! & = & \!\!\! \int_{\R^n} \Big( \int_{\R^n} \widehat{a}(\xi,\mathrm{R} \eta) \lambda_{\mathrm{R}^2\Theta} (\xi/\mathrm{R}) \, d\xi \Big) \widehat{\varphi}(\eta) \lambda_{\mathrm{R}^2\Theta}(\eta) \, d\eta \\
\!\!\! & = & \!\!\! \int_{\R^n} \widetilde{a}_\mathrm{R}(\eta) \widehat{\varphi}(\eta) \lambda_{\mathrm{R}^2\Theta}(\eta) \, d\eta \ = \ \Psi_{\widetilde{a}_\mathrm{R}}(\varphi). \hskip164.5pt \square
\end{eqnarray*}
 
\vskip10pt

\begin{theorem} \label{PSD.CalVaiThm}
If $a \in S_{\rho,\rho}^0(\RR_\Theta)$ and $\rho < 1$, $\Psi_a \hskip-3pt: L_2(\RR_\Theta) \to L_2(\RR_\Theta)$ is bounded.
\end{theorem}

\dem Let $\phi_0 \in \mathcal{C}^\infty(\R^n)$ radial, identically 1 in $\mathrm{B}_1(0)$ and zero outside $\mathrm{B}_2(0)$. Using the partition of unity $\phi_0 + \sum_{j \ge 1} \phi_j \equiv 1$ with $\phi_j(\xi) = \phi_0(2^{-j} \xi) - \phi_0(2^{-(j-1)} \xi)$ we decompose $\Psi_a$ as follows $$\Psi_a = \sum_{j = 0}^\infty \Psi_{a_j}  = \sum_{j = 0}^\infty \Psi_{a_{2 j}} + \sum_{j = 0}^\infty \Psi_{a_{2 j + 1}} = \Psi_\mathrm{even} + \Psi_\mathrm{odd},$$ where $a_j(\xi) = a(\xi)  \phi_j(\xi)$. We shall only bound the even part, since both are treated in a similar way. To do so, we apply Cotlar's lemma as stated above. Given $j,k$ distinct even numbers, we clearly have $\Psi_{a_j} \Psi_{a_k}^\ast = 0$ since $\phi_j$ and $\phi_k$ have disjoint supports. Therefore, it suffices to prove that 
\begin{itemize}
\item[i)] $\displaystyle \hskip2pt \sup_{j \ge 0} \big\| \Psi_{a_j} \big\|_{\mathcal{B}(L_2(\RR_\Theta))} < \infty$,

\vskip3pt

\item[ii)] $\displaystyle \big\| \hskip1pt \Psi_{a_j}^\ast \Psi_{a_k} \hskip1pt \big\|_{\mathcal{B}(L_2(\RR_\Theta))} \le c_jc_k$,
\end{itemize}
for some summable sequence $(c_j)_{j \ge 0} \subset \R_+$ and any pair of distinct even integers $j,k$. The first condition follows from our form of Calder\'on-Vaillancourt theorem in $S_{00}^0(\RR_\Theta)$ and Lemma \ref{LemmaDilation}. Indeed, pick $\mathrm{R}_j = 2^{j\rho}$ and let $$a_{[j]} = \widetilde{(a_j)}_{\mathrm{R}_j} = \mathrm{D}_{\mathrm{R}_j} \big( a_j (\mathrm{R}_j \hskip1pt \cdot \hskip1pt )\big).$$ Then
$$\big\| \Psi_{a_j} \big\|_{\B(L_2(\RR_\Theta))} \le \big\| \Psi_{a_{[j]}} \big\|_{\B(L_2(\RR_{\mathrm{R}_j^2\Theta}))}$$
since $\Psi_{a_j} = \mathrm{D}_{\mathrm{R}_j}^{-1} \Psi_{a_{[j]}} \mathrm{D}_{\mathrm{R}_j}$ and $$\big\| \mathrm{D}_{\mathrm{R}_j}^{-1} \big\|_{\B(L_2(\RR_{\mathrm{R}_j^2\Theta}), L_2(\RR_\Theta))} \big\| \mathrm{D}_{\mathrm{R}_j} \big\|_{\B(L_2(\RR_\Theta), L_2(\RR_{\mathrm{R}_j^2 \Theta}))} = 1.$$ The $L_2$-boundedness of $\Psi_{a_{[j]}}$ follows from Theorem \ref{PSD.CalVaiThm00} since $a_{[j]} \in S_{00}^0(\RR_{\mathrm{R}_j^2 \Theta})$. The proof of this fact follows essentially as in the Euclidean setting. Indeed, write $a_{[j]}$ in terms of $\widehat{a}_j$ |Lemma \ref{LemmaDilation}| and use that $\mathrm{D}_{\mathrm{R}_j}$ is a $*$-homomorphism. In conjunction with the $\xi$-localization of $a_j$ in the annulus of radii $\sim 2^j$, this easily gives $a_{[j]} \in S_{00}^0(\RR_{\mathrm{R}_j^2\Theta})$, we leave details to the reader. It remains to estimate the norm of $\Psi_{a_j}^\ast \Psi_{a_k}$ for even $j \neq k$. After a calculation we obtain that the kernel $k_{jk}$ of such operator is given by 
$$k_{jk} = (id \otimes \tau_\Theta \otimes id) \int_{\R^n} \int_{\R^n}
      (\pi_\Theta(\exp_\xi) \otimes \1) \,
      (\1 \otimes a_j(\xi)^\ast a_k(\eta) \otimes \1) \,
      (\1 \otimes \pi_\Theta(\exp_\eta))
      \, d \xi d \eta,$$
where, in an abuse of notation, the element
$\pi_\Theta(\exp_\xi) = \lambda_\Theta(\xi) \otimes \lambda_\Theta(\xi)^\ast$
is seen as belonging to $\RR_\Theta \bar\otimes \RR_\Theta$ instead of
$\RR_\Theta \bar\otimes \RR_\Theta^\mathrm{op}$. We are also going to shorten
$a \otimes b \otimes \1$ by $(a \otimes b)_{[12]}$ where the leg numbers
just mean that the tensor components are placed in the first and second places respectively. Now, we use 
\begin{eqnarray*}
  \frac{(\1 - \Delta_\Theta)^{\mathrm{N}}} {(1 + 4 \pi^2 |\xi - \eta|^2)^{\mathrm{N}}} \, \lambda_\Theta(\eta) \lambda_\Theta(\xi)^\ast & = & \lambda_\Theta(\eta) \lambda_\Theta(\xi)^\ast,\\
  \pi_\Theta \bigg( \frac{1}{(1 + 4 \pi^2\mathrm{d}^2)^n} \bigg) (\1 - \Delta_\eta)^n \, \lambda_\Theta(\eta) \otimes \lambda_\Theta(\eta)^\ast & = & \lambda_\Theta(\eta) \otimes \lambda_\Theta(\eta)^\ast,\\
  \pi_\Theta \bigg( \frac{1}{(1 + 4 \pi^2 \mathrm{d}^2)^n} \bigg) (\1 - \hskip0.5pt \Delta_\xi)^n \, \lambda_\Theta( \hskip0.5pt \xi) \otimes \lambda_\Theta( \hskip0.5pt \xi)^\ast & = & \lambda_\Theta( \hskip0.5pt \xi) \otimes \lambda_\Theta( \hskip0.5pt \xi)^\ast,
\end{eqnarray*}
where $\mathrm{d}(x) = |x|$ is the Euclidean distance. Integration by parts yields
$$k_{jk} = (id \otimes \tau_\Theta \otimes id) \int_{\R^n} \int_{\R^n} \pi_\Theta (\varphi_\xi)_{[12]} B(\xi,\eta)_{[2]} \pi_\Theta (\varphi_{\eta})_{[23]} \, d\eta d\xi,$$
where
$$B(\xi, \eta) = (\1 - \Delta_\eta)^n (\1 - \Delta_\xi)^n \bigg\{ \frac{(\1 - \Delta_\Theta)^{\mathrm{N}}}{(1 + 4 \pi^2 |\eta - \xi|^2)^{\mathrm{N}}} a_j(\xi)^\ast a_k(\eta) \bigg\}$$
and the function $\varphi_\zeta$ is given by $\exp_\zeta (1 + 4 \pi^2 \mathrm{d}^2)^{-n}$. After expanding the derivatives using the Leibniz rule, we obtain that $B$ is a finite sum of simple terms of the form
$$B_s(\xi, \eta) = \partial_\xi^{\alpha_1} \partial_\eta^{\beta_1} \Big( \frac{1}{(1 + 4 \pi^2 |\xi - \eta|^2)^{\mathrm{N}}} \Big) \, \underbrace{\partial_\xi^{\alpha_2} \partial_\Theta^{\sigma_1} a_j(\xi)^\ast}_{b^s_j(\xi)^\ast} \, \underbrace{\partial_\eta^{\beta_2}  \partial_\Theta^{\sigma_2} a_k(\eta)}_{b^s_k(\eta)},$$
where $\alpha_i, \beta_i, \sigma_i \in \Z_+^n$ satisfy $|\alpha_1 + \alpha_2| \leq 2 n$, $|\beta_1 + \beta_2| \leq 2 n$, $|\sigma_1 + \sigma_2| \leq 2 \mathrm{N}$
and $s$ is the combination of the involved multindices. We can bound each of the above summands in $s$ independently
$$\big\| \Psi_{a_j}^\ast \Psi_{a_k} \big\| \le \summ_s \bigg\| \int_{\R^n} \int_{\R^n} T_{(id \otimes \tau_\Theta \otimes id) \{ \pi_\Theta(\varphi_\xi)_{[12]} B_s(\xi,\eta)_{[2]} \pi_\Theta(\varphi_\eta)_{[23]} \} } d\eta d\xi \bigg\| = \summ_s A_s.$$
Using $\big| \partial_\xi^{\alpha_1} \partial_\eta^{\beta_1} \langle \xi - \eta \rangle^{-2 \mathrm{N}} \big| \lesssim \langle \xi - \eta \rangle^{-2 \mathrm{N}}$ we obtain
$$A_s \lesssim \int_{\R^n} \int_{\R^n} \langle \xi - \eta \rangle^{-2 \mathrm{N}} \big\| \underbrace{T_{(id \otimes \tau_\Theta \otimes id) \{ \pi_\Theta(\varphi_\xi)_{[12]} b_j^s(\xi)^\ast_{[2]} b_k^s(\eta)_{[2]} \pi_\Theta(\varphi_\eta)_{[23]} \}}}_{T_{j\xi k \eta}} \big\| \, d\eta d\xi.$$
$T_{j\xi k \eta}$ can be factorized as $T_{(b_j^s \otimes \1) \bullet \pi_\Theta(\varphi_\xi)} \circ T_{(b_k^s(\eta) \otimes \1) \bullet \pi_\Theta(\varphi_\eta)}$, so that
$$\big\| T_{j \xi k \eta} \big\|_{\B(L_2(\RR_\Theta))} \leq \| b_j^s(\xi) \|_{\RR_\Theta} \|  b_k^s(\eta) \|_{\RR_\Theta} \| T_{\pi_\Theta(\varphi_\xi)} \|_{\B(L_2(\RR_\Theta))} \| T_{\pi_\Theta(\varphi_\eta)} \|_{\B(L_2(\RR_\Theta))}.$$
Recall that $\|T_{\pi_\Theta(\varphi_\zeta)}\|_{\B(L_2(\RR_\Theta))} \le \|\varphi_\zeta\|_{L_1(\R^n)} \lesssim 1$. Moreover, using that $\xi \sim 2^{j}$ and $\eta \sim 2^{k}$ from the supports of $a_j(\xi)$ and $a_k(\eta)$ as well as the H\"ormander condition for $a$, we deduce the following bound
$$A_s \lesssim \int_{\R^n} \int_{\R^n} \frac{\| b_j^s(\xi) \|_{\RR_\Theta} \| b_k^s(\eta) \|_{\RR_\Theta}}{\langle \xi - \eta \rangle^{2 \mathrm{N}}} \, d\xi d\eta  
\lesssim 2^{-2 \mathrm{N} \max\{j,k\}} \, 2^{2 \mathrm{N} \rho \max\{j,k\}} \, 2^{n( j + k)}.$$
Summing all the terms indexed by $s$ we obtain
$$\big\| \Psi_{a_j}^\ast \Psi_{a_k} \big\|_{\B(L_2(\RR_\Theta))} \lesssim 4^{-\max\{j, k\} ((1 - \rho) \mathrm{N} - n)} \le 2^{-j ((1 - \rho) \mathrm{N} - n)} 2^{-k ((1 - \rho) \mathrm{N} - n)} = c_j c_k,$$ which arises from $(c_j)_{j \ge 0}$ summable for $\mathrm{N}$ large enough. The proof is complete. \fin

\begin{remark}
\emph{Let $0 \le \delta \le \rho \le 1$, since $$S_{\rho,\delta}^0(\RR_\Theta) \subset S_{\delta,\delta}^0(\RR_\Theta) \cap S_{\rho,\rho}^0(\RR_\Theta),$$ we deduce $\Psi_a: L_2(\RR_\Theta) \to L_2(\RR_\Theta)$ for $a \in S_{\rho,\delta}^0(\RR_\Theta)$ as long as $(\rho,\delta) \neq (1,1)$.}
\end{remark}

\begin{remark}
\emph{A standard (nonoptimal) proof of Calder\'on-Vaillancourt theorem for $S_{00}^0$ in the Euclidean setting \cite{St}�\ follows from a suitable partition of unity in the variables $(x,\xi) \in \R^n \times \R^n$ with no known analogue for $x_\Theta \in \RR_\Theta$ and $\xi \in \R^n$. An alternative way to proceed is the following. Given $a \in S_{00}^0(\RR_\Theta)$, let $$\mathbf{a}(x,\xi) = \sigma_\Theta^x a(\xi) = \int_{\R^n} \widehat{a}(z,\xi) \sigma_\Theta^x(\lambda_\Theta(z)) \, dz = \int_{\R^n} \widehat{a}(z,\xi) e^{2\pi i \langle x, z \rangle} \lambda_\Theta(z) \, dz.$$ Using the intertwining identity $\Psi_\mathbf{a} \circ \sigma_\Theta = \sigma_\Theta \circ \Psi_a$ and recalling from Appendix B that $\sigma_\Theta: L_2^c(\RR_\Theta) \to L_2^c(\R^n) \bar\otimes \RR_\Theta$ is a complete isometry, it turns our that the $L_2$-boundedness of $\Psi_a$ is equivalent to the boundedness of the operator-valued map $\Psi_\mathbf{a}: L_2^c(\R^n) \bar\otimes \RR_\Theta \to L_2^c(\R^n) \bar\otimes \RR_\Theta$. Now, since $\Psi_\mathbf{a}$ is a right $\RR_\Theta$-module map, it follows from \cite[Remark 2.4]{JMP} that this will hold as long as $\Psi_\mathbf{a}$ is bounded over the Hilbert space $L_2(\R^n;L_2(\RR_\Theta))$. $\Psi_\mathbf{a}$ comes equipped with an operator-valued kernel acting by left multiplication. This kind of maps are generally bad behaved \cite{HLMP} but we know from our proof above that $L_2$-boundedness must hold in this case. Thus this also opens the door to prove Calder\'on-Vaillancourt using a partition of unity in the $x$-component, which mirrors the behavior of its quantum analogue $x_\Theta$.}
\end{remark}

\subsubsection{Bourdaud's condition for forbidden symbols in $\RR_\Theta$}

We have justified that all symbols in $S_{\rho,\delta}^0(\RR_\Theta)$ yield $L_2$-bounded pseudodifferential operators except for the class of so-called forbidden symbols with $\rho=\delta=1$, which is known to fail it even in the Euclidean setting. Bourdaud established a sufficient condition in \cite{Bo} playing the role of the $T(1)$-theorem for pseudodifferential operators and which we now study in $\RR_\Theta$. Given $p \ge 1$ and $s \in \R$, let $\mathrm{W}_{2,s}(\RR_\Theta)$ be the Sobolev space defined as the closure of $\S_\Theta$ with respect to the norm $$\|\varphi\|_{\mathrm{W}_{2,s}(\RR_\Theta)} = \big\| (\1 - \Delta_\Theta)^\frac{s}{2} \varphi \big\|_2.$$

\begin{lemma} \label{PSD.list}
Let $\phi: \R^n \to \R_+$ be a radial smooth function identically $1$ in $\mathrm{B}_1(0)$ and vanishing outside $\mathrm{B}_2(0)$. Let $\psi_j(\xi) = \phi(2^{-j}\xi) - \phi(2^{-j+1} \xi)$ for any integer $j \in \Z$. Then, we have a norm equivalence  
$$\big\| \lambda_\Theta(f) \big\|_{2,s}^2 \sim \sum_{j \in \mathbb{Z}} 2^{2 s j} \big\|\lambda_\Theta(\psi_j f)\big\|_2^2 = \sum_{j \in \mathbb{Z}} 2^{2 s j} \big\| \psi_j f \big\|_2^2.$$ In particular, the following properties hold$\hskip1pt :$ 
\begin{itemize}
\item[\emph{i)}] $\displaystyle \mathrm{W}_{2,-s}(\RR_\Theta)^\ast = \mathrm{W}_{2,s}(\RR_\Theta)$ under the pairing $\langle x, y \rangle = \tau_\Theta(x^\ast \, y)$.

\item[\emph{ii)}] $\displaystyle \big[ \mathrm{W}_{2,s}(\RR_\Theta), \mathrm{W}_{2,-s}(\RR_\Theta) \big]_{\frac{1}{2}} = L_2(\RR_\Theta)$ by complex interpolation.
\end{itemize}
\end{lemma}

The proofs of all assertions above are straightforward. Properties i) and ii) above hold isomorphically from the first assertion, but also isometrically. We shall use the following terminology for the rest of this section. Let us consider a function $\phi: \R^n \to \R_+$ which is identically $1$ in $\mathrm{B}_{1/8}(0)$ and vanishing outside $\mathrm{B}_{1/4}(0)$. Set $\psi_0 = \phi$ and construct $\psi_j(\xi) = \phi(2^{-j}\xi) - \phi(2^{-j+1}\xi)$ for $j \ge 1$. We shall also use 
the partition of unity $\rho_0 = \psi_0 + \psi_1$ and $\rho_j = \psi_{j-1} + \psi_j + \psi_{j+1}$ for $j \ge 1$, so that 
\begin{itemize}
\item $\displaystyle \sum_{j \ge 0} \psi_j(\xi) = 1$,

\item $\displaystyle \sum_{j \ge 0} \rho_j(\xi) = 3 - \phi(\xi)$,

\item $\displaystyle \rho_j(\xi) \psi_j(\xi) = \psi_j(\xi)$ for $j \ge 0$. 
\end{itemize}

\begin{lemma} \label{PSD.symbolDecII}
If $a \in S_{1,1}^0(\RR_\Theta)$ and $\mathrm{N} > n$, we have  
$$a(\xi) = \sum_{k \in \mathbb{Z}^n} \langle k \rangle^{-\mathrm{N}} \sum_{j \geq 0} c_{j,k} \rho_j(\xi) e^{2 \pi i \langle 2^{-j} \xi, k \rangle}$$
where the coefficients $c_{j,k} \in \RR_\Theta$ satisfy the following estimate
$$\sup_{j \ge 0} \sup_{k \in \mathbb{Z}^n} \Big( \| c_{j,k} \|_{\RR_\Theta} + 2^{-j} \Big\| \Big( \sum_{i=1}^n |\partial_\Theta^i c_{j,k} |^2 \Big)^\frac12 \Big\|_{\RR_\Theta} \Big) \, < \, \infty.$$
\end{lemma}

\dem Let $a_j(\xi) = a(\xi) \psi_j(\xi)$ and $b_j(\xi) = a_j(2^j \xi)$, so that $$a(\xi) = \sum_{j \geq 0}^{\null} a_j(\xi) = \sum_{j \ge 0} b_j(2^{-j} \xi) = \sum_{j \ge 0} b_j(2^{-j} \xi) \rho_j(\xi).$$ According to this and recalling that $\mbox{supp} \, \rho_j(2^j \xi) \subset [-\frac12,\frac12]^n = Q$, it suffices to see that $b_j(\xi) = \sum_{k \in \Z^n} \langle k \rangle^{- \mathrm{N}} c_{j,k} e^{2\pi i \langle \xi, k \rangle} \chi_Q(\xi)$ for some $c_{j,k}$ satisfying the estimates in the statement. Now, since $b_j$ is also supported by $Q$, we find that $b_j(\xi) = d_j(\xi) \chi_Q(\xi)$ where $d_j$ is the $\Z^n$-periodization of $b_j$. This gives rise to the identity 
\begin{eqnarray*}
b_j(\xi) & = & \sum_{k \in \Z^n} \widehat{b}_j(k) e^{2\pi i \langle \xi, k\rangle} \chi_Q(\xi) \\ 
& = & \sum_{k \in \Z^n} \Big( \int_{\mathbb{T}^n} b_j(s) e^{-2\pi i \langle s, k \rangle} ds \Big) e^{2\pi i \langle \xi, k\rangle} \chi_Q(\xi).
\end{eqnarray*}
Integrating by parts, we obtain
$$\widehat{b}_j(k) = \frac{1}{(1 + 4 \pi^2 |k|^2)^\frac{\mathrm{N}}{2}} \int_{\mathbb{T}^n} (\1 - \Delta_s)^{\frac{\mathrm{N}}{2}} \big( b_j(s) \big) e^{-2 \pi i \langle s, k \rangle} \, ds = \langle k \rangle^{-\mathrm{N}} c_{j,k}.$$
To estimate $(\1 - \Delta_\xi)^{\frac{\mathrm{N}}{2}}b_j(\xi)$ we notice that $|\xi| \sim 1$, so
$$\big\| \partial_\xi^\alpha b_j(\xi) \big\|_{\mathcal{R}_\Theta} = 2^{j |\alpha|} \, \big\| (\partial_\xi^\alpha a_j)(2^j \xi) \big\|_{\RR_\Theta} \lesssim 2^{j |\alpha|} \, \langle 2^{j} \xi \rangle^{-|\alpha|} \lesssim 1$$ by the H\"ormander condition and therefore $(\1 - \Delta_\xi)^{\frac{\mathrm{N}}{2}}b_j(\xi)$ is uniformly bounded
in norm. The second inequality uses a similar calculation for $\nabla_\Theta (\1 - \Delta_\xi)^{\frac{\mathrm{N}}{2}}b(\xi)$. \fin

\begin{lemma} \label{PSD.gradMult}
We have $$\big\| \lambda_\Theta(\psi_j f) \big\|_{\RR_\Theta} \lesssim 2^{-j} \Big\| \Big( \sum_{k=1}^n \big| \partial_{\Theta}^k \lambda_\Theta(f) \big|^2 \Big)^{\frac12} \Big\|_{\RR_\Theta} \quad \mbox{for} \quad j > 0.$$
\end{lemma}

\dem Given $g \in \S(\R^n)$ 
$$(\sigma_\Theta \big( \lambda_\Theta(f) \big) \ast g)(x) = \int_{\R^n} \sigma_\Theta^{x - y} \big( \lambda_\Theta(f) \big) g(y) \, dy = \sigma_\Theta^x \big( \lambda_\Theta(\widehat{g} f) \big).$$
We also have $\big\| \lambda_\Theta (\psi_j f) \big\|_{\RR_\Theta} = \big\| \sigma_\Theta \lambda_\Theta(\psi_j f) \big\|_{\RR_\Theta \bar\otimes L_\infty(\R^n)}$ and 
\begin{eqnarray*}
\sigma_\Theta \lambda_\Theta(\psi_j f) & = & \sigma_\Theta \lambda_\Theta(f) \ast \widehat{\psi}_j \\
[5pt] & = & \int_{\R^n} \sigma_\Theta^{x - y}\big(\lambda_\Theta(f)\big) \widehat{\psi}_j(y) \, dy \\
& = & \int_{\R^n} \Big( \sigma_\Theta^{x - y} \big( \lambda_\Theta(f) \big) - \sigma_\Theta^x \big( \lambda_\Theta(f) \big) \Big) \widehat{\psi}_j(y) \, dy \\
& = & \int_{\R^n} \Big( \int_0^1 \sum_{k=1}^n y_k \sigma_\Theta^{x-ty} \big(\partial_\Theta^k \lambda_\Theta(f) \big) \, dt \Big) \widehat{\psi}_j(y) \, dy.
\end{eqnarray*}
We have used that the integral of $\widehat{\psi}_j$ is $0$ for any $j > 0$. Taking norms gives 
\begin{eqnarray*}
\lefteqn{\big\| \sigma_\Theta \lambda_\Theta(\psi_j f) \big\|_{\RR_\Theta \bar\otimes L_\infty(\R^n)}} \\ 
\!\!\!\! & \leq & \!\!\!\! \int_{\R^n} \Big\| \sum_{k=1}^n y_k \big(\partial_\Theta^k \lambda_\Theta(f) \big) \Big\|_{\RR_\Theta} |\widehat{\psi}_j(y)| \, dy \\
\!\!\!\! & \le & \!\!\!\! \Big( \int_{\R^n} |y| |\widehat{\psi}_j(y)| \, dy \Big) \Big\| \Big( \sum_{k=1}^n \big| \partial_{\Theta}^k \lambda_\Theta(f) \big|^2 \Big)^{\frac12} \Big\|_{\RR_\Theta} \\
\!\!\!\! & \lesssim & \!\!\! 2^{-j} \! \int_{\R^n} |y| |\widehat{\phi}(y)| \, dy \hskip1pt \Big\| \Big( \sum_{k=1}^n \big| \partial_{\Theta}^k \lambda_\Theta(f) \big|^2 \Big)^{\frac12} \Big\|_{\RR_\Theta} \lesssim 2^{-j} \Big\| \Big( \sum_{k=1}^n \big| \partial_{\Theta}^k \lambda_\Theta(f) \big|^2 \Big)^{\frac12} \Big\|_{\RR_\Theta} \square
\end{eqnarray*} 

\begin{theorem} \label{PSD.SobolevBound}
If $a \in S_{1,1}^0(\RR_\Theta)$ $$\Psi_a: \mathrm{W}_{2,s}(\RR_\Theta) \to \mathrm{W}_{2,s}(\RR_\Theta) \quad \mbox{is bounded for} \quad 0 < s < 1.$$
\end{theorem}

\dem
By Lemma \ref{PSD.symbolDecII} we have that
$$a(\xi) = \sum_{k \in \mathbb{Z}^n} \langle k \rangle^{-\mathrm{N}} \underbrace{\sum_{j \geq 0} c_{j,k} \rho_j(\xi) e^{2 \pi i \langle 2^{-j} \xi, k \rangle}}_{a_k(\xi)}.$$
By taking $\mathrm{N} > n$ we obtain that the symbol $a$ is just a summable combination of terms $a_k$ and we can concentrate on studying such terms. If
$\lambda_\Theta(f) \in \mathrm{W}_{2,s}(\RR_\Theta)$, we have that
$$\Psi_{a_k}(\lambda_\Theta(f)) = \sum_{j \geq 0} c_{j,k} \lambda_\Theta(\rho_j \exp_{2^{-j} k}f) = \sum_{j \geq 0} c_{j,k} b_{j,k}.$$
Taking another partition of unity $(\psi_\ell)_{\ell \geq 0}$ we get
$$c_{j,k} = \sum_{\ell \geq 0} \int_{\R^n} \psi_\ell(\xi) \widehat{c}_{j,k}(\xi) \lambda_\Theta(\xi) \, d\xi = \sum_{\ell \geq 0} c_{j,k}^{\ell}$$
and Lemmas \ref{PSD.symbolDecII} and \ref{PSD.gradMult} give $\|c_{j,k}^\ell \|_{\RR_\Theta} \lesssim 2^{j - \ell}$ for $\ell > 0$. Now decompose 
$$\Psi_{a_k}(\lambda_\Theta(f)) = \sum_{\ell \leq j - 4} c_{j,k}^\ell b_{j,k} + \sum_{j - 4 < \ell < j + 4} c_{j,k}^\ell b_{j,k} + \sum_{\ell \ge j + 4} c_{j,k}^\ell b_{j,k} = \mathrm{L} + \mathrm{D} + \mathrm{U}.$$ Let us begin with the estimate of the upper term $\mathrm{U}$. The Fourier support of $\lambda_\Theta(f) \lambda_\Theta(g) = \lambda_\Theta(f \ast_{\Theta} g)$ is contained in the sum of the Fourier supports of $\lambda_\Theta(f)$ and $\lambda_\Theta(g)$ respectively. In particular, the Fourier support of $c_{j,k}^\ell b_{j,k}$ is contained in $\mbox{supp} \, \psi_\ell + \mbox{supp} \, \rho_j \subset [ \mathrm{B}_{2^{\ell - 2}}(0) \setminus \mathrm{B}_{2^{\ell-4}}(0) ] + \mathrm{B}_{2^{j-1}}(0) \subset \mathrm{B}_{2^{\ell - 1}}(0) \setminus \mathrm{B}_{2^{\ell-5}}(0)$. Now we apply Lemma \ref{PSD.list} to obtain 
\begin{eqnarray*}
\| \mathrm{U} \|_{\mathrm{W}_{2,s}(\RR_\Theta)}^2 \!\!\! & \lesssim & \!\!\! \sum_{\ell \geq 4} 2^{2 \ell s} \Big\| \sum_{j \leq \ell - 4} c_{j,k}^\ell b_{j,k}\Big\|_2^2 \lesssim  \sum_{\ell \geq 4} 2^{2 \ell s} \Big( \sum_{j \leq \ell - 4} \|c_{j,k}^\ell \|_{\RR_\Theta} \|b_{j,k}\|_2 \Big)^2 \\
\!\!\! & \lesssim & \!\!\! \sum_{\ell \geq 4} 2^{2 \ell s} \Big( \sum_{j \leq \ell - 4} 2^{j - \ell} \|b_{j,k}\|_2 \Big)^2 \lesssim \sum_{\ell \geq 4} 2^{2 \ell (s-1)} \Big( \sum_{j \leq \ell - 4} 2^{j} \|b_{j,k}\|_2 \Big)^2.
\end{eqnarray*}
On the other hand, given $0 < \delta < 1-s$ we have that
$$\Big( \sum_{j \leq \ell - 4} 2^{j} \|b_{j,k}\|_2 \Big)^2 \le C_\delta \, 2^{\ell \delta} \sum_{j \leq \ell - 4} 4^{j (1 - \delta)} \|b_{j,k} \|_2^2.$$
In particular, we finally obtain the expected estimate
\begin{eqnarray*}
\| \mathrm{U} \|_{\mathrm{W}_{2,s}(\RR_\Theta)}^2 & \lesssim & \sum_{\ell \geq 0} \sum_{j \leq \ell - 4} 2^{2 \ell (s-1 + \delta)} 4^{j (1 - \delta)} \|b_{j,k}\|_2^2 \\
& = & \sum_{j \geq 0} \Big(\sum_{\ell \ge j + 4} 2^{2 \ell (s-1 + \delta)} \Big) \, 4^{j (1 - \delta)} \|b_{j,k}\|_2^2 \\
& \lesssim & \sum_{j \geq 0} 4^{j (s - 1 + \delta)} 4^{j (1 - \delta)} \|b_{j,k}\|_2^2 \ \lesssim \ \| \lambda_\Theta(f) \|_{\mathrm{W}_{2,s}(\RR_\Theta)}^2.
\end{eqnarray*}
For the lower part $\mathrm{L}$, a similar argument yields that the Fourier support of $c_{j,k}^\ell b_{j,k}$ is contained inside $\mathrm{B}_{2^j}(0) \setminus \mathrm{B}_{2^{j-6}}(0)$. Then we can apply the same principle so that
\begin{eqnarray*}
\| \mathrm{L} \|_{\mathrm{W}_{2,s}(\RR_\Theta)}^2 & \lesssim & \sum_{j \geq 4} 2^{2 j s} \Big\| \sum_{\ell \leq j - 4} c_{j,k}^\ell \, b_{j,k} \Big\|_2^2 \\
& \leq & \sum_{j \geq 4} 2^{2 j s} \Big\| \sum_{\ell \leq j - 4} c_{j,k}^\ell \Big\|_{\RR_\Theta} \| b_{j,k} \|_2^2 \\
& \lesssim & \sum_{j \geq 4} 2^{2 j s} \| b_{j,k} \|_2^2  \ \lesssim \ \| \lambda_\Theta(f) \|_{\mathrm{W}_{2,s}(\RR_\Theta)}^2.
\end{eqnarray*}
In the third inequality we have used that $\sum_\ell c_{j,k}^\ell = c_{j,k}$ and therefore
$$\Big\| \sum_{\ell \leq j - 4} c_{j,k}^\ell \Big\|_{\RR_\Theta} = \Big\| c_{j,k} - \sum_{\ell = j - 3}^{\infty} c_{j,k}^\ell \Big\|_{\RR_\Theta} =
\| c_{j,k} \|_{\RR_\Theta} + \sum_{\ell = j - 3}^{\infty} \| c_{j,k}^\ell \|_{\RR_\Theta} \lesssim 1.$$
The diagonal part $\mathrm{D}$ is easier to bound. Assume for simplicity that $j = \ell$. The Fourier support of $c_{j,k}^j b_{j,k}$ is comparable this time to a fixed dilation of $\mathrm{B}_{2^j}(0)$, not an annulus. Nevertheless, although we do not have a norm equivalence, the    
norm in $\mathrm{W}_{2,s}(\RR_\Theta)$ is still dominated by the corresponding weighted $L_2$-sum, and we get \\ [5pt] \null \hfill $\displaystyle \| \mathrm{D} \|_{\mathrm{W}_{2,s}(\RR_\Theta)}^2 \le \sum_{j \ge 0}^{\null} 2^{2js} \big\| c_{j,k}^j b_{j,k} \big\|_2^2 \lesssim \sum_{j \ge 0} 2^{2js} \| b_{j,k} \|_2^2 \sim \big\| \lambda_{\Theta}(f) \big\|_{\mathrm{W}_{2,s}(\RR_\Theta)}^2. \hskip20pt \square$ 

\vskip5pt

\begin{theorem} \label{PSD.boudreaud}
If $a, a_\dagger^* \in S^0_{1,1}(\RR_\Theta)$, then $\Psi_a: L_2(\RR_\Theta) \to L_2(\RR_\Theta)$ is bounded. 
\end{theorem}

\dem By Theorem \ref{PSD.SobolevBound}, $\Psi_a$ and its adjoint are bounded in $\mathrm{W}_{2,s}(\RR_\Theta)$. By Lemma \ref{PSD.list}, taking duals gives $\Psi_a: \mathrm{W}_{2,-s}(\RR_\Theta) \to \mathrm{W}_{2,-s}(\RR_\Theta)$ and interpolating both inequalities for $\Psi_a$ yields the assertion. \fin

\begin{remark} \label{PSD.FiniteDerivatives}
\emph{A careful examination yields that
$$\big\| \partial_\Theta^\beta \partial_\xi^\alpha a (\xi) \big\|_{\RR_\Theta} \lesssim \langle \xi \rangle^{- |\alpha| + |\beta|} \quad \mbox{for} \quad |\alpha| \leq n+1 \quad \mbox{and} \quad |\beta| \le 1$$ for $a$ and its dual symbol $a_\dagger^*$ suffices to deduce the $L_2$-boundedness of $\Psi_a$.}
\end{remark}

\subsection{$L_p$-boundedness and Sobolev $p$-estimates}
\label{PSD-Lp}

The $L_2$-boundedness results above together with our Calder\'on-Zygmund theory for $\RR_\Theta$ are the tools to find sufficient smoothness conditions on a given symbol for the $L_p$-boundedness of its pseudodifferential operator. As pointed in the Introduction, this naturally requires to work with a different quantum form of the H\"ormander classes, which is more demanding, but still recovers the classical definition for $\Theta=0$. Given $a: \R^n \to \RR_\Theta$ we say that it belongs to $\Sigma_{\rho, \delta}^m(\RR_\Theta)$ when
$$\big| \partial_\Theta^\beta \, \partial_{\Theta, \xi}^{\alpha_1} \, \partial_\xi^{\alpha_2} a(\xi) \big| \le C_{\alpha_1, \alpha_2, \beta} \langle \xi \rangle^{m - \rho |\alpha_1 + \alpha_2| + \delta |\beta|}$$
for all $\alpha_1, \alpha_2, \beta \in \Z_+^n$. Here are some trivial, albeit important, properties:
\begin{itemize}
\item[i)] $\Sigma_{\rho, \delta}^m(\RR_\Theta) \subset S_{\rho, \delta}^m(\RR_\Theta)$ since one condition reduces to the other when $\alpha_1 = 0$.

\vskip5pt

\item[ii)] All of the three derivatives involved in the definition of $\Sigma_{\rho,\delta}^m(\RR_\Theta)$ commute with each other. In particular, the order considered is completely irrelevant.

\vskip3pt

\item[iii)] Fix $(\rho,\delta,m)$ and set $|a|_{\alpha,\beta}^S$ and $|a|^{\Sigma}_{\alpha_1, \alpha_2, \beta}$ for the seminorms given by the optimal constant in the defining inequalities of $S_{\rho,\delta}^m(\RR_\Theta)$ with parameters $(\alpha,\beta)$ or $\Sigma_{\rho,\delta}^m(\RR_\Theta)$ with parameters $(\alpha_1,\alpha_2,\beta)$ respectively. Then, we have $$\lim_{\Theta \to 0} |a|^\Sigma_{\alpha_1, \alpha_2, \beta} = \lim_{\Theta \to 0} |a|_{\alpha_1+\alpha_2, \beta}^S.$$
\end{itemize}

Given $a \in \Sigma_{1,1}^0(\RR_\Theta) \subset S_{1,1}^0(\RR_\Theta)$, we will now prove that the integral kernel $k_a$ associated with $\Psi_a$ satisfies the Calder\'on-Zygmund kernel conditions in Theorem A. In conjunction with our Bourdaud type condition in Theorem \ref{PSD.boudreaud}, it will give the complete $L_p$-boundedness of $\Psi_a$ stated in Theorem B iii). Composition results further yield Sobolev $p$-estimates
$$\big\| \Psi_a: \mathrm{W}_{p,s}(\RR_\Theta) \to \mathrm{W}_{p,s-m}(\RR_\Theta) \big\|_{\mathrm{cb}} < \infty,$$
for many symbols of degree $m$, with $1 < p < \infty$ and $\|\varphi\|_{p,s} = \|(\1 - \Delta_\Theta)^{s/2}\varphi\|_p$. 

\begin{lemma} \label{Mx.Gradiente}
Given $a \in \Sigma_{\rho,\delta}^m(\RR_\Theta)$, let $$k_1 = (\nabla_\Theta \otimes id) (k_a) \quad \mbox{and} \quad k_2 = (id \otimes \nabla_\Theta)(k_a).$$ Then, there exist $b_1, b_2 \in \Sigma_{\rho,\delta}^{m+1}(\RR_\Theta)$ satisfying that $k_{b_1} = k_1$ and $k_{b_2} = k_2$.
\end{lemma}

\dem It is easily checked that $$(\nabla_\Theta \otimes id) (k_a) = \sum_{j=1}^n s(e_j) \otimes \int_{\R^n} \partial_{\Theta}^j [a(\xi) \lambda_\Theta(\xi)] \otimes \lambda_\Theta(\xi)^* \, d\xi = k_{b_1}$$ where $b_1(\xi) = \nabla_\Theta (a) + 2\pi i s(\xi) \otimes a(\xi)$ takes values in $\mathcal{L} (\mathbb{F}_n) \bar\otimes \RR_\Theta$, we omit the extra tensor component just to simplify our notation. A simple calculation also gives that $b_2(\xi) = - 2\pi i s(\xi) \otimes a(\xi)$. It is clear that $\nabla_\Theta(a) \in \Sigma^{m + \delta}_{\rho,\delta}(\RR_\Theta) \subset \Sigma^{m + 1}_{\rho,\delta}(\RR_\Theta)$ while the inclusion for $s(\xi) \otimes a(\xi)$ follows by Leibniz rule and $\partial_{\Theta, \xi}^j \xi = \partial_{\xi}^j \xi$. \fin

\begin{lemma} \label{PSD.lemmaKernelMult}
Given $a \in \S(\R^n; \S_\Theta)$, let $k_a$ be the kernel of $\Psi_a$. Recall that $\pi_\Theta(P)$ is a distribution in $\S_{\Theta \oplus \Theta}$ for any polynomial $P$ in $\xi_1, ..., \xi_n$. Then, the following identities hold in the sense of distributions for all $\alpha \in \Z_+^n$
\begin{eqnarray*}
k_a \bullet \pi_\Theta ( (2 \pi i z )^\alpha) & = & k_{\partial_\xi^\alpha a}, \\
\pi_\Theta \big( (2 \pi i z)^\alpha \big) \bullet k_a & = & k_{\partial_{\Theta, \xi}^\alpha a}.
\end{eqnarray*}
\end{lemma}

\dem Note that 
$$k_a \bullet \pi_\Theta(\exp_\zeta) = \Big( \int_{\R^n} (a(\xi) \otimes \1) \bullet \pi_\Theta(\exp_{\xi + \zeta}) \, d\xi\Big) = k_{a( \, \cdot \, - \, \zeta)}.$$
Taking derivatives formally gives
$$k_a \bullet \pi_\Theta(2 \pi i z_j) = k_a \bullet \frac{d}{ds} {\Big|}_{s = 0}  \pi_\Theta(\exp_{s e_j}) = k_{\partial_\xi^j a}.$$
This symbolic calculation can be justified in the distributional sense. For the second identity, we recall the identity $(\partial_{\Theta,\xi}^ja)(\xi) = \lambda_\Theta(\xi) \partial_\xi^j \big\{ \lambda_\Theta(\xi)^\ast a(\xi) \lambda_\Theta(\xi) \big\} \lambda_\Theta(\xi)^\ast$ and notice that 
$$\int_{\R^n} (a(\xi) \otimes \1) \bullet \pi_\Theta(\exp_{\xi}) \, d\xi = \int_{\R^n} \pi_\Theta(\exp_{\xi}) \bullet (\lambda_\Theta(\xi)^\ast a(\xi) \lambda_\Theta(\xi) \otimes \1) \, d\xi.$$ Therefore, arguing as above, we obtain the identity for left multiplication. \fin

\begin{lemma} \label{PSD.local}
Let $\psi_j(\xi) = \phi(2^{-j}\xi) - \phi(2^{-j+1}\xi)$ be a standard partition of unity in $\R^n$ from a smooth, radial and compactly supported $\phi$. If we let $a_j(\xi) = a(\xi) \psi_j(\xi)$ for $a \in \Sigma_{\rho,\delta}^m(\RR_\Theta)$ and $\ell_1, \ell_2 \ge 0$, we have
$$\Big\| \dTh^{\ell_1} \bullet k_{a_j} \bullet \dTh^{\ell_2} \Big\|_{\RR_\Theta \bar\otimes \RR_\Theta^\mathrm{op}}
\, \le \, C_{\ell_1,\ell_2} \ 2^{j ( n + m - \rho (\ell_1 + \ell_2))}.$$
\end{lemma}

\dem It is clear that 
$$\partial_{\Theta,\xi}^\alpha a_j = \sum_{\beta + \gamma = \alpha} a (\partial^{\beta}_\xi \psi_j) + (\partial_{\Theta,\xi}^{\gamma} a) \psi_j.$$
Since $|\partial_\xi^\beta \psi_j(\xi)| \lesssim 2^{-j |\beta|}$ when $\langle \xi \rangle \sim 2^j$, the $a_j$'s are in $\Sigma_{\rho,\delta}^m(\RR_\Theta)$ with constants independent of $j$. 
Assume first that $\ell_1, \ell_2$ are even numbers, $\ell_k = 2 \mathrm{N}_k$. Then the ${\ell_k}$-th power of $|z|$ is a polynomial of the form
$$|z|^{\ell_k} = \sum_{|\alpha| = \mathrm{N}_k} z^{2\alpha}.$$
Applying Lemma \ref{PSD.lemmaKernelMult} gives that
\begin{equation*}
\Big\| \dTh^{\ell_1} \bullet k_{a_j} \bullet  \dTh^{\ell_2} \Big\|_{\RR_\Theta \bar\otimes \RR_\Theta^\mathrm{op}}
= \frac{1}{(2 \pi)^{\ell_1 + \ell_2}} \Big\| \sum_{|\alpha_1| = \mathrm{N}_1} \sum_{|\alpha_2| = \mathrm{N}_2} k_{\partial_{\Theta,\xi}^{2 \alpha_1} \partial_\xi^{2 \alpha_2} a_j} \Big\|_{\RR_\Theta \bar\otimes \RR_\Theta^\mathrm{op}}
\end{equation*}
and for each of the terms we have the estimate
\begin{eqnarray*}
\big\| k_{\partial_{\Theta,\xi}^{2 \alpha_1} \partial_\xi^{2 \alpha_2} a_j} \big\|_{\RR_\Theta \bar\otimes \RR_\Theta^\mathrm{op}}
& = & \Big\| \int_{\R^n} \big( \partial_{\Theta,\xi}^{2 \alpha_1} \partial_\xi^{2 \alpha_2} a_j(\xi) \otimes \1 \big) \bullet \pi_\Theta(\exp_\xi) \, d\xi  \Big\|_{\RR_\Theta \bar\otimes \RR_\Theta^\mathrm{op}} \\
& \lesssim & 2^{j n} \sup_{|\xi| \lesssim 2^j} \big\| (\partial_\xi^{2 \alpha_1} \partial_{\Theta,\xi}^{2 \alpha_2} a_j)(\xi) \big\|_{\RR_\Theta} \ \lesssim \ 2^{j ( n + m - \rho (\ell_1 + \ell_2))}.
\end{eqnarray*}
For general (noneven) $\ell_1, \ell_2$ we proceed by interpolation. Note that the norm of $k_a$ is not altered under left/right multiplication by $\dTh^{i s}$ for any $s \in \R$. Therefore we have a bounded and holomorphic function $\zeta \mapsto \pi_\Theta(|z|^{\ell_1 + 2\zeta}) \bullet k_{a} \bullet \pi_\Theta(|z|^{\ell_2})$ defined in the band $0 < \Re(\zeta) < 1$. An application of the three lines lemma gives the bound for any $\ell_1$, the same follows for $\ell_2$. This completes the proof. \fin

\begin{proposition} \label{PSD.sizeGrad}
Given $a \in \Sigma_{\rho, \rho}^m(\RR_\Theta)$ and $m_1, m_2 > 0$  
$$\Big\| \dTh^{m_1} \bullet k_a \bullet \dTh^{m_2} \Big\|_{\RR_\Theta \bar\otimes \RR_\Theta^\mathrm{op}} \lesssim 1$$
provided $\rho(m_1 + m_2) > n + m$ for $\rho < 1$ or $\rho(m_1 + m_2) \ge n + m$ for $\rho=1$.
\end{proposition}

\dem Let $$\rho_j(\xi) = \phi(2^{-j} \xi) - \phi(2^{-j+1} \xi)$$ be another partition of unity |this time for $j \in \mathbb{Z}$| and set $b_j = \pi_\Theta(\rho_j)$. Then
$$\dTh^{m_1} \bullet k_a \bullet \dTh^{m_2} = \sum_{j,k \ge 0} + \sum_{j,k \le 0} + \sum_{j \cdot k <0} \big( b_j \bullet \dTh^{m_1} \bullet k_a \bullet \dTh^{m_2} \bullet b_k \big) = A_{+} + A_{-} + A_{\pm}.$$
\noindent {\bf Estimate of $A_+$.} Letting $a_\ell = a \psi_\ell$ as in Lemma \ref{PSD.local}
$$\| A_+ \|_{\RR_\Theta \bar\otimes \RR_\Theta^\mathrm{op}} \le \sum_{j,k \ge 0} \sum_{\ell \ge 0} \Big\| b_j \bullet \dTh^{m_1} \bullet k_{a_\ell} \bullet \dTh^{m_2} \bullet b_k \Big\|_{\RR_\Theta \bar\otimes \RR_\Theta^\mathrm{op}} = \sum_{j,k \ge 0} \sum_{\ell \ge 0} A_+(j,k,\ell).$$
Pick $\ell_1$ and $\ell_2$ large enough (see below) and use Lemma \ref{PSD.local} to estimate $A_+(j,k,\ell)$
\begin{eqnarray*}
A_+(j,k,\ell) & \leq & \big\| \dTh^{\ell_1} \bullet k_{a_\ell} \bullet \dTh^{\ell_2} \big\|_{\RR_\Theta \bar\otimes \RR_\Theta^\mathrm{op}} \\
& \times & \big\| \pi_\Theta(\rho_j |z|^{(m_1 - \ell_1)}) \big\|_{\RR_\Theta \bar\otimes \RR_\Theta^\mathrm{op}} 
\big\| \pi_\Theta(\rho_k |z|^{(m_2 - \ell_2)}) \big\|_{\RR_\Theta \bar\otimes \RR_\Theta^\mathrm{op}} \\
& \lesssim & 2^{\ell (n + m - \rho (\ell_1 + \ell_2))} 2^{j (m_1 - \ell_1)} 2^{k (m_2 - \ell_2)}.
\end{eqnarray*}
Taking $\ell_1 > m_1$, $\ell_2 > m_2$ and $\rho(\ell_1 + \ell_2) > n + m$ we may sum over $j,k,\ell \ge 0$. 

\noindent {\bf Estimate of $A_{-}$.} Letting $a_\ell = a \psi_\ell$ once more, we get 
\begin{eqnarray*}
A_{-} \!\!\!\! & = & \!\!\!\! \sum_{j,k \le 0} \bigg\{ \sum_{\ell > |j|} + \sum_{\ell \leq |j|} \bigg\} \big( b_j \bullet \dTh^{m_1} \bullet k_{a_\ell} \bullet \dTh^{m_2} \bullet b_k \big) = \sum_{j,k \le 0} A_{-}^1(j,k) + A_{-}^2(j,k) \\ 
\\ \!\!\!\! & = & \!\!\!\! \ \underbrace{\sum_{j \le k \le 0} A_{-}^1(j,k)}_{ A_{-}^{11}} \hskip3pt + \underbrace{\sum_{k < j \le 0} A_{-}^1(j,k)}_{A_{-}^{12}} \hskip3pt + \underbrace{\sum_{j \le k \le 0} A_{-}^2(j,k)}_{A_{-}^{21}} \hskip3pt + \underbrace{\sum_{k < j \le 0} A_{-}^2(j,k)}_{A_{-}^{22}}.
\end{eqnarray*}
First, we may bound $A_{-}^1(j,k)$ and $A_{-}^2(j,k)$ in norm via Lemma \ref{PSD.local}
\begin{eqnarray*}
\big\| A_{-}^1(j,k) \big\|_{\RR_\Theta \bar\otimes \RR_\Theta^\mathrm{op}}
& \lesssim & \sum_{\ell > |j|} 2^{j(m_1 - \ell_1)} 2^{\ell((n + m) - \rho(\ell_1 + \ell_2))} 2^{k(m_2 - \ell_2)}, \\
\big\| A_{-}^2(j,k) \big\|_{\RR_\Theta \bar\otimes \RR_\Theta^\mathrm{op}}
& \lesssim & \sum_{\ell \le |j|} 2^{j(m_1 - \ell_1)} 2^{\ell((n + m) - \rho(\ell_1 + \ell_2))} 2^{k(m_2 - \ell_2)}.
\end{eqnarray*}
\begin{itemize}
\item[$A_{-}^{11}$)] Taking $\rho (\ell_1 + \ell_2) > n + m$ and $r = k - j \ge 0$, we get 
\begin{eqnarray*}
\hskip25pt \big\| A_{-}^1(j,k) \big\|_{\RR_\Theta \bar\otimes \RR_\Theta^\mathrm{op}} 
& \lesssim & 2^{j(m_1 - \ell_1)} 2^{|j|(n + m - \rho(\ell_1 + \ell_2))} \, 2^{k(m_2 - \ell_2)} \\ 
& = & 2^{r(m_2 - \ell_2)} \, 2^{|j|(n + m + (1-\rho)(\ell_1 + \ell_2) - (m_1 + m_2))}. 
\end{eqnarray*}
If $\rho < 1$, our condition $\rho(m_1 + m_2) > n + m$ allows us to pick $\ell_1, \ell_2$ satisfying $\ell_1 + \ell_2 = m_1 + m_2$ and $\ell_2 > m_2$. If $\rho=1$ we pick $\ell_j > m_j$ for $j=1,2$. In both cases we get the inequalities $n + m + (1-\rho)(\ell_1 + \ell_2) - (m_1 + m_2) \le 0$ and $m_2 - \ell_2 > 0$. This gives
\begin{eqnarray*}
\hskip30pt \big\| A_{-}^{11} \big\|_{\RR_\Theta \bar\otimes \RR_\Theta^\mathrm{op}} \!\!\! & \le & \!\!\! \sum_{r \geq 0} \Big\| \sum_{j \le 0} A_{-}^1(j,j+r) \Big\|_{\RR_\Theta \bar\otimes \RR_\Theta^\mathrm{op}} \\ \!\!\! & \le & \!\!\! \sum_{r \geq 0} r \sup_{j \le 0} \big\| A_{-}^1(j,j+r) \big\|_{\RR_\Theta \bar\otimes \RR_\Theta^\mathrm{op}} 
\lesssim \sum_{r \geq 0} r 2^{r(m_2 - \ell_2)} \lesssim 1.
\end{eqnarray*}
By almost orthogonality of $b_j$'s the sum inside the norm is a $r$-th diagonal operator, dominated by $r$ times the supremum of the norms of its entries. 

\vskip3pt

\item[$A_{-}^{21}$)] Letting $r = k-j \ge 0$, we get
\begin{eqnarray*}
\hskip35pt \lefteqn{\big\| A_{-}^2(j,k) \big\|_{\RR_\Theta \bar\otimes \RR_\Theta^\mathrm{op}} 
\lesssim 2^{j(m_1 - \ell_1)} \max \Big\{1, 2^{|j|(n + m - \rho(\ell_1 + \ell_2))} \Big\} 2^{k(m_2 - \ell_2)}} \\ 
\!\!\! & = & \!\!\! 2^{r(m_2 - \ell_2)} \max \Big\{ 2^{|j|((\ell_1 + \ell_2) - (m_1 + m_2))}, 2^{|j|(n + m + (1-\rho)(\ell_1 + \ell_2) - (m_1+m_2))} \Big\}.
\end{eqnarray*}
Then any choice with $\ell_1+\ell_2 = m_1 + m_2$ and $\ell_2 > m_2$ gives $\| A_{-}^{21}\| \lesssim 1$. 

\vskip3pt

\item[$A_{-}^{12}$)] We may write $$\hskip30pt A_{-}^{12} = \sum_{k < j \le 0} \Big\{ \sum_{\ell > |k|} + \sum_{|j| < \ell \le |k|} \Big\} \big( b_j \bullet \dTh^{m_1} \bullet k_{a_\ell} \bullet \dTh^{m_2} \bullet b_k \big) = A_{-}^{121} + A_{-}^{122}.$$ Then, $A_{-}^{121}$ is estimated exactly as $A_{-}{11}$. On the other hand, the estimate of $A_{-}^{122}$ is very much similar to that of $A_{-}^{21}$, we leave the details to the reader. 

\vskip3pt

\item[$A_{-}^{22}$)] Interchanging roles of $(j,k)$, \, $A_{-}^{22}$ is estimated as $A_{-}^{21}$ above (even simpler).
\end{itemize} 

\noindent {\bf Estimate of $A_\pm$.} Since the conditions on $m_1, m_2$ are symmetric in the statement and the sum $\sum_{i \cdot j < 0}$ splits into $\sum_{j < 0 < k} + \sum_{k < 0 < j}$, it suffices to estimate one of these two sums. Arguing as above, if we pick $\ell_1, \ell_2$ so that $\rho(\ell_1 + \ell_2) > n+m$, the problem reduces to estimate $$\sum_{j < 0 < k} \Big\{ \sum_{\ell > k} + \sum_{\ell \le k} \Big\} 2^{j(m_1 - \ell_1)} 2^{\ell(n+m - \rho(\ell_1 + \ell_2))} 2^{k(m_2 - \ell_2)} = A+B.$$ Again as above, we pick $r = k-j > 0$ and obtain
\begin{eqnarray*}
A \!\!\! & \lesssim & \!\!\! \sum_{r > 0} \sum_{k > 0} 2^{-r(m_1 - \ell_1)} 2^{k(n+m - (1+\rho)(\ell_1 + \ell_2) + m_1 + m_2)},\\
B \!\!\! & \lesssim & \!\!\! \sum_{r > 0} \sum_{k > 0} 2^{-r(m_1 - \ell_1)} \max \Big\{ 2^{k((m_1+m_2) - (\ell_1+\ell_2))}, 2^{k(n+m - (1+\rho)(\ell_1+\ell_2) + m_1+m_2)} \Big\}. 
\end{eqnarray*}
Since $n+m \le \rho(m_1+m_2)$, it suffices to take $\ell_1 > m_1$ and $\ell_1+\ell_2 > m_1 + m_2$. \fin
 
\begin{theorem} \label{PSD.lpBound}
If $a, a_\dagger^* \in \Sigma_{1, 1}^0(\RR_\Theta)$, we have 
\begin{eqnarray*}
\big\| \Psi_a: \mathrm{H}_1(\RR_\Theta) \to L_1(\RR_\Theta) \big\|_{\mathrm{cb}} < \infty, & \\
\big\| \Psi_a: L_\infty(\RR_\Theta) \to \mathrm{BMO}(\RR_\Theta) \big\|_{\mathrm{cb}} < \infty. &
\end{eqnarray*}
In particular, $\Psi_a: L_p(\RR_\Theta) \to L_p(\RR_\Theta)$ is completely bounded when $1 < p < \infty$. 
\end{theorem}

\dem According to our Calder\'on-Zygmund extrapolation in Theorem \ref{CZExt4}, it suffices to see that $\Psi_a$ is $L_2$-bounded and its kernel $k_a$ satisfies the CZ conditions there. The $L_2$-boundedness follows from the quantum form of Bourdaud's condition in Theorem \ref{PSD.boudreaud}. On the other hand, according to Lemma \ref{Mx.Gradiente}, both $(\nabla_\Theta \otimes id)(k_a)$ and $(id \otimes \nabla_\Theta)(k_a)$ belong to $\Sigma_{1,1}^1(\RR_\Theta)$. In particular, Proposition \ref{PSD.sizeGrad} yields the CZ kernel conditions which we need for $m_1+m_1=n+1$. \fin

\begin{corollary} \label{PSD.lpBound2}
If $a, a_\dagger^* \in \Sigma_{1, 1}^m(\RR_\Theta)$, we have 
$$\big\| \Psi_a: \mathrm{W}_{p,s}(\RR_\Theta) \to \mathrm{W}_{p,s-m}(\RR_\Theta) \big\|_{\mathrm{cb}} < \infty \quad \mbox{for every} \quad 1 < p < \infty.$$
\end{corollary}

\dem We have that 
\begin{equation*}
\xymatrix{
L_p(\RR_\Theta) \ar[rrrr]^{(\1 - \Delta_\Theta)^{\frac{s - m}{2}} \Psi_a (\1 - \Delta_\Theta)^{-\frac{s}{2}}} & & & & L_p(\RR_\Theta) \ar[dl]^{\, \, (\1 - \Delta_\Theta)^{-\frac{s-m}{2}}} \\
& \mathrm{W}_{p,s}(\RR_\Theta) \ar[rr]^{\Psi_a} \ar[ul]^{\, (\1 - \Delta_\Theta)^{\frac{s}{2}} \, } & & \mathrm{W}_{p,s-m}(\RR_\Theta) &
}
\end{equation*}
where $(\1 - \Delta_\Theta)^{u/2}: \mathrm{W}_{p,s}(\RR_\Theta) \to \mathrm{W}_{p,s-u}(\RR_\Theta)$ are complete isometries. On the other hand, the complete $L_p$-boundedness of $(\1 - \Delta_\Theta)^{(s - m)/2} \Psi_a (\1 - \Delta_\Theta)^{-s/2}$ follows from Theorem \ref{PSD.lpBound} once we observe that this pseudodifferential operator and its adjoint are associated to symbols in $\Sigma_{1,1}^0(\RR_\Theta)$, which in turn follows from the composition rules for $\Sigma_{\rho,\delta}^m(\RR_\Theta)$ established in Remark \ref{PSD.remarkCompSigma}. \fin

\begin{remark} \label{RemS1delta}
\emph{The assertions in Theorem \ref{PSD.lpBound} and Corollary \ref{PSD.lpBound2} remain valid for symbols $a \in \Sigma_{1,\delta}^0(\RR_\Theta)$ with $0 \le \delta < 1$. Indeed, $\Sigma_{1,\delta}^0(\RR_\Theta) \subset S_{1,\delta}^0(\RR_\Theta) \subset S_{1,1}^0(\RR_\Theta)$ and the middle class is stable under adjoints. Thus, we may apply our Bourdaud's condition as we did in the proof of Theorem \ref{PSD.lpBound}. In addition, $\Sigma_{1,\delta}^0(\RR_\Theta) \subset \Sigma_{1,1}^0(\RR_\Theta)$ so that Proposition \ref{PSD.sizeGrad} applies. The argument for Sobolev spaces is similar.}
\end{remark}

\begin{remark} \label{RemJustn+2}
\emph{A careful analysis of our proof for Theorem \ref{PSD.lpBound} and Corollary \ref{PSD.lpBound2} yields that the condition $\Sigma_{1,1}^0(\RR_\Theta)$ can be replaced by the weaker condition below $$\big| \partial_\Theta^\beta \partial_{\Theta,\xi}^{\alpha_1} \partial_\xi^{\alpha_2} a(\xi) \big| + \big| \partial_\Theta^\beta \partial_{\Theta,\xi}^{\alpha_1} \partial_\xi^{\alpha_2} a_\dagger^*(\xi) \big| \, \le \, C_{\alpha_1, \alpha_2, \beta} \langle \xi \rangle^{- |\alpha_1 + \alpha_2| + |\beta|}$$ for $|\alpha_1 + \alpha_2| \le n+2$ and $|\beta| \le 1$. Indeed, according to Remark \ref{PSD.FiniteDerivatives}, Bourdaud's condition in Theorem \ref{PSD.boudreaud} can be weakened to $|\alpha_1+\alpha_2| \le n+1$ and $|\beta| \le 1$. Moreover, our proof of Proposition \ref{PSD.sizeGrad} for $m_1 + m_2 = n+1$ requires $|\alpha_1 + \alpha_2| \le n+2$.} 
\end{remark}

The $L_p$-theory for exotic symbols $\Sigma_{\rho,\rho}^m(\RR_\Theta)$ $(\rho < 1)$ is only possible due to the regularizing effect of a negative degree $m$. Fefferman proved in \cite{Fe} the $L_p$ bounds for the critical index $m = - (1-\rho) \frac{n}{2}$. The noncritical range was obtained by Hirschman and Wainger \cite{Hi,Wa} (constant coefficients) and H\"ormander \cite{Ho} (general symbols). Standard interpolation arguments yield even more general statements \cite[VII 5.12]{St}. Now we shall prove (non optimal) inequalities of this kind in $\RR_\Theta$ with applications below for $L_p$-regularity of elliptic PDEs. Namely, in what follows we shall write $\mathrm{N}$ for the best possible constant in Remark \ref{RemOptimalCV}. As explained there we suspect that any $\mathrm{N} > n/4$ is valid and this would be optimal, as it is the case in the Euclidean theory. Consider the index $$\Lambda_{\rho,n} = - (1-\rho) \max \big\{ 2\mathrm{N}, n+2 \big\}.$$ It follows from arguments in \cite{T1} that $2\mathrm{N}$ is (at least) less or equal than $3n+2$.

\begin{corollary} \label{CorollaryRegularizingTerm}
Let $a \in \Sigma_{\rho,\rho}^m(\RR_\Theta)$ be a symbol satisfying $m \le \Lambda_{\rho,n}$ for some $\rho < 1$. Then, the pseudodifferential operator $\Psi_a$ satisfies the following estimates for $1 < p < \infty$  
\begin{eqnarray*}
\hskip29pt \big\| \Psi_a: \mathrm{H}_1(\RR_\Theta) \to L_1(\RR_\Theta) \big\|_{\mathrm{cb}} < \infty, & \\
\big\| \Psi_a: L_p(\RR_\Theta) \hskip1pt \to L_p(\RR_\Theta) \big\|_{\mathrm{cb}} < \infty, & \\
\big\| \Psi_a: L_\infty(\RR_\Theta) \to \mathrm{BMO}(\RR_\Theta) \big\|_{\mathrm{cb}} < \infty. & 
\end{eqnarray*}
Moreover, if $m$ is any real number and $\ell = m - \Lambda_{\rho,n}$  
\begin{eqnarray*}
\big\| \Psi_a: \mathrm{W}_{p,s}(\RR_\Theta) \to \mathrm{W}_{p,s-\ell}(\RR_\Theta) \big\|_{\mathrm{cb}} < \infty.
\end{eqnarray*}
\end{corollary}

\dem Since $a \in \Sigma_{\rho,\rho}^m(\RR_\Theta)$ and $m = \ell + \Lambda_{\rho,n}$ 
\begin{eqnarray*}
\langle \xi \rangle^{m - \rho |\alpha_1 + \alpha_2| + \rho |\beta|} \!\!\! & = & \!\!\! \langle \xi \rangle^{m + (1-\rho) |\alpha_1 + \alpha_2| - |\alpha_1 + \alpha_2| + \rho |\beta|} \\ 
\!\!\! & \le & \!\!\! \langle \xi \rangle^{m + (1-\rho) \max\{2\mathrm{N}, n+2\} - |\alpha_1 + \alpha_2| + \rho |\beta|} \le \langle \xi \rangle^{\ell - |\alpha_1 + \alpha_2| + \rho |\beta|}
\end{eqnarray*} 
as long as $|\alpha_1 + \alpha_2| \le \max\{2\mathrm{N},n+2\}$. This means that $a$ satisfies the H\"ormander condition $\Sigma_{1,\rho}^\ell(\RR_\Theta)$ for $\partial_{\Theta,\xi}$, $\partial_\xi$ of order up to $|\alpha_1 + \alpha_2| \le \max \{2\mathrm{N},n+2\}$. For the first assertion we apply Theorem A. The $L_2$-boundedness is guaranteed by our Calder\'on-Valillancourt theorem since $\Sigma_{1,\rho}^0(\RR_\Theta) \subset S_{\rho,\rho}^0(\RR_\Theta)$ and $2\mathrm{N}$ $\xi$-derivatives suffice, according to Remark \ref{RemOptimalCV}. Next, inclusion $\Sigma_{1,\rho}^0(\RR_\Theta) \subset \Sigma_{1,1}^0(\RR_\Theta)$ together with the fact that  Proposition \ref{PSD.sizeGrad} only requires $|\alpha_1 + \alpha_2| \le n+2$ |see Remark \ref{RemJustn+2}| imply that the CZ kernel conditions also hold. This proves that the first assertion follows from Theorem A. Then, the second assertion follows by adapting the argument in the proof of Corollary \ref{PSD.lpBound2}. Indeed, arguing as in Remark \ref{PSD.remarkCompSigma} we deduce that $$a \in \Sigma_{1,\rho}^\ell(\RR_\Theta) \ \ \Rightarrow \ \ (\1 - \Delta_\Theta)^{\frac{s-\ell}{2}} \Psi_a (\1 - \Delta_\Theta)^{-\frac{s}{2}} = \Psi_b \quad \mbox{for some} \quad b \in \Sigma_{1,\rho}^0(\RR_\Theta).$$ In fact, the same holds limiting the conditions above to a prescribed number of derivatives $\partial_{\Theta,\xi}$ and $\partial_\xi$. Hence, we apply Theorem A once more to conclude. \fin

\begin{remark}
\emph{It is very tempting to claim that Corollary \ref{CorollaryRegularizingTerm} holds for the index $\Lambda_{\rho,n} = - (1-\rho)(n+2)$ since it is reasonable to think that the above result follows from a direct combination of Remarks \ref{RemS1delta} and \ref{RemJustn+2} above. However, at the time of this writing, we are not able to circumvent the adjoint stability used in Remark \ref{RemS1delta}, since we need it for H\"ormander conditions limited to a prescribed number of derivatives. The product stability used above is straightforward instead.}
\end{remark}

\begin{remark} \label{CriticalIndex}
\emph{$L_p$-boundedness up to the critical index $m \hskip-2pt = \hskip-2pt -(1-\rho) \frac{n}{2}$ is still open.}
\end{remark}

\section{\bf $L_p$ regularity for elliptic PDEs}

In this section, we illustrate our results with a basic application to elliptic PDEs in quantum Euclidean spaces. Given $0 \le \delta \le \rho \le 1$ and $m \in \R$, a symbol $a \in S_{\rho,\delta}^m(\RR_\Theta)$ is called elliptic of order $m$ when there exist constants $C,R > 0$ for which the following inequality holds $$\big| a(\xi) \big| \ge C |\xi|^m \quad \mbox{for all} \quad |\xi| \ge R.$$ A prototypical example of elliptic symbol of order $2$ is given by $a(\xi) = \xi^\mathrm{t} A \xi$ for some uniformly positive definite $A \in M_n(\RR_\Theta)$. We shall be interested in the elliptic PDE $$\Psi_a(u) = \varphi$$ with data $\varphi$ in the Sobolev space $\mathrm{W}_{p,s}(\RR_\Theta)$ and $a \in \Sigma_{1,\delta}^m(\RR_\Theta)$. $L_p$-regularity means that, no matter which a priori regularity do we have in a given solution $u$, it must belong at least to the Sobolev space $\mathrm{W}_{p,s+m}(\RR_\Theta)$. When the regularity gained is smaller than $m$ we speak about hypoellipticity. In the Euclidean case, elliptic regularity arises naturally for $(\rho,\delta) = (1,0)$ and still holds for $\rho=1$, whereas the case $\rho < 1$ leads to hypoelliptic scenarios \cite{T1}. As we shall see, this is also the case in the quantum setting. Equipped with our results so far, the main obstruction we shall need to overcome will be to construct suitable parametrices for symbols in $\Sigma$-classes, for which we can not use product stability in that class. Our first step yields Sobolev $p$-estimates $\mathrm{W}_{p,s}(\RR_\Theta) \to \mathrm{W}_{p,s-\ell}(\RR_\Theta)$ for symbols in $S_{\rho,\delta}^m(\RR_\Theta)$ instead of $\Sigma_{\rho,\delta}^\ell(\RR_\Theta)$, provided the order $m$ is small enough. 

\begin{lemma} \label{LemmaNegativeIndex}
Given $s, \ell \in \R$, we have $$\Psi_a: \mathrm{W}_{p,s}(\RR_\Theta) \to \mathrm{W}_{p,s-\ell}(\RR_\Theta)$$ provided $a \in S_{\rho,\delta}^m(\RR_\Theta)$ with degree $m + (1+\delta) \max \big\{ 2 \mathrm{N}, n+2 \big\} \le \ell$.   
\end{lemma}

\dem Arguing as in the proof of Corollary \ref{CorollaryRegularizingTerm}, it suffices to see that $a$ satisfies the $\Sigma_{1,\delta}^\ell(\RR_\Theta)$-condition for $\partial_{\Theta,\xi}$ and $\partial_\xi$ of order up to $\max \{2\mathrm{N},n+2\}$. Then recalling that $$\partial_{\Theta,\xi}^j = \partial_\xi^j + \frac{1}{2\pi i} \sum_{k=1}^n \Theta_{jk} \partial_\Theta^k,$$ we easily get the following estimate 
\begin{eqnarray*}
\big\| \partial_\Theta^\beta \partial_{\Theta,\xi}^{\alpha_1} \partial_\xi^{\alpha_2} a(\xi) \big\|_{\RR_\Theta} & \lesssim & \sum_{\alpha_{11} + \alpha_{12} = \alpha_1} \sum_{|\gamma| = |\alpha_{11}|} \big\| \partial_\Theta^{\beta + \gamma} \partial_\xi^{\alpha_{12} + \alpha_2} a(\xi) \big\|_{\RR_\Theta} \\ & \lesssim & \sum_{\alpha_{11} + \alpha_{12} = \alpha_1} \sum_{|\gamma| = |\alpha_{11}|} \langle \xi \rangle^{m - \rho |\alpha_{12} + \alpha_2| + \delta |\beta + \gamma|}. 
\end{eqnarray*} 
When $|\alpha_1 + \alpha_2| \le \max\{2\mathrm{N},n+2\}$, we use 
\begin{eqnarray*}
\lefteqn{m - \rho |\alpha_{12} + \alpha_2| + \delta |\beta + \gamma|} \\ 
& = & m + (1 - \rho) |\alpha_{12} + \alpha_2| + (1+\delta) |\gamma| - |\alpha_1 + \alpha_2| + \delta |\beta| \\
& \le & m + (1 + \delta) |\alpha_1 + \alpha_2| - |\alpha_1 + \alpha_2| + \delta |\beta| \ \le \ \ell -|\alpha_1 + \alpha_2| + \delta |\beta|
\end{eqnarray*}
since $m + \max \{2\mathrm{N}, n+2\} \le \ell$ by  hypothesis. This completes the proof. \fin

\begin{lemma} \label{LemmaParametrix}
Let $a \in \Sigma_{\rho,\delta}^m(\RR_\Theta)$ be an elliptic symbol for some $0 \le \delta < \rho \le 1$ and degree $m$. Let $\ell = m + \Lambda_{\rho,n}$. Then, for every $k \in \N$, there exist symbols $b_k$ and $c_k$ satisfying the following properties$\hskip1pt :$ 
\begin{itemize}
\item[\emph{i)}] $\Psi_{b_k} \Psi_a = id - \Psi_{c_k}$,

\vskip1pt

\item[\emph{ii)}] $c_k \in S_{\rho,\delta}^{k \gamma}(\RR_\Theta)$ with $\gamma = \delta - \rho < 0$,

\vskip1pt

\item[\emph{iii)}] If $\rho=1$, then $\Psi_{b_k}: \mathrm{W}_{p,s}(\RR_\Theta) \to \mathrm{W}_{p,s+m}(\RR_\Theta)$ for all $s \in \R$.

\vskip1pt

\item[\emph{iv)}] If $\rho<1$, then $\Psi_{b_k}: \mathrm{W}_{2,s}(\RR_\Theta) \to \mathrm{W}_{2,\hskip1pt s+ \hskip1pt \ell} \hskip1pt (\RR_\Theta)$ \hskip1pt for all $s \in \R$.
\end{itemize}
In fact, the last assertion holds under the weaker assumption that $a \in S_{\rho,\delta}^m(\RR_\Theta)$.
\end{lemma}

\dem Let $$b_1(\xi) = \big( 1 - \phi(\xi) \big) a^{-1}(\xi)$$ where $\phi$ is a smooth function which is identically $1$ in $\mathrm{B}_R(0)$ and vanishes outside $\mathrm{B}_{R+1}(0)$. Here $R$ is determined by the ellipticity of $a$, so that $|a(\xi)| \ge C |\xi|^m$ for $|\xi| \ge R$. We claim that  
\begin{itemize}
\item[A)] $b_1 \in \Sigma_{\rho,\delta}^{-m}(\RR_\Theta)$,

\item[B)] $\Psi_{b_1} \Psi_a = id - \Psi_{c_1}$ for some $c_1 \in S_{\rho,\delta}^{\gamma}(\RR_\Theta)$.
\end{itemize}
Assuming the claim, let $b_k$ and $c_k$ be determined by $$\Psi_{b_k} = \sum_{j=0}^{k-1} \Psi_{c_1}^j \Psi_{b_1} \quad \mbox{and} \quad \Psi_{c_k} = \Psi_{c_1}^k.$$ 
\begin{itemize}
\item[i)] $\displaystyle \Psi_{b_k} \Psi_a = \sum_{j=0}^{k-1} \Psi_{c_1}^j \Psi_{b_1} \Psi_a = \sum_{j=0}^{k-1} \Psi_{c_1}^j \big( id - \Psi_{c_1} \big) = id - \Psi_{c_1}^k = id - \Psi_{c_k}$.

\item[ii)] $c_k \in S_{\rho,\delta}^{k \gamma}(\RR_\Theta)$ with $\gamma = \delta - \rho < 0$ follows from Corollary \ref{PSD.stabilityA}, since $\delta < \rho$.

\vskip10pt

\item[iii)] We may not use our results directly since we ignore whether or not $b_k$ belongs to the right $\Sigma$-class, due to the lack (so far) of stability results for the product of symbols in these classes. However, when $\rho=1$ we know from claim A) above and Remark \ref{RemS1delta} that $$\hskip10pt \Psi_{b_1}: \mathrm{W}_{p,s}(\RR_\Theta) \to \mathrm{W}_{p,s+m}(\RR_\Theta).$$ Let us note in passing that Corollary \ref{CorollaryRegularizingTerm} would also do the job here for $\rho<1$ and $\ell$ in place of $m$. Next, it suffices to show that $\Psi_{c_1}^j$ takes the Sobolev space $\mathrm{W}_{p,s+m}(\RR_\Theta)$ to itself. This is clear for the identity map with $j=0$. On the other hand, the boundedness for $j > 0$ trivially follows from the case $j=1$. Since $(a, b_1) \in \Sigma_{1,\delta}^m(\RR_\Theta) \times \in \Sigma_{1,\delta}^{-m}(\RR_\Theta)$, the boundedness of $\Psi_{c_1} = id - \Psi_{b_1} \Psi_a$ follows again from Remark \ref{RemS1delta}. 

\vskip10pt

\item[iv)] By the product stability of $S$-classes from Section \ref{CompAdj} and according to claims A) and B) we know that $b_k \in S_{\rho,\delta}^{-m}(\RR_\Theta)$ and the result follows from the argument in Corollary \ref{PSD.lpBound2}. Indeed, it works in $L_2$ when $\Sigma$-classes are replaced by $S$-classes since we just need to apply Calder\'on-Vaillancourt and composition with the right powers of $\1 - \Delta_\Theta$ in that case.
\end{itemize}
Once we have proved the assertion, it remains to justify our claim. Point A) follows easily once we express the involved derivatives of $a^{-1}(\xi)$ in terms of those for $a$. It is clear that $\partial^j \big( a^{-1}\big)(\xi) + a^{-1}(\xi) \hskip1pt \partial^j(a)(\xi) \hskip1pt  a^{-1}(\xi) = 0$ for the derivations $\partial \in \{\partial_\xi, \partial_{\Theta,\xi}, \partial_\Theta\}$. By ellipticity we obtain the estimates below for $|\xi| \ge R$
\begin{eqnarray*}
\big\| \partial_{\xi}^j \big( a^{-1}\big)(\xi) \big\|_{\RR_\Theta} & \lesssim & \langle \xi \rangle^{-m - \rho}, \\ 
\big\| \partial_{\Theta}^j \big( a^{-1}\big)(\xi) \big\|_{\RR_\Theta} & \lesssim & \langle \xi \rangle^{-m + \delta}, \\
\big\| \partial_{\Theta,\xi}^j \big( a^{-1}\big)(\xi) \big\|_{\RR_\Theta} & \lesssim & \langle \xi \rangle^{-m - \rho}.
\end{eqnarray*}
By Leibniz rule and induction we get $b_1 \in \Sigma_{\rho,\delta}^{-m}(\RR_\Theta)$ which proves A). Then B) follows from the product stability in Corollary \ref{PSD.stabilityA} as in \cite[Theorem III.1.3]{T1}. \fin

\begin{remark}
\emph{The above result for $p=2$ is still open for $1 < p < \infty$. According to point ii) and Lemma \ref{LemmaNegativeIndex}, we know that $\Psi_{c_1}^j$ is bounded on $\mathrm{W}_{p,s+\ell}(\RR_\Theta)$ for $j$ large enough. It would be tempting to deduce the result by complex interpolation with $j=0$. However, imaginary powers of $\Psi_{c_1}$ are generally unbounded in $L_p$ since the same happens for $\Psi_{c_1}$, due to Fefferman's critical index $-(1-\rho)n/2$. Indeed, $\Psi_{c_1}$ will not be bounded in $L_p$ or $\mathrm{W}_{p,s}$ when $|\gamma|$ is small enough and $\rho<1$.}
\end{remark}

\begin{remark} \label{RemParametrix}
\emph{In the absence of stability for products of symbols in $\Sigma$-classes |left open in Section \ref{CompAdj}| Lemmas \ref{LemmaNegativeIndex} and \ref{LemmaParametrix} give together a good substitute for many applications. Lemma \ref{LemmaParametrix} provides a parametrix $\Psi_{b_k}$ which, despite we ignore for the moment whether or not it lives in the right $\Sigma$-class, it does send $\mathrm{W}_{p,s}(\RR_\Theta)$ to the correct Sobolev space. Moreover, we know from Lemma \ref{LemmaNegativeIndex} that the same holds for the error term $\Psi_{c_k}$ provided $k$ is large enough, since $\gamma < 0$.}
\end{remark}

\begin{theorem} \label{Lp-Regularity}
Given $0 \le \delta < \rho \le 1$, consider $a \in \Sigma_{\rho,\delta}^m(\RR_\Theta)$ an elliptic symbol for some $m \in \R$ and let $\ell = m + \Lambda_{\rho,n}$. Given $1 < p < \infty$ and $r,s \in \R$, assume $\varphi \in \mathrm{W}_{p,s}(\RR_\Theta)$ and let $u$ solve $$\Psi_a(u) = \varphi$$ for some $u \in \mathrm{W}_{p,r}(\RR_\Theta)$. Then, the following estimates hold$\hskip1pt:$
\begin{itemize}
\item[\emph{i)}] If $\rho = 1$, we get $\|u\|_{\mathrm{W}_{p,s+m}(\RR_\Theta)} \lesssim \|u\|_{\mathrm{W}_{p,r}(\RR_\Theta)} + \|\varphi\|_{\mathrm{W}_{p,s}(\RR_\Theta)}.$

\item[\emph{ii)}] If $\rho < 1$ and $p=2$, we get $\|u\|_{\mathrm{W}_{2,s+\ell}(\RR_\Theta)} \lesssim \|u\|_{\mathrm{W}_{2,r}(\RR_\Theta)} + \|\varphi\|_{\mathrm{W}_{2,s}(\RR_\Theta)}.$
\end{itemize}
\end{theorem}

\dem According to Lemma \ref{LemmaParametrix} $$u - \Psi_{c_k} (u) = \Psi_{b_k} \Psi_a(u) = \Psi_{b_k}(\varphi)$$ for any $k \ge 0$. This gives in particular 
\begin{itemize}
\item[i)] If $\rho=1$
\begin{eqnarray*}
\|u\|_{p,s+m} & \le & \|\Psi_{b_k}(\varphi)\|_{p,s+m} + \|\Psi_{c_k} (u)\|_{p,s+m} \\
& \lesssim & \|\varphi\|_{p,s} +  \|\Psi_{c_k} (u)\|_{p,s+m}.
\end{eqnarray*}

\item[ii)] If $\rho < 1$ and $p=2$
\begin{eqnarray*} \|u\|_{2,s+\ell} & \le & \|\Psi_{b_k}(\varphi)\|_{2,s+\ell} + \|\Psi_{c_k} (u)\|_{2,s+\ell} 
\\ & \lesssim & \|\varphi\|_{2,s} +  \|\Psi_{c_k} (u)\|_{2,s+\ell}.
\end{eqnarray*}
\end{itemize}
Next, Lemma \ref{LemmaNegativeIndex} gives that $\Psi_{c_k}: \mathrm{W}_{p,r}(\RR_\Theta) \to \mathrm{W}_{p,s+\ell}(\RR_\Theta)$ for $k$ large enough. \fin

\begin{remark}
\emph{As in the Euclidean setting \cite{T1}, Theorem \ref{Lp-Regularity} above gives elliptic $L_p$-regularity in the H\"ormander class $\Sigma_{1,\delta}^m(\RR_\Theta)$ and hypoelliptic $L_2$-regularity in $\Sigma_{\rho,\delta}^m(\RR_\Theta)$ when $\rho < 1$. The latter result remains open for other values of $p \neq 2$. Compared to \cite{T1} our result for $p=2$ quantifies the loss of regularity in terms of $\rho$ and it holds in the larger class $S_{\rho,\delta}^m(\RR_\Theta)$.}
\end{remark}

\section*{{\bf Appendix A.} {\bf Noncommutative tori}} 

Given any $n \times n$ anti-symmetric $\R$-matrix $\Theta$, the subalgebra of $\RR_\Theta$ generated by $w_j = u_j(1)$ is the rotation algebra $\A_{\Theta}$ ---also known as quantum or noncommutative torus--- and we have $$w_j w_k = \exp(2 \pi i \Theta_{jk}) w_k w_j.$$ $\A_{\Theta}$ can also be described as the $\Z^n$-periodic subalgebra $$\A_\Theta \, = \, \big\langle \lambda_\Theta(\k) : \k \in \Z^n \big\rangle'' \, = \, \Big\{ \varphi \in \RR_\Theta \, \big| \, \sigma_\Theta^\k(\varphi) = \varphi \mbox{ for all } \k \in \Z^n \Big\}.$$

The extension of  our results for pseudodifferential operators to noncommutative tori $\A_\Theta$ follows by a combination of well-known transference arguments, which we recall now. Given a symbol $a: \Z^n \to \A_\Theta$ we shall say that
\begin{itemize}
\item $a \in S_{\rho,\delta}^m(\A_\Theta)$ when $$\hskip10pt \big| \partial_\Theta^\beta \partial_\k^\alpha a(\k) \big| \, \le \, C_{\alpha,\beta} \, \langle \k \rangle^{m - \rho |\alpha| + \delta |\beta|}.$$

\item $a \in \Sigma_{\rho,\delta}^m(\A_\Theta)$ when $$\hskip22pt \big| \partial_\Theta^\beta \partial_{\Theta,\k}^{\alpha_1} \partial_\k^{\alpha_2} a(\k) \big| \, \le \, C_{\alpha_1,\alpha_2,\beta} \, \langle \k \rangle^{m - \rho |\alpha_1 + \alpha_2| + \delta |\beta|}.$$
\end{itemize}
In the above definitions, $\partial_\Theta$ remains the same differential operator as in $\RR_\Theta$ whereas $\partial_\k$ is the difference operator $(\partial_\k^j a)(\k) = a(\k + e_j) - a(\k)$. The mixed derivatives $\partial_{\Theta,\k}$ are again $\Theta$-deformations of $\partial_\k$ by $\partial_\Theta$'s 
$$\partial_{\Theta, \k}^j a (\k) \, = \, \partial_\k^j a(\k) + 2 \pi i \big[ x_{\Theta,j}, a(\k) \big] \, = \, \partial_\k^j a(\k) + \frac{1}{2\pi i} \sum_{\ell=1}^n \Theta_{j\ell} \hskip1pt \partial_\Theta^\ell a(\k).$$
The associated pseudodifferential operator is $$\Psi_a(\varphi) \, = \, \sum_{\hskip1pt \k \in \Z^n} a(\k) \widehat{\varphi}(\k) \lambda_\Theta(\k) \quad \mbox{for} \quad \varphi = \sum_{\hskip1pt \k \in \Z^n} \widehat{\varphi}(\k) \lambda_\Theta(\k).$$ We say that $\widetilde{a}: \R^n \to \A_\Theta$ is a Euclidean lifting of $a$ when its restriction to $\Z^n$ coincides with the original symbol $a: \Z^n \to \A_\Theta$. Recall that we impose the lifting to take values in the periodic subalgebra $\A_\Theta$, not just in $\RR_\Theta$. The extension/restriction theorem below provides a useful characterization of the quantum H\"ormander classes in $\A_\Theta$ defined above since it relates them with their siblings in $\RR_\Theta$.  

\begin{Atheorem} \label{ThmExtension}
Assume $\rho > 0 \hskip-1pt :$
\begin{itemize}
\item[\emph{i)}] $a \in S_{\rho,\delta}^m \hskip1pt (\A_\Theta)$ iff it admits a lifting $\widetilde{a} \in S_{\rho,\delta}^m \hskip1pt (\RR_\Theta)$.

\item[\emph{ii)}] $a \in \Sigma_{\rho,\delta}^m(\A_\Theta)$ iff it admits a lifting $\widetilde{a} \in \Sigma_{\rho,\delta}^m(\RR_\Theta)$.
\end{itemize}
In fact, the lifting $\widetilde{a}: \R^n \to \A_\Theta$ has the form $$\widetilde{a}(\xi) \, = \, \sum_{\hskip1pt \k \in \Z^n}^{\null} \phi(\xi-\k) a(\k)$$ for certain Schwartz function $\phi:\R^n \to \R$ satisfying $\phi(\k) = \delta_{\k,0}$ for $\k \in \Z^n$. 
\end{Atheorem}

The proof follows verbatim \cite[Theorem 4.5.3]{RT} since the argument only affects the classical variables $\k \in \Z^n$ and $\xi \in \R^n$. In particular, the exact same argument applies when we take values in $\A_\Theta$. In fact, the same extension procedure applies when the H\"ormander condition is only required for finitely many derivatives in the line of Remark \ref{RemJustn+2}, see \cite[Corollary 4.5.7]{RT}. The equality of the associated pseudodifferential operators is also proved in \cite[Theorem 4.6.12 and Corollary 4.6.13]{RT}. Namely, the class of pseudodifferential operators associated to $S_{\rho,\delta}^m(\A_\Theta)$ or $\Sigma_{\rho,\delta}^m(\A_\Theta)$ can be identified with the corresponding H\"ormander classes in $\RR_\Theta$ for periodic symbols |that is, taking values in $\A_\Theta$| when acting on periodic elements $\varphi = \sigma_\Theta^\k(\varphi)$ for $\k \in \Z^n$. Finally, it is also worth mentioning that the extension above also respects ellipticity, as shown in \cite[Theorem 4.9.15]{RT}. 

\begin{Atheorem} \label{ThmBTorus}
Let $a: \Z^n \to \A_\Theta$ and $1 < p < \infty$$\hskip1pt :$
\begin{itemize}
\item[\emph{i)}] If $a \in S_{\rho,\rho}^0(\A_\Theta)$ with $0 \le \rho < 1$, $\Psi_a: L_2(\A_\Theta) \to L_2(\A_\Theta)$.

\item[\emph{ii)}] If $a \in S_{1,1}^0 \hskip1pt (\A_\Theta) \cap S_{1,1}^0 \hskip1pt (\A_\Theta)^*$, then $\Psi_a: L_2(\A_\Theta) \to L_2(\A_\Theta)$. 

\item[\emph{iii)}] If $a \in \Sigma_{1,1}^0(\A_\Theta) \cap \Sigma_{1,1}^0(\A_\Theta)^*$, then $\Psi_a: L_p(\A_\Theta) \to L_p(\A_\Theta)$. 
\end{itemize}
\end{Atheorem}

\demA Let $b_j = \mathrm{B}_{2^{-j}}(0)$ and $$h_j = \frac{1}{\sqrt{|b_j|}} \lambda_\Theta(1_{b_j}) \quad \mbox{and} \quad \Lambda_j: \lambda_\Theta(\k) \mapsto \lambda_\Theta(\k) h_j.$$ Observe that $\Lambda_j: L_2(\A_\Theta) \to L_2(\RR_\Theta)$ is an isometry for all $j \ge 1$. Indeed
\begin{eqnarray*}
\big\| \Lambda_j(\varphi) \big\|_{L_2(\RR_\Theta)}^2 
& = & \frac{1}{|b_j|} \Big\| \sum_{\hskip1pt \k \in \Z^n} \widehat{\varphi}(\k) \lambda_\Theta(\k) \lambda_\Theta(1_{b_j}) \Big\|_{L_2(\RR_\Theta)}^2 \\
& = & \frac{1}{|b_j|} \Big\| \int_{\R^n} \sum_{\hskip1pt \k \in \Z^n} \widehat{\varphi}(\k) 1_{b_j}(\xi - \k) e^{2\pi i \langle \k, \Theta_\downarrow \xi - \k\rangle} \lambda_\Theta(\xi) \, d\xi \Big\|_{L_2(\RR_\Theta)}^2.
\end{eqnarray*}
By Plancherel theorem and using that $b_j + \k$ are pairwise disjoint, we get
\begin{eqnarray*}
\big\| \Lambda_j(\varphi) \big\|_{L_2(\RR_\Theta)}^2 
\!\!\! & = & \!\!\! \frac{1}{|b_j|} \int_{\R^n} \Big| \sum_{\hskip1pt \k \in \Z^n} \widehat{\varphi}(\k) 1_{b_j}(\xi - \k) e^{2\pi i \langle \k, \Theta_\downarrow \xi - \k\rangle} \Big|^2 \, d\xi \\
\!\!\! & = & \!\!\! \frac{1}{|b_j|} \int_{\R^n} \sum_{\hskip1pt \k \in \Z^n} \big| \widehat{\varphi}(\k) \big|^2 1_{b_j}(\xi - \k) \, d\xi = \sum_{\hskip1pt \k \in \Z^n} \big| \widehat{\varphi}(\k) \big|^2 = \|\varphi\|_{L_2(\A_\Theta)}^2.
\end{eqnarray*}
Then, the assertion follows from the following claim 
$$\lim_{j \to \infty} \Big\| \Lambda_j \big( \Psi_a(\varphi) \big) - \Psi_{\widetilde{a}} \big( \Lambda_j(\varphi) \big) \Big\|_{L_2(\RR_\Theta)} \, = \, 0$$ for any trigonometric polynomial $\varphi$. In other words, for finite linear combinations of the $\lambda_\Theta(\k)$'s. Indeed, assume the limit above vanishes, then $\Psi_a$ is $L_2$-bounded since trigonometric polynomials are dense and 
\begin{eqnarray*}
\big\| \Psi_a(\varphi) \big\|_{L_2(\A_\Theta)} & = & \lim_{j \to \infty} \big\| \Lambda_j \big( \Psi_a(\varphi) \big) \big\|_{L_2(\RR_\Theta)} \\
& = & \lim_{j \to \infty} \big\| \Psi_{\widetilde{a}} \big( \Lambda_j (\varphi) \big) \big\|_{L_2(\RR_\Theta)} \\
& \le & \lim_{j \to \infty} \big\| \Lambda_j (\varphi) \big\|_{L_2(\RR_\Theta)} \ = \ \|\varphi\|_{L_2(\A_\Theta)}.
\end{eqnarray*}
The inequality above follows by application of Theorem \ref{ThmExtension} in conjunction with Theorems \ref{PSD.CalVaiThm00}, \ref{PSD.CalVaiThm} and \ref{PSD.boudreaud}. Let us then justify our claim above. It clearly suffices to prove it with $\varphi = \lambda_\Theta(\k)$ for any $\k \in \Z^n$. Given an arbitrary $\varepsilon > 0$, we shall prove that the quantity $\| \Lambda_j ( \Psi_a(\lambda_\Theta(\k))) - \Psi_{\widetilde{a}} ( \Lambda_j(\lambda_\Theta(\k)) ) \|_2 < C \varepsilon$ for some absolute constant $C$ independent of $(\k,\varepsilon)$ and $j$ large enough. Let $\phi: \R^n \to \R$ be the function used in Theorem \ref{ThmExtension} for the construction of the lifting. Since $\phi$ is a Schwartz function and $\phi(\k) = \delta_{\k,0}$ for $\k \in \Z^n$, there must exists a $\delta > 0$ satisfying 
$$|\xi - \k| < \delta \ \Rightarrow \ \max \Big\{ \big| \phi(\xi - \k) - 1 \big|, \sup_{\mathrm{j} \neq \k} |\phi(\xi - \mathrm{j})| \Big\} < \varepsilon R_\varepsilon^{-n}$$ where $R_\varepsilon$ is large enough to satisfy $$\sum_{|\k| > R_\varepsilon} \frac{1}{|\k|^{n+1}} < \varepsilon.$$ Next, consider the Fourier multiplier $M_{b_{\k\delta}}(\varphi) = \int_{\R^n} 1_{b_{\k\delta}}(\xi) \widehat{\varphi}(\xi) \lambda_\Theta(\xi) d\xi$ where we write $b_{\k\delta}$ for $\mathrm{B}_\delta(\k)$. Then, we decompose the $L_2$-norm into three terms as follows
\begin{eqnarray*}
\lefteqn{\hskip-20pt \big\| \Lambda_j \big( \Psi_a(\lambda_\Theta(\k)) \big) - \Psi_{\widetilde{a}} ( \Lambda_j(\lambda_\Theta(\k)) ) \big\|_2} \\ 
& = & \big\| a(\k) \lambda_\Theta(\k) h_j - \Psi_{\widetilde{a}} \big( \lambda_\Theta(\k) h_j \big) \big\|_2 \\
& \le & \big\| a(\k) \big( \lambda_\Theta(\k) h_j - M_{b_{\k\delta}}(\lambda_\Theta(\k) h_j) \big) \big\|_2 \\
& + & \big\| a(\k) M_{b_{\k\delta}}(\lambda_\Theta(\k) h_j) - \Psi_{\widetilde{a}} \big( M_{b_{\k\delta}}(\lambda_\Theta(\k) h_j)  \big) \big\|_2 \\
& + & \big\| \Psi_{\widetilde{a}} \big( M_{b_{\k\delta}}(\lambda_\Theta(\k) h_j)  \big) - \Psi_{\widetilde{a}} \big( \lambda_\Theta(\k) h_j \big) \big\|_2 \ = \ \mathrm{A} + \mathrm{B} + \mathrm{C}.
\end{eqnarray*}
We recall one more time from Theorem \ref{ThmExtension} and Theorem B in the Introduction that $\Psi_{\widetilde{a}}: L_2(\RR_\Theta)) \to L_2(\RR_\Theta)$ is a bounded map. Moreover, we also know that $a \in \ell_\infty(\Z^n; \A_\Theta)$ since it has degree $0$. In particular
$$\mathrm{A} + \mathrm{C} \le \Big( \sup_{\hskip1pt \k \in \Z^n} \|a(\k)\|_{\A_\Theta} + \big\| \Psi_{\widetilde{a}} \big\|_{L_2(\RR_\Theta) \to L_2(\RR_\Theta)}  \Big) \big\| \lambda_\Theta(\k) h_j -  M_{b_{\k\delta}}(\lambda_\Theta(\k) h_j) \big) \big\|_2.$$ The $L_2$-norm above can be estimated with Plancherel theorem
\begin{eqnarray*}
\lefteqn{\big\| \lambda_\Theta(\k) h_j - M_{b_{\k\delta}}(\lambda_\Theta(\k) h_j) \big) \big\|_2} \\
& = & \frac{1}{\sqrt{|b_j|}} \Big\| \int_{\R^n} (1 - 1_{b_{\k\delta}}(\xi)) 1_{b_j}(\xi - \k) e^{2\pi i \langle \k, \Theta_\downarrow \xi - \k\rangle} \lambda_\Theta(\xi) \, d\xi \Big\|_2 \\
& = & \Big( \frac{1}{|b_j|} \int_{\k + b_j} \big| 1 - 1_{b_{\k\delta}}(\xi) \big|^2 \, d\xi \Big)^{\frac12} \longrightarrow \ \big| 1 - b_{\k\delta}(\k) \big| \ = \ 0 
\end{eqnarray*}
as $j \to \infty$. Therefore, it remains to estimate the term $\mathrm{B}$. Letting 
\begin{eqnarray*}
\mathbf{a}_\k(\xi) & = & \big( \widetilde{a}(\xi) - a(\k) \big) 1_{b_{\k\delta}(\xi)} \\
[12pt] & = & a(\mathrm{k}) \big( \phi(\xi-\mathrm{k}) - 1 \big) 1_{b_{\k\delta}(\xi)} \\
[10pt] & + & \sum_{\begin{subarray}{c} \mathrm{j} \neq \k \\ |\mathrm{j} - \k| \le R_\varepsilon \end{subarray}}^{\null} a(\mathrm{j}) \phi(\xi-\mathrm{j}) 1_{b_{\k\delta}(\xi)} \\
& + & \sum_{|\mathrm{j} - \k| > R_\varepsilon} a(\mathrm{j}) \phi(\xi-\mathrm{j}) 1_{b_{\k\delta}(\xi)} \ = \ \mathbf{a}_{1\k}(\xi) + \mathbf{a}_{2\k}(\xi) + \mathbf{a}_{3\k}(\xi),
\end{eqnarray*}
we clearly have 
$$\mathrm{B} \, = \, \big\| \Psi_{\mathbf{a}_\k}(\lambda_\Theta(\k) h_j) \big\|_2 \, \le \, \sum_{j=1}^3 \big\| \Psi_{\mathbf{a}_{j\k}}: L_2(\RR_\Theta) \to L_2(\RR_\Theta) \big\|$$
since $\lambda_\Theta(\k) h_j$ is a unit vector in $L_2(\RR_\Theta)$. This gives 
\begin{eqnarray*}
\mathrm{B} & \le & \Big( \sup_{\mathrm{j} \in \Z^n} \|a(\mathrm{j})\|_{\A_\Theta} \Big) \Big( \sup_{|\xi - \k| < \delta} \big| \phi(\xi- \k) - 1 \big| \Big) \\
[8pt] & + & \Big( \sup_{\mathrm{j} \in \Z^n} \|a(\mathrm{j})\|_{\A_\Theta} \Big) \Big( \sum_{\begin{subarray}{c} \mathrm{j} \neq \k \\ |\mathrm{j} - \k| \le R_\varepsilon \end{subarray}}^{\null}  \sup_{|\xi - \k| < \delta} |\phi(\xi- \mathrm{j})| \Big) \\
& + & \Big( \sup_{\mathrm{j} \in \Z^n} \|a(\mathrm{j})\|_{\A_\Theta} \Big) \big\| |\xi|^{n+1} \phi(\xi) \big\| \Big( \sum_{|\mathrm{j} - \k| > R_\varepsilon}^{\null}  \frac{1}{|\mathrm{j} -\k|^{n+1}} \Big) \ \lesssim \ 3 \varepsilon.
\end{eqnarray*}
Then, letting $\varepsilon \to 0^+$ this completes the proof of the claim. \fin

\begin{Aremark}
\emph{The above argument was inspired by the proof of \cite[Theorem A]{CPPR}.} 
\end{Aremark}

\demAA The next ingredient we need is the natural BMO space in $\A_\Theta$. Define $\mathrm{BMO}_c(\A_\Theta)$ as the column BMO space associated to the transferred heat semigroup $\varphi \mapsto \sum_\k \widehat{\varphi}(\k) \exp(-t|\k|^2) \lambda_\Theta(\k)$. As in Section \ref{MetricandBMO} it can be regarded as the weak-$*$ closure of $\sigma_\Theta(\A_\Theta)$ with respect to the pair $(\mathrm{H}^c_1(\mathcal{Q}_\Theta), \mathrm{BMO}_c(\mathcal{Q}_\Theta))$. In other words, we find $$\|a\|_{\mathrm{BMO}_c(\A_\Theta)} \, \sim \, \sup_{Q \in \mathcal{Q}} \Big\| \Big( \mean_Q \big| \sigma_\Theta(a) - \sigma_\Theta(a)_Q \big|^2 \, d\mu \Big)^\frac12 \Big\|_{\A_\Theta},$$ where $\mathcal{Q}$ is the set of all Euclidean cubes in $\R^n$ with sides parallel to the axes, $\mu$ stands for the Lebesgue measure and $\sigma_\Theta(a)_Q$ is the average of $\sigma_\Theta(a)$ over the cube $Q$. Up to absolute constants, it is not difficult to recover an equivalent norm when restricting to cubes $Q$ of side length $\ell(Q) \in (0,1) \cup \N$. Moreover, since $a$ is spanned by $\lambda_\Theta(\k)$ for $\k \in \Z^n$, it is clear that $\sigma_\Theta(a)$ is $\Z^n$-periodic. In particular, the quantity above for $\ell(Q) \in \N$ coincides with the same quantity for $Q = [0,1] \times \ldots \times [0,1]$, so that we may assume in addition $\ell(Q) \le 1$. We have proved 
$$\|a\|_{\mathrm{BMO}_c(\A_\Theta)} \, \sim \, \sup_{Q \in \mathbb{T}^n} \Big\| \Big( \mean_Q \big| \sigma_\Theta(a) - \sigma_\Theta(a)_Q \big|^2 \, d\mu \Big)^\frac12 \Big\|_{\A_\Theta}.$$ In other words, $\mathrm{BMO}_c(\A_\Theta)$ embeds into $\mathrm{BMO}_c(\mathbb{T}^n; \A_\Theta)$ using Mei's terminology \cite{Mei}. The interpolation behavior and other natural properties which we explore for $\mathrm{BMO}(\RR_\Theta)$ in Appendix B are well-known in this case \cite{JM}, due to the finiteness of $\A_\Theta$. Note that, according to our definition of $\mathrm{BMO}(\A_\Theta)$, the natural inclusion map $\A_\Theta \to \RR_\Theta$ extends to an embedding $\mathrm{BMO}(\A_\Theta) \to \mathrm{BMO}(\RR_\Theta)$. In other words, $\mathrm{BMO}(\A_\Theta)$ is the subspace of periodic elements in $\mathrm{BMO}(\RR_\Theta)$. Now, recalling that $\Psi_{\widetilde{a}}$ sends periodic elements into periodic elements, this makes the following a commutative diagram       
\begin{equation*}
\xymatrix{
\ar[d]^{\hskip-17pt \Psi_a} L_\infty(\A_\Theta) \ar[rr]^{id} & & L_\infty(\RR_\Theta) \ar[d]^{\Psi_{\widetilde{a}}} \\
\mathrm{BMO}(\A_\Theta) \ar[rr]^{id} & & \mathrm{BMO}(\RR_\Theta) 
}
\end{equation*}
The assertion follows from it and Theorem \ref{ThmBTorus} ii) by interpolation and duality. \fin  

\begin{Aremark}
\emph{Of course, our $L_p$-inequalities also hold in the category of operator spaces and admit the endpoint estimates $\mathrm{H}_1 \to L_1$ and $L_\infty \to \mathrm{BMO}$, as in the quantum Euclidean setting. Besides, the natural analogues of Remarks \ref{RemS1delta} and \ref{RemJustn+2} as well as Corollary \ref{CorollaryRegularizingTerm} concerning $L_p$-estimates still apply. On the other hand, the Sobolev $p$-estimates in Corollary \ref{PSD.lpBound2} and the $L_p$-regularity for elliptic PDEs require in addition analogues of the product stability of H\"ormander classes in Section \ref{CompAdj}, which seems to be straightforward but we shall not generalize it here.}
\end{Aremark}

\section*{{\bf Appendix B.} {\bf BMO space theory in $\RR_\Theta$}}
\label{AppB}

\renewcommand{\theequation}{B.1}
\addtocounter{equation}{-1}

The theory of $\BMO$ spaces was developed as a natural endpoint class for singular integral operators. In particular, the natural requirements for a reasonable BMO space are:  
\begin{enumerate}[label={\arabic*)}]
\item \label{BMO.cond2}
Interpolation endpoint for the $L_p$-scale. 
  
\item \label{BMO.cond3}
John-Nirenberg inequalities and $\Hardy_1-\BMO$ duality. 
  
\item \label{BMO.cond1}
$L_\infty \to \BMO$ boundedness for Calder\'on-Zygmund operators. 
\end{enumerate}
$\BMO$ spaces over von Neumann algebras were introduced by Pisier and Xu in \cite{PX} and have been investigated since then. The theory  when averages over balls or martingale filtrations are replaced by the action of a Markovian semigroup has been addressed for finite von Neumann algebras in \cite{JM}. Interpolation requires a different approach over $\RR_\Theta$ |less intricate than the general semifinite case| which we present here. Duality was developed by Mei \cite{Mei, MeiTent} and endpoint estimates for imaginary powers $A^{is}$ of infinitesimal generators, noncommutative Riesz transforms or more general Fourier multipliers have been studied in \cite{CXY,JunMeiRiesz,JM,JMP,Ri}. In the setting of $\A_\Theta$ and $\RR_\Theta$, Theorems A and B include many more singular integrals.

\subsubsection*{\emph{B.1.} Operator space structures on $\BMO$ and $\Hardy_1$}
\label{AppB.BMOoss}

Let us recall the definitions of several natural operator space structures |o.s.s. in short| for $\BMO(\R^n)$ and its predual. We define the column operator space structure by the family of matrix norms on $f = [f_{i j}] \in M_m[\BMO_c(\R^n)]$ given by
$$\|f\|_{M_m[\BMO_c(\R^n)]} = \sup_{Q \in \Q} \Big\| \mean_Q (f - f_Q)^\ast \, (f - f_Q) \, d\mu \Big\|_{M_m}^{\frac{1}{2}},$$
where $f_Q$ is the average of $f$ over $Q$ and $\Q$ stands for the set of all the Euclidean balls. We will denote the resulting operator space by $\BMO_c(\R^n)$. Similarly, we can define the row o.s.s. by $\|f\|_{M_m[\BMO_r(\R^n)]} = \|f^\ast \|_{M_m[\BMO_c(\R^n)]}$. We shall also denote by $\BMO(\R^n)$ |sometimes $\BMO_{r \wedge c}(\R^n)$ for convenience| the operator space structure  
$$\|f\|_{M_m[\BMO(\R^n)]} = \max \Big\{ \|f\|_{M_m[\BMO_c(\R^n)]}, \|f\|_{M_m[\BMO_r(\R^n)]} \Big\}.$$
These are dual operator spaces, with preduals $\Hardy_1^\dag(\R^n)^\ast = \BMO_\dag(\R^n)$ given by
\begin{eqnarray*}
  \|f\|_{\Hardy^c_1(\R^n)} & = & \Big\| \Big(\int_{\R_+} \big| s (\nabla + \partial_s^2) P_s f \big|^2 \frac{ds}{s} \Big)^\frac12 \Big\|_{L_1(\R^n)}, \\
  \|f\|_{\Hardy^r_1(\R^n)} & = & \Big\| \Big(\int_{\R_+} \big| s (\nabla + \partial_s^2) P_s f^\ast \big|^2 \frac{ds}{s} \Big)^\frac12 \Big\|_{L_1(\R^n)},
\end{eqnarray*}
where $P_s$ is the Poisson semigroup. The quantities above are just pseudonorms. A natural way of turning them into norms is working with $0$-integral functions, in a way dual to the quotient of constants taken in the definition of $\BMO$. Comparable norms can be defined by removing the $\partial_s^2$ inside the square function and by using the semigroup analogue of Lusin area integral, given by
$$\|f\|_{\Hardy^c_1(\R^n)} \, \sim \, \Big\| \Big( \int_{\Gamma_x} \big| (\nabla + \partial_s^2) P_s f(y) \big|^2 \, ds dy \Big)^\frac12 \Big\|_{L_1(\R^n)},$$
where $\Gamma_x = \{ (y,s) \in \R^n \times \R_+ : |y - x| \leq s \}$ is the cone centered at $x$. The row case can be expressed analogously. The o.s.s.
of Hardy spaces can be easily described by taking matrix-valued functions $f = [f_{ij}]$ in the expression above and taking norms in $S_1^m \widehat\otimes L_1(\R^n) = L_1(\R^n; S_1^m)$. That will give a family of matrix norms which describes the operator space structure. Indeed, using \cite[Lemma 1.7]{P2} and the well-known relation 
$$M_m[\Hardy^\dagger_1(\R^n)] = \mathcal{CB}(S_1^m, \Hardy^\dagger_1(\R^n))$$
see e.g. \cite[Theorem 4.1]{P3}, we can easily express the norm of $M_m[\Hardy^\dagger_1(\R^n)]$ in
terms of the known norms. The operator space predual $\Hardy_1(\R^n)$ of
$\BMO(\R^n)$ is given by the sum $\Hardy^c_1(\R^n) + \Hardy^r_1(\R^n)$,
whose norm is
$$\|f\|_{S_1^m \widehat\otimes \Hardy_1(\R^n)} = \inf \Big\{ \|g\|_{S_1^m \widehat\otimes \Hardy_1^r(\R^n)} + \|h\|_{ S_1^m \widehat\otimes \Hardy_1^c(\R^n)}: f = g + h \Big\}.$$

Let us note that, by computations in Section \ref{MetricandBMO} we have that $\sigma_\Theta$ gives an isomorphic embedding $\BMO_\dagger(\RR_\Theta) \to
\BMO_\dagger(\R^n) \bar\otimes \RR_\Theta$. In particular, since $\RR_\Theta$ is hyperfinite, we may equip $\BMO_\dagger(\RR_\Theta)$ with an o.s.s. naturally inherited from $\BMO_\dagger(\R^n)$. Following Mei \cite{MeiTent, MeiH1}, the definition of $\Hardy_1^\dagger(\RR_\Theta)$ will be given by completion on the $0$-trace functions with respect to 
\begin{eqnarray*}
  \| \varphi \|_{\Hardy_1^c(\RR_\Theta)} & = & \Big\| \Big( \int_{\R_+} S_{\Theta,t} \big| \nabla_\Theta S_{\Theta,t} \varphi \big|^2 \, dt \Big)^\frac12 \Big\|_{L_1(\RR_\Theta)},\\
  \| \varphi \|_{\Hardy_1^r(\RR_\Theta)} & = & \Big\| \Big( \int_{\R_+} S_{\Theta,t} \big| \nabla_\Theta S_{\Theta,t} \varphi^\ast \big|^2 \, dt \Big)^\frac12 \Big\|_{L_1(\RR_\Theta)}.
\end{eqnarray*}
The operator space structures of such spaces are defined in the same way as the operator space structures of the classical ones, which could also have been defined with this square function instead of the given one yielding an equivalent norm. Note also that when $\Theta = 0$, the semigroup $S_{\Theta,t}$ behaves (intuitively) like an average over balls of radius $\sqrt{t}$ and a calculation gives that the quantities above are comparable to the Lusin integral and therefore recover the classical $\Hardy_1(\R^n)$. We will write $\Hardy_1(\RR_\Theta)$ or $\Hardy_1^{r + c}(\RR_\Theta)$ for the sum $$\Hardy_1(\RR_\Theta) = \Hardy_1^{r}(\RR_\Theta) + \Hardy_1^{ c}(\RR_\Theta).$$

\subsubsection*{\emph{B.2.} The $\Hardy_1$-$\BMO$ duality}
\label{AppB.BMOPredual}

The von Neumann algebra analogue of the celebrated $\Hardy_1-\BMO$ duality \cite{FS} has been carefully studied in our semigroup setting by Tao Mei. By \cite[Theorem 0.2]{MeiH1}, the duality between $\Hardy_1(\RR_\Theta)$ and $\BMO(\RR_\Theta)$ can be deduced after verifying that the associated heat semigroup $(S_{\Theta,t})_{t \geq 0}$ satisfies the  following conditions:
\begin{enumerate}[label={\rm \roman*)}, ref={\roman*)}]
    \item \label{BMO.MeiDuals0} 
    Bakry's $\Gamma_2 \ge 0$ condition.

\vskip10pt

    \item \label{BMO.MeiDuals1}
    For all $\varepsilon, t > 0$ and $\varphi \in L_1(\RR_\Theta)$
    $$\big\| (S_{\Theta,(1 + \epsilon)t} - S_{\Theta,t}) \varphi \big\|_{L_1(\RR_\Theta)} \, \lesssim \, \varepsilon^r \|\varphi\|_{L_1(\RR_\Theta)}.$$
    
    \item \label{BMO.MeiDuals2}
    For every $t > 0$ and $\varphi \in L_1(\RR_\Theta)$
    $$\sup_{t > 0} \Big\| \Big(\mean_0^{8 t} S_{\Theta,s}(|S_{\Theta,t}(\varphi)|^2) \, ds \Big)^\frac12 \Big\|_{L_1(\RR_\Theta)} \, \lesssim \, \|\varphi\|_{L_1(\RR_\Theta)}.$$
\end{enumerate}
Verifying such identities is relatively easy for the heat semigroup $(S_{\Theta,t})_{t \geq 0}$ after noting that it can be presented as an integrable convolution with respect to the $z$-variable $\sigma_\Theta^z(\varphi)$ and using bounds in $L_{1/2}(\RR_\Theta)$. In particular, we obtain the expected duality theorem. 

\begin{Btheorem}
\label{BMO.duality}
We have  
$$\Hardy^\dagger_1(\RR_\Theta)^\ast = \BMO_\dagger(\RR_\Theta)$$
in the category of operator spaces for $\dagger \in \{r ,c\}$. Also $\Hardy_1(\RR_\Theta)^* = \BMO(\RR_\Theta)$.
\end{Btheorem}


\subsubsection*{\emph{B.3.} Complex interpolation}
\label{AppB.BMOInter}

We are now interested in proving the generalization of the classical interpolation identities between $L_p$, $\BMO$ and $\Hardy_1$. According to Wolff's interpolation theorem \cite{Wolff}, this can be easily reduced to justifying the complete isomorphism $[L_2(\RR_\Theta),\BMO(\RR_\Theta)]_{\theta} = L_p(\RR_\Theta)$ for $p = \frac{2}{1 - \theta}$ which in turn will be reduced, via suitable complemented subspaces, to the same result in $\R^n$ but with operator values in certain hyperfinite von Neumann algebra.

Let us recall a few standard definitions from interpolation theory. Given $\mathrm{X}_0, \mathrm{X}_1$ Banach spaces, assume that they embed inside a topological vector space with dense intersection, so that we can define $\mathrm{X}_0 \cap \mathrm{X}_1$ and $\mathrm{X}_0 + \mathrm{X}_1$ with their natural norms. Let us write $\mathcal{F}(\mathrm{X}_0, \mathrm{X}_1)$ for the space of $(\mathrm{X}_0 + \mathrm{X}_1)$-valued holomorphic functions in the strip $0 < \Re(z) < 1$ which admit a continuous extension to the boundary, with $\mathrm{X}_j$-values at $\partial_j$ for $j=1,2$. Such space is a Banach space with respect to the norm given by
$$\|f\|_{\mathcal{F}(\mathrm{X}_0, \mathrm{X}_1)} = \max \Big\{ \sup_{s \in \R} \| f(i s) \|_{\mathrm{X_0}}, \sup_{s \in \R} \| f(1 + i s) \|_{\mathrm{X_1}} \Big\}.$$
The interpolated space with parameter $0 < \theta < 1$ is
$$\big[ X_0, X_1 \big]_\theta = \mathcal{F}(\mathrm{X}_0, \mathrm{X}_1) / \mathfrak{N}_\theta,$$
where $\mathfrak{N}_\theta$ is the subspace of functions with $f(\theta) = 0$. We can also define a larger interpolation functor $[\mathrm{X}_0, \mathrm{X}_1]^\theta$ that contains $[\mathrm{X}_0, \mathrm{X}_1]_\theta$ isometrically by changing $\mathcal{F}(\mathrm{X}_0, \mathrm{X}_1)$ by a la larger space $\mathcal{F}_\ast(\mathrm{X}_0, \mathrm{X}_1)$ of holomorphic functions in which $f{|}_{\partial_j}$ is a more general $\mathrm{X}_j$-valued distribution. These interpolation functors satisfy that $[\mathrm{X}_0, \mathrm{X}_1]_\theta^\ast = [\mathrm{X}^\ast_0, \mathrm{X}^\ast_1]^\theta$ and both coincide if any of the spaces involved $\mathrm{X}_0, \mathrm{X}_1$ is reflexive \cite[Corollary 4.5.2]{BeLo} and \cite[Theorem 2.7.4]{P3}. If $\mathrm{X}_j$ are operator spaces, the o.s.s. of $[\mathrm{X}_0,\mathrm{X}_1]_\theta$ is given by the identification $$M_m\big( [\mathrm{X}_0,\mathrm{X}_1]_\theta \big) = \big[ M_m(\mathrm{X}_0),M_m(\mathrm{X}_1) \big]_\theta.$$ 

We first need an auxiliary result concerning complex interpolation of tensor products against hyperfinite von Neumann algebras. This result is a consequence of the interpolation identity $[\M_* \widehat{\otimes} \mathrm{X}_0, \M_* \widehat{\otimes} \mathrm{X}_1]_\theta = \M_* \widehat{\otimes} \mathrm{X}_\theta$ which can be found in \cite[page 40]{P2}. We prove it for completeness.

\begin{Blemma} \label{BMO.InterTensor}
We have
$$\big[ \M \bar\otimes \mathrm{X}_0, \M \bar\otimes \mathrm{X}_1 \big]^\theta = \M \bar\otimes [\mathrm{X}_0, \mathrm{X}_1]^\theta$$
for any hyperfinite algebra $\M$ and any pair of dual operator spaces $\mathrm{X}_0$, $\mathrm{X}_1$.
\end{Blemma}

\dem According to \cite{ERHopf} the spaces involved are dual operator spaces. Indeed, von Neumann algebra preduals have the $\mathrm{OAP}$ so $\M \bar\otimes \mathrm{X}^\ast = (\M_\ast \widehat\otimes \mathrm{X})^\ast$. Now, since hyperfiniteness and semidiscreteness are equivalent $id:\M \to \M$ is approximable in the pointwise weak-$\ast$ topology by a net $i_\alpha = \psi_\alpha \phi_\alpha$ where $\phi_\alpha:\M \to M_{m_\alpha}(\C)$ and $\psi_\alpha: M_{m_\alpha}(\C) \to \M$ are ucp. We have 
\begin{equation*}
  \xymatrix{
    \M \bar\otimes \mathrm{X}^\theta \ar[dr]^-{\psi_\alpha \otimes id} \ar[rrr]^-{i_\alpha \otimes id} & & & \big[\M \bar\otimes \mathrm{X}_0, \M \bar\otimes \mathrm{X}_1 \big]^\theta\\
    & M_{m_\alpha}(\mathrm{X}^\theta) \ar@{-}[r]^-{\simeq} & \big[ M_{m_\alpha}(\mathrm{X}_0), M_{m_\alpha}(\mathrm{X}_1) \big]^\theta \ar[ur]^{\phi_\alpha \otimes id} & \\
  }
\end{equation*}
a commutative diagram for $\mathrm{X}^\theta = [\mathrm{X}_0, \mathrm{X}_1]^\theta$. The maps $i_\alpha$ approximate the identity and taking a weak-$\ast$ accumulation point in $\CB(\M \bar\otimes \mathrm{X}\theta, [\M \bar\otimes \mathrm{X}_0, \M \bar\otimes \mathrm{X}_1]^\theta)$, which is a dual space since $\CB(\mathrm{X},\mathrm{Y}^\ast) = (\mathrm{X} \widehat\otimes \mathrm{Y})^\ast$, we obtain a complete isomorphism.
\fin

A key point in our interpolation argument will be to show that the co-action $\sigma_\Theta: \RR_\Theta \to L_\infty(\R^n) \bar\otimes \RR_\Theta$ also carries other $\RR_\Theta$-spaces |$L_p$ and $\BMO$| into their $\RR_\Theta$-valued Euclidean counterparts.

\begin{Bproposition} \label{BMO.IsoLp}
We have complete contractions$\hskip1pt :$ 
\begin{enumerate}[label={\rm \roman*)}]
    \item \label{BMO.Lp1}
    $\sigma_\Theta: L_2^\dagger(\RR_\Theta) \longrightarrow L_2^{\dagger}(\R^n) \bar\otimes \RR_\Theta$ for $\dagger \in \{r,c\}$,

\vskip2pt

    \item \label{BMO.Lp2}
    $\sigma_\Theta: L_p(\RR_\Theta) \longrightarrow L_p(\R^n) \bar\otimes \RR_\Theta$ 
    for any $2 \le p \le \infty$,
    
\vskip2pt

    \item \label{BMO.Lp3}
    $\sigma_\Theta: \BMO_{\dagger}(\RR_\Theta) \longrightarrow \BMO_{\dagger}(\R^n) \bar\otimes \RR_\Theta$ for $\dagger \in \{r, c, r \wedge c\}$.
\end{enumerate} 
\end{Bproposition}

\dem Let us recall the Fubini-type identity  
$$\1 \tau_\Theta(\varphi) = \int_{\R^n} \sigma_\Theta^z(\varphi) \, dz.$$
In particular, given $\varphi = [\varphi_{i j}] \in M_m[L_2^c(\RR_\Theta)]$ we obtain
\begin{eqnarray*}
  \| \varphi \|_{M_m[L_2^c(\RR_\Theta)]}^2
  \!\!\! & = & \!\!\! \big\| (id \otimes \tau_\Theta)\big( \varphi^\ast \varphi \big) \big\|_{M_m}  = \big\| (id \otimes \1 \tau_\Theta)\big( \varphi^\ast \varphi \big) \big\|_{M_m[\RR_\Theta]} \\
  \!\!\! & = & \!\!\! \Big\| \int_{\R^n} \sigma_\Theta^z(\varphi)^\ast \sigma_\Theta^z(\varphi) \, dz \Big\|_{M_m[\RR_\Theta]}
  = \big\| \sigma_\Theta(\varphi) \big\|_{M_m[L_2^c(\R^n) \bar\otimes \RR_\Theta]}^2.
\end{eqnarray*}
The same follows in the row case. In fact, $\sigma_\Theta$ is a complete isometry in case \ref{BMO.Lp1} and also in case  \ref{BMO.Lp3} by construction of $\BMO(\RR_\Theta)$. Assertion \ref{BMO.Lp2} follows by interpolation from Lemma \ref{BMO.InterTensor}. Indeed, since $L_2(\R^n) = [L_2^c(\R^n), L_2^r(\R^n)]_{1/2}$ in the category of operator spaces and all spaces involved are reflexive, we obtain from \ref{BMO.Lp1} that $\sigma_\Theta$ is a complete contraction from $L_2(\RR_\Theta)$ to $L_2(\R^n) \bar\otimes \mathcal{R}_\Theta$. The case $p > 2$ also follows by complex interpolation, using the reflexivity of $L_2$, since the contractivity of the other endpoint for $p = \infty$ was already justified in Corollary \ref{Planc+Corep}. \fin

\begin{Bremark}
\emph{It is interesting to know whether an analogue of Proposition \ref{BMO.IsoLp} holds for $p=1$. Note that $L_1(\R^n)$ is not a dual space, so that we can not use the weak-$*$ closed tensor product. Instead, we shall consider the mixed-norm space $L_\infty(\mathcal{R}_\Theta; L_1(\R^n))$ as introduced in \cite{JDoob,JP}
\[
L_\infty \big(\mathcal{R}_\Theta; L_1(\R^n) \big) \, = \, \big( L_2^r(\R^n) \bar\otimes \mathcal{R}_\Theta \big) \big( L_2^c(\R^n) \bar\otimes \mathcal{R}_\Theta \big),
\]
where the operator space structure for $\omega \in M_m(L_\infty(\mathcal{R}_\Theta; L_1(\R^n)))$ is  
$$\inf \left\{ \Big\| \summ_k \alpha_k \otimes e_{1k} \Big\|_{M_m(L_2^r(\R^n) \bar\otimes \mathcal{R}_\Theta \otimes R)} \Big\| \summ_k \beta_k \otimes e_{k1} \Big\|_{M_m(L_2^c(\R^n) \bar\otimes \mathcal{R}_\Theta \otimes C)} \right\}$$ where the infimum runs over all possible factorizations $\omega = \summ_k \alpha_k \beta_k$.  
Now, the contractivity of $\sigma_\Theta: L_1(\R_\Theta) \to L_\infty(\RR_\Theta; L_1(\R^n))$ follows easily from Proposition \ref{BMO.IsoLp} \ref{BMO.Lp1}. Particular cases of this kind of spaces |over finite von Neumann algebras or discrete $\ell_1$ spaces| have been proved to interpolate in the expected way with the corresponding $L_p$ scale \cite{JDoob,JP}. The lack of an available argument in the literature for the general case has led us to avoid the case $1 < p < 2$ in Proposition \ref{BMO.IsoLp}. This contractivity result is unnecessary for our goals.}
\end{Bremark}

Observe that 
$$\sigma_\Theta(\lambda_\Theta(f)) = \int_{\R^n} f(\xi) \, (\exp_\xi \otimes \lambda_\Theta(\xi)) \, d\, \xi \quad \mbox{for} \quad f \in \S(\R^n).$$
Clearly such element is invariant under the group of trace preserving automorphisms $\beta_z$ given by $\beta_z = \sigma_0^{-z} \otimes \sigma_\Theta^z$.
Let us denote by $(\mathrm{X} \bar\otimes \RR_\Theta)^\beta$ the $\beta$-invariant part of the $\mathrm{X} \bar\otimes \RR_\Theta$ with $\mathrm{X}$ any of the Euclidean function spaces in Proposition \ref{BMO.IsoLp}. We need to see that $(\mathrm{X} \bar\otimes \RR_\Theta)^\beta$
coincides with the image of $\sigma_\Theta$ and that the $\beta$-invariant subspace is complemented. Let us start with the complementation.

\begin{Bproposition} \label{BMO.complemented}
The following subspaces
\begin{enumerate}[label={\rm \roman*)}]
    \item $(L_2^\dagger(\R^n) \bar\otimes \RR_\Theta)^\beta \subset L_2^\dagger(\R^n) \bar\otimes \RR_\Theta$ for $\dagger \in \{r,c\}$,

\vskip2pt

    \item $(L_p(\R^n) \bar\otimes \RR_\Theta)^\beta \subset L_p(\R^n) \bar\otimes \RR_\Theta$ for any $2 \le p \le \infty$,
    
\vskip2pt

    \item $(\BMO_\dagger(\R^n) \bar\otimes \RR_\Theta)^\beta \subset \BMO_\dagger(\R^n) \bar\otimes \RR_\Theta$ for $\dagger \in \{c,r,r \wedge c\}$,
\end{enumerate}
are completely complemented as operator spaces in the respective ambient spaces.
\end{Bproposition}

\dem By amenability of $\R^n$, let $m \in L_\infty(\R^n)^\ast$ be an invariant mean and let $m_\alpha$ be a sequence of probability measures in $L_1(\R^n)$ which approximate $m$. Given $\omega$ in $L_p(\R^n) \bar\otimes \RR_\Theta$, the function $z \mapsto \beta_z \omega$ sits in the space $L_\infty(\R^n) \bar\otimes L_p(\R^n) \bar\otimes \RR_\Theta$, so
\[
  P_\alpha(\omega) = (m_\alpha \otimes id \otimes id)(\beta_z \omega)
\]
defines a family of completely positive operators
\[
  P_\alpha: L_\infty(\R^n) \bar\otimes L_p(\R^n) \bar\otimes \RR_\Theta \to L_p(\R^n) \bar\otimes \RR_\Theta.
\]
Since the image is in a dual space, we use $\CB(\mathrm{X},\mathrm{Y}^\ast) = (\mathrm{X} \widehat{\otimes} \mathrm{Y})^\ast$ and Banach-Alaoglu theorem. Let $P$ be an accumulation point of $(P_\alpha \circ \beta)_\alpha$ in the weak-$*$ topology. $P$ gives a cb-bounded projection into the $\beta$-invariant part. We have only used that $L_p(\R^n)$ is a dual space for weak-$*$ compactness. Therefore, the same proof applies to $\BMO_\dagger(\R^n)$. The projections $P:L_p(\R^n) \bar\otimes \RR_\Theta \to L_p(\R^n) \bar\otimes \RR_\Theta$ form compatible family: they are restrictions of a map defined in the sum of the above spaces. 
\fin

\begin{Bproposition} \label{BMO.EqualSets}
We have
\begin{enumerate}[label={\rm \roman*)}, ref={\rm \roman*)}]
\item $\sigma_\Theta(L_2^\dagger(\RR_\Theta)) = (L_2^\dagger(\R^n) \bar\otimes \RR_\Theta)^\beta$ for $\dagger \in \{r,c\}$,

\vskip2pt

\item $\sigma_\Theta(L_p(\RR_\Theta)) = (L_p(\R^n) \bar\otimes \RR_\Theta)^\beta$ for any $2 \leq p \leq \infty$,

\vskip2pt

\item $\sigma_\Theta(\BMO_\dagger(\RR_\Theta)) = (\BMO(\R^n)_\dagger \bar\otimes \RR_\Theta)^\beta$ for $\dagger \in \{r,c,r \wedge c\}$.
\end{enumerate}
\end{Bproposition}

\dem
Since the spaces $(L_p(\R^n) \bar\otimes \RR_\Theta)^\beta$ are complemented subspaces, it is enough to prove the identity for $p=2$ and $p = \infty$ and interpolation will yield the result for $2 < p < \infty$ since the maps $\sigma_\Theta: L_p(\RR_\Theta) \to L_p(\R^n) \bar\otimes \RR_\Theta$ are compatible. The same argument gives that the $L_2$ case follows by interpolation between $L_2^c(\RR_\Theta)$ and $L_2^r(\RR_\Theta)$. This reduces the proof to the row/column cases, the case $p=\infty$ and $\BMO(\RR_\Theta)$. We shall only prove it for columns and for $p =\infty$, since the argument is similar in $\BMO$. Let us define the map $\mathcal{W}:L_2^c(\R^n) \bar\otimes \RR_\Theta \to L_2^c(\R^n) \bar\otimes \RR_\Theta$ by extension of $\exp_\xi \otimes \lambda_\Theta(\eta) \mapsto \exp_\xi \otimes \lambda_\Theta(\xi)  \lambda_\Theta(\eta)$. A calculation easily yields that $\mathcal{W}$ is a complete isometry. The same follows in the row case if one takes the map $\exp_\xi \otimes \lambda_\Theta(\eta) \mapsto \exp_\xi \otimes \lambda_\Theta(\eta) \lambda_\Theta(\xi)$ instead. We have that $\mathcal{W}$ gives an isomorphism between $L_2^c(\R^n) \otimes \1$ and $\sigma_\Theta[\lambda_\Theta[L_2(\R^n)]]$. We also have that $\mathcal{W}$ intertwines the action $\beta_z$ as follows
\begin{equation*}
  \label{BMO.EqualSets.eq}
  \xymatrix{
    L_2^c(\R^n) \bar\otimes \RR_\Theta \ar[rr]^{\mathcal{W}} \ar[d]^{\sigma_\Theta^z} & & L_2^c(\R^n) \bar\otimes \RR_\Theta \ar[d]^{\beta_z}\\
    L_2^c(\R^n) \bar\otimes \RR_\Theta \ar[rr]^{\mathcal{W}} & & L_2^c(\R^n) \bar\otimes \RR_\Theta.
  }
\end{equation*} 
Therefore, the subspace fixed by $\beta$ corresponds under $\mathcal{W}$ with the subspace fixed by $id \otimes \sigma_\Theta$. But evaluating such space against every $\varphi \otimes id$, with $\varphi \in L_2^r(\R^n)$, gives that the fixed subspace of $id \otimes \sigma_\Theta$ is $L_2^c(\R^n)$ tensored with the subspace fixed by $\sigma_\Theta$. Such subspace is $\C \1$. Indeed, if $\varphi \in \RR_\Theta$ is invariant under $\sigma_\Theta$ we obtain that $\varphi = \lambda_\Theta(\psi)$, where $\psi \in \S(\R^n)'$ is a distribution supported on $\{0\}$. But such distribution is a linear combination of distributions of the form $\langle \psi, f \rangle = f^{(k)}(0)$, where $f \in \S(\R^n)$. The derivatives with $k > 0$ give rise to unbounded elements and so we obtain that $\psi$ has to be a multiple of $\delta_{0}$ or, equivalently, that $\varphi \in \C \1$. 

The case of $p = \infty$ follows similarly. We first define a normal $\ast$-homomorphim $\mathcal{U}: L_\infty(\R^n) \bar\otimes \RR_\Theta \to L_\infty(\R^n) \bar\otimes \RR_\Theta$ by extension of $\exp_\xi \otimes \lambda_\Theta(\eta) \mapsto \exp_{\eta + \xi} \otimes \lambda_\Theta(\xi)$. To prove that such map is a $\ast$-homomorphism we can implement it spatially with techniques analogous to that of Corollary \ref{Planc+Corep}. We have that $\mathcal{U}$ carries $\1 \otimes \RR_\Theta$ in $\sigma_\Theta[\RR_\Theta]$ and that it intertwines the actions in the expected way. Proceeding like in the case $p =2$ we can conclude. 

The case of $\BMO_{\dagger}$ can be deduced from a similar result for mixed spaces. First we note that the result for $\BMO_{r \wedge c}$ follows from the corresponding ones for $\BMO_r$ and $\BMO_c$, we shall only prove it for $\BMO_c$. Fix a Euclidean ball $\mathrm{B} \subset \R^n$ and consider the following operator-valued inner product 
\[
  \langle f, f \rangle_\mathrm{B} = \mean_\mathrm{B} |f_s|^2 \, ds - \Big| \mean_\mathrm{B} f_s \, ds \Big|^2 \quad \mbox{for} \quad f \in L_\infty(\R^n) \otimes_{\mathrm{alg}} \RR_\Theta.
\]
Let $\H^c_\Theta(\mathrm{B})$ denote the corresponding Hilbert module over $\RR_\Theta$ and let $\H^c_\Theta$ be the direct sum, in the $\ell_\infty$-sense, of $\H^c_\Theta(\mathrm{B})$ over all balls $\mathrm{B}$. Clearly $\BMO_c(\RR_\Theta)$ embeds in $\H^c_\Theta$ and we have that 
\begin{equation*}
  \label{BMO.EqualSets.eq2}
  \xymatrix{
    \H^c_\Theta \ar[rr]^{id \otimes \sigma_\Theta} & & \H^c_\Xi\\
    \BMO_c(\RR_\Theta) \ar@{^{(}->}[u] \ar[rr]^{\sigma_\Theta} & & \BMO_c(\R^n) \bar\otimes \RR_\Theta \ar@{^{(}->}[u],
  }
\end{equation*}
where $\Xi$ is the $2n \times 2n$-matrix $\Xi = 0 \otimes \Theta$. Now, we can define a map preserving the operator-valued inner product (and thus an isometry) $\mathcal{W}_c: \H_\Xi^c \to \H_\Xi^c$ by extension of $\exp_{\xi_1} \otimes \exp_{\xi_2} \otimes \lambda_\Theta(\eta) \mapsto \exp_{\xi_1} \otimes \exp_{\xi_2  + \eta} \otimes \lambda_\Theta(\eta)$ for every ball $\mathrm{B}$. Such map carries the copy of $\H_\Theta^c \otimes \1$ that lives in the first and third tensor components into $\sigma_\Theta[\H_\Theta^c]$ and proceeding like in the previous cases we get that $\sigma_\Theta[\H_\Theta^c]$ coincides with the subspace of $\H_\Xi^c$ invariant under the group of automorphisms $\beta_z = id \otimes \sigma_0^{-z} \otimes \sigma_\Theta^{z}$. That result restricts to $\BMO_c$. \fin

\begin{Bproposition} \label{BMO.IsoLp2}
We have complete isometries$\hskip1pt :$ 
\begin{enumerate}[label={\rm \roman*)}]
    \item \label{BMO.Iso2Lp1}
    $\sigma_\Theta: L_2^\dagger(\RR_\Theta) \longrightarrow L_2^{\dagger}(\R^n) \bar\otimes \RR_\Theta$ for $\dagger \in \{r,c\}$,

\vskip2pt

    \item \label{BMO.Iso2Lp2}
    $\sigma_\Theta: L_p(\RR_\Theta) \longrightarrow L_p(\R^n) \bar\otimes \RR_\Theta$ for any $2 \le p \leq \infty$,
    
\vskip2pt 

    \item \label{BMO.Iso2Lp3}
    $\sigma_\Theta: \BMO_{\dagger}(\RR_\Theta) \longrightarrow \BMO_{\dagger}(\RR_\Theta) \bar\otimes \RR_\Theta$ for $\dagger \in \{r, c, r \wedge c\}$.
\end{enumerate} 
\end{Bproposition}

\dem Assertions \ref{BMO.Iso2Lp1} and \ref{BMO.Iso2Lp3} were proved in the proof of Proposition \ref{BMO.IsoLp}. Assertion \ref{BMO.Iso2Lp2} for $p=\infty$ was already justified in Corollary \ref{Planc+Corep}. The rest of the cases trivially follow by complementation and complex interpolation from our results above. \fin
 
\begin{Btheorem} \label{BMO.thmB}
We have
\begin{eqnarray*}
\big[\Hardy_1(\RR_\Theta), \BMO(\RR_\Theta) \big]_{\theta} & = & \big[L_1^\circ(\RR_\Theta), \BMO(\RR_\Theta) \big]_{\theta} \\
& = & \big[\Hardy_1(\RR_\Theta), L_\infty(\RR_\Theta) \big]_{\theta} \ = \ L_p(\RR_\Theta)
\end{eqnarray*}
for $p = \frac{1}{1-\theta}$. All isomorphisms above hold in the category of operator spaces. 
\end{Btheorem}

\dem Since $L_2(\RR_\Theta)$ is reflexive
\begin{eqnarray*}
\big[ L_2(\RR_\Theta), \BMO(\RR_\Theta) \big]_\theta & = & \big[ L_2(\RR_\Theta), \BMO(\RR_\Theta) \big]^\theta \\ 
[-1pt] & = & \big[ \sigma_\Theta(L_2(\RR_\Theta)), \sigma_\Theta(\BMO(\RR_\Theta)) \big]^\theta \\
[-1pt] & = & \big[ P ( L_2(\R^n) \bar\otimes \RR_\Theta ), P( \BMO(\R^n) \bar\otimes \RR_\Theta) \big]^\theta \\
& = & P \big( \big[ L_2(\R^n) \bar\otimes \RR_\Theta, \BMO(\R^n) \bar\otimes \RR_\Theta \big]^\theta \big) \\
& = & P \big( [L_2(\R^n), \BMO(\R^n)]^\theta \bar\otimes \RR_\Theta \big) \\
& = & P \big(  L_p(\R^n) \bar\otimes \RR_\Theta \big) \ = \ \sigma_\Theta \big( L_p(\RR_\Theta) \big) \ = \ L_p(\RR_\Theta)  
\end{eqnarray*}
for $p = \frac{2}{1-\theta}$. Indeed, the second identity follows from Proposition \ref{BMO.IsoLp2}, which gives $\mathrm{X} = \sigma_\Theta(\mathrm{X})$ completely isomorphic for $\mathrm{X} = L_2(\RR_\Theta)$ and $\mathrm{X} = \BMO(\RR_\Theta)$. The third and fourth identities follow from Proposition \ref{BMO.complemented}, which shows that $\sigma_\Theta(\mathrm{X})$ can be identified with $P(\mathrm{Z})$ where $\mathrm{Z}$ is the ambient space of operator-valued functions in $\R^n$ associated to $\mathrm{X}$. Moreover, since $P$ is a bounded projection, it commutes with the complex interpolation functor by complementation. The fifth identity follows from Lemma \ref{BMO.InterTensor} and the sixth one from Mei's interpolation theorem \cite{MeiH1}. The last two identities apply from Propositions \ref{BMO.complemented} and \ref{BMO.IsoLp2} again. Once this is known we use the reflexivity of $L_2(\RR_\Theta)$ and duality $\Hardy_1(\RR_\Theta)^\ast = \BMO(\RR_\Theta)$ to obtain 
$$\big[ \Hardy_1(\RR_\Theta), L_2(\RR_\Theta) \big]_\theta^* \, = \, \big[ \BMO(\RR_\Theta), L_2(\RR_\Theta) \big]^\theta \, = \, L_{q'}(\RR_\Theta)$$ 
for $q = \frac{2}{2-\theta}$. This shows that $[ \Hardy_1(\RR_\Theta), L_2(\RR_\Theta)]_\theta$ must be reflexive and we get 
$$\big[ \Hardy_1(\RR_\Theta), L_2(\RR_\Theta) \big]_\theta \, = \, L_{q}(\RR_\Theta).$$
The interpolation results in the statement follow from Wolff's theorem \cite{Wolff}, which states that if $\mathrm{X}_1$, $\mathrm{X}_2$, $\mathrm{X}_3$, $\mathrm{X}_4$ are spaces with $\mathrm{X}_1 \cap \mathrm{X}_4$ dense inside both $\mathrm{X}_2$ and $\mathrm{X}_3$, then 
$$\mathrm{X}_2 = [\mathrm{X}_1, \mathrm{X}_3]_{\theta_1} \mbox{ and } \mathrm{X}_3 = [\mathrm{X}_2, \mathrm{X}_4]_{\theta_2} \ \Rightarrow \ \mathrm{X}_2 = [\mathrm{X}_1, \mathrm{X}_4]_{\vartheta_1} \mbox{ and } \mathrm{X}_3 = [\mathrm{X}_1, \mathrm{X}_4]_{\vartheta_2}$$
where $\vartheta_1 = \theta_1 \theta_2 / (1 - \theta_1 + \theta_1 \theta_2)$ and $\vartheta_2 = \theta_2 / (1 - \theta_1 + \theta_1 \theta_2)$.
Taking
$$\begin{array}{rcl rcl rcl rcl}
\mathrm{X}_1 \!\!\! & = & \!\!\! \Hardy_1(\RR_\Theta), & \mathrm{X}_2 \!\!\! & = & \!\!\! L_{\frac{4}{3}}(\RR_\Theta), & \mathrm{X}_3 \!\!\! & = & \!\!\! L_2(\RR_\Theta), & \mathrm{X}_4 \!\!\! & = & \!\!\! L_4(\RR_\Theta), \\
\mathrm{Z}_1 \!\!\! & = & \!\!\! \Hardy_1(\RR_\Theta), & \mathrm{Z}_2 \!\!\! & = & \!\!\! L_2 \hskip1pt (\RR_\Theta), & \mathrm{Z}_3 \!\!\! & = & \!\!\! L_4(\RR_\Theta), & \mathrm{Z}_4 \!\!\! & = & \!\!\! \BMO(\RR_\Theta),
\end{array}$$
we first obtain, using $\mathrm{X}_j$-spaces and the interpolation of $\Hardy_1(\RR_\Theta)$ with $L_2(\RR_\Theta)$, that $\Hardy_1(\RR_\Theta)$ and $L_4(\RR_\Theta)$ interpolate in the expected way. Then, using the same procedure with the $\mathrm{Z}_j$-spaces and the interpolation of $L_2(\RR_\Theta)$ with $\BMO(\RR_\Theta)$, we finally get the expected result for the bracket $$\big[ \Hardy_1(\RR_\Theta), \BMO(\RR_\Theta) \big]_\theta.$$The other two brackets in the statement can be treated analogously. \fin

\begin{Bremark} \label{BMO.ClassicalAppli}
\emph{It is worth mentioning that our techniques have at least another potential application in the abelian case. Let $(X,\mu)$ be a $G$-space with a $G$-invariant measure. In that case, we can identify $X$ and $H \times G/H$ as measure spaces, where $H$ is the stabilizer and we have the following Fubini-type identity
$$\int_G f(g^{-1} \, x) \, d\mu_G(g) = \int_{G/H = X} \int_{H} f(h^{-1} \, g^{-1} \, x)  \, d\mu_H(h) \, d\mu(x),$$
see \cite[Chapter 2]{Fo}. If the stabilizer is compact we can exchange integration in $X$ and integration in $G$ in a way analogous to the Fubini-type identity which relates $\tau_\Theta$ and $\sigma_\Theta$. If there is a natural definition of $\BMO(G)$, either with averages over the balls of an invariant measure or with translation-invariant semigroups, and that $\BMO$ interpolates, then we can transfer the interpolation to $\BMO(X)$ provided that $G$ is amenable. This seems to be a very direct approach for proving interpolation of $G$-invariant $\BMO$-spaces over $X= G/K$, where $G$ is a solvable and unimodular Lie group and $K$ is a compact subgroup.}
\end{Bremark}

\subsubsection*{\emph{B.4.} An auxiliary density result}
\label{AppB.Density}

Let us write in what follows $\S_\Theta^\circ$ for the kernel of the trace functional $\tau_\Theta:\S_\Theta \to \C$, which is of course continuous over $\S_\Theta$. It is trivial that $\S_\Theta^\circ \subset \Hardy^\dagger_1(\RR_\Theta)$. We are going to see that it is in fact dense. It will be an easy consequence of the fact that $\sigma_\Theta: \BMO_\dagger(\RR_\Theta) \to \BMO_\dagger(\R^n) \bar\otimes \RR_\Theta$, for $\dagger \in \{r, c, r \wedge c\}$ are normal and complete isometries. Taking preduals we obtain a complete and surjective projection 
\[
  (\sigma_\Theta)_\ast: \Hardy_1^\dagger(\R^n) \widehat\otimes (\RR_\Theta)_\ast \longrightarrow \Hardy^\dagger_1(\RR_\Theta),
\]
for $\dagger \in \{r, c, r + c\}$. We are just going to need that such map carries $\S^\circ(\R^n) \otimes_\pi \S_\Theta$ into $\S^\circ_\Theta$ but indeed much more is true and the map $(\sigma_\Theta)_\ast$ can be explicitly described as a diagonal restriction multiplier. That is, it satisfies the following commutative diagram, where $\S_0(\R^n)$ is the subclass of Schwartz functions with $f(0) = 0$
\begin{equation*}
  \xymatrix{
    \Hardy_1^\dagger(\R^n) \widehat\otimes (\RR_\Theta)_\ast \ar[rr]^{(\sigma_\Theta)_\ast} & &  \Hardy_1^\dagger(\RR_\Theta) \ar[d]^{\lambda_\Theta^{-1}}\\
    \S_0(\R^n) \otimes_\pi \S(\R^n) \ar[u]^{\lambda_0 \otimes \lambda_\Theta} \ar[rr]^-{f \mapsto f{|}_{\Delta}} & &  \S_0(\R^n)
  }
\end{equation*}
Now, the proof of the density is immediate.

\begin{Bcorollary} \label{Bdensity}
$\S_\Theta^\circ$ is dense inside $\Hardy^\dagger_1(\RR_\Theta)$ for $\dagger \in \{r,c, r + c\}$.
\end{Bcorollary}

\dem
We just have to use that $\S^\circ(\R^n) \otimes_\pi \S_\Theta \subset \Hardy_1^\dagger(\R^n) \widehat\otimes (\RR_\Theta)_\ast$ is a dense subset. Since $(\sigma_\Theta)_\ast(\S^\circ(\R^n)) \subset \S^\circ_\Theta$ and the image under a projection of a dense set is a dense set we conclude. \fin

\vskip3pt

\noindent \textbf{Acknowledgement.} A. Gonz\'alez-P\'erez was partially supported by European Research Council Consolidator Grant 614195. M. Junge was partially supported by the NSF DMS-1501103. J. Parcet was partially supported by European Research Council Starting Grant 256997 and CSIC Grant PIE 201650E030. All authors are also supported in part by ICMAT Severo Ochoa Grant SEV-2015-0554 (Spain). M. Junge would like to thank Alain Connes for a very encouraging discussion on diagonals for noncommutative algebras at an Oberwolfach meeting. J. Parcet would like to express his gratitude to the Math Department of the University of Illinois at Urbana-Champaign for their hospitality along his visit in 2015, during which this paper was enormously developed. 

\bibliographystyle{amsplain}

\vskip10pt

\enlargethispage{2cm}

\hfill \noindent \textbf{Adri\'an M. Gonz\'alez-P\'erez} \\
\null \hfill K U Leuven - Departement wiskunde \\ 
\null \hfill 200B Celestijnenlaan, 3001 Leuven. Belgium 
\\ \null \hfill\texttt{adrian.gonzalezperez@kuleuven.be}

\vskip2pt

\hfill \noindent \textbf{Marius Junge} \\
\null \hfill Department of Mathematics
\\ \null \hfill University of Illinois at Urbana-Champaign \\
\null \hfill 1409 W. Green St. Urbana, IL 61891. USA \\
\null \hfill\texttt{junge@math.uiuc.edu}

\vskip2pt

\hfill \noindent \textbf{Javier Parcet} \\
\null \hfill Instituto de Ciencias Matem{\'a}ticas \\ \null \hfill
Consejo Superior de
Investigaciones Cient{\'\i}ficas \\ \null \hfill C/ Nicol\'as Cabrera 13-15.
28049, Madrid. Spain \\ \null \hfill Department of Mathematics
\\ \null \hfill University of Illinois at Urbana-Champaign \\
\null \hfill 1409 W. Green St. Urbana, IL 61891. USA \\
\null \hfill\texttt{javier.parcet@icmat.es}

\end{document}